\documentclass[12pt,a4paper]{amsart}
\usepackage[utf8]{inputenc}
\usepackage[T1]{fontenc}
\usepackage{amssymb}
\usepackage{mathrsfs} 
\usepackage{amsmath}
\usepackage{amsthm}
\usepackage{amscd}
\usepackage{paralist} 
\setdefaultenum{\upshape(i)}{}{}{}
\usepackage[pdftex,colorlinks,linkcolor=black,filecolor=black,citecolor=black,urlcolor=blue]{hyperref}

\allowdisplaybreaks

\newtheoremstyle{mythmstyle}%
{1.5\baselineskip}
{\baselineskip}
{\itshape}
{}
{\bf}
{}
{0pt}
{} 

\newtheoremstyle{mydefstyle}%
{1.5\baselineskip}
{\baselineskip}
{}
{}
{\bf}
{}
{0pt}
{} 

\newtheoremstyle{mypreuvestyle}%
{\baselineskip}
{\baselineskip}
{}
{}
{\em}
{}
{0pt}
{} 

\newif\ifmynonumberenvi\mynonumberenvitrue
\theoremstyle{mythmstyle}
\newtheorem{proclaimmythm}[equation]{} 
\newtheorem*{proclaimmythm*}{}
\newenvironment{proclaim}[2][*]{\ifx*#1\mynonumberenvitrue\begin{proclaimmythm*}{\bf#2.} \ignorespaces\else\mynonumberenvifalse\begin{proclaimmythm}{.\kern0.5em\bf#2.}\label{#1} \ignorespaces\fi}{\ifmynonumberenvi\end{proclaimmythm*}\else\end{proclaimmythm}\fi}
\renewenvironment{proclaim}[2][*]{\ifx*#1\begin{proclaimmythm}{.\kern0.5em\bf#2.} \ignorespaces\else\begin{proclaimmythm}{.\kern0.5em\bf#2.}\label{#1} \ignorespaces\fi}{\end{proclaimmythm}}

\theoremstyle{mydefstyle}
\newtheorem{proclaimmydef}[equation]{} 
\newtheorem*{proclaimmydef*}{}
\newenvironment{definition}[2][*]{\ifx*#1\mynonumberenvitrue\begin{proclaimmydef*}{\bf#2.}\else\mynonumberenvifalse\begin{proclaimmydef}{.\kern0.5em\bf#2.}\label{#1} \ignorespaces\fi}{\ifmynonumberenvi\end{proclaimmydef*}\else\end{proclaimmydef}\fi}
\renewenvironment{definition}[2][*]{\ifx*#1\begin{proclaimmydef}{.\kern0.5em\bf#2.} \ignorespaces\else\begin{proclaimmydef}{.\kern0.5em\bf#2.}\label{#1} \ignorespaces\fi}{\end{proclaimmydef}}

\newenvironment{definition*}[1]{\begin{proclaimmydef*}{\bf#1.} \ignorespaces}{\end{proclaimmydef*}}

\def\QEDbox{\hbox{\lower2.3pt\vbox{\hrule\hbox
   {\vrule\kern1pt\vbox{\kern1.7pt\hbox{$\scriptstyle
   QED$}\kern.6pt}\kern1pt\vrule}\hrule}}}
\def\QED{\hskip0.01em plus 40pt\null{} \null\nobreak\hfill
   \kern3pt\QEDbox} 
\newcommand\QEDici{\\\noalign{\vskip-\baselineskip\smash{\hbox to\linewidth{\vrule width0pt \hfill\global\QEDdejaplacetrue\QEDbox}}\vskip-\baselineskip}}
\newif\ifQEDdejaplace\QEDdejaplacefalse

\theoremstyle{mypreuvestyle}
\newtheorem*{proclaimmypreuve}{}
\newenvironment{preuve}[1][*]{\begin{proclaimmypreuve}{\ifx*#1{\em Proof.}\else{\em#1.}\fi} \ignorespaces}{\ifQEDdejaplace\global\QEDdejaplacefalse\else\QED\fi\end{proclaimmypreuve}}

\newif\ifmynonumberequation\mynonumberequationtrue

\makeatletter
\numberwithin{equation}{section}
\newenvironment{moneq}[1][*]{\ifx*#1\mynonumberequationtrue\begin{equation*}\else\mynonumberequationfalse\begin{equation}\label{#1}\fi}{\ifmynonumberequation\end{equation*}\@ignoretrue\else\end{equation}\@ignoretrue\fi\ignorespaces}
\makeatother

\newcommand{\recalf}[1]{(\ref{#1})}
\newcommand{\recals}[1]{{\S\ref{#1}}}
\newcommand{\recalt}[1]{{[\ref{#1}]}}
\newcommand{\recaltt}[2]{{[\ref{#2}.\ref{#1}]}}

\def\dcap_#1{\mathchoice{%
          {\textstyle\bigcap\limits_{#1}}}%
          {\underset{#1}\cap}%
          {\underset{#1}\cap}%
          {\underset{#1}\cap}}
\def\dcup_#1{\mathchoice{%
          {\textstyle\bigcup\limits_{#1}}}%
          {\underset{#1}\cup}%
          {\underset{#1}\cup}%
          {\underset{#1}\cup}}
\def\ddcap_#1^#2{\mathchoice{%
          {\textstyle\bigcap\limits_{#1}^{#2}}}%
          {\overset{#2}{\underset{#1}\cap}}%
          {\overset{#2}{\underset{#1}\cap}}%
          {\overset{#2}{\underset{#1}\cap}}}
\def\ddcup_#1^#2{\mathchoice{%
          {\textstyle\bigcup\limits_{#1}^{#2}}}%
          {\overset{#2}{\underset{#1}\cup}}%
          {\overset{#2}{\underset{#1}\cup}}%
          {\overset{#2}{\underset{#1}\cup}}}

\newcommand\bigrestricted{{\kern1pt\vrule height3.3ex depth1.7ex width0.6pt\kern1pt}}
\newcommand\caprestricted{{\kern1pt\smash{\vrule height1.7ex depth0.7ex width0.6pt}\relax\vrule width0pt depth0pt\kern2pt}}
\newcommand\midrestricted{{\kern1pt\vrule height1.7ex depth0.9ex width0.6pt\kern1pt}}
\newcommand\bmidrestricted{{\kern1pt\vrule height2.7ex depth1.7ex width0.6pt\kern1pt}}
\newcommand\NN{\mathbf{N}}
\newcommand\restricted{{\kern1pt\vrule height1.3ex depth0.5ex width0.6pt\kern1pt}}

\def\oversetalign#1\to#2{\mathbin{\smash{\overset{{\text{\rlap{\hss#1}}}}{#2}}}}
\def\oversettext#1\to#2{\mathbin{\smash{\overset{{\text{{#1}}}}{#2}}}}

\newcommand\norm[2][1]{\ifcase#1
   \left\Vert#2\right\Vert 
   \or\Vert#2\Vert 
   \or\bigl\Vert#2\bigr\Vert 
   \or\Bigl\Vert#2\Bigr\Vert 
   \or\biggl\Vert#2\biggr\Vert 
   \or\Biggl\Vert#2\Biggr\Vert 
   \else\Vert#2\Vert\fi} 

\newcommand\hominprod[4][1]{\ifcase#1
   \left\langle\homind{#4}#2,#3\right\rangle 
   \or\langle\homind{#4}#2,#3\rangle 
   \or\bigl\langle\homind{#4}#2,#3\bigr\rangle 
   \or\Bigl\langle\homind{#4}#2,#3\Bigr\rangle 
   \or\biggl\langle\homind{#4}#2,#3\biggr\rangle 
   \or\Biggl\langle\homind{#4}#2,#3\Biggr\rangle 
   \else\langle\homind{#4}#2,#3\rangle\fi} 

\newcommand\homsuperinprod[4][1]{\ifcase#1
   \supinprsym{#4}\left(#2,#3\right) 
   \or\supinprsym#4(#2,#3) 
   \or\supinprsym#4\bigl(#2,#3\bigr) 
   \or\supinprsym#4\Bigl(#2,#3\Bigr) 
   \or\supinprsym#4\biggl(#2,#3\biggr) 
   \or\supinprsym#4\Biggl(#2,#3\Biggr) 
   \else\supinprsym#4(#2,#3)\fi} 

\newcommand\inprod[3][1]{\ifcase#1
   \left\langle#2,#3\right\rangle 
   \or\langle#2,#3\rangle 
   \or\bigl\langle#2,#3\bigr\rangle 
   \or\Bigl\langle#2,#3\Bigr\rangle 
   \or\biggl\langle#2,#3\biggr\rangle 
   \or\Biggl\langle#2,#3\Biggr\rangle 
   \else\langle#2,#3\rangle\fi} 
   
\newcommand\inprodsym{\inprod\cdot\cdot}

\newcommand\inprodd[3][1]{\ifcase#1
   \left\langle\mkern-3mu\left\langle#2,#3
        \right\rangle\mkern-3mu\right\rangle 
   \or\langle\mkern-3mu\langle#2,#3
        \rangle\mkern-3mu\rangle 
   \or\bigl\langle\mkern-4mu\bigl\langle#2,#3
        \bigr\rangle\mkern-4mu\bigr\rangle 
   \or\Bigl\langle\mkern-6mu\Bigl\langle#2,#3
        \Bigr\rangle\mkern-6mu\Bigr\rangle 
   \or\biggl\langle\mkern-8mu\biggl\langle#2,#3
        \biggr\rangle\mkern-8mu\biggr\rangle 
   \or\Biggl\langle\mkern-10mu\Biggl\langle#2,#3
        \Biggr\rangle\mkern-10mu\Biggr\rangle 
   \else\langle\mkern-3mu\langle#2,#3
        \rangle\mkern-3mu\rangle\fi} 

\newcommand\supinprsym{\mathscr{S}}

\newcommand\superinprod[3][1]{\ifcase#1
   \supinprsym\left(#2,#3\right) 
   \or\supinprsym(#2,#3) 
   \or\supinprsym\bigl(#2,#3\bigr) 
   \or\supinprsym\Bigl(#2,#3\Bigr) 
   \or\supinprsym\biggl(#2,#3\biggr) 
   \or\supinprsym\Biggl(#2,#3\Biggr) 
   \else\supinprsym(#2,#3)\fi} 

\newcommand\contrfoper{{\iota}}

\newcommand\contrf[3][1]{\ifcase#1
   \contrfoper\left(#2\right)#3 
   \or\contrfoper(#2)#3 
   \or\contrfoper\bigl(#2\bigr)#3 
   \or\contrfoper\Bigl(#2\Bigr)#3 
   \or\contrfoper\biggl(#2\biggr)#3 
   \or\contrfoper\Biggl(#2\Biggr)#3 
   \else\contrfoper(#2)#3\fi} 

\newcount\mysubsectioncount\mysubsectioncount=0
\newcommand\firstofmysubsection{\mysubsectioncount=0}
\newcommand\mysubsection[2]{\vskip2\baselineskip\advance\mysubsectioncount by1\noindent{\textbf{\ref{#1}.\the\mysubsectioncount{} #2}}\medskip}

\newcommand\mysubsectionnonumber[2]{\vskip1.5\baselineskip\noindent{\textbf{#2}}\medskip}

\newcommand\ad{\mathrm{ad}}
\newcommand\Ad{\mathrm{Ad}}

\newcommand\atlas{\mathcal{U}}
\newcommand\Aut{\mathrm{Aut}}
\newcommand\BCH{\mathcal{B}_\text{BCH}}
\newcommand\Ber{\mathrm{Ber}}

\newcommand\body{\mathbf{B}}
\newcommand\CA{\mathcal{A}}
\newcommand\CC{\mathbf{C}}

\newcommand\conjugate{\mathfrak{C}}
\newcommand\Cyl{\mathrm{Cyl}}

\newcommand\Dense{\mathcal{D}}
\newcommand\Density{\mathscr{D}}
\newcommand\Det{\mathrm{Det}}
\newcommand\dfracp[2]{{\displaystyle\fracp{#1}{#2}}}
\newcommand\eexp{\mathrm{e}}

\newcommand\End{\mathrm{End}}
\newcommand\extder{\mathrm{d}}
\newcommand\extpower[1]{{}\kern-4pt\raise1.2ex\hbox{$\scriptstyle#1$}\kern3pt}
\newcommand\FgroupHilbert{\mathbf{F}}

\newcommand\fint{{\mathcal{I}_f}}
\newcommand\Fourier{\mathcal{F}}
\newcommand\Fourierodd{\smash{\overset{\kern2pt\raise-2pt\hbox{$\scriptscriptstyle \mathrm o$}}{\Fourier}}}
\newcommand\fracp[2]{\frac{\partial#1}{\partial#2}}
\newcommand\framebundle{\mathcal{F}}
\newcommand\Gextension{\relax\ifmmode\mathbf{G}\else$\mathbf{G}$-extension\fi}
\newcommand\Gl{\mathrm{Gl}}

\newcommand\Hilbert{\mathcal{H}}
\newcommand\Hilberth{\hat\Hilbert}
\newcommand\Hstarh{{\hat*}}
\newcommand\homind[1]{_{{}_{\scriptstyle#1}}\,}

\newcommand\id{\mathrm{id}}
\newcommand\ie{i.e.}

\newcommand\Inv{\mathcal{I}}

\newcommand\KK{\mathbf{K}}
\newcommand\KKA{\mathcal{K}}
\newcommand\labelIntroSHS{B}
\newcommand\labelIntroSUR{C}
\newcommand\labelIntroSURmod{D}
\newcommand\labelMainSHS{SHS\kern0.1em\relax}
\newcommand\labelMainSUR{SUR\kern0.1em\relax}
\newcommand\labelSHCP{A}
\newcommand\Leb{\lambda}
\newcommand\Liealg[1]{\mathfrak{#1}}

\newcommand\limn{\lim_{n\to\infty}}
\newcommand\mapob{\kern0.75em}
\newcommand\masection[1]{\newpage\section{#1}}
\newcommand\mfdmetric{\mathbf{g}}

\newcommand\mo{^{-1}}
\newcommand\Mob{\textrm{Möb}}
\newcommand\myquote[1]{``#1''}
\newcommand\nasset{{\underline{n}}}
\newcommand\nilpotent{\mathcal{N}}
\newcommand\oddif{^{\mathrm{odd}}}
\newcommand\oddp[1]{#1^{\textrm{odd dim.}}}
\newcommand\oneasmatrix{\mathbf{1}}
\newcommand\oneori[1]{i^{[\kern-0.9pt[#1]\kern-0.9pt]}}
\newcommand\oneormi[1]{(-i)^{[\kern-0.9pt[#1]\kern-0.9pt]}}
\newcommand\ood{o.o.d\,}
\newcommand\OSp{\mathrm{OSp}}
\newcommand\parity[1]{\vert#1\vert}
\newcommand\Phih{\widehat{\Phi}}
\newcommand\piBer{{\Ber_\pi}}

\newcommand\psib{\bar{\psi}}
\newcommand\psih{\hat{\psi}}
\newcommand\Psitriv[1]{(\kern-2.2pt[#1]\kern-2.2pt)}
\newcommand\psur{proto super unitary representation}

\newcommand\refmetnaam[2]{#1\ref{#2}}

\newcommand\rhob{\bar\rho}

\newcommand\rhoh{\hat\rho}
\newcommand\RR{\mathbf{R}}
\newcommand\scirc{\,{\raise 0.8pt\hbox{$\scriptstyle\circ$}}\,}
\def\sdlim_#1{\smash{\displaystyle\lim_{#1}}}

\newcommand\shifttag[1]{\kern#1&\kern-#1}
\newcommand\Sl{\mathrm{Sl}}
\newcommand\smooths{\mathbf{F}}

\newcommand\SScircle{\mathbf{S}}
\newcommand\stresd[1]{\textit{#1}}
\newcommand\stress[1]{\textbf{#1}}
\newcommand\superspmatrix{\mathbf{S}}
\newcommand\supp{\mathrm{supp}}
\newcommand\taub{\bar\tau}

\newcommand\tauh{\hat\tau}

\newcommand\trans{^{t}}

\newcommand\Vol{\mathrm{Vol}}

\newcommand\wod[1]{#1^{\textrm{even dim.}}}

\newcommand\zb{\overline{z}}
\newcommand\Ztwo{\ZZ/2\ZZ}
\newcommand\ZZ{\mathbf{Z}}

\begin{document}
\setdefaultleftmargin{2.5em}{2.1em}{1.87em}{1.7em}{1em}{1em}

\author{Gijs M. Tuynman}

\title{Super unitary representations revisited}

\address{Laboratoire Paul Painlevé, U.M.R. CNRS 8524 et UFR de Mathématiques, Université de Lille I, 59655 Villeneuve d'Ascq Cedex, France}
\email{FirstName[dot]LastName[at]univ-lille1[dot]fr}

\begin{abstract}
With the usual definition of a super Hilbert space and a super unitary representation, it is easy to show that there are lots of super Lie groups for which the left-regular representation is not super unitary. 
I will argue that weakening the definition of a super Hilbert space (by allowing the super scalar product to be non-homogeneous, not just even) will allow the left-regular representation of all (connected) super Lie groups to be super unitary (with an adapted definition). 
Along the way I will introduce a (super) metric on a supermanifold that will allow me to define super and non-super scalar products on function spaces and I will show that the former are intimately related to the Hodge-star operation and the Fermionic Fourier transform. 
The latter also allows me to decompose certain super unitary representations as a direct integral over odd parameters of a family of super unitary representations depending on these odd parameters. 

\end{abstract}

\maketitle

\tableofcontents

\masection{Introduction}
\label{Introductionsection}

The generally accepted definition of a super unitary representation of a super Lie group $G$ is the one that can be found (among others) in \cite[Def. 2, \S2.3]{CCTV:2006} and \cite{AllHilLau:2013}. 
It starts with the interpretation of a super Lie group $G$ as a super Harish-Chandra pair $(G_o,\Liealg g)$, in which $\Liealg g=\Liealg g_0 \oplus \Liealg g_1$ is a super Lie algebra (over $\RR$) and $G_o$ an ordinary Lie group acting on $\Liealg g$ such that:
\begin{enumerate}[{\labelSHCP}1.]
\item
the Lie algebra of $G_o$ is (isomorphic to) $\Liealg g _0$;

\item
the action of $G_o$ preserves each $\Liealg g_\alpha$ (the action is \myquote{even});

\item
the restriction of the $G_o$ action to $\Liealg g_0$ is (isomorphic to) the adjoint action of $G_o$ on it Lie algebra.

\end{enumerate}
The next item to be defined is a super Hilbert space $(\Hilbert,\inprodsym, \supinprsym)$, which is usually taken to be a graded Hilbert space $\Hilbert = \Hilbert_0\oplus \Hilbert_1$ with scalar product $\inprodsym$ and super scalar product $\supinprsym$ (a graded symmetric non-degenerate sesquilinear form) satisfying the following conditions:
\begin{enumerate}[\labelIntroSHS1.]
\item\label{IntroSHS1label}
$\inprod{\Hilbert_0}{\Hilbert_1}=0$;

\item\label{IntroSHS2label}
for all homogeneous $x,y\in \Hilbert$ we have $\superinprod xy = i^{\parity x}\cdot \inprod xy$.

\end{enumerate}
With these ingredients a super unitary representation of $(G_o,\Liealg g)$ on the super Hilbert space $(\Hilbert,\inprodsym, \supinprsym)$ then is a couple $(\rho_o,\tau)$ in which $\rho_o$ is an ordinary unitary representation of $G_o$ on the Hilbert space $\Hilbert$ and $\tau:\Liealg g\to \End\bigl(C^\infty(\rho_o)\bigr)$ an even super Lie algebra representation of $\Liealg g$ on $C^\infty(\rho_o)$, the space of smooth vectors for $\rho_o$ defined by
$$
C^\infty(\rho_o) = \{\,\psi\in \Hilbert \mid g\mapsto \rho(g)\psi \text{ is a smooth map }G\to\Hilbert\,\}
\mapob, 
$$
satisfying the conditions:
\begin{enumerate}[\labelIntroSUR1.]
\item\label{IntroSUR1label}
for each $g\in G_o$ the map $\rho_o(g)$ preserves each $\Hilbert_\alpha$ (the representation is \myquote{even});

\item
for each $X\in \Liealg g_0$ (the Lie algebra of $G_o$!) the map $\tau(X)$ is the restriction of the infinitesimal generator of $\rho_o\bigl(\exp(tX)\bigr)$ to $C^\infty(\rho_o)$;

\item\label{IntroSUR3label}
for each $X\in \Liealg g_\alpha$ the map $\tau(X)$ is graded skew-symmetric with respect to $\supinprsym$;

\item\label{IntroSUR4label}
for all $g\in G_o$ and all $X\in \Liealg g_1$ we have
$$
\tau(g\cdot X) = \rho_o(g)\scirc \tau(X) \scirc \rho_o(g\mo)
\mapob,
$$
where on the left we denote by $g\cdot X$ the action of $G_o$ on $\Liealg g$.

\end{enumerate}

Actually, condition \refmetnaam{\labelIntroSHS}{IntroSHS2label} tells us that the super scalar product $\supinprsym$ can be deduced from the scalar product $\inprodsym$. Moreover, using \refmetnaam{\labelIntroSHS}{IntroSHS2label} one can show that condition \refmetnaam{\labelIntroSUR}{IntroSUR3label} is equivalent to the condition that a suitable multiple of $\tau(X)$ is essentially self-adjoint on $C^\infty(\rho_o)$. 
In other words, one can dispense completely with the notion of the super scalar product $\supinprsym$. 

In the context of non-super Lie groups, any Lie group $G$ admits a canonically defined unitary representation: the left-regular representation on $L^2(G)$, the space of square-integrable functions on $G$ (with respect to an invariant (Haar) measure). 
Unfortunately, this no longer is true for super Lie groups when we adopt the above definition of a super unitary representation. 
One can easily show that there are lots of super Lie groups for which there is no way to turn the natural left-regular representation into a super unitary one according to the above definition. 
The easiest example is the additive (abelian) group of graded dimension $0\vert 1$ usually denoted as $\RR^{0\vert1}$. 

Of course one could accept this fact as a particularity of super Lie groups, but I want to argue that we have another option: changing the definition of a super Hilbert space by allowing the super scalar product to be non-homogeneous (not only even). 
I will show that with this weaker definition (and an adapted version of a super unitary representation) one can turn the left-regular representation of any super Lie group into a super unitary one. 
More precisely, I propose to define a super scalar product to be any (not necessarily even or homogeneous) non-degenerate sesquilinear map $\supinprsym:\Hilbert\times \Hilbert \to \CC$, and then to define a super Hilbert space $(\Hilbert, \inprodsym,\supinprsym)$ by the condition
\begin{enumerate}
\item[]
the  map $\supinprsym:\Hilbert\times \Hilbert \to \CC$ is continuous with respect to the topology induced by the scalar product $\inprodsym$.

\end{enumerate}
I also change the viewpoint of a super unitary representation by defining a super unitary representation of a super Lie group $G$ on a super Hilbert space $(\Hilbert, \inprodsym, \supinprsym)$ to be a couple $(\Dense, \rho)$ with the following properties. 
$\Dense\subset \Hilbert$ is a dense graded subspace and $\rho:G\to \Aut(\Dense\otimes\CA^\KK)$ (notation to be explained in the main text) a group homomorphism satisfying the following conditions.
\begin{enumerate}[\labelIntroSURmod1.]
\item\label{IntroSURmodlabel2}
For all $g\in G$ the map $\rho(g)$ preserves $\supinprsym$.

\item\label{IntroSURmodlabel3}
For all $\psi\in \Dense$ the map $\FgroupHilbert_\psi:G\to \Dense$, $\FgroupHilbert_\psi(g)=\rho(g)\psi$ is smooth.

\item\label{IntroSURmodlabel1}
For all $g\in \body G$ the map $\rho(g)$ preserves $\inprodsym$.

\item\label{IntroSURmodlabel4}
The couple $(\Dense,\rho)$ is maximal with respect to the previous conditions.

\end{enumerate}

Before showing that this definition of a super unitary representation is equivalent to a lookalike of the standard definition given by the conditions \labelIntroSUR{} above, let me first make some comments on this definition. 
The space $\body G$ is the underlying non-super Lie group and condition \refmetnaam{\labelIntroSURmod}{IntroSURmodlabel1} says that it is unitary. 
It thus extends naturally to a unitary representation $\rho_o$ of $\body G$ on $\Hilbert$. 
Looking at condition \refmetnaam{\labelIntroSURmod}{IntroSURmodlabel3} then says in particular that all elements $\psi\in \Dense$ belong to $C^\infty(\rho_o)$, \ie, $\Dense\subset C^\infty(\rho_o)$. 
Moreover, this condition says that $\Dense$ is a dense graded subspace of smooth vectors for the whole representation $\rho$. 
And then condition \refmetnaam{\labelIntroSURmod}{IntroSURmodlabel4} can be interpreted as saying that $\Dense$ is the whole space of smooth vectors for $\rho$. 
On the other hand, as it is natural that the set of smooth vectors for a subgroup is larger than the set of smooth vectors for the whole group, it seems reasonable that we cannot require equality $\Dense=C^\infty(\rho_o)$. 
With these remarks in mind, the following result should not come as a total surprise.

\begin{proclaim}[equivalenceSURintroduction]{Theorem \textit{(a reformulation of \recalt{equivalentDefSuperUnitaryRepNEW})}}
Giving a super unitary representation of a super Lie group $G$ on a super Hilbert space $(\Hilbert, \inprodsym, \supinprsym)$ in the form of a couple $(\Dense,\rho)$ satisfying the conditions \labelIntroSURmod{} above is equivalent, in terms of the super Harish-Chandra pair $(G_o,\Liealg g)\cong G$, to a triple $(\rho_o,\Dense,\tau)$ with the following properties. 
$\rho_o$ is a unitary representation of $G_o\equiv \body G$ on $\Hilbert$, $\Dense\subset C^\infty(\rho_o)$ is a dense graded subspace of $\Hilbert$ and $\tau:\Liealg g\to \End(\Dense)$ an even graded Lie algebra morphism satisfying the conditions:
\begin{enumerate}
\item
for each $g\in G_o$ the map $\rho_o(g)$ preserves each $\Hilbert_\alpha$ (the representation is \myquote{even});

\item
for each $X\in \Liealg g_0$ (the Lie algebra of $G_o$!) the map $\tau(X)$ is the restriction of the infinitesimal generator of $\rho_o\bigl(\exp(tX)\bigr)$ to $\Dense$;

\item
for each $X\in \Liealg g_\alpha$ the map $\tau(X)$ is graded skew-symmetric with respect to $\supinprsym$;

\item
for all $g\in G_o$ and all $X\in \Liealg g_1$ we have
$$
\tau(g\cdot X) = \rho_o(g)\scirc \tau(X) \scirc \rho_o(g\mo)
\mapob;
$$

\item
the couple $(\Dense, \tau)$ is maximal with respect to the above conditions. 

\end{enumerate}
With these ingredients, the maps $\tau(X)$ with $X\in \Liealg g_\alpha$ are the infinitesimal generators of the representation $\rho$.

\end{proclaim}

Comparing this result with the conditions \labelIntroSUR{} shows that they really are a lookalike. 
The one exception is that $\tau$ acts on a dense subspace of $C^\infty(\rho_o)$ and not on the whole of it. 
But the change in the definition of a super Hilbert space has one far more important consequence: the maps $\tau(X)$ with $X\in \Liealg g_1$ are no longer (essentially) self-adjoint (nor a suitable multiple)! 
Now this should not worry the reader, simply because we have no use for it: even if (a multiple) were (essentially) self-adjoint, the $1$-parameter unitary group it generates plays no role, as the $\tau(X)$ with $X\in \Liealg g_1$ are \stress{odd} operators, and thus (as is obvious from the proof of the above equivalence) we are only interested in $\exp\bigl(\xi \tau(X)\bigr) \equiv \oneasmatrix+\xi\tau(X)$ with $\xi\in \CA_1$ an \stress{odd} parameter. 

\bigskip

\noindent\textbf{The main consequence}
\medskip

Let $G$ be a super Lie group in the $\CA$-manifold setting for which we want to show that its left-regular representation is super unitary. 
To do so, we consider the space $C^\infty(G)$ of (super) smooth functions on $G$ (with values in $\CA^\CC$). 
The group $G$ acts on $C^\infty(G)$ by left-translation as
$$
\bigl(\rho(g)\psi\bigr)(h)=\psi(g\mo h)
\mapob.
$$
If $\mfdmetric$ is any left-invariant super metric on $G$, we have a natural left-invariant super density $\nu_\mfdmetric$ (unique up to a scalar multiple) which allows us to integrate elements of $C_c^\infty(G)\subset C^\infty(G)$, the compactly supported elements of $C^\infty(G)$. 
As $\nu_\mfdmetric$ is left-invariant, we obtain a super scalar product $\supinprsym_\mfdmetric$ on $C_c^\infty(G)$ that is preserved by the (left-regular) representation $\rho$. 

But the super metric $\mfdmetric$ also allows us to define a scalar product\slash metric $\inprodsym_\mfdmetric$ on $C_c^\infty(G; \CA^\CC)$. 
This scalar product is obtained via Batchelor's theorem (describing $C^\infty(G)$ as sections of an ordinary (\ie, non-super) vector bundle over the ordinary manifold $G_o=\body G$). 
Moreover, this ordinary scalar product is invariant under $\rho(g)$ for $g\in G_o$. And, even more important, the super scalar product $\supinprsym_\mfdmetric$ defined above is continuous with respect to $\inprodsym_\mfdmetric$. 

We are now in business: we can define $\Hilbert$ as the completion of $C_c^\infty(G)$ with respect to $\inprodsym_\mfdmetric$. The representation $\rho$ of $G_o\equiv \body G$ on $C^\infty_c(G)$ extends in the natural way to a unitary representation $\rho_o$ on $\Hilbert$. 
And because the super scalar product $\supinprsym_\mfdmetric$ is continuous, it also extends to $\Hilbert$ (and happens to remain non-degenerate). 
It then suffices to take the (sic!) maximal (with respect to conditions \refmetnaam{\labelIntroSURmod}{IntroSURmodlabel2}--\ref{IntroSURmodlabel1}) invariant subspace $\Dense_\rho \subset C^\infty(\rho_o)$ to obtain a super unitary representation $(\Dense_\rho,\rho)$ of $G$ on $(\Hilbert, \inprodsym_\mfdmetric, \supinprsym_\mfdmetric)$. 

\medskip

Weakening the definition of a super Hilbert space not only allows us to show that the left-regular representation of any super Lie group is super unitary. 
It also allows us to use direct Berezin-integrals of super unitary representations indexed by odd parameters (over which we take the Berezin integral) and to prove that the result is again a super unitary representation. 
And once we have this notion of a direct Berezin-integral of representations, we can use the Fermionic Fourier transform (the super version of the Fourier transform which, as we will show, is closely related to the Hodge-star operation and which will allow us to provide, in special cases, a close link between the ordinary scalar product $\inprodsym_\mfdmetric$ and the super scalar product $\supinprsym_\mfdmetric$) to decompose super unitary representations as a direct integral over odd parameters, just as one uses the ordinary Fourier transform to decompose unitary representations as a direct integral over even (real) parameters.

\bigskip

\noindent\textbf{An outline}
\medskip

Let me finish this introduction with an outline of this paper. 
I will be working with $\CA$-manifold theory (my version of the geometric $H^\infty$ version of DeWitt supermanifolds), which is equivalent to the theory of graded manifolds of Leites and Kostant (see \cite{DW84}, \cite{Ko77}, \cite{Le80}, \cite{Ro07} \cite{Tu04}, \cite{Va04}).
Any reader using a (slightly) different version of supermanifolds should be able to translate the results to her\slash his version of supermanifolds. 
Readers unfamiliar with my approach will find in the appendix \recals{appendixonAmanifoldssection} some notational conventions as well as an extremely succinct overview of $\CA$-manifold theory. 
We start, in \recals{motivatingexamplesection}, with a motivating example in which we perform heuristic computations that are a lookalike to computations one usually makes for the Heisenberg group: determination of the left-regular representation and decomposition as a direct integral of Fourier modes. 
In order to turn these heuristic computations into valid ones, we start in \recals{superscalarproductssection} with the definition of super scalar products and super metrics. 
With these ingredients we are able to state, in \recals{superunitarydefandequivalencesection}, our weakened definition of a super Hilbert space and to prove our equivalence theorem \recalt{equivalentDefSuperUnitaryRepNEW}. 

In \recals{BatchelorbundlewithHodgestarsection} we recall the definition of the Batchelor bundle associated to an $\CA$-manifold $M$, we introduce the notion of a super metric on an $\CA$-manifold, and we show how to deduce an ordinary metric on $C^\infty_c(M)$. 
This will allow us to obtain a Hilbert space as the completion of this space of compactly supported smooth functions; it will be the analogue of the space of square integrable functions on an ordinary non-super manifold with respect to a metric volume form. 
We also will show that a metric preserving diffeomorphism of $M$ satisfying a \myquote{linearity} condition induces a unitary map of $C^\infty_c(M)$. 
In \recals{Berezinintegrationsection} we will recall some basic facts about Berezin integration on $\CA$-manifolds and we will show that any super metric on an $\CA$-manifold $M$ defines a natural super scalar product on the space $C_c^\infty(M)$ of compactly supported smooth functions on $M$ (just as an ordinary metric does for an ordinary manifold).
In \recals{Integrationalongfiberssection} we prove a \myquote{change of variables} formula for Berezin integration of functions depending upon odd parameters (more precisely, depending upon a point in an $\CA$-manifold). 
This will allow us to justify that the natural super scalar product on $C^\infty_c(M)$ associated to an invariant super metric is invariant under left-translation by the group. 
The next section then treats the Berezin-Fourier transform. 
In it we show that a generalization of the Berezin-Fourier transform to any $\CA$-manifold equipped with a super metric is closely related to the Hodge-star operation and to a natural super scalar product on $C^\infty_c(M)$. 
In \recals{leftregularrepresentationssection} we apply our results to the left-regular representation of a super Lie group and prove our main result \recalt{maintheoremLeftRegularRep}, viz, that the natural (left-regular) representation of $G$ on the space of compactly supported smooth functions (on $G$) admits a unique maximal extension to a super unitary representation $(\Dense_\rho, \rho)$. 
Moreover, the definition of $\Dense_\rho$ is a lookalike of the description of the space of smooth vectors of the left-regular representation of an ordinary Lie group, reinforcing our interpretation of the dense subspace in the definition of a super unitary representation as the set of smooth vectors of this representation. 

Before we return our attention to the motivating example of \recals{motivatingexamplesection}, we first define in \recals{directintegralsofrepssection} the notion of a direct (Berezin) integral of super unitary representations. 
The last three sections then treat three examples: in \recals{backtomotivtingexamplesection} we come back to our motivating example and provide the justifications for all heuristic computations given in \recals{motivatingexamplesection}. 
In \recals{axiplusbetagroupsection} we treat the $a\xi+\beta$ group, which is interesting because it is an example in which we indeed have $\Dense \neq C^\infty(\rho_o)$. 
It also provides us with a second example (after the one of our motivating example) of a decomposition by a Berezin-Fourier transform with an odd parameter (more examples of such decompositions (albeit in a more primitive form) can be found in \cite{Tuynman:2009}). 
And in \recals{OSp12section} we treat the super Lie group $\mathrm{OSp}(1,2)$, which is interesting because its super scalar product is less \myquote{trivial} than in the previous two examples.

\masection{A motivating example}
\label{motivatingexamplesection}

My motivating example will be the analogue of (a part of) the standard decomposition of the left-regular representation of the classical Heisenberg group $\RR^3$ with group multiplication 
$$
(x,y,z)\cdot (x',y',z') = \bigl(x+x'+\tfrac12 (xy'-x'y), y+y', z+z'\bigr)
\mapob.
$$
The computations in this section will be heuristic but \dots{} everybody will recognize that, when done in the ordinary non-super context, only a few cautionary phrases (in particular concerning domains of definition) will suffice to make them rigorous. 
On the other hand, the knowledgeable reader will certainly see lots of problems in the super context, some of which will be pointed out at the end. 
I will not always be very precise and some notions might not be defined here. 
However, all details can be found in the main text.

\mysubsectionnonumber{motivatingexamplesection}{The left-regular representation}

Let $E$ be a graded vector space of dimension $1\vert n$ (with $n\in \NN$) and denote by $(x,\xi_1, \dots, \xi_n)$ global coordinates on $E_0$ (the even part of $E$, with $x$ being even and the $\xi_i$ being odd coordinates). 
We define the group $G$ as being $G=E_0$ equipped with the group multiplication\footnote{In \cite{AllHilLau:2013} this is called a Clifford-Heisenberg group; it is a special case of what I called, in \cite{Tuynman:2010}, a \myquote{Heisenberg-like group.}}
$$
(x,\xi)\cdot (x',\xi') = (\,x+x'+\tfrac12\,\inprod\xi{\xi'}, \xi+\xi' \,)
\mapob,
$$
where $\inprodsym$ denotes the standard scalar product:
$$
\inprod\xi{\xi'} = \sum_{i=1}^n \xi_i \cdot \xi'_i
\mapob.
$$
Any smooth function $f$ on $G$ is of the form
$$
f(x,\xi) = \sum_{I\subset \{1, \dots, n\}} f_I(x)\cdot \xi^I
\mapob,
$$
where the $f_I$ are smooth functions of a single real\slash even variable, and where $\xi^I$ is defined as
\begin{moneq}[ConventionxitothepowersetI]
\xi^\emptyset = 1
\qquad\text{and}\qquad
I=\{i_1< \dots< i_k\}
\quad\Rightarrow\quad
\xi^I = \xi_{i_1}\cdots \xi_{i_k}
\mapob.
\end{moneq}
It thus seems reasonable that the Hilbert space on which the left-regular representation is defined is given as
$$
\Hilbert = \Bigl\{\ \psi(x,\xi) = \sum_{I\subset \{1, \dots, n\}} \psi_I(x)\cdot \xi^I\ \Bigm\vert\ \psi_I\in L^2(\RR, \Leb)\ \Bigr\}
\cong
\bigl(  L^2(\RR, \Leb) \bigr)^{2^n}
\mapob,
$$
equipped with the natural super scalar product $\supinprsym$ defined by
$$
\superinprod\chi\psi = \int_\RR\extder\Leb(x)\ \int \extder\xi_1\cdots \extder\xi_n\  \overline{\chi(x,\xi)}\cdot \psi(x,\xi)
\mapob.
$$
The left-regular representation $\rho$ itself is defined in the usual way by
$$
\bigl(\rho{(y,\eta)}\psi\bigr)(x,\xi)
=
\psi\bigl((y,\eta)\mo\cdot (x,\xi)\bigr)
=
\psi(x-y-\tfrac12\,\inprod\eta\xi,\xi-\eta)
\mapob.
$$
As the Lebesgue measure and Berezin integration are translation invariant, it is immediate that $\rho(y,\eta)$ preserves the super scalar product:
\begin{align*}
\superinprod[2]{\rho(y,\eta)\chi}{\rho(y,\eta)\psi}
&
=
\int_\RR \extder\Leb(x)\ \int \extder\xi_1\cdots \extder\xi_n\  \overline{\bigl(\rho(y,\eta)\chi\bigr)(x,\xi)}
\cdot \bigl(\rho(y,\eta)\psi\bigr)(x,\xi)
\\&
=
\int_\RR \extder\Leb(x)\ \int \extder\xi_1\cdots \extder\xi_n\  \overline{\chi(x-y-\tfrac12\,\inprod\eta\xi,\xi-\eta)}
\\&
\kern14em
\cdot \psi(x-y-\tfrac12\,\inprod\eta\xi,\xi-\eta)
\\
\text{\tiny ($z=x-y-\tfrac12\,\inprod\eta\xi\,$)}\quad
&
=
\int_\RR \extder\Leb(z)\ \int \extder\xi_1\cdots \extder\xi_n\  \overline{\chi(z,\xi-\eta)}
\cdot \psi(z,\xi-\eta)
\\
\text{\tiny ($\zeta = \xi-\eta\,$)}\quad
&
=
\int_\RR \extder\Leb(z)\ \int \extder\zeta_1\cdots \extder\zeta_n\  \overline{\chi(z,\zeta)}
\cdot \psi(z,\zeta)
\\&
=
\superinprod\chi\psi
\mapob.
\end{align*}
It thus seems reasonable to say that $\rho$ is a super unitary representation. 

\mysubsectionnonumber{motivatingexamplesection}{Fourier decomposition via the central coordinate}

But we can go (at least) one step further. 
We start with the observation that functions of the form $\psi(x,\xi) = \psi_k(\xi)\,\eexp^{ikx}$ with $k\in\RR$ are invariant under the action of $\rho{(y,\eta)}$:
\begin{align*}
\bigl(\rho{(y,\eta)}\psi\bigr)(x,\xi)
&
=
\psi(x-u - \tfrac12\inprod\eta\xi, \xi-\eta)
=
\psi_k(\xi-\eta)\,\eexp^{ik(x-y-\frac12\inprod\eta\xi)}
\\&
=
\psi_k(\xi-\eta)\,\eexp^{-ik(y+\frac12\inprod\eta\xi)}\,\eexp^{ikx}
\mapob.
\end{align*}
This suggests the family of representation $\rhoh_k$ of $G$ on the space of (complex valued) functions depending on the odd coordinates $\xi$ only, a space we denote as $C^\infty(\CA_1^n; \CA^\CC)$, given by
$$
\bigl(\rhoh_k(y,\eta)\chi\bigr)(\xi) = \chi(\xi-\eta)\,\eexp^{-ik(y+\frac12\inprod\eta\xi)}
\mapob.
$$
It then is immediate that a partial Fourier transform $\Fourier$ with respect to the (real) variable $x$ intertwines the representation $\rho$ with this family $\rhoh_k$, $k\in\RR$:  
\begin{align*}
\bigl(\Fourier(\rho{(y,\eta)}\psi)\bigr)(k,\xi)
&
=
\frac1{\sqrt{2\pi}} \int_\RR \eexp^{-ikx}\,\bigl(\rho{(y,\eta)}\psi\bigr)(x,\xi)\ \extder x
\\&
=
\frac1{\sqrt{2\pi}} \int_\RR \eexp^{-ikx}\,\psi(x-y-\tfrac12\inprod\eta\xi,\xi-\eta)\ \extder x
\\
\text{\tiny ($z=x-y-\tfrac12\,\inprod\eta\xi\,$)}\quad
&
=
\frac1{\sqrt{2\pi}} \int_\RR \eexp^{-ikz}\,\eexp^{-ik(y+\frac12\inprod\eta\xi)}\, \psi(z,\xi-\eta)\ \extder z
\\&
=
\eexp^{-ik(y+\frac12\inprod\eta\xi)}\cdot (\Fourier\psi)(k,\xi-\eta)
=
\bigl(\rhoh_k(y,\eta)(\Fourier\psi)\bigr)(k,\xi)
\mapob.
\end{align*}
The super unitary representation $\rho$ thus decomposes as a direct integral of the family of representations $\rhoh_k$. 
Moreover, these representations preserve the natural super scalar product on $C^\infty(\CA_1^n; \CA^\CC)$ given by
$$
\superinprod\chi\psi = \int \extder\xi_1\cdots \extder \xi_n \ \overline{\chi(\xi)}\cdot\psi(\xi)
\mapob,
$$
simply because Berezin integration is translation invariant:
\begin{align*}
\superinprod{\rhoh_k(y,\eta)\chi}{\rhoh_k(y,\eta)\psi}
&
=
\int \extder\xi_1\cdots \extder\xi_n\  \overline{\bigl(\rhoh_k(y,\eta)\chi\bigr)(\xi)}
\cdot \bigl(\rhoh_k(y,\eta)\psi\bigr)(\xi)
\\&
=
\int \extder\xi_1\cdots \extder\xi_n\  \overline{\chi(\xi-\eta)\,\eexp^{-ik(y+\frac12\inprod\eta\xi)}}
\\&
\kern12em
\cdot \psi(\xi-\eta)\,\eexp^{-ik(y+\frac12\inprod\eta\xi)}
\\&
=
\int \extder\xi_1\cdots \extder\xi_n\  \overline{\chi(\xi-\eta)}
\cdot \psi(\xi-\eta)
\\
\text{\tiny $(\zeta=\xi-\eta\,$)}\quad
&
=
\int \extder\zeta_1\cdots \extder\zeta_n\  \overline{\chi(\zeta)}
\cdot \psi(\zeta)
=
\superinprod\chi\psi
\mapob.
\end{align*} 
Again it thus is reasonable to say that the representations $\rhoh_k$ are super unitary.

\mysubsectionnonumber{motivatingexamplesection}{Decomposing the $0$-Fourier mode}

We now first concentrate on the super unitary representation $\rhoh_{k=0}$, which is given by
$$
\bigl(\rhoh_{k=0}(y,\eta)\psi\bigr)(\xi) = \psi(\xi-\eta)
\mapob.
$$
Analogously to the case of the representation $\rho$, we first note that functions of the form $\psi(\xi) = c\cdot\eexp^{i\,\inprod\xi\kappa} = c\cdot\eexp^{-i\,\inprod\kappa\xi}$ with $\kappa=(\kappa_1, \dots, \kappa_n)$ $n$ odd parameters (we write $\kappa\in \CA_1^n$) are invariant under the action of $\rhoh_0(y,\eta)$:
$$
\bigl(\rhoh_{k=0}(y,\eta)\psi\bigr)(\xi) = \psi(\xi-\eta)
=
c\cdot \eexp^{i\,\inprod{\xi-\eta}\kappa}
=
c\cdot \eexp^{-i\,\inprod\eta\kappa}\cdot \eexp^{i\,\inprod\xi\kappa}
\mapob.
$$
This suggests the family of representations $\rhob_\kappa$, $\kappa\in \CA_1^n$ of $G$ on the $1$-dimensional space $\CA^\CC$ of constant (functions) $c$ given by
$$
\rhob_\kappa(y,\eta)c = \eexp^{-i\,\inprod\eta\kappa} \cdot c \equiv \eexp^{i\,\inprod\kappa\eta} \cdot c
\mapob.
$$
We now \myquote{recall} that the standard Berezin-Fourier transform $\Fourierodd$ on functions in $C^\infty(\CA_1^n; \CA^\CC)$ is defined as 
$$
(\Fourierodd\psi)(\kappa) = \int \extder \xi_1\cdots \extder \xi_n \ 
\eexp^{-i\,\inprod\xi\kappa}\cdot \psi(\xi)
\mapob.
$$
It then is immediate that this Berezin-Fourier transform intertwines the representation $\rhoh_{k=0}$ with the family $\rhob_\kappa$:
\begin{align*}
\bigl(\Fourierodd\rhoh_{k=0}(y,\eta)\psi\bigr)(\kappa)
&
=
\int\extder\xi_1 \cdots \extder \xi_n\  \eexp^{-i\inprod\xi\kappa}\,\bigl(\rhoh_{k=0}(y,\eta)\psi\bigr)(\xi) 
\\&
=
\int\extder\xi_1 \cdots \extder \xi_n\  \eexp^{-i\inprod\xi\kappa}\,\psi(\xi-\eta) 
\\
\text{\tiny $(\zeta=\xi-\eta\,$)}\quad
&
=
\int \extder\zeta_1 \cdots \extder \zeta_n\  \eexp^{i\inprod\kappa{\eta}}\,\eexp^{-i\inprod{\zeta}\kappa}\,\psi(\zeta) 
=
\eexp^{i\inprod\kappa{\eta}}\, (\Fourierodd\psi)(\kappa)
\\&
=
\bigl( \rhob_\kappa(y,\eta)(\Fourierodd\psi)\bigr)(\kappa)
\mapob.
\end{align*}
The super unitary representation $\rhoh_{k=0}$ thus decomposes as a direct Berezin integral  of the family  of representations $\rhob_\kappa$. 
Moreover, it is immediate that these representations preserve the natural (super) scalar product $\superinprod c{c'} = \overline{c}\cdot c'$ on these $1$-dimensional spaces, so these representations are super unitary.

\mysubsectionnonumber{motivatingexamplesection}{Decomposing the $\mathbf{k}$-Fourier mode for $\mathbf{k\neq0}$}

In order to keep this motivating example sufficiently simple, we will not investigate the general representation $\rhoh_k$ with $k\neq0$, but we will concentrate on the special cases $n=1$ and $n=2$. 
For $n=1$, any (smooth) function of a single odd variable is determined by two (complex) constants according to
$$
\psi(\xi) = \psi_0 + \xi\cdot \psi_1
\mapob.
$$
In terms of these constants the super scalar product is given by
$$
\superinprod\chi\psi = \overline{\psi_0}\cdot \chi_1 + \overline{\psi_1}\cdot \chi_0
$$
and the representation $\rhoh_k$ is given by
\begin{align*}
\bigl(\rhoh_k(y,\eta)\psi\bigr)(\xi) 
&
= 
\psi(\xi-\eta)\,\eexp^{-ik(y+\frac12\,\eta\,\xi)}
=
\bigl(\psi_0+(\xi-\eta)\psi_1\bigr)\cdot \eexp^{-iky}\cdot \bigl(1-\tfrac12\,ik\eta\xi)\bigr)
\\&
=
\eexp^{-iky}\cdot\bigl( (\psi_0 - \eta\psi_1) + \xi\,( \psi_1+ \tfrac12\,ik\eta\psi_0  )\bigr)
\end{align*}
or in terms of the couple $(\psi_0,\psi_1)$:
$$
\rhoh_k(y,\eta)\cdot
\begin{pmatrix} \psi_0\\ \psi_1 \end{pmatrix}
=
\eexp^{-iky}\cdot
\begin{pmatrix} 1 & -\eta \\ \tfrac12\,ik\eta & 1 \end{pmatrix}
\cdot
\begin{pmatrix} \psi_0\\ \psi_1 \end{pmatrix}
\mapob.
$$
For $k\neq0$ it is easy to show that there do not exist invariant graded subspaces, so the representations $\rhoh_k$ are irreducible for $k\neq0$ and $n=1$.

We next turn our attention to the case $n=2$, where a function of two odd coordinates is determined by four complex constants $\psi_0, \psi_1, \psi_2, \psi_{12}$ by
$$
\psi(\xi_1,\xi_2) = \psi_0+\xi_1\,\psi_1 + \xi_2\,\psi_2+\xi_1\xi_2\,\psi_{12}
\mapob.
$$
In terms of these four components, the representation $\rhoh_k$ is given by
\begin{align*}
(\rhoh_k(y,\eta)\psi)(\xi_1,\xi_2)
&
=
\eexp^{-ik(y+\frac12\eta_1\xi_1+\frac12\eta_2\xi_2)}\cdot\bigl(\psi_0+(\xi_1-\eta_1)\,\psi_1 
\\&
\kern4em
+ (\xi_2-\eta_2)\,\psi_2+(\xi_1-\eta_1)(\xi_2-\eta_2)\,\psi_{12}\bigr)
\mapob.
\end{align*}
The (infinitesimal) operators $\tauh_k(f_j)$ (essentially the derivative of $\rhoh_k(y,\eta)$ with respect to $\eta_j$ at $(y=0,\eta=0)\,$)  are given by
\begin{align*}
\tauh_k(f_1)\psi(\xi_1,\xi_2)
&
=
-\tfrac12\,ik\,\xi_1\,(\psi_0 + \xi_2\psi_2) -(\psi_1 + \xi_2\psi_{12})
\\
\tauh_k(f_2)\psi(\xi_1,\xi_2)
&
=
-\tfrac12\,ik\,\xi_2\,(\psi_0 + \xi_1\psi_1) -(\psi_2 - \xi_1\psi_{12})
\mapob,
\end{align*}
or in matrix form:
$$
\tauh_k(f_1)
=
\begin{pmatrix}
0 & -1 & 0 & 0
\\
-\tfrac12\,ik & 0 & 0 & 0 
\\
0 & 0 & 0 & -1
\\
0 & 0 & -\tfrac12\,ik & 0
\end{pmatrix}
\qquad,\qquad
\tauh_k(f_2)
=
\begin{pmatrix}
0 & 0 & -1 & 0
\\
0 & 0 & 0 & 1 
\\
-\tfrac12\,ik & 0 & 0 & 0
\\
0 & \tfrac12\,ik & 0 & 0
\end{pmatrix}
\mapob.
$$

A direct computation shows that the representation space $C^\infty(\CA_1^2; \CA^\CC)$ for $\rhoh_k$ splits as the direct sum of two invariant graded subspaces $C^\infty(\CA_1^2; \CA^\CC) = \Hilberth_{-1} \oplus \Hilberth_{1}$ where $\Hilberth_\varepsilon$, $\varepsilon=\pm1$ (of dimension $1\vert1$) is generated by the (even and odd) functions $\chi_{\varepsilon,0}$ and $\chi_{\varepsilon,1}$ given by
$$
\chi_{\varepsilon,0}(\xi_1,\xi_2) = \eexp^{-\frac12\varepsilon k\xi_1\xi_2}
\qquad\text{and}\qquad
\chi_{\varepsilon,1}(\xi_1,\xi_2) = \eexp^{-\frac12\varepsilon k\xi_1\xi_2}\cdot (\xi_1 + i\varepsilon \,\xi_2)
$$
or equivalently:
$$
\psi_\varepsilon\in \Hilberth_\varepsilon
\qquad\Longleftrightarrow\qquad
\psi_\varepsilon(\xi_1,\xi_2) = \eexp^{-\frac12\varepsilon k\xi_1\xi_2}\cdot\bigl(\lambda+(\xi_1 + i\varepsilon \,\xi_2)\,\mu\bigr)\ ,\ \lambda,\mu\in \CC
\mapob.
$$
On these two subspaces the representation is given by
\begin{multline*}
\qquad
\bigl(\rhoh_k(y,\eta)\psi_\varepsilon\bigr)(\xi_1,\xi_2)
=
\eexp^{-iky-\frac12\varepsilon k\xi_1\xi_2}
\cdot\eexp^{-\frac12\varepsilon k\eta_1\eta_2}
\cdot\eexp^{-\frac12ik (\eta_1 -i\varepsilon\eta_2)(\xi_1+i\varepsilon\xi_2)}
\\
\cdot\bigl(\lambda-(\eta_1 + i\varepsilon \,\eta_2)\,\mu+(\xi_1 + i\varepsilon \,\xi_2)\,\mu\bigr)
\mapob,
\qquad
\end{multline*}
or in terms of the basis $\chi_0,\chi_1$:
\begin{align*}
\rhoh_k(y,\eta)\chi_{\varepsilon,0}
&
=
\eexp^{-iky-\frac12\varepsilon k\eta_1\eta_2}
\cdot\bigl(\chi_{\varepsilon,0} - \tfrac12 ik(\eta_1 - i\varepsilon\,\eta_2)\chi_{\varepsilon,1}\bigr)
\\
\rhoh_k(y,\eta)\chi_{\varepsilon,1}
&
=
\eexp^{-iky-\frac12\varepsilon k\eta_1\eta_2}
\cdot\bigl(-(\eta_1+i\varepsilon\,\eta_2)\chi_{\varepsilon,0} + \eexp^{k\varepsilon\eta_1\eta_2}\chi_{\varepsilon,1}\bigr)
\mapob.
\end{align*}
A slightly longer computation shows that these two subspaces $\Hilberth_\varepsilon$ are the only graded subspaces that are invariant under the representation $\rhoh_k$. 
We thus have decomposed the representation $\rhoh_k$ for $k\neq0$ and $n=2$ as a direct sum of two irreducible representations.

\mysubsectionnonumber{motivatingexamplesection}{Does it work?}

Now even if the above computations look reasonable, just a slightly better look at them reveals quite a lot of problems. 
The first question one should ask is how we justify our choice for the Hilbert space $\Hilbert$ as consisting of $L^2$-functions? We did not obtain it as a completion of a vector space equipped with a scalar product. More precisely, we only introduced a super scalar product, which is not suitable for the business of completion. 
But even if we accept that $\Hilbert$ should be the right Hilbert space, the computation that intertwines $\rho$ with the family $\rhoh_k$ is fraught with problems. 
In the first place, the notion of nilpotent translation of the Lebesgue measure is not defined, \ie, the term $\inprod\eta\xi$ in the change of coordinates $z=x-y-\inprod\eta\xi/2$. 
And even if it were, such a translation is not well defined on $L^2$-functions: the extension of a function of an even variable to nilpotent arguments uses derivatives. 
More precisely, if $n$ is any (even) nilpotent element, we have
$$
f(x+n) = \sum_{j=0}^\infty \frac{n^j}{j!}\,f^{(j)}(x)
=
f(x) + n\,f'(x) + \frac{n^2}{2!}\,f''(x) + \cdots
\mapob,
$$
which is actually a finite sum because of the nilpotency of $n$. 
The problem of course is that $L^2$-functions do not have derivatives in general. 

But let us assume that this problem has been solved, then the next problem is that the family of representations $\rhoh_k$ is not well defined on the $L^2$-direct integral space $L^2\bigl(\RR; C^\infty(\CA_1^n; \CA^\CC)\bigr)$ (the image of our Hilbert space $\Hilbert$ under the Fourier transform). 
The problem is that it contains the factor $\eexp^{-ik\inprod\eta\xi/2}$, which is multiplication by $\prod_j(1-\tfrac12ik\eta_j\xi_j)$, and multiplication by $k$ is not globally defined on $L^2(\RR; \CC)$. 
Now this is of course the other side of the same coin that says that derivation is not globally defined on $L^2$-functions.
So one could hope that the solution of the previous problem solves this problem at the same time.

\masection{Super scalar products and related stuff}
\label{superscalarproductssection}

\begin{definition}{Conventions and definitions}
$\bullet$
In this paper we will deal with the two fields $\RR$ and $\CC$ and two associated graded rings $\CA$ (over $\RR$) and its complexification $\CA^\CC \equiv \CA\otimes_\RR \CC$ (over $\CC$). Some results will be the same for all four objects, some only for two of them. 
In order not to have to repeat any statement twice (or sometimes even four times), we will introduce the symbol $\KK$ for either $\RR$ or $\CC$, $\CA^\KK$ for either $\CA\equiv\CA^\RR$ or $\CA^\CC$ and $\KKA$ for one of the four objects $\RR$, $\CC$, $\CA$ or $\CA^\CC$. 

$\bullet$
More in particular, an $\CA^\KK$-vector space denotes a graded bi-module $E$ over $\CA^\KK$ of the form $E=V\otimes_\KK \CA^\KK$ with $V$ a graded vector space over $\KK$ (see also \recals{appendixonAmanifoldssection}). 
The body map $\body : E\to V$ is defined by $\body(v\otimes \lambda) = v\cdot \body\lambda$. 
We will in particular identify $V$ with $\body E$, thus writing $E=\body E\otimes_\KK \CA^\KK$. In the same spirit we will identify $v\otimes \lambda$ with $v\cdot \lambda$, which is the same as saying that we identify $v\in \body E$ with $v\otimes 1\in E$. 

$\bullet$
For a $\KKA$-vector space $E$ we denote by $\End(E)$ the space of (right-) linear endomorphisms:
$$
\End(E) = \{\,f:E\to E \mid f \text{ is right-linear over } \KKA\,\}
$$
and by $\Aut(E)$ the space of automorphisms of $E$, \ie the even invertible elements of $\End(E)$:
$$
\Aut(E) = \{\,f\in \End(E) \mid f \text{ even and bijective} \,\}
\mapob.
$$

$\bullet$ 
In the non-super context there is a difference in terminology when dealing with real or complex vector spaces concerning maps that preserve a metric: in the real case one talks about orthogonal maps, and in the complex case of unitary ones. 
As most of our results apply to both cases, it would be rather cumbersome to distinguish these two cases, the more so when we use $\KK$ to denote either $\RR$ of $\CC$. 
We thus will stick exclusively to the name \myquote{super unitary}, even when we are considering graded vector spaces over $\CA^\RR$ (where one would expect the name super orthogonal). 
A similar remark applies to the terminology symmetric versus hermitian: the first is used in the real context, the second in the complex case. Being totally inconsistent, we stick to the name \myquote{symmetric} in both cases (we will mainly be concerned with graded skew-symmetric maps, which thus will include the ones that appear in the complex setting and should be called graded skew-hermitian). 

$\bullet$
Let $x$ be any graded object, which we can decompose as $x=x_0+x_1$ with $x_\alpha$ homogeneous of degree $\alpha$. 
We then define the conjugation operator $\conjugate$ by
$$
\conjugate(x) \equiv \conjugate(x_0+x_1) = x_0-x_1
\mapob.
$$
For any $n\in \ZZ$ or $\Ztwo$ we define (as usual) the operator $\conjugate^n$ to be either the identity (when $n$ is even) or $\conjugate$ (when $n$ is odd). 

$\bullet$
For any graded homogeneous object $x$ we will denote by $\parity x$ its parity, \ie, $\parity x$ is either $0$ (when $x$ is even) or $1$ (when $x$ is odd). 
However, we will use the same \myquote{operation} on finite sets $I$ to denote their cardinal: $\parity{\{a,b,c\}}=3$. But we will use this nearly exclusively in combination with other parities, in which case we interpret $\parity I$ as the cardinal of $I$ modulo $2$. 

\end{definition}

We start with a reminder of the classical definition of scalar products and metrics, not because we think the reader needs them, but in order to show the analogy with our definitions in the super context. 
Now already in ordinary differential geometry we have two kinds of metrics: normal, euclidean or Riemannian metrics on the one hand, and Lorentzian or pseudo-Riemannian metrics on the other, where the difference lies in the fact that the former is positive definite, while the latter is not. 
In this paper we will adopt the following naming convention. 
By a (hermitian) scalar product we will mean a non-degenerate symmetric sesquilinear form on an ordinary (ungraded) vector space over $\KK$ (a Lorentzian metric), and by a metric we will mean an ordinary (euclidean) positive definite hermitian scalar product, again on an ungraded vector space over $\KK$. 
When we add the qualification \myquote{super}, we get a super scalar product and a super metric, both of which will be defined on graded vector spaces of any kind. The super scalar product has \myquote{only} the natural non-degeneracy condition, whereas a super metric has an additional condition reflecting the positive definiteness condition for euclidean metrics.

\begin{definition}[defsofmetricsinnongradedcase]{Classical Definitions}
Let $V$ be a vector space over $\CC$. 

\noindent$\bullet$
A \stresd{hermitian scalar product on $V$} is a non-degenerate symmetric sesquilinear form $\inprod\ \ :V\times V\to \CC$. 
More precisely, the map $\inprod\ \ :V\times V \to \CC$ should satisfy the three conditions.
\begin{enumerate}
\item
$\forall x,y,z\in V$ $\forall r\in \CC$: $\inprod{x}{y+rz} = \inprod xy + r\,\inprod xz$ (linearity in the second variable).

\item
$\forall x,y,z\in V$ $\forall r\in \CC$: $\inprod{x+ry}z = \inprod xz + \overline{r}\ \inprod yz$ (anti-linearity in the first variable).

\item
$\forall x,y\in V$: $\inprod xy = \overline{\inprod yx}$ (symmetry).

\item
$\forall x\in V: (\,\forall y\in V:\inprod xy=0\,)\ \Rightarrow\  x=0$ (non-degeneracy).

\end{enumerate}

\noindent$\bullet$
A \stresd{metric on $V$} is a positive definite hermitian scalar product $\inprod\ \ $ on $V$. 
More precisely, the hermitian scalar product should satisfy the additional condition
\begin{enumerate}
\setcounter{enumi}{4}
\item
$\forall x\in V : \inprod xx\ge0$.

\end{enumerate}

\noindent$\bullet$
If $V$ is a vector space over $\RR$, we can define a \stresd{scalar product} as a non-degenerate hermitian symmetric sesquilinear map $\inprod\ \ :V\times V\to \CC$ (\ie, satisfying exactly the conditions (i)-(iv) above, restricting the scalar $r$ to $r\in \RR$) satisfying the additional (reality) condition
\begin{enumerate}
\setcounter{enumi}{5}
\item
$\forall x,y\in V$: $\inprod xy = \overline{\inprod xy}$.
\label{realitycondonungradedrealvectorspaces}

\end{enumerate}
This amounts to requiring that $\inprod\ \ $ takes its values in $\RR$. 
And as for a vector space over $\CC$, a \stresd{metric} on a vector space over $\RR$ is a positive definite scalar product. 

\noindent$\bullet$
A \stresd{Hilbert space} is a vector space $V$ (over $\RR$ or $\CC$) equipped with a metric $\inprod\ \ $ such that $V$ is complete with respect to the norm $\Vert x\Vert = \sqrt{\inprod xx}$. 

\end{definition}

\begin{definition}[defofgradedhermitianforms]{Definition}
Let $V$ be a graded vector space over $\KKA$.

\medskip

\noindent
$\bullet$
In the case $\KKA=\KK$, a \stresd{(graded, right-) sesquilinear form} on $V$ is a map $\supinprsym : V\times V \to \CC$ satisfying the conditions (i) and (ii) below; in the case $\KKA=\CA^\KK$, it is a map $\supinprsym : V\times V \to \CA^\CC$ satisfying these conditions.
\begin{enumerate}
\item
\stresd{Right-linearity} in the second variable: for all $v,w_1,w_2\in V$ and all $\lambda\in \KKA$ we have
$$
\superinprod v{w_1+w_2 \lambda} = \superinprod v{w_1} + \superinprod v{w_2}\cdot \lambda
\mapob.
$$

\item\label{defofantirightlinearity}
\stresd{Anti-right-linearity} in the first variable (see also \recalt{leftlinearityofrightbilinear} and \recalt{remarkparityhermformonleft}): for all $v_1,v_2, w\in V$ and all $\lambda\in \KKA$ we have
$$
\superinprod{v_1+v_2\lambda}w = \superinprod{v_1}w + \superinprod{v_2}{\overline\lambda w}
\mapob.
$$

\end{enumerate}

\noindent
$\bullet$ If $\supinprsym$ is a graded sesquilinear form on a graded vector space $V$ over $\CA^\KK$, then the map $\body \supinprsym:\body V \times \body V\to \CC$ defined by
$$
(\body \supinprsym)(v,w) = \body\bigl( \supinprsym(v,w)\bigr)
$$
is a graded sesquilinear form on the graded vector space $\body V$ over $\KK$.

\medskip

\noindent
$\bullet$ 
In the case $\KKA=\CA^\KK$, a graded sesquilinear form $\supinprsym$ on $V$ is said to be \stresd{smooth} if it satisfies the condition
\begin{enumerate}
\setcounter{enumi}{2}

\item
\stresd{Smoothness}: 
$\superinprod{\body V}{\body V}\subset \CC$.

\end{enumerate}
Another way to state this condition is that $\body \supinprsym$ equals the restriction of $\supinprsym$ to $\body V$.

\medskip

\noindent
$\bullet$ 
A graded sesquilinear form on $V$ is said to be \stresd{graded symmetric} if it satisfies the condition
\begin{enumerate}
\setcounter{enumi}{3}
\item
\stresd{Graded symmetry}: for all homogeneous $w,v\in V$ we have
$$
\superinprod w v = (-1)^{\parity w \parity v}\,\overline{\superinprod vw}
\mapob.
$$

\end{enumerate}

\noindent
$\bullet$ 
A graded symmetric graded sesquilinear form on $V$ is said to be \stresd{non-degenerate} if it satisfies the condition
\begin{enumerate}
\setcounter{enumi}{4}
\item
\stresd{Non-degeneracy}: 
$$
\forall v\in \body V
\quad:\quad
v\neq0
\quad\Rightarrow\quad
\exists w\in \body V : (\body\supinprsym)(v,w)\neq0
\mapob,
$$
where, in the case of $\KKA=\KK$ the use of the body map should be ignored.

\end{enumerate}

$\bullet$
A \stresd{super scalar product} on $V$ is a graded symmetric non-degenerate graded sesquilinear form on $V$.

\medskip

$\bullet$
A graded sesquilinear form $\supinprsym$ is said to be \stresd{homogeneous of degree $\alpha\in \Ztwo$} if for all $\beta,\gamma\in \Ztwo$ we have the inclusion
$$
\superinprod{V_\beta}{V_\gamma} \subset \CA^\CC_{\alpha+\beta+\gamma}
\mapob,
$$
where in the case $\KKA=\KK$ we use the inclusion $\CC\subset \CA^\CC$ to interpret $\supinprsym$ as taking values in $\CA^\CC$. 
Any hermitian form decomposes as the sum of two homogeneous hermitian forms $\supinprsym = \supinprsym_0 + \supinprsym_1$ when we define $\supinprsym_\alpha$ by
$$
\homsuperinprod vw{_\alpha}
=
\superinprod{v_0}{w_0}_\alpha
+
\superinprod{v_0}{w_1}_{\alpha+1}
+
\superinprod{v_1}{w_0}_{\alpha+1}
+
\superinprod{v_1}{w_1}_\alpha
\mapob,
$$
where $\superinprod xy_\beta$ denotes the homogeneous part of degree $\beta$ of $\superinprod xy\in \CA^\CC$.

\medskip

$\bullet$
A right-linear map $A:V\to V$ is said to be \stresd{graded skew-symmetric} with respect to a graded sesquilinear form $\supinprsym$ on $V$ if we have, for any two $v,w\in V$:
$$
\superinprod{Av}w + \superinprod v{A_0w} + \superinprod[2]{\conjugate(v)}{A_1w}
=
0
\mapob,
$$
where $A=A_0+A_1$ is the decomposition of $A$ into its homogeneous parts.

\end{definition}

\begin{definition}{Nota Bene}
We explicitly do \stress{not} require a super scalar product to be homogeneous!

\end{definition}

\begin{definition}[smoothnessremarkoninfinitedimensions]{Remarks}
$\bullet$ 
We have defined the notion of a hermitian form on a graded vector space over $\RR$ or $\CA$ as taking values in $\CC$ respectively $\CA^\CC$. 
The reader might have expected that in those two cases it should take values in $\RR$ respectively $\CA$ only. 
Our definition does not exclude that possibility, but we will need later on the possibility (necessity!) that a super scalar product on a graded vector space over $\RR$ or $\CA$ is allowed to take values in these larger sets (see in particular \recalt{supscalprodonbodyneedsi} and \recalt{discussioniinrealsupermetric}).

\medskip

$\bullet$
We have formulated the non-degeneracy condition in the case $\KKA=\CA^\KK$ in terms of $\body \supinprsym$. For smooth $\supinprsym$ one can show that this condition is equivalent to the same condition without the body maps. But we will need it this way because when we define a metric on an $\CA$-manifold $M$, it will not be smooth in the sense given above on every tangent space.

\medskip

$\bullet$
The name \myquote{smoothness} we used above for graded sesquilinear forms might look strange to some. 
However, when $V$ is finite dimensional, it is the necessary and sufficient condition for such a map to be smooth in the category of $\CA$-manifolds (see also \recals{appendixonAmanifoldssection}). 
Now, looking carefully at the definitions of smooth maps defined on open sets of finite dimensional $\CA$-vector spaces as explained in \recals{appendixonAmanifoldssection} shows that it should be rather straightforward to generalize this notion to smooth maps defined on open sets of arbitrary graded vector spaces of the form $V\otimes_\KK \CA^\KK$ with $V$ a {normed} vector space over $\KK$. 
If that indeed works, only an additional continuity condition is needed with our definition of smoothness to make these graded hermitian forms smooth maps over $\CA$. 
(As any (multi-) linear map defined on a finite dimensional vector space is automatically continuous (with respect to the euclidean topology), no continuity condition is needed for the finite dimensional case.) 
On the other hand, careful checking of all the details is outside the scope of this paper, so we only retain our \myquote{smoothness} condition without claiming (or using) any kind of smoothness property of a graded hermitian form in the category of $\CA$-manifolds (but we will add a continuity condition in our definition of a super Hilbert space \recalt{defofmainSHS}). 

\end{definition}

\begin{proclaim}[leftlinearityofrightbilinear]{Lemma}
Let $V$ be a graded vector space over $\KKA$ and let $\supinprsym$ be a graded sesquilinear form on $V$. 
Then we have the following \myquote{anti-linearity} property (for $v,w\in V$ and $\lambda\in \KKA$):
\begin{align*}
\superinprod{\lambda v}w
&
=
\overline{\lambda}\cdot\homsuperinprod vw{_0} 
+
\overline{\conjugate(\lambda)}\cdot \homsuperinprod vw{_1}
\\&
=
\overline{\lambda_0}\cdot\superinprod vw 
+ 
\overline{\lambda_1}\cdot\bigl(\,
\homsuperinprod vw{_0} - \homsuperinprod vw{_1}\,\bigr) 
\mapob.
\end{align*}

\end{proclaim}

\begin{definition}[remarkparityhermformonleft]{Remark}
Let $\supinprsym$ be a graded sesquilinear form on $V$ and let $\supinprsym_\alpha$ be its homogeneous parts. 
It then follows from \recalt{leftlinearityofrightbilinear} that the (anti-) linearity conditions are such that $\supinprsym_\alpha$ acts as an object of parity $\alpha$, meaning that interchanging a super object with $\supinprsym_\alpha$ gives an additional sign. 
As we write the symbol $\supinprsym$ to the left of the vectors on which we evaluate, it follows that the parity of a graded sesquilinear form is concentrated on the left. 
The way we formulated the condition of \myquote{anti-right-linearity} in \recalt{defofgradedhermitianforms} avoids all mention of parity.
\end{definition}

\begin{proclaim}[superscalarbackandforthbodyextension]{Lemma}
Let $V$ be a graded vector space over $\CA^\KK$.
\begin{enumerate}
\item
If $\supinprsym$ is a graded sesquilinear form on $V$, then $\body\supinprsym$ is a graded sesquilinear form on $\body V$. 
And if $\supinprsym$ is homogeneous, graded symmetric or non-degenerate, so is $\body \supinprsym$.

\item\label{extensiontoAvectorspaceofgradedsesquilinearform}
Any graded sesquilinear form $\supinprsym^r$ on $\body V$ extends in a unique way to a smooth graded sesquilinear form $\supinprsym$ on $V$ (\ie, $\body \supinprsym=\supinprsym^r$) by defining
\begin{align*}
\superinprod{v\otimes \lambda}{w\otimes\mu}
&
=
\homsuperinprod{v}{w_0}{^r}\cdot \lambda_0\cdot \mu
+
\homsuperinprod{v}{w_1}{^r}\cdot \overline{\lambda}_0\cdot \mu
\\&
\kern5em
+
\homsuperinprod{v}{w_0}{^r}\cdot \lambda_1\cdot \mu
-
\homsuperinprod{v}{w_1}{^r}\cdot \overline{\lambda}_1\cdot \mu
\\&
\equiv
\homsuperinprod{v}{w_0}{^r}\cdot \lambda\cdot \mu
+
\homsuperinprod{v}{w_1}{^r}\cdot \conjugate(\,\overline{\lambda}\,)\cdot \mu
\mapob,
\end{align*}
for $v,w\in V$ and $\lambda,\mu\in \KKA$. 
If $\supinprsym^r$ is homogeneous, graded symmetric or non-degenerate on $\body V$, then so is its extension $\supinprsym$ to $V$.

\end{enumerate}
\end{proclaim}

\begin{definition}{Nota Bene}
Contrary to a normal scalar product, for a super scalar product we do not have the implication $\superinprod vw = 0 \Rightarrow \superinprod wv = 0$. A simple counter example is given by a graded vector space of dimension $1\vert 2$ with basis $e,f_1,f_2$ and the (smooth, even) super scalar product defined by $\superinprod ee = 1$, $\superinprod{f_1}{f_2} = 1 = -\superinprod{f_2}{f_1}$, all other scalar products between basis elements being zero. Then $\superinprod{e+f_1}{e+f_2} = 2$ and $\superinprod{e+f_2}{e+f_1} = 0$. 

We can modify the above even super scalar product to a non-homogeneous super scalar product by \myquote{adding} the \myquote{products} $\superinprod e{f_i} = \superinprod{f_i}e = 1$, leaving the $\superinprod{f_i}{f_i}$ to be the only super scalar products between basis elements that are zero. 
In this case we have $\superinprod{f_2}{e+f_1}=0$ and $\superinprod{e+f_1}{f_2}=2$.

\end{definition}

\begin{proclaim}[supscalprodonbodyneedsi]{Lemma}
Let $\supinprsym$ be an even super scalar product on a graded vector space $V$ over $\CA^\KK$. 
Then:
\begin{enumerate}
\item
$(\body\supinprsym)({\body V_0},{\body V_1})=0$;

\item
the restriction of $\body\supinprsym$ to $\body V_0$ is a hermitian scalar product on the ungraded vector space $\body V_0$ over $\KK$ in the sense of \recalt{defsofmetricsinnongradedcase};

\item
the restriction of $-i\cdot\body\supinprsym$ to $\body V_1$ is a hermitian scalar product on the ungraded vector space $\body V_1$ over $\KK$ in the sense of \recalt{defsofmetricsinnongradedcase}.

\end{enumerate}

\end{proclaim}

\begin{definition}[defsupermeticoverC]{Definition}
A super scalar product $\supinprsym$ on a graded vector space $V$ over $\CC$ or $\CA^\CC$ will be called a \stresd{super metric} if it is \textbf{even} and satisfies the two conditions
\begin{enumerate}
\item
$\forall v\in \body V_0: (\body\supinprsym)(v,v)\ge0$ and

\item
$\forall v\in \body V_1: -i\cdot(\body\supinprsym)(v,v)\ge0$.

\end{enumerate}

\end{definition}

\begin{definition}[remarkonchoicesignforsupermetric]{Remark}
In the definition of a super metric we have made an arbitrary choice by requiring that $-i\cdot(\body\supinprsym)(v,v)\ge0$ for all $v\in \body V_1$. We could as well have required that $+\,i\cdot(\body\supinprsym)(v,v)\ge0$. 
However, there are two other instances where we can make an arbitrary choice of a sign: in the definition of an equivalence of super Hilbert spaces \recalt{remarksonequivalenceofSHSandSUR} and in the definition of the Berezin-Fourier transform \recalt{remarksondefofBerezinFourier}. 
If one wishes to have consistent results, these choices are not independent: fixing one choice imposes the other two.

\end{definition}

\begin{proclaim}[supermetricinducesmetricsonevenandodd]{Corollary}
If $\supinprsym$ is a super metric on a graded vector space $V$ over $\CC$ or $\CA^\CC$, then $\body\supinprsym\caprestricted_{\body V_0}$ and $-i\cdot\body\supinprsym\caprestricted_{\body V_1}$ are metrics (in the sense of \recalt{defsofmetricsinnongradedcase}) on the ungraded vector spaces (over $\CC$) $\body V_0$ and $\body V_1$ respectively.

\end{proclaim}

\begin{definition}[discussioniinrealsupermetric]{Discussion}
In \recalt{defsupermeticoverC} we have defined the notion of a super metric for graded vector spaces over $\CC$ or $\CA^\CC$. 
For graded vector spaces over $\RR$ or $\CA$ the definition of a super metric is similar, but with a slightly unexpected twist: they should still take values in $\CC$ or $\CA^\CC$ respectively! 
There are several arguments in favor of this definition. 

In the first place, defining a metric on a graded vector space over $\RR$ as a non-degenerate graded symmetric  bilinear form with values in $\RR$ implies automatically that the odd dimension of this vector space should be even (as on the odd-odd part such a form defines a skew-symmetric form and thus non-degeneracy implies even dimension). 
In the second place, even when we restrict attention to graded vector spaces with an even dimension of the odd part, we will miss an important property, needed when we want to apply partition of unity arguments: a linear combination with positive coefficients is not guaranteed to remain non-degenerate. 
And in the third place, we have seen that in the complex setting, the restriction of a super metric to the odd part of the graded vector space is $i$ times an ordinary (hermitian!) metric. 

\end{definition}

\begin{definition}[defsupermetricoverR]{Definition}
Let $V$ be a graded vector space over $\RR$ or $\CA$.
A \stresd{super metric on $V$} is a super metric $\supinprsym$ on $V$ in the sense of \recalt{defsupermeticoverC} with the additional \stresd{reality condition}
$$
\forall v,w\in \body V \text{ homogeneous } : (\body\supinprsym)(v,w) = (-1)^{\parity v\parity w} \ \overline{(\body\supinprsym)(v,w)}
\mapob.
$$
This condition should be compared with condition (\ref{realitycondonungradedrealvectorspaces}) in \recalt{defsofmetricsinnongradedcase} for metrics on ungraded vector spaces over $\RR$ (see also \cite[\S IV.7]{Tu04}). 

\end{definition}

\begin{proclaim}[extensionandrestrictionrealsupermetrics]{Lemma}
We collect here some interesting properties of super metrics on a graded vector space $V$ over $\RR$ or $\CA$.
\begin{enumerate}
\item
\label{restrictingsupermetrictobody}
Let $V$ be a graded vector space over $\CA$ and let $\supinprsym$ be a super metric on $V$, then $\body\supinprsym$ is a super metric on $\body V$.

\item
\label{extensionsupermetricfrombody}
Let $V$ be a graded vector space over $\CA$ and let $\supinprsym$ be a super metric on $\body V$. Then there exists a unique smooth super metric $\supinprsym^{\CA}$ on $V$ with $\body \supinprsym^\CA=\supinprsym$.

\item
If $V$ is a graded vector space over $\CA$ and $\supinprsym$ a super metric on $V$, then $\body\supinprsym\caprestricted_{\body V_0}$ and $-i\,\body\supinprsym\caprestricted_{\body V_1}$ are (ordinary) metrics in the sense of \recalt{defsofmetricsinnongradedcase} on the ungraded vector spaces (over $\RR$) $\body V_0$ and $\body V_1$ respectively.

\item
Let $\supinprsym^a$, $a=1, \dots, \ell$ be super metrics on $V$, a graded vector space over $\RR$ or $\CA$, and let $r_1, \dots, r_\ell$ be elements in $\RR$ respectively $\CA_0$ such that $\body r_a>0$. 
Then $\sum_{a=1}^\ell r_a\cdot \supinprsym^a$ is a super metric on $V$.

\end{enumerate}

\end{proclaim}

\masection{Super Hilbert spaces and super unitary representations}
\label{superunitarydefandequivalencesection}

\begin{definition}[defofmainSHS]{Definitions}
As said in the introduction, I propose to weaken the definition of a super Hilbert space. 
Let us start with the definition of a \stresd{proto super Hilbert space}, by which I will mean a triple $(E,\inprodsym, \supinprsym)$ in which $E$ is a graded vector space over $\KK$, $\inprodsym$ a metric on $E$ seen as an ungraded vector space and $\supinprsym$ a super scalar product on $E$ satisfying the following two conditions.
\begin{enumerate}[\labelMainSHS 1.]
\item\label{mainSHSlabel2}
$\inprod{E_0}{E_1}=0$. 

\item\label{mainSHSlabel3}
$\supinprsym$ is continuous with respect to the topology induced by the metric $\inprodsym$.

\end{enumerate}
And then a \stresd{super Hilbert space} will be a proto super Hilbert space $(\Hilbert,\inprodsym, \supinprsym)$ such that $\Hilbert$ is complete with respect to the topology induced by the metric $\inprodsym$, \ie, the couple $(\Hilbert,\inprodsym)$ is a Hilbert space in the usual sense. 
When $(E,\inprodsym, \supinprsym)$ is a proto super Hilbert space, the couple $(E,\inprodsym)$ is a usual inner product space, one which we can complete to a Hilbert space $\Hilbert$. 
And as $\supinprsym$ is continuous on $E$, it extends to a continuous sesquilinear map on $\Hilbert$. 
However, non-degeneracy of the extension of $\supinprsym$ to $\Hilbert$ is not automatic and has to be checked separately before one can conclude that we have a super Hilbert space. 
The basic example of a super Hilbert space is $\KK$ itself, with $\inprodsym$ and $\supinprsym$ given by
$$
\inprod xy = \superinprod xy = \overline{x}\,y
\mapob.
$$

With this notion of a super Hilbert space, I now define a \stresd{\psur{}} of a super Lie group $G$ on a super Hilbert space $(\Hilbert, \inprodsym, \supinprsym)$ as a couple $(\Dense, \rho)$ with the following properties. 
$\Dense\subset \Hilbert$ is a dense graded subspace and $\rho:G\to \Aut(\Dense\otimes\CA^\KK)$ a group homomorphism satisfying the following conditions.
\begin{enumerate}[\labelMainSUR1.]
\item\label{MainSURlabelrhopreservessupinpr}
for all $g\in  G$, $\rho(g)$ preserves $\supinprsym$ (or more accurately: preserves the restriction to $\Dense \otimes \CA^\KK$ of the extension of $\supinprsym$ to $\Hilbert \otimes \CA^\KK$);

\item\label{MainSURlabelrhopsiissmooth}
For all $\psi\in \Dense$ the map $\FgroupHilbert_\psi:G\to \Dense \otimes \CA^\KK$, $\FgroupHilbert_\psi(g) = \rho(g)\psi$ is smooth.

\item\label{MainSURlabelrhoonbodyisunitary}
For all $g\in \body G$ the map $\rho(g)$ preserves $\inprodsym$.

\end{enumerate}
If $(\Dense,\rho)$ and $(\Dense',\rho')$ are two \psur{}s of $G$ on the super Hilbert space $(\Hilbert,\inprodsym, \supinprsym)$, then $(\Dense',\rho')$ \stresd{extends} $(\Dense,\rho)$ if $\Dense\subset \Dense'$ and for all $g\in G$ we have $\rho'(g)\restricted_{\Dense} = \rho(g)$. 
We will say that a \psur{} $(\Dense, \rho)$ is a \stresd{super unitary representation of $G$ on $\Hilbert$} if it does not admit a non-trivial extension. 

\end{definition}

\begin{proclaim}[DenseinsideCinftyrhoo]{Lemma}
If the couple $(\Dense, \rho)$ is a \psur{} of $G$ on $(\Hilbert, \inprodsym, \supinprsym)$, then: 
\begin{enumerate}
\item
$\forall g\in \body G$ the map $\rho(g)$ preserves $\Dense$ (thus justifying \refmetnaam{\labelMainSUR}{MainSURlabelrhoonbodyisunitary});

\item
the group homomorphism $\rho:\body G \to \Aut(\Dense)$ extends to a unitary representation $\rho_o$ of $\body G$ on $(\Hilbert, \inprodsym)$;

\item
we have the inclusion $\Dense \subset C^\infty(\rho_o)$. 

\end{enumerate}

\end{proclaim}

\begin{preuve}
If $(\Dense,\rho)$ satisfies \refmetnaam{\labelMainSUR}{MainSURlabelrhopsiissmooth}, the map $\FgroupHilbert_\psi:G\to \Dense \otimes \CA^\KK$ is smooth. But smoothness implies that when applied to an element of $\body G$ the result lies in $\body(\Dense \otimes \CA^\KK)= \Dense$. 
It thus makes sense to require that $\rho(g)$ with $g\in \body G$ preserves $\inprodsym$ (and not its extension to $\Dense\otimes \CA^\KK$). 

Once we know that $\rho(g)$ with $g\in \body G$ preserves $\Dense$, condition \refmetnaam{\labelMainSUR}{MainSURlabelrhoonbodyisunitary} says that it is unitary. 
It thus extends by continuity to a unitary representation $\rho_o$ on the whole of $\Hilbert$. 
We then recall that the set $C^\infty(\rho_o)$ of smooth vectors for $\rho_o$ is defined by
$$
C^\infty(\rho_o) = \{\,\psi\in \Hilbert \mid \FgroupHilbert_\psi\caprestricted_{\body G} \text{ is smooth}\,\}
\mapob.
$$
Condition \refmetnaam{\labelMainSUR}{MainSURlabelrhopsiissmooth}, applied to the restriction of $\FgroupHilbert_\psi$ to $\body G$, then tells us that we have the inclusion $\Dense \subset C^\infty(\rho_o)$. 
\end{preuve}

In the non-super context, a unitary representation is a group homomorphism $\rho:G\to \Aut(\Hilbert)$ such that each $\rho(g)$ preserves the metric on the Hilbert space $\Hilbert$. 
Associated we then have the space of smooth vectors $C^\infty(\rho) \subset \Hilbert$; they form a common dense domain of essential self-adjointness for the infinitesimal generators $\tau(X)$ defined as 
$$
\tau(X)\psi = \frac{\extder}{\extder t} \bigrestricted_{t=0}\ \rho\bigl(\exp(Xt)\bigr)\psi
\mapob.
$$
In the super context the interplay between a super unitary representation and its set of smooth vectors is a bit more complicated.  
For a super unitary representation $(\Dense,\rho)$ it is the set $\Dense$ that is\slash plays the role of the set of smooth vectors of this representation according to condition \refmetnaam{\labelMainSUR}{MainSURlabelrhopsiissmooth} and the fact that there is no non-trivial extension. 
In \recalt{DenseinsideCinftyrhoo} we have shown that we can extend the restriction of $\rho$ to $\body G$ to a unitary representation $\rho_o$ of $\body G$ on $\Hilbert$. 
Now there is a very simple reason why we cannot extend $\rho$ to $\Hilbert$ on the whole of $G$: if $\Dense$ is the domain of definition of an odd generator $\tau(X)$, $X\in \body \Liealg g_1$ (and nobody will be surprised that the domain of a generator is not the whole of $\Hilbert$), then its flow is given by $\exp\bigl(\xi \tau(X)\bigr)$ for an odd parameter $\xi$. But then the smoothness condition tells us that we have 
$$
\exp\bigl(\xi \tau(X)\bigr) = \oneasmatrix + \xi\tau(X)
$$
and thus the domain of definition of this flow is the same as the domain of $\tau(X)$. 
In the conventional approach to a super unitary representation, the infinitesimal generators $\tau(X)$ for $X$ odd are supposed to be defined on $C^\infty(\rho_o)$ and they (or a suitable multiple) are supposed to be essentially self-adjoint. 
However, the condition to be (essentially) self-adjoint is equivalent to the condition that it generates a $1$-parameter unitary group of transformations, and for odd generators, this $1$-parameter group has no rôle in the super representation, as we are interested only in $\exp\bigl(\xi \tau(X)\bigr)$ with $\xi$ an odd parameter. 
Moreover, nobody will be surprised that the set of smooth vectors for a subgroup is larger than the set of smooth vectors for the whole group. 
And thus it should not come as a surprise that the alleged set of smooth vectors $\Dense$ for $\rho$ might be smaller than the set of smooth vectors $C^\infty(\rho_o)$ for the subgroup $\body G \subset G$. 
In our motivating example it turns out (see \recals{backtomotivtingexamplesection}) that we have equality $\Dense=C^\infty(\rho_o)$, but in the \myquote{$a\xi+\beta\,$}-group (see \recals{axiplusbetagroupsection}) we have an example with $\Dense \neq C^\infty(\rho_o)$.

\begin{definition}{Definition (notions of equivalence)}
$\bullet$
Let $(E, \inprodsym_E, \supinprsym_E)$ and $(F, \inprodsym_F, \supinprsym_F)$ be two proto super Hilbert spaces. 
We will say that these two spaces are \stresd{equivalent} if there exists a linear bijection $A:E\to F$ satisfying the following conditions (where $A=A_0+A_1$ is the decomposition of $A$ into homogeneous parts).
\begin{enumerate}
\item\label{conditiondirectsuminequivalenceSHS}
$\ker(A_0) + \ker(A_1) = E$ (and thus $E=\ker(A_0) \oplus \ker(A_1)\,$).

\item\label{conditionequivalenceofSHSisunitary}
$\forall x,y\in E$: $\inprod{Ax}{Ay}_F = \inprod xy_E$.

\item\label{conditionoddisoforequivalenceofSHS}
$$
\homsuperinprod{Ax}{Ay}{_F}
=
\begin{cases}
\homsuperinprod{x}{y}{_E}
&\qquad x \in E\ ,\ y\in \ker(A_1)
\\
\homsuperinprod{\conjugate x}{y}{_E}
&\qquad x\in \ker(A_1) \ , \ y\in \ker(A_0)
\\
i\,\homsuperinprod{\conjugate x}{y}{_E}
&\qquad x,y\in \ker(A_0) 
\mapob.
\end{cases}
$$

\end{enumerate}
Alternatively we will say that a map $A$ satisfying these conditions is an \stresd{equivalence} between the two super Hilbert spaces. 

\medskip

$\bullet$
Let $G$ be an $\CA$-Lie group, let $(\Dense_E, \rho_E)$ be a \psur{} of $G$ on $(\Hilbert_E, \inprodsym_E, \supinprsym_E)$ and let $(\Dense_F, \rho_F)$ be a \psur{} of (the same) $G$ on $(\Hilbert_F, \inprodsym_F, \supinprsym_F)$. 
We will say that these two representations are \stresd{equivalent} is there exists an equivalence of proto super Hilbert spaces $A:\Dense_E\to \Dense_F$ intertwining the two representations:
$$
\forall g\in G
\quad:\quad
\rho_F(g) \scirc A = A \scirc \rho_E(g)
\mapob.
$$

\end{definition}

\begin{definition}[remarksonequivalenceofSHSandSUR]{Remarks}
$\bullet$
One might have expected that an equivalence of two proto super Hilbert spaces should be an \stress{even} linear bijection. 
However, the following argument should convince the reader that at least odd linear bijections are completely natural as equivalences for super unitary representations. 
Let $E$ be a graded vector space over $\RR$ of graded dimension $p\vert q$ and let $F=\prod E$ be the parity reversal of $E$ of graded dimension $q\vert p$. 
It is immediate that the two groups $\Aut(E)$ and $\Aut(F)$ are isomorphic: it suffices to look at their matrix representation. 
It thus is quite natural to say that the two tautological representations on $E$ and $F$ are equivalent. 
But the equivalence is given by the parity reversal map, and thus is odd. 
And once we have accepted the idea of an odd equivalence, non-homogeneous ones appear naturally by taking direct sums. 
Condition (\ref{conditiondirectsuminequivalenceSHS}) might be too restrictive, but it is compatible with the procedure of taking direct sums.

\smallskip

$\bullet$
The reader might have expected that our (non-homogeneous) linear bijection $A:E\to F$ should not only be compatible with the metrics as in condition (\ref{conditionequivalenceofSHSisunitary}), but also with the super scalar products in the form of the condition
$$
\forall x,y\in E
\quad:\quad
\homsuperinprod{Ax}{Ay}{_F} = \homsuperinprod xy{_E}
\mapob.
$$
However, just looking at an odd map $A$ shows that such a condition is not compatible with the properties of a super scalar product. 
In particular graded symmetry (and reality in the case of $\KK=\RR$) are incompatible with this formula. 
On the other hand, when taking condition (\ref{conditiondirectsuminequivalenceSHS}) into account, some elementary computations show that the given formula in condition (\ref{conditionoddisoforequivalenceofSHS}) is compatible with these properties (and that it nearly is imposed by them). 
Moreover, when we extend to $\CA^\KK$ vector spaces, it will be compatible with anti-right-linearity in the first variable. 
There is, though, one arbitrary choice I made: the use of the factor $i$ in the case $x,y\in \ker(A_0)$. 
Choosing $-i$ instead will also produce a formula compatible with the properties of a super scalar product. 
However, the choice for $+i$ is imposed when we want that a super metric \recalt{defsupermeticoverC} is equivalent to a super metric (the opposite choice will turn the positive definite condition into negative definite, see also \recalt{remarkonchoicesignforsupermetric}). 

\smallskip

$\bullet$
If $(\Hilbert,\inprodsym,\supinprsym)$ is a super Hilbert space and $\Dense\subset \Hilbert$ a dense graded subspace, then the triple $(\Dense, \inprodsym\caprestricted_\Dense , \supinprsym\caprestricted_\Dense)$ is again a super Hilbert space and in particular the restriction of $\supinprsym$ to $\Dense$ remains non-degenerate (something that is not true for an arbitrary graded subspace!). 
This implies that the definition of equivalence of two \psur{}s makes sense.

\end{definition}

\begin{proclaim}{Lemma}
Let $(\Dense_E, \rho_E)$ be a \psur{} of $G$ on the super Hilbert space $(\Hilbert_E, \inprodsym_E, \supinprsym_E)$ and let $(\Dense_F, \rho_F)$ be a \psur{} of $G$ on $(\Hilbert_F, \inprodsym_F, \supinprsym_F)$. 
Then the following properties are equivalent.
\begin{enumerate}
\item
There exists a homogeneous equivalence of super Hilbert spaces $A:\Dense_E\to \Dense_F$ intertwining $\rho_E$ and $\rho_F$.

\item
There exists a homogeneous equivalence of super Hilbert spaces $A : \Hilbert_E\to \Hilbert_F$ such that $A(\Dense_E) = \Dense_F$ and intertwining $\rho_E$ and $\rho_F$. 

\end{enumerate}

\end{proclaim}

\begin{proclaim}[equivalentDefSuperUnitaryRepNEW]{Theorem (see also \cite[Thm 4.18]{deGoursacMichel:2015}\footnote{After completion of this paper (which had a long gestation period), I became aware of the paper \cite{deGoursacMichel:2015} and realized that I could and should have known of its  existence. Their theorem 4.18 is, apart from inessential details, the same as my theorem \recalt{equivalentDefSuperUnitaryRepNEW}. Moreover, their proof is much tighter than mine. Credit for this result thus is theirs.}}
Let $(\Hilbert, \inprodsym, \supinprsym)$ be a super Hilbert space, let $G$ be a connected $\CA$-Lie group and let $(\Dense, \rho)$ be a \psur{} of $G$ on $\Hilbert$. 
If we denote by $\Liealg g \cong T_eG$ the graded Lie algebra of $G$, then the map $\tau : \body \Liealg g \to \End(\Dense)$ defined by (see \refmetnaam{\labelMainSUR}{MainSURlabelrhopsiissmooth})
\begin{moneq}[defofthemaptauforSUR]
\tau(X)\psi = X\FgroupHilbert_\psi
\end{moneq}
is an even graded Lie algebra morphism with the following properties. 
\begin{enumerate}
\item\label{alternateSUR1}
For each $X\in \body \Liealg g_0$ the map $\tau(X)$ is the restriction of the infinitesimal generator of $\rho_o\bigl( \exp(tX)\bigr)$ to $\Dense$ (see \recalt{DenseinsideCinftyrhoo}).

\item\label{alternateSUR2}
For all $X\in \body \Liealg g$ the map $\tau(X)$ is graded skew-symmetric with respect to $\supinprsym$.

\item\label{alternateSUR3}
For all $g\in \body G$ and all $X\in \body \Liealg g_1$ we have
$$
\tau\bigl(\Ad(g)X\bigr) = \rho_o(g)\scirc \tau(X) \scirc \rho_o(g\mo)
\mapob.
$$

\end{enumerate}
The triple $(\rho_o, \Dense, \tau)$ is called the infinitesimal form of $(\Dense, \rho)$.

\smallskip

Conversely, if a triple $(\rho_o, \Dense,\tau)$ satisfies the conditions (\ref{alternateSUR1})--(\ref{alternateSUR3}) with $\rho_o$ an even unitary representation of $\body G$ on $\Hilbert$ (meaning that $\rho_o(g)$ preserves $\Hilbert_\alpha$), $\Dense\subset C^\infty(\rho_o)$ a dense graded subspace of $\Hilbert$ and invariant under the maps $\rho_o(g)$, $g\in \body G$ and $\tau:\body \Liealg g\to \End( \Dense )$ an even graded Lie algebra morphism, then there exists a unique group homomorphism $\rho:G\to \Aut(\Dense \otimes \CA^\KK)$ satisfying the conditions \refmetnaam{\labelMainSUR}{MainSURlabelrhopreservessupinpr}--\ref{MainSURlabelrhoonbodyisunitary} such that the map $\tau$ is given by \recalf{defofthemaptauforSUR} and such that the restriction of $\rho$ to $\body G$ is $\rho_o$. 

As a consequence, $(\Dense, \rho)$ is a super unitary representation (no non-trivial extension) if and only if the couple $(\Dense, \tau)$ does not admit a non-trivial extension in the sense that if $(\Dense',\tau')$ satisfies the same conditions with $\Dense\subset \Dense'$ and $\tau'(X)\restricted_{\Dense} = \tau(X)$ for all $X\in \body \Liealg g$, then $\Dense'=\Dense$ (and thus $\tau'=\tau)$.

\end{proclaim}

\begin{preuve}[Proof of \recalt{equivalentDefSuperUnitaryRepNEW}, direct part]
Suppose that $(\Dense, \rho)$ is a \psur{} of $G$. 
Our first task is to show that the map $\tau$ is well defined. 
According to \refmetnaam{\labelMainSUR}{MainSURlabelrhopsiissmooth} the map $\FgroupHilbert_\psi:G\to \Dense \otimes \CA^\KK$ defined as $\FgroupHilbert_\psi(g) = \rho(g)\psi$ is smooth for all fixed $\psi\in \Dense$. 
As such we thus can apply the tangent vector $X\in \Liealg g \cong \Liealg g \body T_eG$ to this function to obtain an element of $\Dense$, which shows that the element $\tau(X)\psi)$ is well defined. 
Because $\rho(g)$ in linear, it is immediate that $\tau(X)\psi$ is right-linear in $\psi$ for fixed $X\in \body \Liealg g$, and thus we have indeed a map $\tau:\body \Liealg g\to \End(\Dense)$. 
Moreover, applying a tangent vector is a linear operation, from which it follows immediately that $\tau$ is left-linear in $X$. 
Now if $X$ and $\psi$ are homogeneous, $\FgroupHilbert_\psi$ is homogeneous of the same parity as $\psi$ and hence $\tau(X)\psi=X\FgroupHilbert_\psi$ is homogeneous of parity $\parity X + \parity \psi$. 
This shows that $\tau$ is even.

Remains to prove that this $\tau$ is a super Lie algebra morphism satisfying the properties (\ref{alternateSUR1})--(\ref{alternateSUR3}). 
To start with property (\ref{alternateSUR1}) we take $X\in \body\Liealg g_0$ and $\psi\in \Dense$ and we note that (more or less by definition of the action of a tangent vector on a (smooth) function) $\tau(X)\psi$ can be defined as the derivative at $t=0$ of the smooth map $t\mapsto \rho\bigl(\exp(tX)\bigr)\psi$ from $\CA_0$ to $\Dense\otimes \CA^\KK$. 
But this is a smooth map of a single \stress{even} coordinate, so it must be the \Gextension{} of \myquote{the same} smooth map from $\RR\to \Dense$. 
Now, as $\exp(tX)\in \body G$ for $t\in \RR$ (and $X\in \body \Liealg g_0$), we have (because $\rho_o$ is the extension of $\rho\restricted_{\body G}$ to $\Hilbert$, \recalt{DenseinsideCinftyrhoo}) the equality
$$
\rho\bigl(\exp(tX)\bigr)\psi = \rho_o\bigl(\exp(tX)\bigr)\psi
\mapob.
$$
Property (\ref{alternateSUR1}) follows immediately. 

\medskip

For property (\ref{alternateSUR2}) we take two homogeneous $\psi,\chi\in \Dense$ and we note that (by additivity) we only have to prove the graded skew-symmetry property for homogeneous $X$. 
We thus start with $X\in \body \Liealg g_0$ and we consider the function $f:\RR\to \CC$ defined by
$$
f(t)
= 
\superinprod[3]{\rho\bigl( \exp(tX)\bigr)\chi}{\rho\bigl( \exp(tX)\bigr)\psi}
\mapob.
$$
According to \refmetnaam{\labelMainSUR}{MainSURlabelrhopreservessupinpr} this is a constant function. 
But according to \refmetnaam{\labelMainSHS}{mainSHSlabel3} the map $\supinprsym$ is sesquilinear continuous, hence smooth (we are here working in the category of ordinary normed vector spaces). 
And thus, using \refmetnaam{\labelMainSUR}{MainSURlabelrhopsiissmooth}, the map $f$ is smooth as composition of smooth maps. 
We thus can compute its derivative at $t=0$, which should be zero, giving
\begin{align*}
0
&
=
f'(0)
= 
\frac{\extder}{\extder t}\bigrestricted_{t=0} \superinprod[3]{\rho\bigl( \exp(tX)\bigr)\chi}{\rho\bigl( \exp(tX)\bigr)\psi}
\\[2\jot]
&
=
\superinprod{\tau(X)\chi}{\psi} 
+ 
\superinprod{\chi}{\tau(X)\psi}
\mapob.
\end{align*}
This is graded skew-symmetry with respect to $\supinprsym$ for $X\in \body\Liealg g_0$. 

For $X\in \body\Liealg g_1$ we start with the observation that for any smooth map $f:\CA_1\to E\otimes \CA^\KK$ for any graded vector space $E$ (over $\KK$), we necessarily have the equality
$$
f(\xi) = f(0) + \xi\cdot f'(0)
\mapob,
$$
with $f(0), f'(0)\in E$. 
When we apply this to the smooth (by \refmetnaam{\labelMainSUR}{MainSURlabelrhopsiissmooth}) map $f(\xi) = \rho\bigl(\exp(\xi X)\bigr)\psi \equiv \FgroupHilbert_\psi\bigl(\exp(\xi X)\bigr)$ for any $X\in \body\Liealg g_1$ and $\psi\in \Dense$, we immediately find the equality
\begin{moneq}[deftauforoddbyexpansion]
\rho\bigl(\exp(\xi X)\bigr)\psi
=
\psi + \xi\cdot\tau(X)\psi
\mapob.
\end{moneq}
We then use \refmetnaam{\labelMainSUR}{MainSURlabelrhopreservessupinpr} and the properties of a super scalar product to compute:
\begin{align*}
\superinprod[2]{\chi}{\psi}
&
=
\superinprod[2]{\rho\bigl(\exp(\xi X)\bigr)\chi}{\rho\bigl(\exp(\xi X)\bigr)\psi}
\\&
=
\superinprod[2]{\chi+\xi\cdot\tau( X)\chi}{\psi+\xi\cdot\tau( X)\psi}
\\&
=
\superinprod[2]{\chi}{\psi}
+
\superinprod[2]{\xi\cdot\tau( X)\chi}{\psi}
+
\superinprod[2]{\chi}{\xi\cdot\tau( X)\psi}
\\&
=
\superinprod[2]{\chi}{\psi}
\\&
\kern2em
+
\Bigl( \superinprod[2]{\tau( X)\chi}{\psi}
+
(-1)^{\parity \chi} \cdot
\superinprod[2]{\chi}{\tau( X)\psi}\Bigr)\cdot\xi \cdot (-1)^{1+\parity \chi+\parity \psi}
\mapob.
\end{align*}
As this is valid for all $\xi\in \CA_1$, it follows (it is one of the essential features of our graded ring $\CA$) that we have
$$
\superinprod[2]{\tau( X)\chi}{\psi}
+
(-1)^{\parity \chi} \cdot
\superinprod[2]{\chi}{\xi\cdot\tau( X)\psi}
=
0
\mapob,
$$
which is graded skew-symmetry with respect to $\supinprsym$ for $X\in \body\Liealg g_1$, finishing the proof of property (\ref{alternateSUR2}).

\medskip

For property (\ref{alternateSUR3}) we take $g\in \body G$, which guarantees that $\Ad(g)X\in \body\Liealg g_1$ and $\rho(g\mo)\psi \in \Dense$, not in $\Dense\otimes\CA^\KK\,$, so we can use \recalf{deftauforoddbyexpansion} to make the computation
\begin{align*}
\psi + \xi\cdot\tau\bigl(\Ad(g)X\bigr)\psi
&
=
\rho\Bigl( \exp\bigl(\xi\Ad(g)X\bigr)\Bigr)\psi
=
\rho\bigl( g\cdot \exp(\xi X)\cdot g\mo \bigr)\psi
\\&
=
\rho(g) \Bigl( \rho\bigl(\exp(\xi X)\bigr) \bigl( \rho(g\mo)\psi \bigr) \Bigr)
\\&
=
\rho(g) \Bigl(  \rho(g\mo)\psi + \xi\cdot \tau(X)\bigl( \rho(g\mo)\psi \bigr) \Bigr)
\\&
=
\psi + \xi\cdot \bigl( \rho(g)\scirc \tau(X)\scirc \rho(g\mo) \bigr)\psi
\mapob.
\end{align*}
Again because this is valid for all $\xi\in \CA_1$, it follows that we have
$$
\tau\bigl(\Ad(g)X\bigr)\psi
=
\bigl( \rho(g)\scirc \tau(X)\scirc \rho(g\mo) \bigr)\psi
\mapob.
$$
As $\rho_o$ is the extension of $\rho\restricted_{\body G}$ to $\Hilbert$, it follows that for $g\in \body G$ and $\psi\in \Dense$ we have the equality $\rho(g)\psi = \rho_o(g)\psi$, proving property (\ref{alternateSUR3}). 

\medskip

And thus it only remains to show 
that $\tau$ is a graded Lie algebra morphism. 
In order to do so, we will need the following preliminary result. Let $A:\Dense\to \Dense$ be a homogeneous linear map, graded skew-symmetric with respect to $\supinprsym$ and let $X\in \body\Liealg g_0$. 
Then we have the equality
\begin{moneq}[preliminarypropertyequivSURmapA]
\fracp{}{x}\bigrestricted_{x=0} A\Bigl(\rho\bigl(\exp(xX)\bigr)\psi\Bigr)
=
A\bigl(\tau(X)\psi\bigr)
\mapob.
\end{moneq}
The reason that this is not as obvious as it looks, is because we do not assume that $A$ is continuous. 
To circumvent that problem, we use the graded skew-symmetry of $A$ and the continuity of $\supinprsym$ to compute, with an arbitrary $\chi\in \Dense$:
\begin{align*}
\shifttag{7em}
\superinprod[4]{\fracp{}{x}\bigrestricted_{x=0} A\Bigl(\rho\bigl(\exp(xX)\bigr)\psi\Bigr)}{\chi}
=
\fracp{}{x}\bigrestricted_{x=0}\superinprod[4]{ A\Bigl(\rho\bigl(\exp(xX)\bigr)\psi\Bigr)}{\chi}
\\&
=
-\fracp{}{x}\bigrestricted_{x=0}\superinprod[3]{ \rho\bigl(\exp(xX)\bigr)\conjugate^{\parity A}(\psi)}{A\chi}
\\&
=
-\superinprod[3]{ \tau(X)\bigl(\conjugate^A(\psi)\bigr)}{A\chi}
=
\superinprod[3]{ A\bigl(\tau(X)\psi\bigr)}{\chi}
\mapob.
\end{align*}
The result then follows from the non-degeneracy property of $\supinprsym$.

Next we recall that the commutator of two vector fields can be realized as a second order derivative involving their flows: if the vector fields $X$ and $Y$ have flows $\phi^X_t$ and $\phi^Y_s$, then we have
$$
[X,Y]_m = 
\fracp{}{t}\bigrestricted_{t=0}\ \fracp{}{s}\bigrestricted_{s=0} 
(\phi^X_{-t}\scirc \phi^Y_s \scirc \phi^X_t)(m)
=
\fracp{}{t}\bigrestricted_{t=0}
(T\phi^X_{-t})( Y\caprestricted_{\phi^X_t(m)})
\mapob.
$$
This classical result remains valid for (smooth) homogeneous vector fields in the super context \cite[V.5.15]{Tu04}. 
We thus can apply it to left-invariant vector fields associated to homogeneous $X, Y\in \body \Liealg g$, whose flows are right-multiplication by $\exp(tX)$ respectively $\exp(sY)$ with $t$ and $s$ of the appropriate parity. 

Having said this, we start with $\psi\in \Dense$, $X\in \body\Liealg g_0$ and $Y\in \body \Liealg g$ homogeneous and we compute:
\begin{align*}
\tau([X,Y])\psi
&
=
\fracp{}{x}\bigrestricted_{x=0}\ \fracp{}{y}\bigrestricted_{y=0} 
\rho\bigl( \exp(xX)\,\exp(yY)\,\exp(-xX) \bigr)\psi
\\&
=
\fracp{}{x}\bigrestricted_{x=0}\ \fracp{}{y}\bigrestricted_{y=0} 
\rho\bigl( \exp(xX)\bigr)\rho\bigl(\exp(yY)\bigr)\rho\bigl(\exp(-xX) \bigr)\psi
\\
\text{{\tiny $x\in \RR$ and $\rho_o=\rho\restricted_{\body G}$}}
\quad&
=
\fracp{}{x}\bigrestricted_{x=0}\ \fracp{}{y}\bigrestricted_{y=0} 
\rho_o\bigl( \exp(xX)\bigr)\rho\bigl(\exp(yY)\bigr)\rho\bigl(\exp(-xX) \bigr)\psi
\\
\text{{\tiny $\rho_o(g)$ is continuous}}
\quad&
=
\fracp{}{x}\bigrestricted_{x=0} 
\rho_o\bigl( \exp(xX)\bigr)\tau(Y)\rho\bigl(\exp(-xX) \bigr)\psi
\\
\text{{\tiny \recalf{preliminarypropertyequivSURmapA} with $A=\tau(Y)$}}\quad&
=
\tau(X)\bigl(\tau(Y)\psi\bigr) + \tau(Y)\bigl(\tau(-X)\psi\bigr)
=
\bigl[\tau(X),\tau(Y)\bigr]\psi
\mapob.
\end{align*}
This proves the morphism property for $X\in \body\Liealg g_0$ and $Y\in \body \Liealg g$ homogeneous. 

By bilinearity and graded skew-symmetry of the graded Lie bracket, it now suffices to show the morphism property for $X,Y\in \body \Liealg g_1$, for which the computations start similarly, but now with odd coordinates.
\begin{align*}
\tau([X,Y])\psi
&
=
\fracp{}{\xi}\bigrestricted_{\xi=0}\fracp{}{\eta}\bigrestricted_{\eta=0} 
\rho\bigl( \exp(\xi X)\bigr)\rho\bigl(\exp(\eta Y)\bigr)\rho\bigl(\exp(-\xi X) \bigr)\psi
\\&
\oversetalign{\recalf{deftauforoddbyexpansion}}\to=
\kern0.5em
\fracp{}{\xi}\bigrestricted_{\xi=0}\fracp{}{\eta}\bigrestricted_{\eta=0} 
\bigl(\oneasmatrix+\xi\cdot\tau(X)\bigr)\bigl(\oneasmatrix+\eta\cdot\tau(Y)\bigr)\bigl(\oneasmatrix+\xi\cdot\tau(-X)\bigr)\psi
\\&
=
\fracp{}{\xi}\bigrestricted_{\xi=0}\fracp{}{\eta}\bigrestricted_{\eta=0} 
\psi + \eta\cdot \tau(Y)\psi + \eta\xi\cdot\bigl( \tau(X)\tau(Y) + \tau(Y)\tau(X) \bigr)\psi
\\&
=
\bigl( \tau(X)\tau(Y) + \tau(Y)\tau(X) \bigr)\psi
=
\bigl[\tau(X),\tau(Y)\bigr]\psi
\mapob,
\end{align*}
which shows the morphism property for $X,Y\in \body\Liealg g_1$, finishing the proof that $\tau$ is a graded Lie algebra morphism. 
\end{preuve}

In order to prove the \myquote{converse-part} of \recalt{equivalentDefSuperUnitaryRepNEW}, we need a preliminary result, which in turn needs several preliminary definitions.

\begin{definition}[defofwodandoddpofagradedvectorspace]{Definition}
Let $E$ be a finite dimensional graded vector space of dimension $d\vert n$ over $\CA^\KK$ with homogeneous basis $e_1, \dots, e_d, f_1, \dots, f_n$, \ie, $e_1, \dots, e_d$ is a basis of $\body E_0$ and $f_1, \dots, f_n$ is a basis of $\body E_1$. 
It is immediate that an element $x$ of $E_0$ retains information about the odd basis vectors because $x$ is determined by $d$ elements of $\CA^\KK_0$ (its coordinates with respect to the even basis vectors $e_i$) and $n$ elements of $\CA^\KK_1$ (its coordinates with respect to the odd basis vectors $f_j$). 
But sometimes we wish to retain only the even or only the odd basis vectors. 
We thus define the subspaces $\oddp E$ and $\wod E$ by
$$
\wod E = \Bigl\{\ \sum_{i=1}^d x_i\,e_i  \Bigm\vert  x_i\in \CA^\KK\,\Bigr\}
\quad\text{and}\quad
\oddp E = \Bigl\{\ \sum_{j=1}^n y_j\,f_j \Bigm\vert  y_j\in \CA^\KK\,\Bigr\}
\mapob.
$$
A more intrinsic definition would be
\begin{align*}
\wod E 
&
= \body E_0 \otimes \CA^\KK \equiv (\body E_0 \oplus \{0\})\otimes \CA^\KK
\\
\oddp E 
&
= \body E_1 \otimes \CA^\KK \equiv (\{0\}\oplus \body E_1)\otimes \CA^\KK
\mapob,
\end{align*}
but the previous one is easier to understand. 
Both $\wod E$ and $\oddp E$ are graded subspaces of $E$ and we have in particular the following equalities:
$$
E=\wod E \oplus \oddp E
\qquad\text{and}\qquad
E_0 = \wod E_0 \oplus \oddp E_0
\mapob.
$$

\end{definition}

\begin{definition}[forgettingaboutoddcoordinates]{Forgetting the odd coordinates}
Let $f$ be a (local) smooth function of $d$ even coordinates and $n$ odd ones. 
There thus exist ordinary smooth functions $f_I$ of $d$ real coordinates such that (see \recalt{structureofsmoothfunctions})
$$
f(x,\xi) = \sum_{I\subset \nasset} \xi^I\,(\Gextension f_I)(x)
\mapob.
$$
The body map $\body$ is defined on smooth maps and then yields what remains when we \myquote{kill} all nilpotent elements. 
In particular we have
$$
\body f = f_\emptyset
\mapob.
$$
When we now apply $\Gextension$, we get a function of $d$ even coordinates, which is given by
$$
(\Gextension\body f)(x) = (\Gextension f_\emptyset)(x) = f(x,0)
\mapob.
$$
It follows that the combination of operations $\body$ and $\Gextension$ first kills all nilpotent elements (and in particular the odd coordinates) and then reestablishes the nilpotent parts in the even coordinates. 
In other words, we get a smooth function of even coordinates only, which is the functions $f$ that has forgotten about the dependence on the odd coordinates. 
The important observation is that these two operations are \myquote{intrinsic} and thus preserve in particular equalities involving composition of maps.

Now let $M$ be any ordinary manifold of dimension $d$. We then can take an atlas with charts $U_a$, $a\in I$ and extend all coordinate changes $\varphi^{ab}$ to functions $\Gextension \varphi^{ab}$ between open sets of $\CA_0^d$. 
This defines an $\CA$-manifold $\Gextension M$ of dimension $d\vert0$. 
Conversely we can start with an $\CA$-manifold $M$ of dimension $d\vert n$, take its body to obtain an ordinary manifold $\body M$ of dimension $d$ on which we can apply the above construction to obtain an $\CA$-manifold $\Gextension\body M$ of dimension $d\vert0$. 
Looking carefully at the construction of $\Gextension\body M$ will show that it boils down to start with an atlas with coordinate charts and then to replace the transition functions $\varphi^{ab}$ by $\Gextension\body \varphi^{ab}$. 
Said differently, we take an atlas and then forget about the odd coordinates (setting them to zero) and we retain only the even coordinates. 
What should be obvious is that if the $\CA$-manifold $M$ is modelled on the graded vector space $E$, then $\Gextension \body M$ is modelled on $\wod E$. 
For that reason we also denote the $\CA$-manifold $\Gextension\body M$ as $\wod M$. 

An important application of this notion is to a super Lie group $G$. 
A super Lie group $G$ is modelled on its super Lie algebra $\Liealg g$ and one can show that the map $\Phi:\wod G \times \oddp{\Liealg g}_0 \to G$ defined by 
\begin{moneq}[diffeoGwithGwodtimesoddLiealgg1]
\Phi(g,X) = g\cdot \exp(X)
\end{moneq}
is a diffeomorphism of $\CA$-manifolds \cite[VI.1.7]{Tu04}. 
This means that, as an $\CA$-manifold, a super Lie group $G$ is completely determined by the underlying ordinary Lie group $\body G$ and the odd part $\body \Liealg g_1$ of its super Lie algebra $\Liealg g$, simply because we have the equalities
$$
\wod G = \Gextension\body G
\qquad\text{and}\qquad
\oddp{\Liealg g} = \body \Liealg g_1 \otimes \CA
\mapob.
$$
We will use the identification\slash diffeomorphism \recalf{diffeoGwithGwodtimesoddLiealgg1} quite often in the sequel. 
As such it will turn out to be important to know how the product of two group elements is given in this interpretation of $G$ as a direct product. 
It is immediate that we have, for $g,h\in \wod G$ and $X,Y\in \oddp{\Liealg g}_0$:
$$
\Phi(g,x)\cdot\Phi(h,Y) = g\cdot \exp(X)\cdot h\cdot \exp(Y)
=
gh\cdot \exp\bigl(\Ad(h\mo)X\bigr)\cdot \exp(Y)
\mapob.
$$
Now $\Ad(h\mo)X$ and $Y$ belong to $\oddp{\Liealg g}_0$, but their product does not necessarily lie in the image of $\oddp{\Liealg g}_0$ under $\exp$. 
The purpose of \recalt{separatingevenandoddinproductofexponentials} is to provide, for $X,Y\in \oddp{\Liealg g}_0$, two elements $B^{(0)}(X,Y)\in \wod{\Liealg g}_0$ and $B^{(1)}(X,Y)\in \oddp{\Liealg g}_0$ such that we have
$$
\exp(X)\cdot \exp(Y) = \exp\bigl(B^{(0)}(X,Y)\bigr)\cdot \exp\bigl(B^{(1)}(X,Y)\bigr)
\mapob.
$$
And then we can finish our computation of the group product in terms of the diffeomorphism $\Phi$:
$$
\Phi(g,X) \cdot \Phi(h,Y)
=
\Phi\Bigl(gh\,\exp\bigl(B^{(0)}(\Ad(h\mo)X,Y)\bigr)\ ,\ \exp\bigl(B^{(1)}(\Ad(h\mo)X,Y)\bigr)\Bigr)
$$
or equivalently, omitting the map $\Phi$:
$$
(g,X) \cdot (h,Y)
=
\Bigl(gh\,\exp\bigl(B^{(0)}(\Ad(h\mo)X,Y)\bigr)\ ,\ \exp\bigl(B^{(1)}(\Ad(h\mo)X,Y)\bigr)\Bigr)
\mapob.
$$

\end{definition}

\begin{definition}{Definition}
Let $\Liealg g$ be a Lie algebra (super or not) and let $B:\Liealg g\times \Liealg g\to \Liealg g$ be a map. We will say that $B$ is a \stresd{homogeneous algebraic expression in repeated commutators of degree $k\ge2$} if there exist constants $c_{i_1, \dots, i_k}\in \RR$, $i_j=1,2$ such that for all $X_1,X_2\in \Liealg g$ we have
\begin{align}
B(X_1,X_2) 
&
= \sum_{i_1, \dots, i_k=1}^2 c_{i_1, \dots, i_k}\cdot \bigl(\ad(X_{i_{k}}) \scirc \cdots \scirc \ad(X_{i_2})\bigr)(X_{i_1})
\notag
\\&
=
\sum_{i_1, \dots, i_{k}=1}^2 c_{i_1, \dots, i_k}\cdot 
\biggl[X_{i_k},\Bigl[X_{i_{k-1}},\dots\bigl[X_{i_3},[X_{i_2},X_{i_1}]\bigr]\dots\Bigr]\biggr]
\mapob.
\label{defformularepeatedcommutators}
\end{align}
We will say that $B$ is an \stresd{algebraic expression in repeated commutators} (abbreviated as AERC) if there exist a sequence $B_k$, $k\ge2$ of homogeneous algebraic expressions in repeated commutators $B_k$ of degree $k$ such that we have $B(X,Y) = \sum_{k\ge2} B_k(X,Y)$. 

The typical example of an AERC is the Baker-Campbell-Hausdorff formula. 
This is the map $\BCH:\Liealg g \times \Liealg g\to \Liealg g$ such that
$$
\exp(X)\cdot\exp(Y) = \exp\bigl(X+Y+\BCH(X,Y)\bigr)
\mapob.
$$
It starts as (see \cite[p31]{DK00})
\begin{align*}
\BCH(X,Y) 
&
= 
\tfrac12 \,[X,Y] + \tfrac1{12}\, \Bigl(\,\bigl[X,[X,Y]\bigr] + \bigl[Y,[Y,X]\bigr]\,\Bigr) 
\\&
\kern15em
+ \tfrac1{24}\,\Bigl[Y,\bigl[X,[Y,X]\bigr]\Bigr]+\cdots
\mapob.
\end{align*}

\end{definition}

\begin{proclaim}[separatingevenandoddinproductofexponentials]{Lemma}
Let $\Liealg g$ be a Lie superalgebra of finite dimension $d\vert n$. 
Then there exist AERC's $B^{(0)}(X,Y)$ and $B^{(1)}(X,Y)$ with the following properties.
\begin{enumerate}
\item\label{firstconditionEvenOddBCH}
Each term in $B^{(0)}$ contains only an even number of commutators (only even $k$ in \recalf{defformularepeatedcommutators}) and each term in $B^{(1)}$ contains only an odd number of commutators (only odd $k$ in \recalf{defformularepeatedcommutators}).

\item
For $X,Y\in \oddp{\Liealg g}_0$ we have the equality
$$
\exp(X) \cdot \exp(Y) = \exp\bigl(B^{(0)}(X,Y)\bigr) \cdot \exp\bigl(X+Y+B^{(1)}(X,Y)\bigr)
\mapob.
$$

\end{enumerate}

\end{proclaim}

\begin{preuve}
The proof is a \myquote{simple} application of the Baker-Campbell-Hausdorff formula. 
We write $B^{(0)}(X,Y) = \sum_{k=1}^\infty B_{2k}(X,Y)$ and $B^{(1)}(X,Y) = \sum_{k=1}^\infty B_{2k+1}(X,Y)$ and we apply the Baker-Campbell-Hausdorff formula on both sides of the wanted equality
$$
\exp(X) \cdot \exp(Y) = \exp\bigl(B^{(0)}(X,Y)\bigr) \cdot \exp\bigl(X+Y+B^{(1)}(X,Y)\bigr)
\mapob.
$$
This gives us the equality
\begin{align}
\exp\bigl(X+Y+\BCH(X,Y)\bigr) 
&
= 
\exp\Bigl( B^{(0)}(X,Y) + X+Y+B^{(1)}(X,Y) 
\notag
\\&
\kern3.5em
+ \BCH\bigl(B^{(0)}(X,Y)  \,,\,X+Y+B^{(1)}(X,Y)\bigr)   \Bigr)
\mapob.
\label{equalityofexponentialsforBCHinoddproduct}
\end{align}
For an ordinary (non-super) Lie algebra we should worry about convergence and we will certainly not be allowed to deduce (without careful justification) the equality
\begin{multline}
\quad
\BCH(X,Y)
= 
B^{(0)}(X,Y)+B^{(1)}(X,Y) 
\\
\label{equalityofexponentialsforBCHinoddproductwithoutexp}
+ \BCH\bigl(B^{(0)}(X,Y)  \,,\,X+Y+B^{(1)}(X,Y)\bigr)   
\mapob.
\qquad
\end{multline}
However, $X$ and $Y$ belong to $\oddp{\Liealg g}_0$ and thus are of the form $\sum_i \xi_i \,X_i$ with $X_i\in \body \Liealg g_1$ and $\xi_i\in \CA_1$. In particular these coefficients are nilpotent. It follows immediately that any homogeneous AERC $f(X,Y)$ of degree $k$ is identically zero for any $k>2n$. 
The argument is that with respect to a basis $f_1, \dots, f_n$ of $\Liealg g_1$ the two vectors $X$ and $Y$ are expressed with (at most) $2n$ different odd coefficients. And thus any term involving more than $2n$ terms must contain one of these odd coefficients twice (multiplicatively), proving that such a term is zero. 

This has two consequences: first that all AERC's in $X$ and $Y$ only contain a finite number of homogeneous terms, so no convergence problems arise. And secondly that the arguments of the exponential function in \recalf{equalityofexponentialsforBCHinoddproduct} have zero body and thus lie in the neighborhood where the exponential map is a diffeomorphism. 
We thus are allowed to deduce the equality \recalf{equalityofexponentialsforBCHinoddproductwithoutexp}, which we (re)write as
\begin{multline}
\quad
\BCH(X,Y)
= 
\sum_{k\ge2}B_k(X,Y) 
\\
+ \BCH\Bigl(\ \sum_{k\ge1}B_{2k}(X,Y)  \,,\,X+Y+\sum_{k\ge1}B_{2k+1}(X,Y)\Bigr)   
\mapob.
\qquad
\label{iterativederterminationB0andB1}
\end{multline}

We then note that if $B_k$ and $B_\ell$ are homogeneous AERC's of degree $k$ and $\ell$ respectively, then 
$$
\bigl[ B_k(X,Y),B_\ell(X,Y)\bigr]
$$
is a homogeneous AERC of degree $k+\ell$. 
It follows that we can determine the homogeneous AERC's $B_k$ iteratively from \recalf{iterativederterminationB0andB1}, starting with $B_2$. 
The reason this works is because if we equate the homogeneous AERC's of degree $k$ in \recalf{iterativederterminationB0andB1}, we find on the left hand side the homogeneous AERC of degree $k$ of $\BCH(X,Y)$, which is known. 
On the right hand side we first find $B_k(X,Y)$, followed by all homogeneous AERC's of degree $k$ coming from the $\BCH$ term. But as each term of $\BCH$ contains at least one commutator, it follows that the terms of degree $k$ can only involve $B_i$'s with $i<k$. 

To see how this works, let us determine the first few terms.
Looking at the term $\BCH$ on the right hand side of \recalf{iterativederterminationB0andB1}, we note that its first argument starts at degree $2$ and its second argument at degree $1$ (the $X+Y$ term). As $\BCH$ itself starts at degree $2$, it follows that this term starts at degree $3$. 
It follows that we can equate $B_2(X,Y)$ with the degree $2$ term of $\BCH(X,Y)$ on the left hand side of \recalf{iterativederterminationB0andB1}, \ie, 
$$
B_2(X,Y) = \tfrac12\,[X,Y]
\mapob.
$$
The only term of degree $3$ appearing in the $\BCH$ term on the right hand side of \recalf{iterativederterminationB0andB1} comes from the degree $2$ term of $\BCH$ and is given by 
$$
\tfrac12\, \bigl[B_2(X,Y), X+Y  \bigr]
\mapob,
$$
which gives us the equality
$$
\tfrac1{12}\, \Bigl(\,\bigl[X,[X,Y]\bigr] + \bigl[Y,[Y,X]\bigr]\,\Bigr) 
=
B_3(X,Y) + \tfrac12\, \bigl[B_2(X,Y), X+Y  \bigr]
\mapob.
$$
As we already know $B_2(X,Y)$, this determines $B_3(X,Y)$ as 
$$
B_3(X,Y) = \tfrac13\,\bigl[X,[X,Y]\bigr] - \tfrac16\,\bigl[Y,[Y,X]\bigr]
\mapob.
$$

Careful analysis shows that the only term of degree $4$ appearing in the $\BCH$ term on the right hand side of \recalf{iterativederterminationB0andB1} comes from the degree $3$ term of $\BCH$ applied to $B_2(X,Y)$ and $X+Y$; it is given by
$$
\tfrac1{12}\,  \bigl[X+Y,[X+Y,B_2(X,Y)]\bigr] 
\mapob.
$$
This gives us the equality
$$
\tfrac1{24}\,\Bigl[Y,\bigl[\,X,[Y,X]\,\bigr]\Bigr]
=
B_4(X,Y) + \tfrac1{12}\,  \bigl[X+Y,[X+Y,B_2(X,Y)]\,\bigr]
\mapob,
$$
which tells us that $B_4(X,Y)$ is given by
$$
B_4(X,Y) = \tfrac18\,\Bigl[Y,\bigl[\,X,[Y,X]\,\bigr]\Bigr] - \tfrac1{24}\,\Bigl[X,\bigl[\,X,[X,Y]\,\bigr]\Bigr] + \tfrac1{24}\, \Bigl[Y,\bigl[\,Y,[Y,X]\,\bigr]\Bigr]
\mapob.
$$

As $B^{(0)}$ and $B^{(1)}$ satisfy by construction the properties (i) and (ii), this finishes the proof.
\end{preuve}

\begin{preuve}[Proof of \recalt{equivalentDefSuperUnitaryRepNEW}, converse part]
Given a triple $(\rho_o, \Dense,\tau)$ satisfying conditions (\ref{alternateSUR1})--(\ref{alternateSUR3}), we  have to define a representation $\rho$ of $G$ on $\Dense \otimes \CA^\KK$ satisfying the conditions \refmetnaam{\labelMainSUR}{MainSURlabelrhopreservessupinpr}--\ref{MainSURlabelrhoonbodyisunitary} and to show that it is unique when $\tau$ is given by \recalf{defofthemaptauforSUR}. 
The idea is to extend the representation $\rho_o$ of $\body G$ in two stages to the sought-for representation $\rho\,$: first to $\wod G$ and then to $G$, and to prove at each stage that the extension preserves $\supinprsym$. 
In that way \refmetnaam{\labelMainSUR}{MainSURlabelrhopreservessupinpr} and \refmetnaam{\labelMainSUR}{MainSURlabelrhoonbodyisunitary} are automatically satisfied. 
Proving \refmetnaam{\labelMainSUR}{MainSURlabelrhopsiissmooth} then turns out not to be very hard, only the homomorphism property requires some work.

We start by proving that $\rho_o$ preserves $\supinprsym$. 
For $\chi,\psi\in \Dense$ and $X\in \body\Liealg g_0$ we consider the function $f:\RR\to \CC$ defined by
$$
f(t) = \superinprod[3]{\rho_o\bigl(\exp(tX)\bigr)\chi}{\rho_o\bigl(\exp(tX)\bigr)\psi}
\mapob.
$$
As $\supinprsym$ is sesquilinear continuous, it is smooth (because, as in the proof of direct part, we are here working in the category of ordinary normed vector spaces), so $f$ is smooth by composition of smooth maps and its derivative is given by (using property (\ref{alternateSUR1}) of $\tau$)
\begin{align*}
f'(t) 
&
= 
\superinprod[3]{\tau(X)\rho_o\bigl(\exp(tX)\bigr)\chi}{\rho_o\bigl(\exp(tX)\bigr)\psi} 
\\&
\kern3em
+ \superinprod[3]{\rho_o\bigl(\exp(tX)\bigr)\chi}{\tau(X)\rho_o\bigl(\exp(tX)\bigr)\psi}
\mapob.
\end{align*}
By graded skew-symmetry of the $\tau(X)$ with respect to $\supinprsym$ (condition (\ref{alternateSUR2})) this is zero, and hence $f$ is constant. But $\exp$ defines a diffeomorphism between a neighborhood of $0\in \Liealg g_0$ and $\body G$. And thus for all $g$ in this neighborhood the map $\rho_o(g)$ preserves $\supinprsym$. By the representation property and connectedness of $G$ it follows that it is true for all $g\in \body G$.

We next extend the representation $\rho_o$ of $\body G$ to a map (also denoted by $\rho_o$) $\rho_o:\wod G  \equiv \Gextension\body G\to \End\bigl(\Dense\otimes \CA^\KK\bigr)$ as follows. 
For each $\psi\in \Dense$ the map $\FgroupHilbert_\psi:\body G\to \Dense$ is smooth by the assumption $\Dense\subset C^\infty(\rho_o)$. 
Hence its \Gextension{} (denoted by the same symbol) $\FgroupHilbert_\psi:\wod G\to \Dense\otimes \CA^\KK$ is a smooth map in the category of $\CA$-manifolds. 
It is elementary to show that for any $g\in \wod G$, the map $\psi\mapsto \FgroupHilbert_\psi(g)$ is linear in $\psi$. 
For each $g\in \wod G$ we thus have a linear map (over $\KK$) $\rho_o(g):\Dense\to \Dense\otimes \CA^\KK$. 
Moreover, as $\rho_o(g)$ preserves the subspaces $\Hilbert_\alpha$, so does its \Gextension{}. 
Extending this $\rho_o(g)$ by $\CA^\KK$-linearity to $\Dense\otimes \CA^\KK$, we obtain an even right-linear map (over $\CA^\KK$)
$$
\rho_o(g):\Dense\otimes \CA^\KK \to \Dense\otimes \CA^\KK
\mapob.
$$
As we have the equality $\rho_o(gh)\psi = \rho_o(g)\bigl(\rho_o(h)\psi\bigr)$ for all $g,h\in \body G$, this property extends to the \Gextension{}, proving that we have a morphism of groups
$$
\rho_o:\wod G \to \Aut\bigl(\Dense \otimes \CA^\KK\bigr)
\mapob.
$$
In other words, we have a representation of $\wod G$ on $\Dense\otimes \CA^\KK$. 

In order to prove that this extension preserves $\supinprsym$ (or rather its extension to $\Dense\otimes \CA^\KK\,$), we argue as follows. 
For $\chi,\psi\in \Dense$ we have shown that the function $f:\body G\to \CC$ defined as
$$
f(g) = \superinprod[2]{\rho_o(g)\chi}{\rho_o(g)\psi}
$$
is constant $\superinprod{\chi}{\psi}$. 
This implies that in any local coordinate system $(r_1, \dots, r_d)$ on $\body G$ all $k$-th order partial derivatives $\partial_{j_1}\dots\partial_{j_k}f$ are zero. 
But we can also compute $\Gextension f$, the \Gextension{} of $f$, directly from its definition in terms of the \Gextension{} of $\FgroupHilbert_\chi$ and $\FgroupHilbert_\psi$ and the extension of $\supinprsym$ to $\Dense\otimes \CA^\KK$. 
When we do this at a point with local coordinates $x_i\in \CA_0$ with $x_i=r_i+n_i$, $r_i\in \RR$ and $n_i$ nilpotent (and even), we obtain:
\begin{align*}
(\Gextension f)(x_1, \dots, x_d)
&
=
\supinprsym\biggl(\,{\sum_{k=0}^\infty \frac1{k!} \sum_{i_1, \dots, i_k=1}^d \partial_{i_1}\dots\partial_{i_k}\FgroupHilbert_\chi(r)\,n_{i_1}\cdots n_{i_k}}
\ ,\ 
\\&
\kern7em
{\sum_{\ell=0}^\infty \frac1{\ell!}\sum_{j_1, \dots, j_\ell=1}^d \partial_{j_1}\dots\partial_{j_\ell}\FgroupHilbert_\psi(r)\,n_{j_1}\cdots n_{j_\ell}}
\,\biggr)
\\&
=
\sum_{p=0}^\infty \frac1{p!} \,\partial_{i_1}\dots\partial_{i_p}
\Bigl(\superinprod[2]{\FgroupHilbert_\chi(r)}{\FgroupHilbert_\psi(r)}\Bigr)\,n_{i_1}\cdots n_{i_p}
\\&
\equiv
\sum_{p=0}^\infty \frac1{p!} \,(\partial_{i_1}\dots\partial_{i_p}
f)(r)\,n_{i_1}\cdots n_{i_p}
=
\superinprod[2]{\FgroupHilbert_\chi(r)}{\FgroupHilbert_\psi(r)}
\\&
=
\superinprod[2]{\chi}{\psi}
\mapob,
\end{align*}
which shows that indeed $\rho_o(g)$ preserves $\supinprsym$ for all $g\in \wod G$ as wanted.\footnote{The second equality is a combinatorial fact using that $\supinprsym$ is sesquilinear continuous (condition \refmetnaam{\labelMainSHS}{mainSHSlabel3}) and thus that we have, for any partial derivative, the equality $\partial_i\bigl(\superinprod[2]{f(x)}{g(x)}\,\bigr)= \superinprod[2]{(\partial_if)(x)}{g(t)}+\superinprod[2]{f(x)}{(\partial_ig)(x)}$. Writing out the details is a bit cumbersome, so I hope the reader will forgive me for omitting them.}

We now recall that the map $\Phi(g,X)=g\cdot \exp(X)$ \recalf{diffeoGwithGwodtimesoddLiealgg1} defines a diffeomorphism $\wod G \times \oddp{\Liealg g}_0 \to G$.
We will use (mostly implicitly) this diffeomorphism to extend the representation $\rho_o$ of $\wod G$ to a representation $\rho$ of $G$ on $\Dense\otimes \CA^\KK$ by defining $\rho$ on such an element by
$$
\rho\bigl(g\cdot \exp(X)\bigr)\psi
=
\rho_o(g) \scirc \exp\bigl(\tau(X)\bigr)\psi
\mapob,
$$
where $\exp\bigl(\tau(X)\bigr)$ is defined by its power series:
$$
\exp\bigl(\tau(X)\bigr)\psi=\sum_{k=0}^\infty \frac{1}{k!}\,\tau(X)^k(\psi)
\mapob.
$$
More precisely, any $X\in \oddp{\Liealg g}_0$ is necessarily of the form $X=\sum_{i=1}^n f_i \otimes \xi_i \equiv \sum_{i=1}^n f_i\cdot \xi_i$ with $f_1, \dots, f_n$ a basis of $\body \Liealg g_1$ and $\xi_i\in \CA_1$. 
The map $\tau(X)\in \End\bigl(\Dense\otimes \CA^\KK\bigr)$ then is defined by linearity as
$$
\tau(X) = \sum_{i=1}^n \tau(f_i)\cdot \xi_i
\quad\text{or equivalently}\quad
\tau(X)(\psi\otimes \lambda) = -\sum_{i=1}^n \xi_i\cdot \bigl(\tau(f_i)\psi\bigr)\otimes \lambda
\mapob,
$$
where we used that $\tau$ is even and $f_i$ odd, which explains the minus sign in the second formula. 
Because $X$ contains $n$ odd elements, it follows immediately that the power series for $\exp\bigl(\tau(X)\bigr)$ actually stops at $k=n$, any higher order containing at least the square of one of these odd elements, which yields zero. 
And finally we note that the map $\exp\bigl(\tau(X)\bigr)$ is even, simply because each map $\tau(f_i)$ is odd, as is its coefficient $\xi_i$. 

In order to show that the maps $\FgroupHilbert_\psi:G\to \Dense\otimes\CA^\KK$,  $g\mapsto \rho(g)\psi$ are smooth for $\psi\in \Dense$ we argue as follows. 
Using the diffeomorphism $\wod G \times \oddp{\Liealg g}_0\cong G$ we have to prove smoothness in the couple $(g,X)\in \wod G \times \oddp{\Liealg g}_0\cong G$. 
Now the expression $\exp\bigl(\tau(X)\bigr)\psi$ is given by
$$
\exp\bigl(\tau(X)\bigr)\psi
=
\sum_{k=0}^n \frac1{k!}\,\tau(X)^k\psi
=
\sum_{I\subset \{1, \dots, n\}} \xi^I \psi_I
\mapob,
$$
with $\psi_I\in \Dense$ an element of the form (with $c_I\in \RR$)
$$
\psi_I = c_I\, \tau(f_{i_1})\scirc \cdots \scirc \tau(f_{i_k})\psi
\mapob.
$$
It follows that $\FgroupHilbert_\psi(g,X)=\rho(g,X)\psi$ (with $(g,X)\in \wod G \times \oddp{\Liealg g}_0\cong G\,$) is given by
$$
\FgroupHilbert_\psi(g,X) = \sum_{I\subset \{1, \dots, n\}} \xi^I \rho_o(g)\psi_I
\mapob.
$$
As $\rho_o(g)\psi$ is smooth in $g\in \wod G$, this shows that $\FgroupHilbert_\psi$ is indeed a smooth function as required for \refmetnaam{\labelMainSUR}{MainSURlabelrhopsiissmooth}.

To show that this $\rho$ is a homomorphism, we take $g,h\in \wod G$ and $X,Y\in \oddp{\Liealg g}_0$ and we compute (writing $Z=\Ad(h\mo)X\in \oddp{\Liealg g}_0$ in the second line):
\begin{align*}
\shifttag{5em}
\rho\bigl( g\cdot \exp(X)\cdot h\cdot \exp(Y)\bigr)
=
\rho\Bigl(gh\cdot\exp\bigl(\Ad(h\mo)X\bigr)\cdot\exp(Y)\Bigr)
\\&
\equiv
\rho\bigl(gh\cdot\exp(Z)\cdot\exp(Y)\bigr)
\\&
\oversetalign{\recalt{separatingevenandoddinproductofexponentials}}\to=
\quad 
\rho\Bigl(gh\cdot\exp\bigl(B^{(0)}(Z,Y)\bigr)\cdot\exp\bigl(Z+Y+B^{(1)}(Z,Y)\bigr)\Bigr)
\\&
\oversetalign{def.}\to=
\quad
\rho_o\Bigl(gh\cdot\exp\bigl(B^{(0)}(Z,Y)\bigr)\Bigr)\scirc\exp\Bigl(\tau\bigl(Z+Y+B^{(1)}(Z,Y)\bigr)\Bigr)
\\&
=
\quad
\rho_o(gh)\scirc\rho_o\Bigl(\exp\bigl(B^{(0)}(Z,Y)\bigr)\Bigr)\scirc\exp\Bigl(\tau\bigl(Z+Y+B^{(1)}(Z,Y)\bigr)\Bigr)
\\&
=
\quad
\rho_o(gh)\scirc\exp\Bigl(\tau\bigl(B^{(0)}(Z,Y)\bigr)\Bigr)\scirc\exp\Bigl(\tau\bigl(Z+Y+B^{(1)}(Z,Y)\bigr)\Bigr)
\\&
=
\quad
\rho_o(gh)\scirc\exp\Bigl(B^{(0)}\bigl(\tau(Z),\tau(Y)\bigr)\Bigr)
\\&
\kern8em
\scirc\exp\Bigl(\tau(Z)+\tau(Y)+B^{(1)}\bigl(\tau(Z),\tau(Y)\bigr)\Bigr)
\\&
\oversetalign{\recalt{separatingevenandoddinproductofexponentials}}\to=
\quad
\rho_o(gh)\scirc\exp\bigl(\tau(Z)\bigr)
\scirc\exp\bigl(\tau(Y)\bigr)
\\&
=
\quad
\rho_o(g)\scirc\rho_o(h)\scirc\exp\bigl(\tau(Z)\bigr)
\scirc\exp\bigl(\tau(Y)\bigr)
\\&
=
\quad
\rho_o(g)\scirc\exp\bigl(\rho_o(h)\scirc\tau(Z)\scirc\rho_o(h)\mo\bigr) \scirc \rho_o(h)
\scirc\exp\bigl(\tau(Y)\bigr)
\\&
=
\quad
\rho_o(g)\scirc\exp\Bigl(\tau\bigl(\Ad(h)Z\bigr)\Bigr) \scirc \rho_o(h)
\scirc\exp\bigl(\tau(Y)\bigr)
\\&
=
\quad
\rho_o(g)\scirc\exp\bigl(\tau(X)\bigr) \scirc \rho_o(h)
\scirc\exp\bigl(\tau(Y)\bigr)
\\&
\oversetalign{def.}\to=
\quad
\rho\bigl(g\cdot\exp(X)\bigr) \scirc \rho\bigl(h \cdot \exp(Y)\bigr)
\mapob.
\end{align*}
Let us explain this computation line by line.
\begin{enumerate}[{line }1:]
\item
We apply the definition of the Adjoint representation.

\item
We abbreviate $\Ad(h\mo)X$ to $Z$.

\item
We apply \recalt{separatingevenandoddinproductofexponentials}.

\item
We first recall that for any Lie superalgebra we have $[\body\Liealg g_\alpha,\body\Liealg g_1] \subset \body\Liealg g_{1-\alpha}$. 
It follows directly that necessarily $B^{(0)}(Z,Y)\in \wod{\Liealg g}_0 $ and $B^{(1)}(Z,Y)\in \oddp{\Liealg g}_0$ because $B^{(0)}$ only has repeated commutators of an even number of elements in $\oddp{\Liealg g}_0\equiv(\body\Liealg g_1 \otimes\CA)_0$ and $B^{(1)}$ only has repeated commutators of an odd number of such elements. 
In particular $\exp\bigl(B^{(0)}(Z,Y)\bigr)$ belongs to $\wod G$, so we can apply our definition of $\rho$. 

\item
We use that $\rho_o$ is a homomorphism.

\item
We use the fact that $\tau$ is the infinitesimal form of $\rho_o$ on $\Dense$ by hypothesis (\ref{alternateSUR1}) and that we have extended $\rho_o$ to $\wod G$ by \Gextension{}. 

\item
We use the fact that $\tau$ is a graded Lie algebra homomorphism and that $B^{(0)}$ and $B^{(1)}$ are AERC's, so we can commute these two operations. 

\item
We again use \recalt{separatingevenandoddinproductofexponentials}, which is allowed because this result is independent of the particular $\CA$-Lie group. 

\item
We use again that $\rho_o$ is a homomorphism.

\item
We use the definition of $\exp\bigl(\tau(Z)\bigr)$ as a power series (finite!).

\item
This is our hypothesis (\ref{alternateSUR3}).

\item
By definition of $Z$ we have the equality $\Ad(h)Z = X$.

\item
And finally we apply our definition of $\rho$ backwards.

\end{enumerate} 

Now that we have the representation $\rho$, we have to prove that it preserves $\supinprsym$. As we already know this to be true for $g\in \wod G$, it suffices to prove that we have
$$
\superinprod[3]{\rho\bigl(\exp(X)\bigr)\chi}{\rho\bigl(\exp(X)\bigr)\psi}
=
\superinprod{\chi}{\psi}
$$
for any $X\in \oddp{\Liealg g}_0$. 
To do so, we start with the observation that the condition that $\tau(X)$ with $X\in \body \Liealg g_1$ is graded skew-symmetric with respect to $\supinprsym$ (condition (\ref{alternateSUR2})) implies that, for $\lambda\in \CA_1$, the even map $\tau(\lambda X)$ is skew-symmetric with respect to $\supinprsym$ and thus for any $X\in \oddp{\Liealg g}_0$ the (even) map $\tau(X)$ is skew-symmetric.
We then compute:
\begin{align*}
\shifttag{5em}
\superinprod[2]{\rho\bigl(\exp(X)\bigr)\chi}{\rho\bigl(\exp(X)\bigr)\psi}
=
\superinprod[2]{\exp\bigl(\tau(X)\bigr)\chi}{\exp\bigl(\tau(X)\bigr)\psi}
\\&
=
\superinprod[4]{\ \sum_{k=0}^\infty \frac1{k!}\,\tau(X)^k\chi}{\sum_{\ell=0}^\infty \frac1{\ell!}\,\tau(X)^\ell\psi}
\\&
=
\sum_{k=0}^\infty\sum_{\ell=0}^\infty \,\frac1{k!\,\ell!}\cdot\superinprod{\tau(X)^k\chi}{\tau(X)^\ell\psi}
\\&
\llap{\text{\tiny write $p=k+\ell$}\quad}=
\sum_{p=0}^\infty\ \sum_{k=0}^p \frac{1}{k!\,(p-k)!}\,(-1)^k\,\superinprod{\chi}{\tau(X)^p\psi}
\\&
=
\superinprod{\chi}{\psi} + \sum_{p=1}^\infty \frac{\superinprod{\chi}{\tau(X)^p\psi}}{p!} \cdot \biggl(\ \sum_{k=0}^p \frac{p!}{k!\,(p-k)!}\,(-1)^k\cdot 1^{p-k}\biggr)
\\&
=
\superinprod{\chi}{\psi}
\mapob.
\end{align*}

We now only have to prove uniqueness of $\rho$ to finish the proof of \recalt{equivalentDefSuperUnitaryRepNEW}. 
So let $\rhoh : G\to \Aut(\Dense\otimes \CA^\KK)$ be another map satisfying the conditions \refmetnaam{\labelMainSUR}{MainSURlabelrhopreservessupinpr}--\ref{MainSURlabelrhoonbodyisunitary} with infinitesimal form $\tau$ and whose restriction to $\body G$ equals $\rho_o$. 
Hence $\rhoh$ and $\rho$ coincide on $\body G$. 
But then by \refmetnaam{\labelMainSUR}{MainSURlabelrhopsiissmooth} the restriction of $\rhoh$ to $\wod G$ must be the \Gextension{} of the restriction to $\body G$, hence $\rhoh$ coincides with $\rho$ on $\wod G$. 
And finally, for $X\in \body \Liealg g_1$ and $\xi\in \CA_1$ we must have by \refmetnaam{\labelMainSUR}{MainSURlabelrhopsiissmooth} and the fact that $\tau$ is the infinitesimal form of $\rhoh$:
$$
\rhoh\bigl(\exp(\xi X)\bigr)\psi = \psi + \xi \tau(X)\psi
\oversettext{def.}\to=
\rho\bigl(\exp(\xi X)\bigr)\psi
\mapob.
$$
As the elements of the form $\xi X$ with $X\in \body \Liealg g_1$ and $\xi\in \CA_1$ generate $\exp(\oddp{\Liealg g}_0)$, the homomorphism property shows that $\rhoh$ and $\rho$ coincide at all elements $(g,X)\in \wod G \times \oddp{\Liealg g}_0 \cong G$. 
\end{preuve}

\begin{definition}{Remark}
The condition in \recalt{equivalentDefSuperUnitaryRepNEW} that $G$ should be connected is needed only to prove, in the converse part, that $\rho$ preserves $\supinprsym$ if $\tau$ satisfies condition (\ref{alternateSUR2}). 
If one knows by other means that the maps $\rho_o(g)$, $g\in \body G$ preserve $\supinprsym$, then one can drop the connectedness condition on $G$. 

\end{definition}

\masection{The Batchelor bundle and metrics}
\label{BatchelorbundlewithHodgestarsection}

In this section we make the first step in our program to turn the left-regular representation of any super Lie group into a super unitary one. 
We start by recalling the definition of the Batchelor bundle of an $\CA$-manifold $M$. 
We then define the notion of a super metric $\mfdmetric$ on $M$ and we show how this defines an ordinary metric $\inprodsym_\mfdmetric$ on the space $C^\infty_c(M;F\otimes \CA^\KK)$ of compactly supported smooth functions with values in an arbitrary graded vector space $F\otimes \CA^\KK$ equipped with a metric. 
And we terminate this section by showing that the pull-back action associated to a particular kind of diffeomorphism of $M$ will preserve the metric $\inprodsym_\mfdmetric$, and thus defines a unitary map on $C^\infty_c(M;F\otimes \CA^\KK)$. 
But we start with a notational convention that will be used quite often in the sequel.

\begin{definition}[defnotationsetinexponent]{Definition\slash Notation}
$\bullet$
Let $X_1, \dots, X_n$ be $n$ objects that can be multiplied\slash composed. Then for any $I\subset \{1, \dots, n\}$ we define
$$
X^I = X_{i_1} \scirc \cdots \scirc X_{i_k}
\mapob,
$$
when $I=\{i_1, \dots, i_k\}$ with $i_1<i_2<\cdots<i_k$, with the convention that $X^\emptyset=1$ (or $X^\emptyset=\id$ in case of maps). 
This generalizes the convention \recalf{ConventionxitothepowersetI}.
When the set $\{1, \dots, n\}$ itself appears in an exponent or subscript, we will abbreviate it as $\nasset$, \ie, we define $\nasset$ to be shorthand for
$$
\nasset = \{1, \dots, n\}
\mapob.
$$

\end{definition}

\firstofmysubsection
\mysubsection{BatchelorbundlewithHodgestarsection}{We start with an $\CA$-manifold}

Let $M$ be an $\CA$-manifold modelled on the even part of an $\CA$-vector space $E$ of graded dimension $d\vert n$ and let $\mathcal{U} =\{\,U_a\mid a\in I\,\}$ be an atlas of coordinate charts $\varphi_a:U_a\to O_a\subset E_0$. 
It is immediate that the maps $\body \varphi_a:\body U_a \to \body O_a\subset \RR^p$ form an atlas for $\body M$. 
Denoting by $(x^a,\xi^a)$ the even and odd coordinates on a chart $U_a$ (\ie, the coordinates on $\varphi_a(U_a)=O_a\subset E_0$), the change of coordinates $\varphi^{ba}\equiv \varphi_b \scirc\varphi_a\mo$ can be written as 
$$
(x^b,\xi^b) = (\varphi_b\scirc \varphi_a\mo)(x^a,\xi^a) \equiv \varphi^{ba}(x^a,\xi^a)
\mapob.
$$
In \cite{Batchelor:1979} (see also \cite[\S IV.8]{Tu04}) it is shown that for any $\CA$-manifold there always exists an atlas such that the change of coordinate functions $\varphi^{ba}$ have the following special form:
\begin{moneq}[Batcheoratlasformcoordinatechange]
(x^b,\xi^b) = \varphi^{ba}(x^a,\xi^a)
\qquad\text{with}\qquad
x^b=\varphi^{ba}_0(x^a)
\text{ and }
\xi_i^b = A^{ba}_{ij}(x^a)\cdot \xi_j^a
\mapob,
\end{moneq}
where $\varphi^{ba}_0$ and $A^{ba}_{ij}$ are smooth functions of the even coordinates $x^a$ only.
We will call such an atlas a \stresd{Batchelor atlas}. 

Now let $\mathcal U$ be a Batchelor atlas (for the moment any atlas will do, but later on we will need a Batchelor one). 
Then the maps $\body\varphi_a:\body U_a\to \body O_a\subset \RR^d$ form an atlas for the underlying ordinary manifold $\body M$. 
Moreover, we can interpret the functions $\body A^{ba}_{ij}$ as smooth maps $\body A^{ab}:\body U_a\cap \body U_b\to \Gl(n,\RR)$. 
As such they satisfy the cocycle condition 
$$
\sum_{j=1}^n (\body A_{ij}^{cb})(m)\,(\body A_{jk}^{ba})(m)= (\body A_{ik}^{ca})(m)
\mapob,
$$
simply because the $\varphi^{ab}$ do. 
It follows that the functions $T^{ab}:\body U_a\cap \body U_b\to \Gl(n,\RR)$ defined as $T^{ab}_{ij}(m) = \bigl(\body A^{ab}(m)\bigr)^{-1}_{ji}$, \ie, $T^{ba}$ is the inverse-transpose of $\body A^{ba}$, also satisfy the cocycle condition. 
We thus can use them to define a vector bundle $V\!M\to \body M$ with typical fiber $\RR^n$. 
More precisely, $\pi:V\!M\to \body M$ is a vector bundle with local trivializing charts $\Psi_a:\pi\mo(\body U_a)\to \body U_a\times \RR^n$ and transition functions 
$$
\Psi^{ba} = \Psi_b \scirc \Psi_a\mo : (\body U_a \cap \body U_b) \times \RR^n \to (\body U_a\cap \body U_b)\times \RR^n
$$ 
given by
$$
\Psi^{ba}(m,v) = (m,w)
\qquad\text{with}\qquad
w_i = \sum_{j=1}^n T^{ba}_{ij}(m)\,v_j\in \RR
\mapob.
$$

We will say that \stresd{$M$ is oriented in the odd directions}, abbreviated as \stresd{\ood}, when the vector bundle $V\!M\to \body M$ is oriented. 
This is possible only if there exists an atlas for which all the matrices $\body A^{ab}$ have positive determinant, in which case there are two possible choices for such an orientation. When $M$ is \ood, we will say that a Batchelor atlas is \ood{} when it satisfies the condition that all matrices $\body A^{ab}$ have positive determinant. 

Now let $F$ be any graded vector space over $\KK$ equipped with a topology. We then can consider the space $C^\infty(M; F\otimes\CA^\KK)$ of smooth functions on $M$ with values in $F\otimes \CA^\KK$, which is a graded vector space over $\KK$ (and when $F=\KK$, it even is a graded $\KK$-algebra). 
We will now argue that the choice of a Batchelor atlas allows us to create an isomorphism (of graded vector spaces or graded $\KK$-algebras) 
\begin{moneq}[BatcheloridentificationsmoothfunctionswithsectionsNEW]
\sigma:
C^\infty(M;F\otimes\CA^\KK) \to \Gamma^\infty\bigl(\,\bigwedge V\!M \otimes F\to \body M\,\bigr)
\mapob,
\end{moneq}
\ie, between $C^\infty(M; F\otimes\CA^\KK)$ and the space of smooth sections of the (non-super) exterior algebra bundle (tensored with $F$) $\bigwedge V\!M \otimes F\to \body M$. 
To do so, let $f:M\to F\otimes\CA^\KK$ be a smooth function. 
Then, for any coordinate chart $U_a$ in the Batchelor atlas there exist smooth functions $f_I^a:\body O_a\to F$ such that we have 
\begin{moneq}[localexpressionsmoothsuperfunctionsNEW]
(f\scirc\varphi_a\mo)(x^a,\xi^a)
=
\sum_{I\subset \{1, \dots, n\}} (\xi^a)^I\,(\Gextension f_I^a)(x^a)
\mapob.
\end{moneq}
Using the functions $f_I^a$ we then define the local section $\sigma(f):\body U_a\to \bigwedge V\!M \otimes F$ by
\begin{moneq}[tempdefofsigmasubflinkfuncwithBatchNEW]
(\Psi_a\scirc\sigma(f)\scirc \body\varphi_a\mo)(x^a)
=
\Bigl(\ x^a
\ ,\ 
\sum_{I\subset \{1,\dots, n\}} e^I \otimes f_I^a(x^a)
\ \Bigr)
\mapob,
\end{moneq}
where $e_1, \dots, e_n$ denotes the canonical basis of $\RR^n$ and where $e^I$ is defined as $\xi^I$ by replacing the product in $\CA$ by the wedge-product, see \recalt{defnotationsetinexponent}. 
The definition of the transition functions for the bundle $V\!M$ together with the fact that our atlas is a Batchelor atlas guarantee that these local sections glue together to form a global smooth section. 
It is then routine to show that $\sigma$ is an isomorphism of graded vector spaces over $\KK$ (or in case $F=\KK$ of graded $\KK$-algebras). 
Moreover, it follows directly from the definition of the DeWitt topology on $M$ that $\sigma$ preserves compact support, \ie, that $\sigma$ also is an isomorphism between compactly supported smooth functions on $M$ with values in $F\otimes \CA^\KK$ and compactly supported smooth sections of $\bigwedge V\!M\otimes F$.

\begin{definition}[NBontransportingNNgradingtosmoothfunctions]{Nota Bene}
The Batchelor bundle (as well as the notion of \myquote{being oriented in the odd directions}) can be defined without the use of a Batchelor atlas by noting that we have the equality
$$
\body A_{ij}^{ba}(x^a) = 
\body\fracp{\xi^b_i}{\xi^a_j}(x^a)
\mapob.
$$
However, the identification between $C^\infty(M;F\otimes\CA^\KK)$ and $\Gamma^\infty\bigl(\,\bigwedge V\!M \otimes F\to \body M\,\bigr)$ depends upon the choice of a Batchelor atlas. 

\end{definition}

\begin{definition}{Remark}
The Batchelor atlas allows us to do more than make the identification $C^\infty(M;F\otimes\CA^\KK) \cong \Gamma^\infty\bigl(\,\bigwedge V\!M \otimes F\to \body M\,\bigr)$, it also allows us to transport the $\NN$-grading of $\Gamma^\infty\bigl(\,\bigwedge V\!M\otimes F\to \body M\,\bigr) = \bigoplus_{k\in\NN}\Gamma^\infty\bigl(\,\bigwedge^k V\!M\otimes F\to \body M\,\bigr)$ to $C^\infty(M;F\otimes\CA^\KK)$. 
This is a slightly stronger statement that the fact that $\sigma$ preserves parity, which \myquote{only} says that $\sigma$ provides an isomorphism
$$
\sigma:
\bigl(C_{(c)}^\infty(M;F\otimes\CA^\KK)\bigr)_\alpha \to \bigoplus_{k\in\NN}\Gamma_{(c)}^\infty\bigl(\,\bigwedge{}\!\!^{2k+\alpha} V\!M\otimes F\to \body M\,\bigr)
\mapob.
$$
The functions of degree $k\in\NN$ are those functions whose expression in local coordinates \recalf{localexpressionsmoothsuperfunctionsNEW} \stress{in the given Batchelor atlas} always have products of $k$ odd coordinate functions $\xi^I$ (\ie, with cardinality of $I$ being $k$).

\end{definition}

\begin{definition}{Nota Bene\slash Convention}
From now on we will always assume, unless stated explicitly otherwise, that we have chosen a Batchelor atlas. 
If needed, we could extend it to a maximal atlas by including all coordinate charts that have the change of coordinate form \recalf{Batcheoratlasformcoordinatechange} with all charts in the given Batchelor atlas. 

\end{definition}

An additional advantage of using a Batchelor atlas is that it allows us to extend a partition of unity for $\body M$ to a partition of unity of $M$. 
More precisely, let $\atlas=\{\,U_a\mid a\in A\,\}$ be a Batchelor atlas for $M$ and let the family $\rho_a:\body U_a\to [0,1]$ be a (smooth) partition of unity subordinated to the cover $\{\,\body U_a\mid a\in A\,\}$ of $\body M$. 
Then the family $\Gextension \rho_a:U_a\to \CA_0$ is a partition of unity subordinated to the cover $\atlas$. 
This works because when expressing these functions in the same chart, we never get an additional dependence on odd coordinates via a change of coordinates in the Batchelor atlas.

\mysubsection{BatchelorbundlewithHodgestarsection}{An intermezzo on matrices}

Let $E$ be a finite dimensional graded vector space over $\KKA$ of graded dimension $d\vert n$, let $e_1, \dots, e_d, e_{d+1}, \dots, e_{d+n}$ be a homogeneous basis of $E$ (recall that this means that $e_1, \dots, e_d$ is a basis of $\body E_0$ (these are thus even vectors) and $e_{d+1}, \dots, e_{d+n}$ a basis of $\body E_1$ (and thus odd vectors), see also \recals{appendixonAmanifoldssection}). 
If $f:E\to E$ is a right-linear map, we can associate to $f$ its matrix $f_{ij}\in  \Gl(d\vert n,\KKA)$ (the set of square matrices of size $d+n$ with values in $\KKA$) defined by
$$
f(e_i) = \sum_{j=1}^{d+n} e_j \,f_{ji}
\mapob.
$$
As we know that the first $p$ basis elements are even, we can decompose this matrix into four by writing
\begin{moneq}[decomposingmatrixintofour]
\begin{pmatrix}
f_{11} & \dots & f_{1,p+q}
\\
\vdots &&\vdots
\\
f_{p+q,1} & \dots & f_{p+q,p+q}
\end{pmatrix}
=
\begin{pmatrix} A & B \\ C & D \end{pmatrix}
\end{moneq}
with $A$ a square matrix of size $p\times p$ and the other three of appropriate corresponding sizes. 
If $f$ is even, $A$ and $D$ have even entries, whereas $B$ and $C$ have odd ones. 

For a left-linear map $f:E\to E$ we can also define its matrix, but now its definition is given by
$$
\contrf{e_i}f = \sum_{i=1}^{p+q} f_{ij}\, e_j
\mapob.
$$
But here again we can decompose this matrix in four submatrices as in \recalf{decomposingmatrixintofour}, and again if $f$ is even, $A$ and $D$ will have even entries, whereas $B$ and $C$ will have odd entries. 

Yet another way to obtain a matrix is when we have a sesquilinear map $S:E\times E \to \CA^\CC$. In that case we define the matrix $S_{ij}\in \Gl(d\vert n,\CA^\CC)$ by
\begin{moneq}[matrixforsesquilinearmaps]
S_{ij} = S(e_i,e_j)
\mapob.
\end{moneq}
And again we can decompose this matrix in four submatrices as in \recalf{decomposingmatrixintofour}, and once again: if $S$ is even, $A$ and $D$ will have even entries, whereas $B$ and $C$ will have odd entries. 

This decomposition of a matrix into four submatrices becomes quite natural in the setting of tangent maps. 
Consider a smooth map $\Phi:E_0\to E_0$ (or only defined on an open subset of $E_0$) and denote by $(x,\xi)$ the even and odd coordinates on $E$ according to $v=\sum_{i=1}^d x_i\, e_i + \sum_{j=1}^n \xi_j\,e_{d+j}$. 
We the can write
$$
\Phi(x,\xi) = \bigl(F(x,\xi), \phi(x,\xi)\bigr)
\mapob,
$$
with $F$ an even function (with $d$ components) and $\phi$ an odd function (with $n$ components). 
With these notations, the matrix of its tangent map (the Jacobian matrix) at a point $v\in E_0$ is given by
\begin{moneq}[matrixrepoftangentmap]
\mathrm{matrix}(T_v\Phi)
=
\begin{pmatrix}
\displaystyle \fracp{F}{x}(v) 
& \displaystyle\fracp{\phi}{x}(v) 
\\[3\jot]
\displaystyle\fracp{F}{\xi}(v) 
& \displaystyle\fracp{\phi}{\xi}(v)
\end{pmatrix}
\end{moneq}
or explicitly in terms of a (tangent) vector:
\begin{align}
\notag
\shifttag{6em}
\contrf[4]{\sum_{i=1}^d t_i\,e_i+\sum_{j=1}^n \tau_j\,e_{d+j}}{T_v\Phi}
=
\bigl(\ t\quad \tau\ \bigr) \cdot
\begin{pmatrix}
\displaystyle \fracp{F}{x}(v) 
& \displaystyle\fracp{\phi}{x}(v) 
\\[3\jot]
\displaystyle\fracp{F}{\xi}(v) 
& \displaystyle\fracp{\phi}{\xi}(v)
\end{pmatrix}
\cdot
\begin{pmatrix}
e_{1,\dots,d}
\\
e_{d+1,\dots,d+n}
\end{pmatrix}
\label{useofjacobianmatrixofdiffeo}
\\&
=
\sum_{p=1}^d\biggl(
\sum_{i=1}^d t_i\cdot \fracp{F_p}{x_i}(v) + \sum_{j=1}^n\tau_j \cdot \fracp{F_p}{\xi_j}(v)\, \biggr)\cdot e_p
\\&
\kern3em
+
\sum_{q=1}^n\biggl(
\sum_{i=1}^d t_i\cdot \fracp{\phi_q}{x_i}(v) + \sum_{j=1}^n\tau_j \cdot \fracp{\phi_q}{\xi_j}(v)\, \biggr)\cdot e_{d+q}
\mapob.
\notag
\end{align}
When $\Phi$ represents the change of (local) coordinates on an $\CA$-manifold $M$, it is customary to denote the coordinates on the source by $(x,\xi)$ and to change the name at the target, say to $(y,\eta)$. 
In particular one writes $\Phi(x,\xi) = (y,\eta) = \bigl(y(x,\xi),\eta(x,\xi)\bigr)$. 
In that case the homogeneous basis is given by the tangent vectors $(\partial_{ x_i}, \partial_{ \xi_j})$ or $(\partial_{ y_i}, \partial_{ \eta_j})$. 
And then \recalf{useofjacobianmatrixofdiffeo} takes the form of the associated change of basis in the tangent space $T_mM$:
\begin{moneq}[changeofbasisfortangentvectorsbyjacobian]
\begin{pmatrix}
\displaystyle\fracp{}{x}\bigrestricted_m
\\[3\jot]
\displaystyle\fracp{}{\xi}\bigrestricted_m
\end{pmatrix}
=
\begin{pmatrix}
\displaystyle \fracp{y}{x}(x,\xi) 
& \displaystyle\fracp{\eta}{x}(x,\xi) 
\\[3\jot]
\displaystyle\fracp{y}{\xi}(x,\xi) 
& \displaystyle\fracp{\eta}{\xi}(x,\xi)
\end{pmatrix}
\cdot
\begin{pmatrix}
\displaystyle\fracp{}{y}\bigrestricted_m
\\[3\jot]
\displaystyle\fracp{}{\eta}\bigrestricted_m
\end{pmatrix}
\mapob.
\end{moneq}
When both local coordinate systems belong to a Batchelor atlas, we have in particular $\partial_\xi y = 0$.

\mysubsection{BatchelorbundlewithHodgestarsection}{We add a metric to obtain a Hilbert space}

We define a \stresd{super metric $\mfdmetric\,$}\footnote{As is usual, I use the letter $g$ to denote a metric on a manifold. But as a generic element of a group $G$ is also denoted by the letter $g$, this might be confusing when both appear in the same formula (even when most of the time it will be obvious who is what). I thus write a metric in bold face, providing a visual distinction for these two objects.} on the $\CA$-manifold $M$ to be a smooth assignment of a super metric $\mfdmetric(m)$ in the sense of \recalt{defsupermetricoverR} at $T_mM$ for all $m\in M$ (see also \cite[\S IV.7]{Tu04}).\footnote{Note that this definition of a super metric differs significantly from the more standard definition of a super metric as given for instance in \cite{GarnierWurzbacher:2012} or \cite{LRT:2013}} 
Here the smoothness condition means that in terms of local coordinates (even and odd together) $(x_1, \dots, x_{d+n})$ the matrix elements
\begin{moneq}[defofmetricmatrixinaglobalbasis]
\mfdmetric_{ij}(m) =
\mfdmetric(m)\bigl(\, \partial_{x_i}\caprestricted_m\,,\, \partial_{x_j}\caprestricted_m   \, \bigr)
\qquad,\quad
i,j=1, \dots, n+d
\end{moneq}
are smooth functions of $m$.
In any local chart with coordinates $(x,\xi)$ we then can define the four matrices $\mfdmetric^{\alpha\beta}(x,\xi)$, $\alpha, \beta\in \Ztwo$ (see the previous subsection) by
\begin{align*}
\mfdmetric^{00}_{pq}(x,\xi) 
&=
\mfdmetric\Bigl(\,
\fracp{}{x_p}\bigrestricted_{(x,\xi)}\ ,\ 
\fracp{}{x_q}\bigrestricted_{(x,\xi)}
\, \Bigr)
\quad ,\quad 
\mfdmetric^{01}_{pq}(x,\xi)
=
\mfdmetric\Bigl(\, 
\fracp{}{x_p}\bigrestricted_{(x,\xi)}\ ,\ 
\fracp{}{\xi_q}\bigrestricted_{(x,\xi)}
\, \Bigr)
\\
\mfdmetric^{10}_{pq}(x,\xi) 
&= 
\mfdmetric\Bigl(\, 
\fracp{}{\xi_p}\bigrestricted_{(x,\xi)}\ ,\ 
\fracp{}{x_q}\bigrestricted_{(x,\xi)}
\, \Bigr)
\quad ,\quad 
\mfdmetric^{11}_{pq}(x,\xi) = 
\mfdmetric\Bigl(\, 
\fracp{}{\xi_p}\bigrestricted_{(x,\xi)}\ ,\ 
\fracp{}{\xi_q}\bigrestricted_{(x,\xi)}
\, \Bigr)
\mapob.
\end{align*}
The matrices $\mfdmetric^{00}$ and $\mfdmetric^{11}$ are even of size $d\times d$ and $n\times n$ respectively, whereas the matrices $\mfdmetric^{10}$ and $\mfdmetric^{01}$ are odd and of size $n\times d$ and $d\times n$ respectively. 
Now if $(x^a, \xi^a)$ and $(x^b,\xi^b)$ are two coordinate systems, we get two sets of matrices $\mfdmetric^{a,\alpha\beta}$ and $\mfdmetric^{b,\alpha\beta}$. 
Using \recalf{changeofbasisfortangentvectorsbyjacobian} it is not hard to show that these two sets are related by
\begin{moneq}[changeofbasisforametricformula]
\begin{pmatrix}
\mfdmetric^{a,00} & \mfdmetric^{a,01}
\\[2\jot]
\mfdmetric^{a,10} & \mfdmetric^{a,11}
\end{pmatrix}
=
\begin{pmatrix}\displaystyle 
\fracp{x^b}{x^a} & \displaystyle \fracp{\xi^b}{x^a}
\\[4\jot]
\displaystyle \fracp{x^b}{\xi^a} & \displaystyle \fracp{\xi^b}{\xi^a}
\end{pmatrix}
\cdot
\begin{pmatrix}
\mfdmetric^{b,00} & \mfdmetric^{b,01}
\\[2\jot]
\mfdmetric^{b,10} & \mfdmetric^{b,11}
\end{pmatrix}
\cdot
\begin{pmatrix}\displaystyle 
\Bigl(\fracp{x^b}{x^a}\Bigr)^t & \displaystyle \Bigl(\fracp{x^b}{\xi^a}\Bigr)^t
\\[4\jot]
\displaystyle -\Bigl(\fracp{\xi^b}{x^a}\Bigr)^t & \displaystyle \Bigl(\fracp{\xi^b}{\xi^a}\Bigr)^t
\end{pmatrix}
\mapob.
\end{moneq}
Now if we realize that $\partial_{x^a}\xi^b$ and $\partial_{\xi^a} x^b$ are odd matrices, taking the body map of the above equality gives us the two equalities
\begin{moneq}[changemetricbodycoordinates]
\body \mfdmetric^{a,00} = \body\fracp{x^b}{x^a} \cdot \body \mfdmetric^{b,00} \cdot \body \Bigl(\fracp{x^b}{x^a}\Bigr)^t
\quad\text{and}\quad
\body \mfdmetric^{a,11} = \body\fracp{\xi^b}{\xi^a} \cdot \body \mfdmetric^{b,11} \cdot \body \Bigl(\fracp{\xi^b}{\xi^a}\Bigr)^t
\mapob.
\end{moneq}
This shows in the first place that $\body \mfdmetric^{a,00}$ defines a global metric $\mfdmetric_{\body M}$ on $\body M$. 
We then recall that the matrices $\body A^{ba}=\body \partial_{\xi^a}\xi^b$ are the inverse transpose of the transition functions of the Batchelor vector bundle $V\!M\to \body M$. 
Taking the definition of a super metric into account, we conclude that the matrices $( -i\,\body \mfdmetric^{a,11})\mo$ define an ordinary metric $\mfdmetric_{V\!M}$ on (the fibers of) $V\!M$. 
More precisely, in the trivializing chart $U_a$ the metric $\mfdmetric_{V\!M}(m)$ is defined by the equation 
\begin{moneq}[defofinducedmetriconBatchelorBundle]
\sum_{k=1}^n \mfdmetric_{V\!M}(m)({e_j},{e_k})\cdot \body \bigl(-i\,\mfdmetric^{a,11}_{k\ell}(m)\bigr) = \delta_{j\ell}
\mapob,
\end{moneq}
where, as before, the $e_i$ denote the canonical basis of the typical fiber $\RR^n$.

Following \cite[p79]{Wa83} we extend this scalar product to the exterior algebra bundle $\bigwedge V\!M$ by setting the scalar product of homogeneous elements of different degrees to zero, and by setting, in the trivializing chart $U_a$,
$$
\mfdmetric_{V\!M}(m)({v_1\wedge \dots \wedge v_k},{w_1\wedge \dots\wedge w_k}) = \det\bigl(\,\mfdmetric_{V\!M}(m)({v_p},{w_q})\,\bigr)
\mapob.
$$
It then follows that if $v_1, \dots, v_n$ is an orthonormal basis of the fiber $\pi\mo(m)\cong \RR^n$ above $m\in \body M$ of the bundle $V\!M$, then the monomials $v^I\equiv v_{i_1}\wedge \dots\wedge v_{i_k}$, $I\subset \{1, \dots, n\}$ (see \recalt{defnotationsetinexponent}) form an orthonormal basis of $\bigwedge \pi\mo(m)\subset \bigwedge V\!M$. 

Now let $F$ be a graded vector space over $\KK$ and $\inprodsym_F$ an ordinary metric on $F$ (seen without grading), thus providing a topology on $E$. 
We then can extend the metric $\mfdmetric_{V\!M}$ on the fibers of $\bigwedge V\!M$ to an ordinary metric $\mfdmetric_{V\!M}^F$ on the fibers of $\bigwedge V\!M\otimes F$ by defining 
$$
\mfdmetric_{V\!M}^F(\omega\otimes e, \omega' \otimes e') = \mfdmetric_{V\!M}(\omega, \omega')\cdot \inprod e{e'}_F
\mapob.
$$
Using this metric $\mfdmetric_{V\!M}^F$ on the fibers of $\bigwedge V\!M\otimes F$ and the metric $\mfdmetric_{\body M}$ on $\body M$, we can define in a natural way an (ordinary) metric $\inprodsym$ on the space $\Gamma_c^\infty\bigl(\,\bigwedge V\!M\otimes F\to M\,\bigr)$ of compactly supported smooth sections of $\bigwedge V\!M\otimes F$ as follows. 
Let us denote by $\Vol_{\mfdmetric_{\body M}}$ the metric volume density associated to the metric $\mfdmetric_{\body M}$ given in a local coordinate system $x_1^a, \dots, x_d^a$ by
$$
\Vol_{\mfdmetric_{\body M}}(x^a) = \sqrt{\bigl\vert\det\bigl(\body\mfdmetric^{a,00}(x^a)\bigr)\bigr\vert}\ \extder\Leb^{(d)}(x^a)
\mapob.
$$
For any two sections $\chi,\psi\in \Gamma_c^\infty\bigl(\,\bigwedge V\!M \otimes F\to M\,\bigr)$ we then define $\inprod\chi\psi_\mfdmetric$ by
\begin{moneq}[firstdefofordinarymetriconsuperfunctionsinE]
\inprod\chi\psi_\mfdmetric
=
\int_{\body M} \mfdmetric_{V\!M}^F(m)\bigl( \chi(m) ,\psi(m)  \bigr)\ \Vol_{\mfdmetric_{\body M}}(m)
\mapob.
\end{moneq}
It is then routine to show that this is indeed a metric in the usual sense \recalt{defsofmetricsinnongradedcase}.

\begin{definition}{Remark}
Just as for the Batchelor bundle we do not need a Batchelor atlas to define the metrics $\mfdmetric_{\body M}$ and $\mfdmetric_{V\!M}$. 

\end{definition}

\begin{definition}[defofordinarmetriconCinftyMbymfdmetric]{Definition}
With these preparations we now use the isomorphism $\sigma$ to transport the (ordinary) metric $\inprodsym_\mfdmetric$ on $\Gamma_c^\infty\bigl(\,\bigwedge V\!M \otimes F\to M\,\bigr)$ to $C_c^\infty(M; F\otimes\CA^\KK)$, the space of compactly supported smooth functions on $M$ (with values in $F\otimes\CA^\KK$). 
By abuse of notation we will denote this metric on $C_c^\infty(M; F\otimes\CA^\KK)$ also by $\inprodsym_\mfdmetric$. 
This would result in the confusing definition for $f,g\in C_c^\infty(M; F\otimes\CA^\KK)$:
$$
\inprod fg_\mfdmetric = \inprod{\sigma(f)}{\sigma(g)}_\mfdmetric
\mapob,
$$
where on the right hand side $\inprodsym_\mfdmetric$ denotes the metric defined in \recalf{firstdefofordinarymetriconsuperfunctionsinE}, and on the left hand side the induced metric on $C_c^\infty(M; F\otimes\CA^\KK)$. 

We then define the (ordinary) graded Hilbert space $L^2(M; F; \mfdmetric)$ as the completion of the metric space $\bigl( C_c^\infty(M; F\otimes\CA^\KK), \inprodsym_\mfdmetric \bigr)$. 
The extended metric (scalar product) on $L^2(M; F; \mfdmetric)$ will still be denoted by $\inprodsym_\mfdmetric$, while the metric on $F$ is understood implicitly. 

\end{definition}

\begin{proclaim}{Lemma}
If the homogeneous parts of $F$ are orthogonal with respect to $\inprodsym_F$, \ie, $\inprod{F_0}{F_1}_F=0$, then the homogeneous parts of $C^\infty_c(M;F\otimes \CA^\KK)$ are orthogonal with respect to $\inprodsym_\mfdmetric$. 

\end{proclaim}

\begin{preuve}
Let $f,g\in C^\infty_c(M;F\otimes \CA^\KK)$ be such that $f$ is even and $g$ odd. 
In any local coordinate chart (in the Batchelor atlas) we thus have
$$
f(x,\xi) = \sum_{I\subset \{1, \dots, n\} } \xi^I\,f_I(x)
\qquad\text{and}\qquad
g(x,\xi) = \sum_{J\subset \{1, \dots, n\} } \xi^J\,g_J(x)
\mapob.
$$
The corresponding local expressions for $\chi=\sigma(f)$ and $\psi=\sigma(g)$ are given by \recalf{tempdefofsigmasubflinkfuncwithBatchNEW} as
$$
\chi(x)= \Bigl(x,\sum_{I\subset \{1, \dots, n\} } e^I\otimes f_I(x) \ \Bigr)
\qquad\text{and}\qquad
\psi(x)= \Bigl(x,\sum_{J\subset \{1, \dots, n\} } e^J\otimes g_J(x) \ \Bigr)
\mapob.
$$
The fact that $f$ is even implies that the parity of ${f_I(x)}$ equals $\parity I$, the parity of $I$; similarly for $g$ odd the parity of $\parity{g_J(x)}$ equals $\parity J +1$. 
Now we have
$$
\mfdmetric_{V\!M}^F(\chi(x), \psi(x)) = \sum_{I,J\subset \{1, \dots, n\} }\mfdmetric_{V\!M}(e^I, e^J)\cdot \inprod{f_I(x)}{g_J(x)}_F
\mapob.
$$
But for $\parity I\neq \parity J$ we have $\mfdmetric_{V\!M}(e^I, e^J)=0$ (homogeneous elements of different degree are orthogonal), and for $\parity I = \parity J$ we have $\inprod{f_I(x)}{g_J(x)}_F=0$ because then $f_I(x)$ and $g_J(x)$ are of opposite parity in $F$. 
The conclusion is that for even $f$ and odd $g$ we have $\mfdmetric_{V\!M}^F(\chi(x), \psi(x)) =0$, proving that $\inprod fg_\mfdmetric=0$ as wanted. 
\end{preuve}

\begin{proclaim}[Batcheloratlasadaptedtometric]{Lemma}
Let $\mfdmetric$ be a super metric on $M$ and suppose we are given an 
atlas. 
Then there exists a cover $\atlas=\{\,U_b\mid b\in A\,\}$ by local coordinate charts $U_b$ in (or better: compatible with) the given 
atlas on which the matrices $-i\,\body \mfdmetric^{b,11}$ are constant the identity matrix. 
We will say that the resulting  
atlas $\atlas$ is adapted to $\mfdmetric$. 
If the initial atlas is a Batchelor atlas and\slash or \ood{}, then so will be the resulting atlas. 

\end{proclaim}

\begin{preuve}
Let $U_a$ be any local chart in the given  
atlas with local coordinates $(x^a,\xi^a)$ and let $H^a(x)$ be the matrix defined by 
$$
H^a_{pq}(x^a) = -i\,\body \mfdmetric^{a,11}
\biggl(\,
\fracp{}{\xi^a_p}\bigrestricted_{(x^a,\xi^a)}
\ ,\ 
\fracp{}{\xi^a_q}\bigrestricted_{(x^a,\xi^a)}
\,\biggr)
\mapob.
$$
As $\mfdmetric$ is a super metric, it follows that this is a positive definite symmetric bilinear form. 
By a simple Gram-Schmid orthogonalization process we can find a smooth matrix valued function $S:\body U_a \to \Gl(n,\RR)$ with positive determinant (it will actually be a triangular matrix) such that 
\begin{moneq}
\sum_{j,k=1}^n S_{jp}(x) \cdot 
H^a_{jk}(x^a) \cdot 
S_{kq}(x)
=
\delta_{pq}
\mapob.
\end{moneq}
This suggests we introduce the chart $U_b$ with the change of coordinates $\varphi^{ba}$ defined as $U_b=U_a$ and $(x^b,\xi^b) = \varphi^{ba}(x^a,\xi^a)$ with $x^b=x^a$ and 
\begin{moneq}
\xi^a_p = \sum_q (\Gextension S)(x^a)_{pq}\,\xi^b_q \equiv  \sum_q (\Gextension S)(x^b)_{pq}\,\xi^b_q
\mapob.
\end{moneq}
As we have $\body \partial_{\xi^b} \xi^a = S(x^a)$, it follows from \recalf{changemetricbodycoordinates} that we have 
$$
H^b_{pq}(x^b) = \delta_{pq}
\mapob,
$$
\ie, the matrix valued function $H^b$ is constant the identity matrix on $U_b$. 
As the change of coordinates $\varphi^{ba}$ is linear in the odd coordinates, it is immediate that if the initial atlas is a Batchelor atlas, then so is the new one. 
And because the determinant of $S$ is positive, it follows immediately that if the initial atlas is \ood, then so is the new one. 
\end{preuve}

\begin{proclaim}[metriconsuperfunctionsinEinlocalcoordinates]{Lemma}
Let $\atlas=\{\,U_a\mid a\in A\,\}$ be a Batchelor atlas adapted to the super metric $\mfdmetric$ on $M$ \recalt{Batcheloratlasadaptedtometric} and let $\chi,\psi\in C^\infty_c(M;F\otimes \CA^\KK)$ be smooth functions. 
Let furthermore $\rho_a:U_a\to [0,1]$ be a smooth partition of unity subordinated to $\body \atlas$. 
Then the metric $\inprodsym_\mfdmetric$ on $C^\infty_c(M;F\otimes \CA^\KK)$ is given by 
\begin{align*}
\inprod \chi\psi_\mfdmetric
&
=
\sum_{a\in A} 
\int_{\body U_a} \extder\Leb^{(d)}(x^a)\ 
\rho_a(x^a)\cdot
\sqrt{\bigl\vert\det\bigl(\body\mfdmetric^{a,00}(x^a)\bigr)\bigr\vert}
\\&
\kern7em
\cdot
\sum_{I\subset \{1, \dots, n\} } \inprod[2]{\chi_{a,I}(x^a)}{\psi_{a,I}(x^a)}_F
\mapob,
\end{align*}
where $(x^a, \xi^a)$ denotes a local system of coordinates on $U_a$ and where the local expression of $\chi$ (respectively $\psi$) is given by
$$
\chi(x^a,\xi^a)
=
\sum_{I\subset \{1, \dots, n\} } (\xi^a)^I\,\chi_{a,I}(x^a)
$$
for smooth functions $\chi_{a,I}:\body U_a\to F$.

\end{proclaim}

\begin{preuve}
As the atlas is adapted to the super metric, the matrices $-i\,\body \mfdmetric^{a,11}$ are constant the identity. 
This implies that the canonical basis $e_1, \dots, e_n$ of the typical fiber of $V\!M$ is an orthonormal basis. 
We thus have
\begin{align*}
\inprod \chi\psi_\mfdmetric
&
=
\int_{\body M} \mfdmetric_{V\!M}^F(m)\bigl( \chi(m) ,\psi(m)  \bigr)\ \Vol_{\mfdmetric_{\body M}}(m)
\\&
=
\sum_{a\in A}
\int_{\body U_a} \rho_a(m)\cdot \mfdmetric_{V\!M}^F(m)\bigl( \chi(m) ,\psi(m)  \bigr)\ \Vol_{\mfdmetric_{\body M}}(m)
\\&
=
\sum_{a\in A} 
\int_{\body U_a} \extder\Leb^{(d)}(x^a)\ 
\rho_a(x^a)\cdot
\sqrt{\bigl\vert\det\bigl(\body\mfdmetric^{a,00}(x^a)\bigr)\bigr\vert}
\\&
\kern6em
\cdot\sum_{I,J\subset \{1, \dots, n\} }
\mfdmetric_{V\!M}(e^I,e^J) \cdot \inprod{\chi_{a,I}(x^a)}{\psi_{a,J}(x^a)}_F
\\&
=
\sum_{a\in A} 
\int_{\body U_a} \extder\Leb^{(d)}(x^a)\ 
\rho_a(x^a)\cdot
\sqrt{\bigl\vert\det\bigl(\body\mfdmetric^{a,00}(x^a)\bigr)\bigr\vert}
\\&
\kern6em
\cdot\sum_{I\subset \{1, \dots, n\} }
\inprod{\chi_{a,I}(x^a)}{\psi_{a,I}(x^a)}_F
\mapob,
\end{align*}
as claimed. 
\end{preuve}

\mysubsection{BatchelorbundlewithHodgestarsection}{And then we add a diffeomorphism}

\begin{definition}{Definitions}
Let $\Phi:M\to M$ be a diffeomorphism. 
We will say that \stresd{$\Phi$ is linear in the odd coordinates (with respect to a Batchelor atlas)}, if for any $m\in M$, any coordinate system $(x,\xi)$ around $m$ and any coordinate system $(y,\eta)$ around $\Phi(m)$ (coordinate systems in the fixed Batchelor atlas) this diffeomorphism is of the form
\begin{moneq}[diffeolinearinoddcoordinates]
\Phi(x,\xi) = (y,\eta) 
\qquad\text{with}\qquad
y=\Phi_0(x)
\quad,\quad
\eta_i = \sum_{j=1}^n\Phi_{1,ij}(x)\,\xi_j
\end{moneq}
for (local) smooth functions $\Phi_0$ and $\Phi_{1,ij}$ of the even coordinates $x$ only (the condition thus is slightly stronger than the name suggests, as the even coordinates $y$ do not depend upon the odd coordinates $\xi$). 
It is immediate that $\Phi_0:\body M\to \body M$ is a diffeomorphism of the underlying ordinary (\ie, non-super) manifold $\body M$ and that $\Phi_1(x)$ is an invertible $n\times n$ matrix depending smoothly on $x$. 

We will say that \stresd{$\Phi$ preserves the super metric $\mfdmetric$} if for any $m\in M$ and any pair of tangent vectors $v,w\in T_mM$ we have the equality
\begin{moneq}[invariantmetricunderdiffeo]
\mfdmetric_{\Phi(m)}\bigl(T_m\Phi(v), T_m\Phi(w)\bigr)
=
\mfdmetric_m(v,w)
\mapob.
\end{moneq}

\end{definition}

\begin{definition}{Construction}
If a diffeomorphism $\Phi:M\to M$ is linear in the odd coordinates, it defines a bundle isomorphism $\Phih:V\!M\to V\!M$ as follows. 
Let $(x,\xi)$ and $(y,\eta)$ be as in \recalf{diffeolinearinoddcoordinates} and let $(x,v)$ with $v\in \RR^n$ be a point in $V\!M$ above $\body x\cong \body m\in \body M$ in the local trivialization associated to the coordinate chart $(x,\xi)$. 
Then we define $\Phih(x,v)$ by
$$
\Phih(x,v) = (y,w)
\qquad\text{with}\qquad
y=\Phi_0(x)
\quad,\quad
w_i = \sum_{j=1}^n \bigl(\Phi_1(x)\mo\bigr)_{ji}\,v_j
\mapob.
$$
Note that we use the inverse-transpose of the matrix $\Phi_1(x)$ (just as we did in the definition of $V\!M$), which makes this definition a valid definition on $V\!M$. 
Associated to this bundle isomorphism $\Phih:V\!M\to V\!M$ we have the induced bundle isomorphism (denoted by the same symbol, taking the trivial action on $F$) $\Phih:\bigwedge V\!M \otimes F\to \bigwedge V\!M \otimes F$. 
And this bundle isomorphism induces an isomorphism $\Phih^*$ of the space of smooth sections $\Gamma\bigl(\,\bigwedge V\!M \otimes F\to \body M\,\bigr)$ as follows. 
For any smooth section $s:\body M\to \bigwedge V\!M \otimes F$ we define
$$
(\Phih^*s)(m)
=
\Phih\mo\Bigl(s\bigl(\Phi_0(m)\bigr)\Bigr)
\mapob.
$$
It is then routine to show that the identification $\sigma$ \recalf{BatcheloridentificationsmoothfunctionswithsectionsNEW} between super smooth functions on $M$ with values in $F\otimes \CA^\KK$ and smooth sections of $\bigwedge V\!M\otimes F \to \body M$ intertwines the usual pull-back $\Phi^*$ of super smooth functions with $\Phih^*$: for any $f\in C^\infty(M;F\otimes\CA^\KK)$ we have the equality
$$
\sigma(\Phi^*f) = \Phih^*\bigl(\sigma(f)\bigr)
\mapob.
$$

\end{definition}

\begin{proclaim}[linearinoddcoordinateshencepullbackunitary]{Proposition}
Let $\Phi:M\to M$ be a diffeomorphism that preserves the super metric $\mfdmetric$ and is linear in the odd coordinates (with respect to a given Batchelor atlas). 
Then we have the following properties.
\begin{enumerate}
\item
$\Phi:\body M \to \body M$ preserves the ordinary metric $\mfdmetric_{\body M}$.

\item
The induced vector bundles isomorphism $\Phih:V\!M \to V\!M$ preserves the ordinary metric $\mfdmetric_{V\!M}$.

\item
The pull-back operation $\Phi^*: C^\infty_c(M; F\otimes \CA^\KK)$ is an even map preserving the metric $\inprodsym_\mfdmetric$. 
It thus induces an even unitary map on $L^2(M; F; \mfdmetric)$.

\end{enumerate}

\end{proclaim}

\begin{preuve}
Using coordinate systems $(x,\xi)$ and $(y,\eta)$ as in \recalf{diffeolinearinoddcoordinates} and using that $\Phi$ is linear in the odd coordinates, we have (see \recalf{useofjacobianmatrixofdiffeo})
$$
\contrfoper(
\begin{pmatrix}\displaystyle 
\fracp{}{x} \\[4\jot] 
\displaystyle \fracp{}{\xi}
\end{pmatrix}
)T_m\Phi
=
\begin{pmatrix}\displaystyle 
\fracp{\Phi_0}{x}(x) & \displaystyle \fracp{\Phi_{1}}{x}(x)\cdot \xi
\\[4\jot]
\mathbf0 & \displaystyle {\Phi_{1}}(x)
\end{pmatrix}
\cdot
\begin{pmatrix}\displaystyle 
\fracp{}{y} 
\\[4\jot] 
\displaystyle \fracp{}{\eta}
\end{pmatrix}
\mapob.
$$
Substituting this in the invariance condition of $\mfdmetric$ and then taking the body map gives us (just as \recalf{changemetricbodycoordinates}) the equalities
\begin{align*}
\body \mfdmetric^{x,00}(x) 
&
= \body\fracp{\Phi_0}{x} \cdot \body \mfdmetric^{y,00}\bigl(\Phi_0(x)\bigr) \cdot  \bigl(\body\fracp{\Phi_0}{x}\bigr)^t
\\[2\jot]
\body \mfdmetric^{x,11}(x) 
&
= \body\Phi_1(x) \cdot \body \mfdmetric^{y,11}\bigl(\Phi_0(x)\bigr) \cdot \body \bigl(\Phi_1(x)\bigr)^t
\mapob.
\end{align*}
This means, taking the definition of the bundle isomorphism $\Phih$ into account, that the diffeomorphism $\Phi_0:\body M\to \body M$ preserves the metric $\mfdmetric_{\body M}$ on $\body M$ and that the bundle isomorphism $\Phih$ preserves the metric $\mfdmetric_{V\!M}$ on the fibers of $V\!M$. 
It is immediate that the induced\slash extended action on $\bigwedge V\!M \otimes F$ (with trivial action on $F$) preserves the metric $\mfdmetric_{V\!M}^F$. 
The following routine computation then shows that $\Phih^*$ preserves the metric $\inprodsym_\mfdmetric$ on $\Gamma_c^\infty\bigl(\,\bigwedge V\!M \otimes F \to M\,\bigr)$.
\begin{align*}
\inprod[2]{\Phih^*\chi}{\Phih^*\psi}_\mfdmetric
&
=
\int_{\body M} \mfdmetric_{V\!M}^F(m)\bigl( (\Phih^*\chi)(m) ,(\Phih^*\psi)(m)  \bigr)\ \Vol_{\mfdmetric_{\body M}}(m)
\\&
=
\int_{\body M} \mfdmetric_{V\!M}^F(m)\biggl( \Phih\mo\Bigl(\chi\bigl(\Phi_0(m)\bigr)\Bigr) ,\Phih\mo\Bigl(\psi\bigl(\Phi_0(m)\bigr)\Bigr)  \biggr)\ \Vol_{\mfdmetric_{\body M}}(m)
\\
\text{\tiny inv. of $\mfdmetric_{V\!M}^F$}\quad&
=
\int_{\body M} \mfdmetric_{V\!M}^F\bigl(\Phi_0(m)\bigr)\Bigl( \chi\bigl(\Phi_0(m)\bigr) , \psi\bigl(\Phi_0(m)\bigr)  \Bigr)\ \Vol_{\mfdmetric_{\body M}}(m)
\\
\text{\tiny inv. of $\mfdmetric_{\body M}$}\quad&
=
\int_{\body M} \mfdmetric_{V\!M}^F\bigl(\Phi_0(m)\bigr)\Bigl( \chi\bigl(\Phi_0(m)\bigr) , \psi\bigl(\Phi_0(m)\bigr)  \Bigr)\ \Vol_{\mfdmetric_{\body M}}\bigl(\Phi_0(m)\bigr)
\\&
=
\int_{\body M} \mfdmetric_{V\!M}^F(m)\bigl( \chi(m) , \psi(m)  \bigr)\ \Vol_{\mfdmetric_{\body M}}(m)
=
\inprod{\chi}{\psi}_\mfdmetric
\mapob.
\end{align*}
As $\sigma$ intertwines $\Phih^*$ with $\Phi^*$ and as $\Phi^*$ preserves the parity of functions, it is immediate that $\Phi^*$ preserves parity and the metric $\inprodsym_\mfdmetric$.  
\end{preuve}

\masection{Densities, integration on \texorpdfstring{$\CA$}{A}-manifolds and super scalar products}
\label{Berezinintegrationsection}

The purpose of this section is to start with a super metric $\mfdmetric$ on an $\CA$-manifold $M$ and to use it to construct (via a trivializing density) a super scalar product $\supinprsym_\mfdmetric$ on the space $C^\infty_c(M; E\otimes \CA^\KK)$ of compactly supported smooth functions with values in a proto super Hilbert space $E$. 
In order to do so, we first recall some well know properties of the Berezinian and Berezin integration on $\CA$-manifolds. We then recall\slash introduce the notion of a (super) density and we show how to integrate densities over $\CA$-manifolds. 
If $M$ is \ood{}, a super metric $\mfdmetric$ defines two nearly identical trivializing densities $\nu_\mfdmetric$ and $\nu_{B,\mfdmetric}$, which allows us to define two super scalar product $\supinprsym_\mfdmetric$ and $\supinprsym_{B,\mfdmetric}$. 
We finish by showing that the super scalar product $\supinprsym_{B,\mfdmetric}$ always is continuous with respect to the (ordinary) metric $\inprodsym_\mfdmetric$ defined in \recalt{defofordinarmetriconCinftyMbymfdmetric}.

\begin{definition}{Definition}
Let $E$ be a finite dimensional graded vector space over $\KKA$ and let $f:E\to E$ be an automorphism, \ie, an even bijective  (right- and left-) linear map. 
In order to define $\Ber(f)$, the \stresd{Berezinian of $f$}, we decompose its matrix with respect to a homogeneous basis of $E$ into four submatrices as in \recalf{decomposingmatrixintofour}:
$$
\mathrm{matrix}(f) = 
\begin{pmatrix} A & B \\ C & D \end{pmatrix}
\mapob.
$$
And then $\Ber(f)$ is defined by the formula
$$
\Ber(f) = \Det(A-BD\mo C)\cdot \Det(D)\mo
\mapob,
$$
which is independent of the choice of the chosen homogeneous basis. 
We thus obtain a map $\Ber:\Aut(E)\to \KKA_0^*$ defined on the group of all automorphisms of $E$ and taking values in $\KKA_0$ (the invertible elements in the even part of $\KKA$). It is a group homomorphism. 

We also introduce $\piBer(f)$, a variant of the Berezinian defined as
$$
\piBer(f) = \vert \Det(A-BD\mo C)\vert \cdot \Det(D)\mo
\mapob.
$$
It is, as $\Ber$, a homomorphism $\Aut(E)\to \KKA_0^*$.
Note that in case $\KKA=\CA^\KK$, the absolute value is only defined for elements $x\in \KKA_0$, \ie, for elements with $\body x\neq0$, in which case its definition reads
$$
\vert x\vert = \vert \body x\vert \cdot \frac{x}{\body x}
\mapob.
$$
This is the standard \Gextension{} of the smooth function \myquote{absolute value} defined on $\KK^*\equiv \KK\setminus\{0\}$.

\end{definition}

\begin{definition}{Nota Bene}
The Berezinian is also defined for the matrix of a metric (which is in particular an even sesquilinear map, see \recalf{matrixforsesquilinearmaps}). 
However, in that case its value depends upon the choice of the homogeneous basis (as we will see in the proof of \recalt{canonicaltrivdensonOOD}).

\end{definition}

\begin{proclaim}{Lemma}
$M$ is \ood{} only if there exists an atlas $\{U_a\mid a\in A\}$ such that for all $a,b\in A$ we have
$$
\body\piBer(T\varphi^{ba})>0
\mapob.
$$

\end{proclaim}

\begin{definition}[defofelementaryBerezinintegration]{Berezin integration}
Let $E$ be a finite dimensional graded vector space over $\CA$, let $U\subset E_0$ be an open set and let $f: U\to \CA^\KK$ be a smooth map. 
This implies that there exist smooth functions $f_I:U\to \CA^\KK$, $I\subset \{1, \dots, n\}$ such that we have 
$$
f(x,\xi) = \sum_{I\subset \{1, \dots, n\}} \xi^I \,(\Gextension f_I)(x)
\mapob,
$$
where $(x,\xi)$ are the even and odd coordinates on $E_0$. 
We then define the \stresd{Berezin integral of $f$ over $U$} as the element in $\KK$ defined by
$$
\int_U \extder(x,\xi)\ f(x,\xi) = \int_{\body U} f_\nasset(x)\ \extder\Leb^{(d)}(x)
\mapob,
$$
provided of course that the function $f_\nasset:\body U\to \KK$ is Lebesgue integrable over $\body U$ (with $\Leb^{(d)}$ denoting the Lebesgue measure on $\RR^d$). 

\end{definition}

\begin{proclaim}[changeofcoordinatesinsuperint]{Proposition (\cite{Le80}, see also \cite{Roths:1987} and \cite[Thm. 4.6.1]{Va04})}
Let $E$ be a graded vector space of graded dimension $d\vert n$, let $\varphi:U\to V$ be a (super) diffeomorphism between two open subsets of $E_0$, and let $(x, \xi) = (x_1, \dots, x_d, \xi_1, \dots, \xi_n)$ denote the even and odd coordinates on $E_0$. 
Then for any function $f:V\to \CA^\KK$ with compact support we have the equality
\begin{moneq}[superchangeofcoordinates]
\int_U \extder(x,\xi)\ 
\piBer\bigl( T_{(x,\xi)}\varphi \bigr) \cdot f\bigl(\varphi(x,\xi)\bigr)
=
\int_{V} \extder(y,\eta)\ 
f(y,\eta)
\mapob.
\end{moneq}

\end{proclaim}

\begin{definition}{Remark}
The use of $\piBer$ instead of the ordinary Berezinian $\Ber$ in \recalf{superchangeofcoordinates} allows us to ignore whether $\varphi$ is orientation preserving (in the even coordinates) or not.

\end{definition}

\begin{definition}{Definitions}
Let $M$ be an $\CA$-manifold of graded dimension $d\vert n$ modelled on the graded vector space $E$ (over $\CA$).

\noindent
$\bullet$
The \stresd{frame bundle $\framebundle M\to M$} is the bundle whose fibers $\framebundle_m$ consist of all homogeneous bases of the tangent space $T_mM$ (see \recalt{defofwodandoddpofagradedvectorspace} or just after \recalt{defoftheDeWittTopologyinAppendix}). It is a principal $\Gl(d\vert n,\CA)$ bundle over $M$. 
The right-action of $\Gl(d\vert n,\CA)$ on $\framebundle_mM$ is defined as follows. 
For $(v_1, \dots, v_{d+n})\in \framebundle_mM$, \ie, a basis of $T_mM$, and $(A_{ij})\in \Gl(d\vert n,\CA)$ we have
$$
(v_i)\cdot (A_{jk}) = (w_\ell)
\qquad\text{with}\qquad
w_i = \sum_{j=1}^{d+n} v_j \,A_{ji}
\mapob.
$$

\medskip

\noindent
$\bullet$
A \stresd{density} on $M$ is a map $\nu:\framebundle M\to \CA^\KK$ satisfying the condition that for all $v=(v_i)\in \framebundle_mM$ and for all $A=(A_{jk}) \in \Gl(d\vert n,\CA)$ we have the equality
$$
\nu\bigl( v\cdot A) = \piBer(A)\cdot \nu(v)
\mapob.
$$
A density can be seen as a section of the \stresd{density \myquote{line} bundle $\Density(M)$} associated to the principal $\Gl(d\vert n,\CA)$ bundle $\framebundle M$ by the representation $\Gl(d\vert n,\CA)\to\Aut(\CA^\KK)$, $A\mapsto \piBer(A)\mo$ on $\CA^\KK$. It thus has $\CA^\KK$ as typical fiber.

\end{definition}

\begin{definition}[integrationofdensityonsupermanifold]{Construction}
Let $\nu$ be a smooth density with compact support on an $\CA$-manifold $M$ of graded dimension $d\vert n$. 
We then can define the integral of $\nu$ over $M$, denoted as $\int_M \nu$ as follows. 
We choose an open cover $\mathcal{U}=\{\,U_a\mid a\in A\,\}$ of coordinate charts with (local, even and odd) coordinates $x^a, \xi^a$ and a (smooth) partition of unity $\rho_a$ associated to this cover. 
We then define $\int_M\nu$ by
$$
\int_M\nu = 
\sum_{a\in A} 
\int_{U_a} \extder(x^{a} , \xi^{a})\  
\rho_a(x^a,\xi^a)\cdot 
\nu\Bigl( \fracp{}{x^a}\bigrestricted_{(x^a,\xi^a)}, \fracp{}{\xi^a}\bigrestricted_{(x^a,\xi^a)} \Bigr)
\mapob.
$$

\end{definition}

\begin{definition}[partitionofunitywheninBatcheloratlas]{Remark}
When we use a Batchelor atlas, we can \myquote{simplify} the partition of unity used in the definition of $\int_M \nu$. 
It then suffices to choose first a partition of unity $\rhoh_a$ associated to the cover $\{\,\body U_a \mid a\in A\,\}$ of $\body M$, and then to extend these maps by \Gextension{} to a partition of unity $\rho_a$ defined as
$$
\rho_a(x^a,\xi^a) = (\Gextension\rhoh)(x^a)
$$
independent of the odd coordinates. 

\end{definition}

\begin{proclaim}{Proposition}
The value of $\int_M\nu$ is independent of the choices of a coordinate cover and the partition of unity.

\end{proclaim}

\begin{proclaim}[parityofBerezinintegration]{Lemma}
The map $\nu \mapsto \int_M\nu$ defined on $\Gamma^\infty_c\bigl(\Density(M)\bigr)$, the set of smooth densities with compact support on $M$ (smooth sections with compact support of the density bundle $\Density(M)$) is a homogeneous (right-) linear map of parity $n$ (when $d\vert n$ is the dimension of $M$).

\end{proclaim}

\begin{proclaim}[superintinvariantdiffeo]{Lemma}
Let $\varphi:M\to M$ be a diffeomorphism, let $\nu$ be a density with compact support and let $\varphi^*\nu$ be the density defined by
$$
(\varphi^*\nu)(v_m) = \nu\bigl(\varphi_*(v_m)\bigr)
\equiv
\nu\bigl( (T_m\varphi)(v_m)  \bigr)
\mapob,
$$
where $v_m$ denotes any frame at $m\in M$ (and thus $\varphi_*v_m$ a frame at $\varphi(m)$). 
Then $\int_M \nu = \int_M \varphi^*\nu$. 

\end{proclaim}

\begin{definition}[Mobiusexampledensities]{Examples}
We consider two $2$-dimensional manifolds: the cylinder $\Cyl=\SScircle^1\times \RR$ and the Möbius bundle $\Mob$, both (real) line bundles over the circle $\SScircle^1$, as well as their graded counter parts $\Cyl\oddif$ and $\Mob\oddif$ of graded dimension $1\vert1$. 
The Möbius bundle $\Mob$ is defined using three charts $U_i\times \RR$ with $U_1=(0,\pi)\subset \RR$, $U_2=(\tfrac23\pi, \tfrac53\pi)$ and $U_3=(\tfrac43\pi,\tfrac73\pi)$. 
We \myquote{thus} have the maps $\varphi^{ij}$ giving the change of coordinates from (a part of) $U_j$ to (a part of) $U_i$. 
For $x\in (\tfrac23\pi,\pi)$, the coordinate change  $\varphi^{21}$ is given by
$$
\varphi^{21}(x,y) = (x,y)
\mapob,
$$
for $x\in (\tfrac43\pi,\tfrac53\pi)$ the coordinate change $\varphi^{32}$ is given by 
$$
\varphi^{32}(x,y) = (x,y)
\mapob,
$$
and for $x\in (2\pi,\tfrac73\pi)$ the change of charts $\varphi^{13}$ is given by
$$
\varphi^{13}(x,y) = (x-2\pi,-y)
\mapob.
$$
The graded version $\Cyl\oddif$ is just the direct product of the (graded version of the) circle with $\CA_1$, whereas the graded version $\Mob\oddif$ is given by three charts as above by replacing the direct product with $\RR$ by the direct product with $\CA_1$ and by replacing in the coordinate changes the (real) coordinate $y$ by the odd coordinate $\xi$. 
It is immediate that sections of the line bundle $\Cyl\to \SScircle^1$ can be identified with (real valued) functions on $\SScircle^1$. 
On the other hand, sections of the line bundle $\Mob\to \SScircle^1$ necessarily have a zero somewhere. 
The $\CA$-manifold $\Mob\oddif$ is particularly interesting because it is not \ood. 

We first investigate the notion of a density on $\Cyl\oddif$. 
It is fairly easy to show that any density $\nu$ on $\Cyl\oddif$ determines two functions $f_0,f_1:\SScircle^1\to \KK$ by
\begin{moneq}[densityonCyloddbytwofunctions]
\nu\Bigl( \fracp{}{x}\bigrestricted_{(x,\xi)}, \fracp{}{\xi}\bigrestricted_{(x,\xi)}\Bigr) = f_0(x) + \xi\,f_1(x)
\mapob.
\end{moneq}
It follows immediately that the density $\nu_o$ defined by $f_0(x) \equiv 1$ and $f_1(x) \equiv 0$ is a trivializing density in the sense that any smooth density $\nu$ is of the form $\nu = f\cdot \nu_o$ for some function $f\in C^\infty(\Cyl\oddif)$. 
As a (smooth) function $f:\Cyl\oddif\to \CA^\KK$ on $\Cyl\oddif$ is also determined by two smooth functions on the circle $f_0,f_1:\SScircle^1\to \KK$ by $f(x,\xi) = f_0(x)+\xi\,f_1(x)$, it is elementary to show that we have
$$
\int_{\Cyl\oddif}\nu = \int_{\SScircle^1} f_1(x)\ \extder x
\equiv \int_0^{2\pi} f(x)\ \extder x
\mapob,
$$
both for the density $\nu$ determined by the functions $f_0,f_1$ as in \recalf{densityonCyloddbytwofunctions}, or for the density $\nu=f\cdot \nu_o$, product of the trivializing density $\nu_o$ and the function $f$. 

\smallskip

We next investigate the notion of a density on $\Mob\oddif$. On any one of the three charts $U_i$ with coordinates $(x_i,\xi_i)$, a density $\nu$ determines two functions $f_0^{(i)}, f_1^{(i)}$ by
$$
\nu\Bigl( \fracp{}{x_i}\bigrestricted_{(x_i,\xi_i)}, \fracp{}{\xi_i}\bigrestricted_{(x_i,\xi_i)}\Bigr) = f_0^{(i)}(x_i) + \xi_i\,f_1^{(i)}(x_i)
\mapob.
$$
As we have the Berezinians $\piBer(T\varphi^{21})=\piBer(T\varphi^{32}) = -\piBer(T\varphi^{13})$, we must have
$$
f_0^{(1)}(x_1) + \xi_1\,f_1^{(1)}(x_1)
=
f_0^{(2)}(x_2) + \xi_2\,f_1^{(2)}(x_2)
=
f_0^{(2)}(x_1) + \xi_1\,f_1^{(2)}(x_1)
\mapob,
$$
for $x_1\in (\tfrac23\pi,\pi)$. 
For $x_2\in (\tfrac43\pi,\tfrac53\pi)$ we must have
$$
f_0^{(2)}(x_2) + \xi_2\,f_1^{(2)}(x_2)
=
f_0^{(3)}(x_3) + \xi_3\,f_1^{(3)}(x_3)
=
f_0^{(3)}(x_2) + \xi_2\,f_1^{(3)}(x_2)
\mapob,
$$
whereas for $x_3\in (2\pi,\tfrac73\pi)$ we must have
\begin{align*}
f_0^{(3)}(x_3) + \xi_3\,f_1^{(3)}(x_3)
&
=
\piBer(\varphi'_{31})\cdot\bigl(f_0^{(1)}(x_1) + \xi_1\,f_1^{(1)}\bigr)(x_1)
\\&
=
-f_0^{(1)}(x_3-2\pi) + \xi_3\,f_1^{(3)}(x_3-2\pi)
\mapob.
\end{align*}
We now have to distinguish the real case from the complex case. 
If $\nu$ is a density with values in $\CA$, our computations show that the functions $f_0^{(i)}$ determine a section $s_0$ of $\Mob\to \SScircle^1$ and that the functions $f_1^{(i)}$ determine a (real valued) function $f_1$ on $\SScircle^1$ (a section of $\Cyl$). 
The inevitable conclusion is that in this case the bundle of densities $\Density(\Mob\oddif)$ cannot be trivial as any section of the Möbius bundle necessarily has a zero. 

On the other hand, if $\nu$ is a density with values in $\CA^\CC$, the conclusion changes. The functions $f_1^{(i)}$ still determine a function $f_1$, now complex valued, on $\SScircle^1$, but the functions $f_0^{(i)}$ no longer determine a section of the Möbius bundle, simply because they now are complex valued. 
Moreover, it is easy to show that the functions $f_0^{(j)}$ defined as
$$
f_0^{(j)}(x) = \eexp^{ix/2}
$$
satisfy the given conditions. It follows that the ($\CA^\CC$-valued) density $\nu_o$ defined as
$$
\nu_o\Bigl( \fracp{}{x_j}\bigrestricted_{(x_j,\xi_j)}, \fracp{}{\xi_j}\bigrestricted_{(x_j,\xi_j)}\Bigr) = \eexp^{ix_j/2} 
$$
is a trivializing density, \ie, any $\CA^\CC$-valued density $\nu$ on $\Mob\oddif$ can be written as $\nu = f\cdot\nu$ for some $\CA^\CC$-valued function on $\Mob\oddif$. 

However, in both cases a density $\nu$ determines a well defined function $f_1:\SScircle^1\to \KK$ and it is straightforward to show that we have
$$
\int_{\Mob\oddif} \nu = \int_{\SScircle^1} f_1(x) \ \extder x
\equiv \int_0^{2\pi} f_1(x) \extder x
\mapob.
$$

\end{definition}

\begin{proclaim}[canonicaltrivdensonOOD]{Proposition}
Let $M$ be \ood{} and $\mfdmetric$ a metric on $M$.  
\begin{enumerate}
\item
\label{denstrivwithmetricandOOD}
The density $\nu_\mfdmetric$ defined on any coordinate chart $U_a$ with local coordinates $(x^a, \xi^a)$ by
$$
\nu_\mfdmetric\Bigl( \fracp{}{x^a}\bigrestricted_{m}, \fracp{}{\xi^a}\bigrestricted_{m} \Bigr)
=
{\sqrt{\bigl\vert\Ber \bigl( \mfdmetric_{ij}^a(m) \bigr)\bigr\vert}} 
$$
is a globally defined density which is everywhere invertible, \ie, trivializing.

\item
\label{BatdenstrivwithmetricandOOD}
Given a Batchelor atlas, the density $\nu_{B,\mfdmetric}$ defined on any coordinate chart $U_a$ in the Batchelor atlas with local coordinates $(x^a, \xi^a)$ by
$$
\nu_{B,\mfdmetric}\Bigl( \fracp{}{x^a}\bigrestricted_{m}, \fracp{}{\xi^a}\bigrestricted_{m} \Bigr)
=
{\sqrt{\bigl\vert\Ber \bigl( \Gextension\body\mfdmetric_{ij}^a(m) \bigr)\bigr\vert}} 
$$
is a globally defined trivializing density (see also \recalt{forgettingaboutoddcoordinates}).

\end{enumerate}

\end{proclaim}

\begin{preuve}
$\bullet$ (\ref{denstrivwithmetricandOOD}):
If $U_a$ and $U_b$ are two local coordinates charts, then we have seen in \recalf{changeofbasisfortangentvectorsbyjacobian} and \recalf{changeofbasisforametricformula} that the bases $\partial_{x^a}, \partial_{\xi^a}$ and $\partial_{x^b}, \partial_{\xi^b}$ are related by
\begin{moneq}
\begin{pmatrix}
\displaystyle\fracp{}{x^a}\bigrestricted_m
\\[3\jot]
\displaystyle\fracp{}{\xi^a}\bigrestricted_m
\end{pmatrix}
=
\begin{pmatrix}\displaystyle 
\fracp{x^b}{x^a} & \displaystyle \fracp{\xi^b}{x^a}
\\[4\jot]
\displaystyle \fracp{x^b}{\xi^a} & \displaystyle \fracp{\xi^b}{\xi^a}
\end{pmatrix}
\cdot
\begin{pmatrix}
\displaystyle\fracp{}{x^b}\bigrestricted_m
\\[3\jot]
\displaystyle\fracp{}{\xi^b}\bigrestricted_m
\end{pmatrix}
\quad\text{with}\quad
\begin{pmatrix}\displaystyle 
\fracp{x^b}{x^a} & \displaystyle \fracp{\xi^b}{x^a}
\\[4\jot]
\displaystyle \fracp{x^b}{\xi^a} & \displaystyle \fracp{\xi^b}{\xi^a}
\end{pmatrix}
=
\mathrm{matrix}(T_m\varphi^{ba})
\end{moneq}
and that the matrices $\mfdmetric^a_{ij}$ and $\mfdmetric^b_{ij}$ are related by
\begin{moneq}
\begin{pmatrix}
\mfdmetric^{a,00} & \mfdmetric^{a,01}
\\[2\jot]
\mfdmetric^{a,10} & \mfdmetric^{a,11}
\end{pmatrix}
=
\begin{pmatrix}\displaystyle 
\fracp{x^b}{x^a} & \displaystyle \fracp{\xi^b}{x^a}
\\[4\jot]
\displaystyle \fracp{x^b}{\xi^a} & \displaystyle \fracp{\xi^b}{\xi^a}
\end{pmatrix}
\cdot
\begin{pmatrix}
\mfdmetric^{b,00} & \mfdmetric^{b,01}
\\[2\jot]
\mfdmetric^{b,10} & \mfdmetric^{b,11}
\end{pmatrix}
\cdot
\begin{pmatrix}\displaystyle 
\Bigl(\fracp{x^b}{x^a}\Bigr)^t & \displaystyle \Bigl(\fracp{x^b}{\xi^a}\Bigr)^t
\\[4\jot]
\displaystyle -\Bigl(\fracp{\xi^b}{x^a}\Bigr)^t & \displaystyle \Bigl(\fracp{\xi^b}{\xi^a}\Bigr)^t
\end{pmatrix}
\mapob.
\end{moneq}
In terms of the action of $\Gl(d\vert n, \CA)$ on the frame bundle $\framebundle M$ (the action is on the right rather than on the left) the bases $\partial_{x^a}, \partial_{\xi^a}$ and $\partial_{x^b}, \partial_{\xi^b}$ are related by the matrix
$$
\begin{pmatrix}\displaystyle 
\Bigl(\fracp{x^b}{x^a}\Bigr)^t & \displaystyle \Bigl(\fracp{x^b}{\xi^a}\Bigr)^t
\\[4\jot]
\displaystyle -\Bigl(\fracp{\xi^b}{x^a}\Bigr)^t & \displaystyle \Bigl(\fracp{\xi^b}{\xi^a}\Bigr)^t
\end{pmatrix}
\mapob.
$$
As we have (for even matrices) the equality
\begin{align*}
\Ber(\begin{pmatrix} A^t & C^t \\ - B^t & D^t\end{pmatrix})
&
=
\det(A^t + C^t\,D^{t-1}\,B^t) \cdot \det(D^t)\mo
\\&
=
\det(A - B\,D\mo\,C) \cdot \det(D)\mo 
=
\Ber(\begin{pmatrix} A & B \\ C & D\end{pmatrix})
\mapob,
\end{align*}
it follows that we have the equality
\begin{moneq}[changeofBergfortrivdensity]
\Ber\bigl(\mfdmetric^{a}_{ij}(m)\bigr) = \Ber(T_m\varphi^{ba})^2\cdot \Ber\bigl(\mfdmetric^{b}_{ij}(m)\bigr)
\mapob.
\end{moneq}
As we assume that $M$ and the chosen Batchelor atlas are \ood, we have in particular that $\Ber(T_m\varphi^{ba})^2 = \piBer(T_m\varphi^{ba})^2$ and $\body\piBer(T_m\varphi^{ba})>0$. It follows that we have
\begin{align*}
\nu_\mfdmetric\Bigl( \fracp{}{x^a}\bigrestricted_{m}, \fracp{}{\xi^a}\bigrestricted_{m} \Bigr)
&
=
\sqrt{\bigl\vert\Ber\bigl(\mfdmetric^{a}_{ij}(m)\bigr)\bigr\vert}
=
\piBer(T_m\varphi^{ba})\cdot \sqrt{\bigl\vert\Ber\bigl(\mfdmetric^{b}_{ij}(m)\bigr)\bigr\vert}
\\&
=
\piBer(T_m\varphi^{ba}) \cdot
\nu_\mfdmetric\Bigl( \fracp{}{x^b}\bigrestricted_{m}, \fracp{}{\xi^b}\bigrestricted_{m} \Bigr)
\mapob.
\end{align*}
As this is the right behavior, $\nu_\mfdmetric$ is indeed a well defined global density on $M$. 
As we obviously have $\body\nu_\mfdmetric\neq0$, it is a trivializing section.

\medskip

$\bullet$ (\ref{BatdenstrivwithmetricandOOD}):
The case of the density $\nu_{B,\mfdmetric}$ is completely similar; we only have to change \recalf{changeofBergfortrivdensity} to
\begin{align*}
\Ber\bigl(\Gextension\body\mfdmetric^{a}_{ij}(m)\bigr) 
&
= \Gextension\body\bigl(\Ber(T_m\varphi^{ba})^2\cdot \Ber\bigl(\Gextension\body\mfdmetric^{b}_{ij}(m)\bigr)
\\&
=
\Ber(T_m\varphi^{ba})^2\cdot \Ber\bigl(\Gextension\body\mfdmetric^{b}(m)\bigr)
\mapob,
\end{align*}
where the last equality follows from the fact that for two charts in a Batchelor atlas we have $\Ber(T_m\varphi^{ba}) = \Gextension\body\Ber(T_m\varphi^{ba})$ (as $T_m\varphi^{ba}$ does not depend upon the odd coordinates).
\end{preuve}

\begin{definition}{Remark}
We have seen in \recalt{Mobiusexampledensities} that the bundle of (\myquote{real}) densities $\Density(\Mob\oddif)$ on the (super) Möbius bundle is not trivial. 
As any $\CA$-manifold admits a metric (on the Möbius bundle we could take $\mfdmetric=\extder x\otimes \extder x + i\,\extder\xi\otimes \extder \xi$), this shows that the condition that $M$ is \ood{} in \recaltt{denstrivwithmetricandOOD}{canonicaltrivdensonOOD} is not superfluous.

\end{definition}

\begin{definition}{Notation\slash Definition}
For any subset $I\subset \{1, \dots, n\}$ we define $I^c\subset \{1, \dots, n\}$ as 
$$
I^c = \{1, \dots, n\}\setminus I
\mapob.
$$ 
Now let $I,J\subset \{1, \dots, n\}$ be any two disjoint subsets for which we write
$$
I = \{i_1, \dots, i_k\}
\qquad\text{and}\qquad
J = \{j_1, \dots, j_\ell\}
$$
with
$$
i_1<i_2< \cdots < i_k
\qquad\text{and}\qquad
j_1<\cdots < j_\ell
\mapob.
$$
We then define the integer $\varepsilon(I, J)$ as the number of transpositions needed to put the sequence 
$$
i_1, i_2, \dots, i_k, j_1, \dots, j_\ell
$$
in increasing order.

\end{definition}

\begin{proclaim}[superscalarprforHilbertvaluedfunctions]{Proposition}
Let $(E, \inprodsym_E, \supinprsym_E)$ be a proto super Hilbert over $\KK$ and let $M$ be an $\CA$-manifold of graded dimension $d\vert n$ equipped with a trivializing density $\nu$. 
Then the map $\supinprsym_{\nu} : C^\infty_c(M,E\otimes \CA^\KK) \times C^\infty_c(M,E\otimes \CA^\KK)\to \CC$ given by
$$
\homsuperinprod{\chi}{\psi}{_{\nu}} = \int_M \nu(m) \cdot \homsuperinprod[2]{\chi(m)}{\psi(m)}{_E}
$$
defines a super scalar product $\supinprsym_{\nu}$ on $C^\infty_c(M,E\otimes \CA^\KK)$, the space of compactly supported smooth functions on $M$ with values in $E\otimes\CA^\KK$. 
If $\supinprsym_E$ is homogeneous of parity $\alpha$, $\supinprsym_{\nu}$ is homogeneous of parity $\alpha+n$. 

\end{proclaim}

\begin{preuve}
That it defines a graded symmetric sesquilinear form is immediate. 
That it indeed takes values in $\CC$ is a direct consequence of the fact that Berezin integration of a $\CA^\KK$-valued density on $M$ yields an element of $\KK$. 
As $\nu$ is even, the conclusion concerning the homogeneity\slash parity of $\supinprsym_\nu$ follows immediately from the fact that the parity of $\int_M$ is $n \,\mathrm{mod}\,2$, see \recalt{parityofBerezinintegration}. 
Remains non-degeneracy. 

Let $\chi:M\to E\otimes \CA^\KK$ be a smooth function with compact support. If $\chi$ is not (identically) zero, there exists $m\in M$ such that $\chi(m)\neq0$. Now choose a local coordinate chart $m\in U\subset M$ with local coordinates $x_i,\xi_j$. 
On $U$ there exist smooth functions $\chi_I:\body U\to E$ such that $\chi$ is given by
$$
\chi(x,\xi) = \sum_{I\subset \{1, \dots, n\}} \xi^I \cdot \Gextension\chi_I(x)
\mapob.
$$
As there exists $m\in U$ with $\chi(m)\neq0$, not all functions $\chi_I$ can be identically zero, say $\chi_I$ (for some $I\subset \{1, \dots, n\}$ which we now consider fixed). 
It follows that there exists $m_o\in \body U$ with $0\neq \chi_I(m_o)\in E$. 
Because $\supinprsym_E$ is non-degenerate, there exists $v\in E$ such that $\homsuperinprod[2]{\chi_I(m_o)}{v}{_E}\neq0$. 
Taking a suitable multiple of $v$, we may assume $\bigl\vert\homsuperinprod[2]{\chi_I(m_o)}{v}{_E}\bigr\vert=1$. 

Now $\supinprsym_E$ is continuous, so the map $m\mapsto \homsuperinprod[2]{\chi_I(m)}{v}{_E}$ (on $\body U$) is smooth and non-zero at $m_o\in \body U$. 
We thus can choose a smooth function $\rho:\body U\to [0,1]\subset \RR$ with compact support such that $\rho(m_o)=1$ and such that $\rho(m)\neq0$ implies $\bigl\vert\homsuperinprod[2]{\chi_I(m_o)}{v}{_E}\bigr\vert >\tfrac12$. 
We then define the (smooth) function $\psi:M\to E\otimes \CA^\KK$ by
$$
\psi(m) = 
\begin{cases}
\displaystyle
\frac{(\Gextension\rho)(x)}{\bigl(\Gextension\homsuperinprod{\chi_I}{v}{_E}\bigr)(x)\cdot \nu(\partial_x\caprestricted_m,\partial_\xi\caprestricted_m)}
\cdot \conjugate^{\parity I}(v)\cdot \xi^{I^c}
&\quad\text{for } (x,\xi)\cong m\in U
\\[5\jot]
0
&\quad\text{otherwise.}
\end{cases}
$$
It follows immediately that we have (no longer writing the \Gextension{})
\begin{align*}
\homsuperinprod{\chi}{\psi}{_{\nu}}
&
=
\int_M \nu(m) \cdot \homsuperinprod[2]{\chi(m)}{\psi(m)}{_E}
\\&
=
\int_{ U} \extder(x,\xi)\ 
\nu(\partial_x\caprestricted_m,\partial_\xi\caprestricted_m)\cdot
\frac{\displaystyle\homsuperinprod[4]{\sum_{J\subset \{1, \dots, n\}}\chi_J(x)\,\xi^J}{\conjugate^{\parity I}(v)\, \xi^{I^c}}{_E} \cdot \rho(x)}{\homsuperinprod[2]{\chi_I(x)}{v}{_E}\cdot \nu(\partial_x\caprestricted_m,\partial_\xi\caprestricted_m)}
\\&
=
\int_{ U} \extder(x,\xi)\ 
\rho(x)\cdot \xi^I \cdot \xi^{I^c}
\\&
= (-1)^{\varepsilon(I,I^c)}\cdot
\int_{\body U} \extder\Leb^{(p)}(x) \ \rho(x) \neq0
\mapob,
\end{align*}
where for the third equality we used that $J=I$ is the only term contributing to the integral. 
This proves non-degeneracy. 
\end{preuve}

\begin{definition}[defofsuperscalarprodforfunctionsonoodM]{Definition}
Let $M$ be \ood{} and $\mfdmetric$ a metric on $M$. Then we define the super scalar products $\supinprsym_\mfdmetric$ and $\supinprsym_{B,\mfdmetric}$ on $C^\infty_c(M,E\otimes\CA^\KK)$ as those associated to the trivializing densities $\nu_\mfdmetric$ and $\nu_{B,\mfdmetric}$ (defined in \recalt{canonicaltrivdensonOOD}) according to \recalt{superscalarprforHilbertvaluedfunctions}. 

\end{definition}

\begin{proclaim}[descriptionsuperscalarproductBginadaptedBatlas]{Lemma}
Let $\atlas=\{\,U_a\mid a\in A\,\}$ be an \ood{} Batchelor atlas adapted to the super metric $\mfdmetric$ on $M$ \recalt{Batcheloratlasadaptedtometric} and let $\chi,\psi\in C^\infty_c(M;E\otimes \CA^\KK)$ be smooth functions. 
Let furthermore $\rho_a:\body U_a\to [0,1]$ be a smooth partition of unity subordinated to $\body \,\atlas$. Then the super scalar product $\supinprsym_{B,\mfdmetric}$ on $C^\infty_c(M;E\otimes \CA^\KK)$ is given by
\begin{align*}
\homsuperinprod \chi\psi{_{B,\mfdmetric}}
&
=
\sum_{a\in A} 
\int_{\body U_a} \extder\Leb^{(d)}(x^a)\ 
\rho_a(x^a)\cdot
\sqrt{\bigl\vert\det\bigl(\body\mfdmetric^{a,00}_{ij}(x^a)\bigr)\bigr\vert}
\\&
\kern5em
\cdot
\sum_{I\subset \{1, \dots, n\} } (-1)^{\varepsilon(I,I^c)}\, \homsuperinprod[3]{\conjugate^{\parity I}\bigl(\chi_{a,I}(x^a)\bigr)}{\conjugate^n\bigl(\psi_{a,I^c}(x^a)\bigr)}{_E}
\mapob,
\end{align*}
where $(x^a, \xi^a)$ denotes a local system of coordinates on $U_a$ and where the local expression of $\chi$ (respectively $\psi$) is given by
$$
\chi(x^a,\xi^a)
=
\sum_{I\subset \{1, \dots, n\} } (\xi^a)^I\,\Gextension\chi_{a,I}(x^a)
$$
for smooth functions $\chi_{a,I}:\body U_a\to E$.

\end{proclaim}

\begin{preuve}
As the atlas is adapted to the super metric, the matrices $-i\,\body \mfdmetric^{a,11}$ are constant the identity. 
This implies that the canonical basis $e_1, \dots, e_n$ of the typical fiber of $V\!M$ is an orthonormal basis. 
We thus have
\begin{align*}
\homsuperinprod \chi\psi{_{B,\mfdmetric}}
&
=
\int_M \nu_{B,\mfdmetric}(m) \cdot \homsuperinprod[2]{\chi(m)}{\psi(m)}{_E}
\\&
\oversetalign{\recalt{integrationofdensityonsupermanifold}}\to=
\quad
\sum_{a\in A} 
\int_{ U_a} \extder(x^{a},\xi^{a})\  
(\Gextension\rho_a)(x^a)\cdot 
\nu_{B,\mfdmetric}\Bigl( \fracp{}{x^a}\bigrestricted_{(x^a,\xi^a)}, \fracp{}{\xi^a}\bigrestricted_{(x^a,\xi^a)} \Bigr)
\\&
\kern6em
\cdot \sum_{I,J\subset \{1, \dots, n\} }
\homsuperinprod[2]{(\xi^a)^I\,\chi_{a,I}(x^a)}{(\xi^A)^J\,\psi_{a,J}(x^a)}{_E}
\\&
\oversetalign{\recalt{canonicaltrivdensonOOD}}\to=
\quad
\sum_{a\in A} 
\int_{ U_a} \extder(x^{a}, \xi^{a})\  
(\Gextension\rho_a)(x^a)\cdot 
{\sqrt{\bigl\vert\Ber \bigl( \Gextension\body\mfdmetric_{ij}^a(x^a) \bigr)\bigr\vert}} 
\\&
\kern6em
\cdot \sum_{I,J\subset \{1, \dots, n\} }(\xi^a)^I\, (\xi^a)^J\,
\homsuperinprod[2]{\conjugate^{\parity I}\,\chi_{a,I}(x^a)}{\conjugate^{I+ J}\,\psi_{a,J}(x^a)}{_E}
\\
&
=
\sum_{a\in A} 
\int_{\body U_a} \extder \Leb^{(d)}(x^a)\ 
\rho_a(x^a)\cdot 
{\sqrt{\bigl\vert\det \bigl( \body\mfdmetric_{ij}^{a,00}(x^a) \bigr)\bigr\vert}} 
\\&
\kern7em
\cdot \sum_{I\subset \{1, \dots, n\} }(-1)^{\varepsilon(I,I^c)}
\homsuperinprod[2]{\conjugate^{\parity I}\,\chi_{a,I}(x^a)}{\conjugate^{n}\,\psi_{a,I^c}(x^a)}{_E}
\mapob,
\end{align*}
where for the third equality we first used the extension of $\supinprsym_E$ to $E\otimes \CA^\KK$ by right-linearity \recaltt{extensiontoAvectorspaceofgradedsesquilinearform}{superscalarbackandforthbodyextension} to extract the powers of $\xi^a$, giving us the conjugations on $\chi_{a,I}$ and $\psi_{a,J}$, and then that $\supinprsym_E$ takes values in $\CC$, which is even, allowing to put the powers of $\xi^a$ to the left. 
For the fourth equality we used $-i\,\body \mfdmetric^{a,11}=\oneasmatrix$ to simplify the Berezinian and we applied the definition of Berezin integration \recalt{defofelementaryBerezinintegration}. 
\end{preuve}

\begin{proclaim}[superscalarBgiscontinuous]{Proposition}
The super scalar product $\supinprsym_{B,\mfdmetric}$ \recalt{defofsuperscalarprodforfunctionsonoodM}, \recalt{superscalarprforHilbertvaluedfunctions} is continuous with respect to the (ordinary) metric $\inprodsym_\mfdmetric$ \recalf{firstdefofordinarymetriconsuperfunctionsinE}, \recalt{metriconsuperfunctionsinEinlocalcoordinates}.

\end{proclaim}

\begin{preuve}
We use the descriptions un terms of an adapted Batchelor atlas given in \recalt{metriconsuperfunctionsinEinlocalcoordinates} and \recalt{descriptionsuperscalarproductBginadaptedBatlas}. 
We first note that according to \recalt{metriconsuperfunctionsinEinlocalcoordinates} the norm-squared of an element $\chi\in C^\infty_c(M;E\otimes \CA^\KK)$ is given by
\begin{align*}
\norm \chi^2_\mfdmetric
&
=
\sum_{a\in A} 
\int_{\body U_a} \extder\Leb^{(d)}(x^a)\ 
\rho_a(x^a)\cdot
\sqrt{\bigl\vert\det\bigl(\body\mfdmetric^{a,00}(x^a)\bigr)\bigr\vert}
\cdot
\sum_{I\subset \{1, \dots, n\} } \norm{\chi_{a,I}(x^a)}^2_E 
\mapob,
\end{align*}
As $E$ is a proto super Hilbert space, we have $\inprod{E_0}{E_1}_E=0$, which implies in particular that we have (for $\chi,\psi\in E$) the equality
$$
\inprod{\conjugate \chi}{\conjugate\psi}_E \equiv 
\inprod{\chi_0-\chi_1}{\psi_0-\psi_1}_E
=
\inprod{\chi_0}{\psi_0}_E + \inprod{\chi_1}{\psi_1}_E
=
\inprod{\chi}{\psi}_E
\mapob,
$$
and thus in particular $\norm{\conjugate\chi}_E = \norm\chi_E$. 
Moreover, $\supinprsym_E$ is continuous with respect to the topology induced by the norm, which means that there exists a constant $C\ge0$ such that
$$
\vert \homsuperinprod\chi\psi{_E} \vert \le
C\cdot\norm\chi_E \cdot\norm\psi_E
\mapob.
$$

With these preparations we now turn our attention to $\vert\homsuperinprod \chi\psi{_{B,\mfdmetric}}\vert$ (with $\chi,\psi\in C^\infty_c(M;E\otimes \CA^\KK)$), for which we obtain, using \recalt{descriptionsuperscalarproductBginadaptedBatlas}, the following estimates:
\begin{align*}
\vert 
\homsuperinprod \chi\psi{_{B,\mfdmetric}}
\vert
&
\le
\sum_{a\in A} 
\int_{\body U_a} \extder \Leb^{(d)}(x^a)\ 
\rho_a(x^a)\cdot 
{\sqrt{\bigl\vert\det \bigl( \body\mfdmetric_{ij}^{a,00}(x^a) \bigr)\bigr\vert}} 
\\&
\kern7em
\cdot \sum_{I\subset \{1, \dots, n\} }
\bigl\vert
\homsuperinprod[2]{\conjugate^I\,\chi_{a,I}(x^a)}{\conjugate^{n}\,\psi_{a,I^c}(x^a)}{_E}
\bigr\vert
\\&
\le
\sum_{a\in A} 
\int_{\body U_a} \extder \Leb^{(d)}(x^a)\ 
\rho_a(x^a)\cdot 
{\sqrt{\bigl\vert\det \bigl( \body\mfdmetric_{ij}^{a,00}(x^a) \bigr)\bigr\vert}} 
\\&
\kern7em
\cdot \sum_{I\subset \{1, \dots, n\} }
C\cdot\norm{\conjugate^I\,\chi_{a,I}(x^a)}_E \cdot\norm{\conjugate^{n}\,\psi_{a,I^c}(x^a)}_E
\\&
=
\sum_{a\in A} 
\int_{\body U_a} \extder \Leb^{(d)}(x^a)\ 
\rho_a(x^a)\cdot 
{\sqrt{\bigl\vert\det \bigl( \body\mfdmetric_{ij}^{a,00}(x^a) \bigr)\bigr\vert}} 
\\&
\kern7em
\cdot \sum_{I\subset \{1, \dots, n\} }
C\cdot\norm{\chi_{a,I}(x^a)}_E \cdot\norm{\psi_{a,I^c}(x^a)}_E
\mapob.
\end{align*}
We now use the standard argument for the Cauchy-Schwarz inequality, \ie, the positivity of the quadratic function $g(t)$ defined as
\begin{align*}
g(t) 
&
= 
\sum_{a\in A} 
\int_{\body U_a} \extder \Leb^{(d)}(x^a)\ 
\rho_a(x^a)\cdot 
{\sqrt{\bigl\vert\det \bigl( \body\mfdmetric_{ij}^{a,00}(x^a) \bigr)\bigr\vert}} 
\\&
\kern7em
\cdot \sum_{I\subset \{1, \dots, n\} }
\bigl(\,\norm{\chi_{a,I}(x^a)}_E +t\cdot\norm{\psi_{a,I^c}(x^a)}_E\bigr)^2
\\&
=
\norm\chi^2_\mfdmetric + t^2\cdot \norm\psi^2_\mfdmetric
+2t
\sum_{a\in A} 
\int_{\body U_a} \extder \Leb^{(d)}(x^a)\ 
\rho_a(x^a)\cdot 
{\sqrt{\bigl\vert\det \bigl( \body\mfdmetric_{ij}^{a,00}(x^a) \bigr)\bigr\vert}} 
\\&
\kern7em
\cdot \sum_{I\subset \{1, \dots, n\} }
\norm{\chi_{a,I}(x^a)}_E \cdot\norm{\psi_{a,I^c}(x^a)}_E
\mapob,
\end{align*}
to finally obtain the estimate
$$
\vert 
\homsuperinprod \chi\psi{_{B,\mfdmetric}}
\vert
\le
C\cdot \norm\chi_\mfdmetric \cdot \norm\psi_\mfdmetric
\mapob,
$$
which shows that $\supinprsym_{B,\mfdmetric}$ is continuous with respect to the topology on $C^\infty_c(M;E\otimes \CA^\KK)$ induced by the metric $\inprodsym_\mfdmetric$. 
\end{preuve}

\begin{proclaim}{Corollary}
If $(E, \inprodsym_E, \supinprsym_E)$ is a proto super Hilbert space (over $\KK$), then $C^\infty_c(M; E\otimes \CA^\KK)$ equipped with the metric $\inprodsym_\mfdmetric$ and the super scalar product $\supinprsym_{B,\mfdmetric}$ also is a proto super Hilbert space.

\end{proclaim}

\masection{Integrating along fibers}
\label{Integrationalongfiberssection}

It frequently happens that we have a function $f$ on a product space $E\times F$ and that we want to integrate only over the $F$-part in such a way that the result is a function on $E$. 
In the ordinary context one usually speaks about this as integrals depending upon parameters and several results as to when the resulting function on $E$ is continuous or smooth are available. 
In the super context we have to be a bit careful because Berezin integration breaks down a smooth function of even and odd variables into smooth functions of real variables only (that then are integrated). 
Without care, we thus would not get a (smooth) function of the remaining even and odd variables on $E$. 
The purpose of this section is to provide a rigorous framework for such integrals depending upon even and odd parameters and to show that we have a change-of-variables formula for such integrals. 
We then use it to show that if a super metric $\mfdmetric$ on an $\CA$-manifold is invariant under the action of an $\CA$-Lie group, then the super scalar product $\supinprsym_{\mfdmetric}$ is invariant under the induced group action on the space of functions.

\begin{definition}{Discussion}
Let $E$ be a finite dimensional graded vector space of dimension $p\vert q$, let  $F$ be a finite dimensional graded vector space of dimension $d\vert n$ and let $U\subset E_0\times F_0$ be an open set. 
When we use (even and odd) coordinates $(s,\sigma)$ on $E_0$ and $(x,\xi)$ on $F_0$, any smooth function $f:U\to \CA^\KK$ is determined by $2^q\times 2^n=2^{q+n}$  ordinary smooth functions $f_{I,J}:\body U\to \KK$ according to (note the order of the products of coordinates $\xi$ and $\sigma$)
$$
f(s,\sigma,x,\xi) = \sum_{I\subset \{1, \dots, n\}}\sum_{J\subset \{1, \dots, q\}} \xi^I\,\sigma^J\,\bigl(\Gextension(f_{I,J})\bigr)(s,x)
\mapob.
$$
Associated to $U$ we define the subset $U_{(s,\sigma)}\subset E_0$ by
\begin{align*}
U_{(s,\sigma)}
&
= 
\{\,(x,\xi)\in E_0\mid (s,\sigma,x,\xi)\in U\,\}
\mapob.
\end{align*}
We now assume that the function $f$ is such that its restriction to $U_{(s,\sigma)}$ has compact support for all $(s,\sigma)$.  
With these ingredients we then define the function $\fint:\pi_E(U)\to \CA^\KK$ (where $\pi_E:E_0\times F_0 \to E_0$ denotes the canonical projection) as integration of $f$ over the fibers $F_0$ of $U$ by
\begin{align*}
\fint(s,\sigma) 
&
= 
\int_{U_{(s,\sigma)}} \extder( x,\xi)\ 
f(s,\sigma,x,\xi)
\\&
\oversetalign{def.}\to=
\sum_{J\subset \{1, \dots, q\}} \sigma^J\int_{\body U_{(s,\sigma)}}\extder\Leb^{(d)}(x)\ f_{\nasset,J}(s,x)
\mapob.
\end{align*}
A slightly better but less clear definition would be the following (which clarifies some obscure points in the above formula). For $s\in \body E_0\cong \RR^p$ we define the set $\body U_s = \{\,x\in \body F_0\cong \RR^d\mid (s,x)\in \body U\,\}$ (which is the same as $\body U_{(s,\sigma)}$, as the map $\body$ discards any dependence on nilpotent and thus odd coordinates). 
We then define the (ordinary!) smooth functions $\fint\kern-1pt_J:\body \pi_E(U)\to \KK$ by
$$
\fint\kern-1pt_J(s) = \int_{\body U_s} f_{\nasset,J}(s,x) \ \extder\Leb^{(d)}(x)
\mapob.
$$
That these functions are smooth is because of our assumption that the restriction of $f$ to all $U_{(s,\sigma)}$ has compact support (and of course because the functions $f_{I,J}$ are smooth). 
With these functions we then define the (super) smooth function $\fint$ by
$$
\fint(s,\sigma) = \sum_{J\subset \{1, \dots, q\}} \sigma^J\, \bigl(\Gextension( F_J)\bigr)(s)
\mapob.
$$
We now want to show that this partial integration (not to be confused with integration by parts) or fiberwise integration has a nice transformation property under a special kind of (local) diffeomorphisms of $E_0\times F_0$: those that preserve the \myquote{fibers.}

\end{definition}

\begin{definition}{Definition}
Let $U,V\subset E_0\times F_0$ be two open sets and let $\Phi:U\to V$ be a (local) diffeomorphism. 
Using the product structure $E_0\times F_0$ we thus can \myquote{decompose} $\Phi$ into two components $\Phi_E:U\to E_0$ and $\Phi_F:U\to F_0$ according to
$$
\Phi(u,v) = \bigl(\Phi_E(u,v),\Phi_F(u,v)\bigr)
\mapob.
$$
We will say that $\Phi$ is a \stresd{fiber-diffeomorphism} if $\Phi_E$ is independent of $v$. 
A more sophisticated way would be to say that $\Phi$ is a fiber-diffeomorphism if there exists a diffeomorphism $\Phi_E:\pi_E(U) \to \pi_E(V)$ such that the following diagram is commutative:
\begin{moneq}[fiberdiffeoascommutativediagram]
\begin{CD}
U @>\Phi>> V
\\ 
@V{\pi_E}VV @VV{\pi_E}V
\\
\pi_E(U) @>>\Phi_E> \pi_E(V)\rlap{\mapob.}
\end{CD}
\end{moneq}
When we use coordinates (even and odd) on $E_0\times F_0$ and when we denote those on $U$ by $\bigl((s,\sigma),(x,\xi)\bigr)\cong (s,\sigma,x,\xi)$ ($(s,\sigma)$ being coordinates on $E_0$ and $(x,\xi)$ on $F_0$) and those on $V$ by $(t,\tau,y,\eta)$ (they are the same coordinates, but we use different symbols to differentiate between source and target space), then a fiber-diffeomorphism has the form
$$
\Phi(s,\sigma,x,\xi) = (t,\tau,y,\eta)
\mapob,
$$
with
$$
t=t(s,\sigma)
\quad,\quad
\tau = \tau(s,\sigma)
\qquad,\qquad
y=y(s,\sigma,x,\xi)
\quad,\quad
\eta = \eta(s,\sigma,x,\xi)
\mapob.
$$

\end{definition}

If $\Phi:U\to V$ is a fiber-diffeomorphism, the restriction $\Phi_{(s,\sigma)}$ of $\Phi$ to a slice $U_{(s,\sigma)}$ is a \myquote{diffeomorphism} $\Phi_{(s,\sigma)}:U_{(s,\sigma)}\to V_{(t,\tau)}$, where $(t,\tau)\in E_0$ is determined by $(t,\tau) = \Phi_E(s,\sigma)$. 
Of course it is not a real diffeomorphism, it only is a homeomorphism with inverse $(\Phi\mo)_{(t,\tau)}$, the restriction of $\Phi\mo$ to the slice $V_{(t,\tau)}$. 
However, in many ways it behaves\slash can behave like a diffeomorphism. 
In particular we can attribute a tangent map to it (this will be a particular case of a generalized tangent map as described in \cite[V.3.19]{Tu04}). 
If $\Phi:U\to V$ is a fiber-diffeomorphism, then its tangent map is given by (see \recalf{matrixrepoftangentmap} and \recalf{changeofbasisfortangentvectorsbyjacobian})
\begin{moneq}
T\Phi(
\begin{pmatrix} \dfracp{}g 
\\[4\jot] 
\dfracp{}\chi 
\\[4\jot] 
\dfracp{}x 
\\[4\jot] 
\dfracp{}\xi 
\end{pmatrix}
)
=
\begin{pmatrix}
\dfracp tg & \dfracp \tau g & \dfracp yg & \dfracp \eta g
\\[4\jot]
\dfracp t\chi & \dfracp \tau\chi & \dfracp y\chi & \dfracp \eta\chi
\\[4\jot]
0 & 0 & \dfracp yx & \dfracp \eta x
\\[4\jot]
0 & 0 & \dfracp y\xi & \dfracp \eta\xi
\end{pmatrix}
\cdot
\begin{pmatrix} \dfracp{}t 
\\[4\jot] 
\dfracp{}\tau 
\\[4\jot] 
\dfracp{}y 
\\[4\jot] 
\dfracp{}\eta 
\end{pmatrix}
\mapob.
\end{moneq}
Obviously the matrix elements are smooth functions on $U$. 
According to \cite[V.3.19]{Tu04} the tangent map of $\Phi_{(s,\sigma)}$ is described by the matrix
\begin{moneq}
T\Phi_{(s,\sigma)}(
\begin{pmatrix} 
\dfracp{}x 
\\[4\jot] 
\dfracp{}\xi 
\end{pmatrix}
)
=
\begin{pmatrix}
\dfracp yx & \dfracp \eta x
\\[4\jot]
\dfracp y\xi & \dfracp \eta\xi
\end{pmatrix}
\cdot
\begin{pmatrix}
\dfracp{}y 
\\[4\jot] 
\dfracp{}\eta 
\end{pmatrix}
\mapob.
\end{moneq}
For fixed $(s,\sigma)$ the matrix elements of $\Phi_{(s,\sigma)}$ are not (necessarily) smooth functions, but they are smooth in the complete set of variables $(s,\sigma,x,\xi)$.

\begin{proclaim}[compositionfiberdiffeos]{Lemma}
Let $U,V,W\subset E_0\times F_0$ be three open sets and $\Phi:U\to V$ and $\Psi:V\to W$ two fiber-diffeomorphisms. 
Then $\Psi\scirc \Phi:U\to W$ is a fiber diffeomorphism and we have the following properties concerning the generalized tangent maps:
$$
T(\Psi\scirc\Phi)_{(s,\sigma)} = T\Psi_{\Phi_E(s,\sigma)} \scirc T\Phi_{(s,\sigma)}
$$
and
$$
\Ber_{(\pi)}\bigl(T\Psi\scirc\Phi)_{(s,\sigma)}\bigr) = \Ber_{(\pi)}\bigl(T\Psi_{\Phi_E(s,\sigma)}\bigr) \cdot \Ber_{(\pi)}\bigl(T\Phi_{(s,\sigma)}\bigr)
\mapob.
$$

\end{proclaim}

\begin{preuve}
It is immediate that the composition of two fiber-diffeomorphisms again is a fiber-diffeomorphism. 
When we denote the coordinates on $W$ by $(g, \chi, z,\zeta)$ (still the same, just changing the names again), we have for the full tangent maps
$$
T(\Psi\scirc \Phi) = (T\Psi) \scirc (T\Phi)
\mapob,
$$
which in terms of the matrix representation gives
\begin{moneq}
T(\Psi\scirc\Phi)(
\begin{pmatrix} \dfracp{}s 
\\[4\jot] 
\dfracp{}\sigma 
\\[4\jot] 
\dfracp{}x 
\\[4\jot] 
\dfracp{}\xi 
\end{pmatrix}
)
=
\begin{pmatrix}
\dfracp ts & \dfracp \tau s & \dfracp ys & \dfracp \eta s
\\[4\jot]
\dfracp t\sigma & \dfracp \tau\sigma & \dfracp y\sigma & \dfracp \eta\sigma
\\[4\jot]
0 & 0 & \dfracp yx & \dfracp \eta x
\\[4\jot]
0 & 0 & \dfracp y\xi & \dfracp \eta\xi
\end{pmatrix}
\cdot
\begin{pmatrix}
\dfracp gt & \dfracp \chi t & \dfracp zt & \dfracp \zeta t
\\[4\jot]
\dfracp g\tau & \dfracp \chi\tau & \dfracp z\tau & \dfracp \zeta\tau
\\[4\jot]
0 & 0 & \dfracp zy & \dfracp \zeta y
\\[4\jot]
0 & 0 & \dfracp z\eta  & \dfracp \zeta\eta
\end{pmatrix}
\cdot
\begin{pmatrix} \dfracp{}g
\\[4\jot] 
\dfracp{}\chi
\\[4\jot] 
\dfracp{}z
\\[4\jot] 
\dfracp{}\zeta 
\end{pmatrix}
\mapob.
\end{moneq}
It follows immediately that we have
$$
\begin{pmatrix}
\dfracp zx & \dfracp \zeta x
\\[4\jot]
\dfracp z\xi & \dfracp \zeta\xi
\end{pmatrix}
=
\begin{pmatrix}
\dfracp yx & \dfracp \eta x
\\[4\jot]
\dfracp y\xi & \dfracp \eta\xi
\end{pmatrix}
\cdot
\begin{pmatrix}
\dfracp zy & \dfracp \zeta y
\\[4\jot]
\dfracp z\eta & \dfracp \zeta\eta
\end{pmatrix}
\mapob.
$$
This, together with the multiplicative property of the Berezinian $\Ber_{(\pi)}$, proves the second part of the lemma. 
\end{preuve}

\begin{proclaim}[thefiberedchangeofvariablestheorem]{Proposition}
Let $\Phi:U\to V$ be a fiber-diffeomorphism with
$$
\Phi(s,\sigma,x,\xi)
=
\bigl(\Phi_E(s,\sigma),\Phi_F(s,\sigma,x,\xi)\bigr)
\equiv
(t,\tau,y,\eta)
$$
and let $f:V\to \CA^\KK$ be a smooth function with compact support on all slices $V_{(t,\tau)}$. 
Then we have the equality
\begin{align}
\notag
\shifttag{6em}
\int_{V_{(t,\tau)}}\kern-1em \extder(y , \eta)\ 
f(t,\tau,y,\eta)
\equiv
\int_{V_{(t,\tau)}}\kern-1em \extder(y , \eta )\ 
f\bigl(\Phi_E(s,\sigma),y,\eta\bigr)
\\&
=
\int_{U_{(s,\sigma)}} \kern-1em \extder( x, \xi)\ 
f\bigl(\Phi(s,\sigma,x,\xi)\bigr) \cdot \piBer(T_{(x,\xi)}\Phi_{(s,\sigma)})
\mapob.
\label{thechangeofvariablesformula}
\end{align}

\end{proclaim}

\begin{preuve}
We follow the scheme of the proof as given in \cite[Thm 4.6.1]{Va04}. 
It is an elementary consequence of \recalt{compositionfiberdiffeos} that we can decompose $\Phi$ into elementary parts and prove this property for each part separately. 

\medskip

\noindent
$\bullet$
Our first \myquote{simplification} is to get rid of the change of coordinates $(s,\sigma)\mapsto (t,\tau)$. 
To do so, we note that the map $\Phi_E:\pi_E(U)\to \pi_E(V)$, $(s,\sigma)\mapsto (t,\tau)=\Phi_E(s,\sigma)$ is a diffeomorphism. 
It follows that we can define the fiber-diffeomorphisms $\Phi_1$ by
$$
\Phi_1(s,\sigma,x,\xi) = \bigl(\Phi_E(s,\sigma), (x,\xi)\bigr) = (t,\tau,x,\xi)
$$
and the fiber-diffeomorphism $\Phi_2 = \Phi\scirc \Phi_1\mo$ by
$$
\Phi_2(t,\tau,x,\xi) = \Bigl(t,\tau, \Phi_F\bigl(\Phi_E\mo(t,\tau), x,\xi\bigr)\,\Bigr)
\mapob.
$$
Said differently, we have decomposed the map $(s,\sigma,x,\xi)\mapsto (t,\tau,y,\eta)$ as
$$
(s,\sigma,x,\xi)\mapsto (t,\tau,x,\xi)\mapsto (t,\tau,y,\eta)
\mapob.
$$
As the integration does not concern the coordinates $(s,\sigma)$, the result is true for $\Phi_1$. 
We thus may assume that $\Phi$ is of the form
$$
\Phi(s,\sigma,x,\xi) = \bigl(s,\sigma, \Phi_F(s,\sigma,x,\xi)\bigr) \equiv (s,\sigma,y,\eta)
\mapob.
$$
Said differently, we may assume that the diffeomorphism $\Phi_E:\pi_E(U)\to \pi_E(V)$ in \recalf{fiberdiffeoascommutativediagram} is the identity. 
Nevertheless, the open sets $U$ and $V$ need not be the same and even when they are, we will keep denoting them differently in order to distinguish which set of coordinates we are using. 

\medskip

\noindent
$\bullet$
The next step is to dissociate the odd and even variables $\xi$ and $x$ as follows. 
According to the inverse function theorem a map $\Psi$ between open subsets of $F_0$ is a local diffeomorphism around a point $v$ if and only if the matrices $\body (\partial y/\partial x)(v)$ and $\body (\partial \eta/\partial\xi)(v)$ are invertible. 
When we apply this to the map
$$
\Phi(s,\sigma,x,\xi) = (s,\sigma,y,\eta) = \bigl(s,\sigma,\Phi_y(s,\sigma,x,\xi), \Phi_\eta(s,\sigma,x,\xi)\bigr)
\mapob,
$$
it follows easily (see also \cite[Cor. 4.4.2]{Va04}) that the maps
$$
\Phi_1(s,\sigma,x,\xi) = (s,\sigma,x,\eta) \equiv \bigl(s,\sigma,x,\Phi_\eta(s,\sigma,x,\xi)\bigr)
$$
and
$$
\Phi_2 = \Phi\scirc \Phi_1\mo
\qquad,\qquad
\Phi_2(s,\sigma,x,\eta) = (s,\sigma,y,\eta)
$$
are fiber-diffeomorphisms. 
The conclusion is that we may assume that either only the $x$ variables change or only the $\xi$ variables. 

\medskip

\noindent
$\bullet$
We first assume that only the $\xi$ variables change, \ie, that $\Phi$ is of the form
$$
\Phi(s,\sigma,x,\xi) = \bigl(s,\sigma,x,\eta(s,\sigma,x,\xi)\bigr)
\mapob.
$$
(Note that this implies that we have $U=V$ because open sets are saturated with respect to nilpotent elements or, said differently, open sets $U$ satisfy $U=\body\mo\body U$. But as said, we will keep the names $U$ and $V$ for source and target space to indicate which set of coordinates we are using.)  
As this $\Phi$ is a diffeomorphism, the matrix of its tangent map (the Jacobian) is invertible, which implies in particular that the matrix of size $n\times n$
$$
\body\fracp{\eta_i}{\xi_j}(s,\sigma,x,\xi)
$$
is invertible. 
Another way to obtain this matrix is to note that the functions $\eta_i(s,\sigma,x,\xi)$ are smooth, which implies that there are $2^{q+n}$ ordinary smooth functions $\eta_{i,I,J}(s,x)$ such that
$$
\eta_i(s,\sigma,x,\xi)
=
\sum_{I\subset \{1, \dots, n\}} \sum_{J\subset\{1, \dots, q\}} \xi^I\,\sigma^J\, \bigl(\Gextension(\eta_{i,I,J})\bigr)(s,x)
\mapob.
$$
As $\eta_i$ is odd, the functions $\eta_{i,I,J}$ with $\parity I + \parity J$ even are identically zero. 
And in this vision, the matrix element $\body\partial \eta_i/\partial\xi_j$ is the function $\eta_{i,\{j\},\emptyset}(s,x)$. 
Using this invertible matrix, we define the fiber-diffeomorphism $\Phi_1$ by the formula
$$
\Phi_1(s,\sigma,x,\xi) = (s,\sigma,x,\zeta)
\qquad,\qquad
\zeta_i = \sum_{j=1}^n \xi_j\,\bigl(\Gextension(\eta_{i,\{j\},\emptyset})\bigr)(s,x)
\mapob.
$$
The idea behind this fiber-diffeomorphism is that $\Phi_1$ is simply linear in the odd variables and that the remaining fiber-diffeomorphism $\Phi_2 = \Phi\scirc \Phi_1\mo$ is of the form
\begin{align}
\notag
\eta_i(s,\sigma,x,\xi) 
&
= \xi_i + \text{terms with a single $\sigma$ variable}
\\&
\kern3em
+\text{terms with at least three odd variables $\sigma$ and\slash or $\xi$}
\mapob.
\label{changeofvarinoddwithidentityatlevetl1}
\end{align}
With this decomposition $\Phi=\Phi_2\scirc \Phi_1$ we finally can attack the proof of our change-of-variables formula \recalf{thechangeofvariablesformula}. 

\medskip

\noindent
$\bullet$
So suppose $\Phi$ is simply linear in the odd variables, \ie, of the form
$$
\Phi(s,\sigma,x,\xi) = (s,\sigma,x,\eta)
\qquad,\qquad
\eta_i = \sum_{j=1}^n \xi_j\,\bigl(\Gextension(\Phi_{ji})\bigr)(s,x)
\mapob.
$$
It then is immediate that we have
$$
\piBer(T_{(x,\xi)}\Phi_{(s,\sigma)}) = \Bigl(\det\bigl(\Gextension(\Phi_{ji})\bigr)(s,x)\Bigr)\mo
\mapob,
$$
a function that is independent of the odd variables. 

To see that \recalf{thechangeofvariablesformula} holds, consider a function $f:V\to \CA^\KK$ with
$$
f(s,\sigma,x,\eta) 
= 
\sum_{I\subset \{1, \dots, n\}}\sum_{J\subset \{1, \dots, q\}} \eta^I\,\sigma^J\,\bigl(\Gextension(f_{I,J})\bigr)(s,x)
\mapob.
$$
According to the definition we thus have (using that $U=V$ and thus in particular $\body V_{(s,\sigma)} = \body U_{(s,\sigma)}$) 
\begin{moneq}[firstequalitylinearoddchangeofvariables]
\int_{V_{(s,\sigma)}}\kern-1em \extder(x , \eta )\ 
f(s, \sigma, x,\eta)
=
\sum_{J\subset \{1, \dots, q\}} \sigma^J\int_{\body U_{(s,\sigma)}}\extder\Leb^{(d)}(x)\ f_{\nasset,J}(s,x)
\mapob.
\end{moneq}
On the other hand, the function $f\bigl(\Phi(s,\sigma,x,\xi)\bigr) \cdot \piBer(T_{(x,\xi)}\Phi_{(s,\sigma)})$ is described by
\begin{align}
\notag
\shifttag{5em}
f\bigl(\Phi(s,\sigma,x,\xi)\bigr) \cdot \piBer(T_{(x,\xi)}\Phi_{(s,\sigma)})
=
\\&
\sum_{I\subset \{1, \dots, n\}}\sum_{J\subset \{1, \dots, q\}} (\xi_j\Phi_{ji})^I\,\sigma^J\,\Bigl(\Gextension(f_{I,J}\cdot\det\bigl(\Phi_{ji})\mo\bigr)\Bigr)(s,x)
\mapob.
\label{fcomposedwithPhiandBerforintegration}
\end{align}
To compute
$$
\int_{U_{(s,\sigma)}} \kern-1em \extder( x ,\xi )\ 
f\bigl(\Phi(s,\sigma,x,\xi)\bigr) \cdot \piBer(T_{(x,\xi)}\Phi_{(s,\sigma)})
$$
we thus have to determine in \recalf{fcomposedwithPhiandBerforintegration} those terms with a factor $\xi^\nasset$. 
This will happen only when $I$ contains $n$ elements, \ie, when $I=\{1, \dots, n\}$. And in that case we have
$$
(\xi_j\Phi_{ji})^I
=
\det\bigl(\Phi_{ji}(s,x)\bigr)\cdot\xi^\nasset
\mapob.
$$
It follows that we have
\begin{align*}
\shifttag{3.5em}
\int_{U_{(s,\sigma)}} \kern-1em \extder( x , \xi )\ 
f\bigl(\Phi(s,\sigma,x,\xi)\bigr) \cdot \piBer(T_{(x,\xi)}\Phi_{(s,\sigma)})
\\&
=
\sum_{J\subset \{1, \dots, q\}} \sigma^J\int_{\body U_{(s,\sigma)}}\extder\Leb^{(d)}(x)\ f_{\nasset,J}(s,x)\cdot \det\bigl(\Phi_{ji}(s,x)\bigr) \cdot \det\bigl(\Phi_{ji}(s,x)\bigr)\mo
\\&
\oversetalign{\recalf{firstequalitylinearoddchangeofvariables}}\to=
\kern1em
\int_{V_{(s,\sigma)}}\kern-1em \extder( x , \eta )\ 
f(s,\sigma,x,\eta)
\mapob,
\end{align*}
proving \recalf{thechangeofvariablesformula} in this case. 

\medskip

\noindent
$\bullet$
Now suppose $\Phi$ is of the form \recalf{changeofvarinoddwithidentityatlevetl1}, which is equivalent to saying that we have $\body\partial \eta_i/\partial\xi_j = \delta_{ij}$. We then can decompose $\Phi$ as a product of $n$ fiber-diffeomorphisms, each of which changes only one odd coordinate. 
More precisely, we decompose the map $(\xi_1, \dots, \xi_n)\mapsto (\eta_1, \dots, \eta_n)$ as
\begin{multline*}
\quad
(\xi_1, \dots, \xi_n) \mapsto (\eta_1, \xi_2, \dots, \xi_n) \mapsto
(\eta_1,\eta_2, \xi_3, \dots, \xi_n) \mapsto 
\\
\cdots \mapsto (\eta_1, \dots, \eta_{n-1}, \xi_n) \mapsto (\eta_1, \dots, \eta_n)
\mapob.
\quad
\end{multline*}
And even more precisely, knowing that the $\eta_i$ are functions of all coordinates $(s,\sigma,x,\xi)$, we define the fiber-diffeomorphism $\Phi_1$ by
$$
\Phi_1(s,\sigma,x,\xi) = \bigl(s,\sigma,x,(\eta_1(s,\sigma,x,\xi), \xi_2, \dots, \xi_n)\bigr)
\mapob.
$$
As before, it is an immediate consequence of the inverse function theorem (because $\body\partial \eta_i/\partial\xi_j = \delta_{ij}$) and the fact that it is injective at the body level that $\Phi_1$ is a fiber-diffeomorphism. 
Its inverse thus is of the form
$$
\Phi_1\mo(s,\sigma,x,\zeta) = \bigl(s,\sigma,x,(\psi_1(s,\sigma,x,\zeta), \zeta_2, \dots, \zeta_n)\bigr)
\mapob.
$$
We then can define the fiber-diffeomorphism $\Phi_2$ by
$$
\Phi_2(s,\sigma,x,\zeta) = \Bigl( s,\sigma,x,\bigl(\zeta_1, \eta_2\bigl(s,\sigma,x,(\psi_1(s,\sigma,x,\zeta), \zeta_2, \dots, \zeta_n)\bigr), \zeta_3, \dots, \zeta_n\bigr)\Bigr)
\mapob,
$$
which is a complicated way to write the map $(\eta_1, \xi_2, \dots, \xi_2)\mapsto (\eta_1, \eta_2, \xi_3, \dots, \xi_n)$. 
A direct computation shows that we have
$$
(\Phi_2\scirc\Phi_1)(s,\sigma,x,\xi)
=
\Bigl(s,\sigma,x,\bigl(\eta_1(s,\sigma,x,\xi), \eta_2(s,\sigma,x,\xi), \xi_3, \dots, \xi_n\bigr)\Bigr)
\mapob.
$$
The remaining fiber-diffeomorphisms $\Phi_3, \dots, \Phi_n$ are defined similarly. 
We thus have the decomposition $\Phi=\Phi_n\scirc \cdots\scirc \Phi_1$, in which each $\Phi_i$ is of the specific form (only changing odd variables, $\body\partial \eta_i/\partial\xi_j = \delta_{ij}$) and moreover, $\Phi_i$ changes only the $i$th odd coordinate. 
As a permutation of the odd coordinates is simply linear (in the odd variables), it follows from the previous step that we only have to show that $\Phi_1$ satisfies \recalf{thechangeofvariablesformula}. 

By definition we can write
$$
\eta_1(s,\sigma,x,\xi) 
=
\xi_1\cdot\bigl(1+A(s,\sigma,x,\xi_2, \dots, \xi_n)\bigr) + B(s,\sigma,x,\xi_2, \dots, \xi_n)
$$
with
$$
\fracp{\eta_1}{\xi_1}(s,\sigma,x,\xi) = 1+A(s,\sigma,x,\xi_2, \dots, \xi_n)
\mapob,
$$
where $A$ is an even nilpotent function (no terms without odd coordinates).
There are also some restrictions on $B$, but those are not important for what follows. 
It is immediate that we have
$$
\piBer(T_{(x,\xi)}\Phi_{1(s,\sigma)})=\bigl(1+A(s,\sigma,x,\xi_2, \dots, \xi_n)\bigr)\mo
\mapob.
$$
With notation as before, we thus have to determine those terms with a factor $\xi^\nasset$ in the expression
\begin{align*}
\shifttag{5em}
f\bigl(\Phi_{1}(s,\sigma,x,\xi)\bigr) \cdot \piBer(T_{(x,\xi)}\Phi_{1(s,\sigma)})
=
\\&
\sum_{I\subset \{1, \dots, n\}}\sum_{J\subset \{1, \dots, q\}} (\eta_1,\xi_2, \dots\xi_n)^I\,\sigma^J\,\bigl(\Gextension(f_{I,J})\Bigr)(s,x) \cdot (1+A)\mo
\mapob.
\end{align*}
Now (because $A$ does not depend upon $\xi_1$) if $1\notin I$, the product $(\eta_1,\xi_2, \dots\xi_n)^I\cdot (1+A)\mo$ will only contain $\xi_2, \dots, \xi_2$, so can never give rise to $\xi^\nasset$. 
To have a term with $\xi^\nasset$, we thus must have $1\in I$, and then
\begin{align*}
\shifttag{5em}
(\eta_1,\xi_2, \dots\xi_n)^I \cdot (1+A)\mo
= 
\eta_1\cdot (\xi_2, \dots, \xi_n)^{I\setminus \{1\}} \cdot (1+A)\mo
\\&
=
\xi_1\cdot (\xi_2, \dots, \xi_n)^{I\setminus \{1\}} +  B\cdot (\xi_2, \dots, \xi_n)^{I\setminus \{1\}}\cdot (1+A)\mo
\mapob.
\end{align*}
As the functions $A$ and $B$ do not depend upon $\xi_1$, it follows that the only term giving rise to a factor $\xi^\nasset$ comes from $I=\nasset\equiv \{1, \dots, n\}$, and in that case we get exactly
$$
\sum_{J\subset \{1, \dots, q\}} \xi^\nasset\,\sigma^J\,\bigl(\Gextension(f_{\nasset,J})\Bigr)(s,x) 
\mapob.
$$
We thus can compute
\begin{align*}
\shifttag{3.5em}
\int_{U_{(s,\sigma)}} \kern-1em \extder( x , \xi )\ 
f\bigl(\Phi_{1}(s,\sigma,x,\xi)\bigr) \cdot \piBer(T_{(x,\xi)}\Phi_{1(s,\sigma)})
\\&
=
\sum_{J\subset \{1, \dots, q\}} \sigma^J\int_{\body U_{(s,\sigma)}}\extder\Leb^{(d)}(x)\ f_{\nasset,J}(s,x)
\\&
=
\int_{V_{(s,\sigma)}}\kern-1em \extder( x , \eta )\ 
f(s,\sigma,x,\eta)
\mapob,
\end{align*}
as wanted. 

\medskip

\noindent
$\bullet$
Having proved \recalf{thechangeofvariablesformula} when only the odd coordinates $\xi$ change, we now turn our attention to the case when only the even coordinates change, \ie, when $\Phi$ is of the form
$$
\Phi(s,\sigma,x,\xi) = \bigl(s,\sigma,y(s,\sigma,x,\xi), \xi\bigr)
\mapob.
$$
And, as for the case when only odd coordinates change, we first decompose such a fiber-diffeomorphism as a product of two, in which the first has a coordinate change $y = y(s,x)$ independent of all odd variables, and a second one in which we have 
\begin{moneq}[yequalxplusnilpotentinchangeofvariables]
y_i(s,\sigma,x,\xi) = x_i + A_i(s,\sigma,x,\xi)
\end{moneq}
with $A_i$ an even nilpotent function. 
Developing the function $y$ with respect to the odd variables, we have
$$
y(s,\sigma,x,\xi) = 
\sum_{I\subset \{1, \dots, n\}}\sum_{J\subset \{1, \dots, q\}} \xi^I\,\sigma^J\,\bigl(\Gextension(y_{I,J})\Bigr)(s,x) 
\mapob.
$$
Because $\Phi$ is a diffeomorphism, its body is a diffeomorphism of the body sets, implying that the transformation $(s,x)\mapsto \bigl(s,y_{\emptyset,\emptyset}(s,x)\bigr)$ is a diffeomorphism. 
We thus can define the fiber-diffeomorphism $\Phi_1$ by
$$
\Phi_1(s,\sigma,x,\xi) = \bigl(s,\sigma,\Gextension y_{\emptyset,\emptyset}(s,x), \xi\bigr)
$$
and the fiber-diffeomorphism $\Phi_2=\Phi\scirc \Phi_1\mo$. 
By construction $\Phi_2$ satisfies the requirement \recalf{yequalxplusnilpotentinchangeofvariables}.

\medskip

\noindent
$\bullet$
We thus first concentrate on a fiber-diffeomorphism $\Phi$ of the form 
$$
\Phi(s,\sigma,x,\xi) = \bigl(s,\sigma,\Gextension y(s,x), \xi\bigr)
\mapob,
$$
in which case we have
$$
\piBer(T_{(x,\xi)}\Phi_{(s,\sigma)})=\Bigl\vert \Gextension \det\Bigl(\fracp{y_i}{x_j}\Bigr)(s,x)\Bigr\vert
\mapob.
$$
And then we compute:
\begin{align*}
\shifttag{3.5em}
\int_{U_{(s,\sigma)}} \kern-1em \extder( x , \xi )\ 
f\bigl(\Phi(s,\sigma,x,\xi)\bigr) \cdot \piBer(T_{(x,\xi)}\Phi_{(s,\sigma)})
\\&
=
\int_{U_{(s,\sigma)}} \kern-1em \extder( x , \xi )\ 
f\bigl(s,\sigma,y(s,x),\xi\bigr) \cdot \Bigl\vert \Gextension \det\Bigl(\fracp{y_i}{x_j}\Bigr)(s,x)\Bigr\vert
\\&
=
\sum_{J\subset \{1, \dots, q\}} \sigma^J\int_{\body U_{(s,\sigma)}}\extder\Leb^{(d)}(x)\ f_{\nasset,J}\bigl(s,y(s,x)\bigr)\cdot \Bigl\vert \det\Bigl(\fracp{y_i}{x_j}\Bigr)(s,x)\Bigr\vert
\\&
=
\int_{V_{(s,\sigma)}}\kern-1em \extder( y , \xi )\ 
f(s,\sigma,y,\xi)
\mapob,
\end{align*}
where the last equality follows directly from the classical change of coordinates formula for Lebesgue integration. 
We thus have shown \recalf{thechangeofvariablesformula} in this case. 

\medskip

\noindent
$\bullet$
Remains the last case in which the fiber-diffeomorphism is of the form
$$
\Phi(s,\sigma,x,\xi) = \bigl(s,\sigma,y(s,\sigma,x,\xi), \xi\bigr)
$$
with $y_i$ satisfying \recalf{yequalxplusnilpotentinchangeofvariables}. 
Note that this implies that we have the equality $U=V$. 
With the same arguments as used in the case of the odd variables, we can decompose such a fiber-diffeomorphism as a $n$-fold composition of fiber-diffeomorphisms that change a single even variable:
$$
\Phi=\Phi_p \scirc \cdots \scirc \Phi_1
$$
with 
$$
\Phi_i\bigl(s,\sigma,(y_1, \dots, y_{i-1}, x_i, \dots, x_p), \xi\bigr)
=
\bigl(s,\sigma,(y_1, \dots, y_i, x_{i+1}, \dots, x_p), \xi\bigr)
\mapob.
$$
And as in the odd case, a permutation of the even coordinates falls under the previous step, so we only have to show \recalf{thechangeofvariablesformula} for a fiber-diffeomorphism of the form
$$
\Phi(s,\sigma,x,\xi) = \Bigl(s,\sigma,\bigl(y_1(s,\sigma,x,\xi), x_2, \dots, x_p\bigr), \xi\Bigr)
$$
with $y_1$ of the form 
\begin{moneq}
y_1(s,\sigma,x,\xi) = x_1 + A_1(s,\sigma,x,\xi)
\mapob.
\end{moneq}
We thus have
$$
\fracp{y_1}{x_1}(s,\sigma,x,\xi) = 1+\fracp{A_1}{x_1}(s,\sigma,x,\xi) 
=
\piBer(T_{(x,\xi)}\Phi_{(s,\sigma)})
\mapob,
$$
where we note that no absolute values are needed because $A_1$ is nilpotent.
In order to compute
\begin{align*}
\shifttag{3.5em}
\int_{U_{(s,\sigma)}} \kern-1em \extder( x , \xi )\ 
f\bigl(\Phi(s,\sigma,x,\xi)\bigr) \cdot \piBer(T_{(x,\xi)}\Phi_{(s,\sigma)})
\end{align*}
we thus have to determine the terms containing $\xi^\nasset$ in
\begin{align*}
\shifttag{4em}
f\bigl(\Phi(s,\sigma,x,\xi)\bigr) \cdot \piBer(T_{(x,\xi)}\Phi_{(s,\sigma)})
\\&
=
f\bigl(s,\sigma,(x_1+A_1, x_2, \dots, x_p),\xi\bigr)
\cdot
\Bigl(1+\fracp{A_1}{x_1}(s,\sigma,x,\xi)\Bigr)
\mapob.
\end{align*}
As $A_1$ is nilpotent, we can use the \Gextension{} property of a smooth function to write (the infinite sum is actually finite because $A_1$ is nilpotent of maximal order $p+n$):
\begin{align*}
\shifttag{4em}
f\bigl(s,\sigma,(x_1+A_1, x_2, \dots, x_p),\xi\bigr)
\cdot
\Bigl(1+\fracp{A_1}{x_1}(s,\sigma,x,\xi)\Bigr)
\\&
=
\sum_{k=0}^\infty \frac{\partial^k f}{\partial x_1^k}(s,\sigma,x,\xi)\cdot \frac{A_1(s,\sigma,x,\xi)^k}{k!}\cdot \Bigl(1+\fracp{A_1}{x_1}(s,\sigma,x,\xi)\Bigr)
\\&
=
f(s,\sigma,x,\xi) + \fracp{}{x_1} \biggl(\ \sum_{k=0}^\infty \frac{\partial^k f}{\partial x_1^k}\cdot \frac{A_1^{k+1}}{(k+1)!}\ \biggr)
\mapob.
\end{align*}
As $f$ and $A_1$ are smooth functions, there exist smooth functions (of ordinary real variables) $F_{I,J}$ such that we have
$$
\sum_{k=0}^\infty \frac{\partial^k f}{\partial x_1^k}\cdot \frac{A_1^{k+1}}{(k+1)!}
=
\sum_{I\subset \{1, \dots, n\}}\sum_{J\subset \{1, \dots, q\}} \xi^I\,\sigma^J\,\bigl(\Gextension(F_{I,J})\bigr)(s,x)
\mapob.
$$
It follows immediately that we have
\begin{align*}
\shifttag{3.5em}
\int_{U_{(s,\sigma)}} \kern-1em \extder( x , \xi )\ 
f\bigl(\Phi(s,\sigma,x,\xi)\bigr) \cdot \piBer(T_{(x,\xi)}\Phi_{(s,\sigma)})
\\&
=
\sum_{J\subset \{1, \dots, q\}} \sigma^J\int_{\body U_{(s,\sigma)}}\extder\Leb^{(d)}(x)\ \Bigl(f_{\nasset,J}(s,x)+ \fracp{F_{\nasset,J}}{x_1}(s,x)\Bigr)
\\&
=
\int_{V_{(s,\sigma)}}\kern-1em \extder( y , \xi )\ 
f(s,\sigma,y,\xi)
+
\sum_{J\subset \{1, \dots, q\}} \sigma^J
\int_{\body U_{(s,\sigma)}}\extder\Leb^{(d)}(x)\  \fracp{F_{\nasset,J}}{x_1}(s,x)
\mapob.
\end{align*}
But $f$ and thus $F_{I,J}$ have compact support in any slice $\body U_{(s,\sigma)}$, so we have
$$
\int_{\body U_{(s,\sigma)}}\extder\Leb^{(d)}(x)\  \fracp{F_{\nasset,J}}{x_1}(s,x)
=
0
\mapob,
$$
because $F_{\nasset,J}$ is zero at the boundary (being a compactly supported function defined on an open set). 
This finishes the last case. 
\end{preuve}

Once we have this result, we are in business. 
We can define \myquote{partial} densities on locally trivial fiber bundles, integration of this kind of density along the fibers and their behavior under fiber-diffeomorphisms. 
We will do so only for direct products, but the generalization to locally trivial fiber bundles should be obvious.

\begin{definition}{Definition}
Let $G$ and $M$ be two $\CA$-manifolds of graded dimensions $p\vert q$ and $d\vert n$ respectively.

\noindent
$\bullet$
An \stresd{$M$-frame at $(g,m)\in G\times M$} is an ordered homogeneous basis of $\{0\}\times T_mM
=
\ker(T_{(g,m)}\pi_G)\subset T_{(g,m)}(G\times M)$, where $\pi_G:G\times M\to G$ denotes the canonical projection. 
The \stresd{$M$-frame bundle $\framebundle^M (G\times M)\to G\times M$} is the bundle whose fibers $\framebundle^M_{(g,m)} (G\times M)$ consist of all $M$-frames at $(g,m)$. 
Said differently, we have an \myquote{isomorphism}
$$
\framebundle^M_{(g,m)} \cong \framebundle_mM
\mapob.
$$
It is a principal $\Gl(d\vert n,\CA)$ bundle over $G\times M$. 
The right-action of $\Gl(d\vert n,\CA)$ on $\framebundle^M(G\times M)$ is defined as follows. 
For $(v_1, \dots, v_{d+n})\in \framebundle_{(g,m)}^M(G\times M)$, \ie, a basis of $T_mM$, and $(A_{ij})\in \Gl(d\vert n,\CA)$ we have
$$
(v_i)\cdot (A_{jk}) = (w_\ell)
\qquad\text{with}\qquad
w_i = \sum_{j=1}^{d+n} v_j \,A_{ji}
\mapob.
$$

\medskip

\noindent
$\bullet$
An \stresd{$M$-fiber-density} on $G\times M$ is a map $\nu:\framebundle^M (G\times M)\to \CA^\KK$ satisfying the condition that for all $v=(v_i)\in \framebundle_{(g,m)}^M (G\times M)$ and for all $A=(A_{jk}) \in \Gl(d\vert n,\CA)$ we have the equality
$$
\nu\bigl( v\cdot A) = \piBer(A)\cdot \nu(v)
\mapob.
$$
An $M$-fiber-density can be seen as a section of the \stresd{$M$-fiber-density \myquote{line} bundle $\Density^M(G\times M)$} associated to the principal $\Gl(d\vert n,\CA)$ bundle $\framebundle^M (G\times M)$ by the representation $\Gl(d\vert n,\CA)\to\Aut(\CA^\KK)$, $A\mapsto \piBer(A)\mo$ on $\CA^\KK$. It thus has $\CA^\KK$ as typical fiber.

\end{definition}

\begin{definition}{Construction}
Let $\nu$ be a smooth $M$-fiber-density on $G\times M$ that has compact support on each slice $\{g\}\times M$, \ie, for all $g\in G$ the set
$$
\{\,m\in M\mid \nu(g,m)\neq0\,\} \subset M
$$
has compact closure. 
Associated to such an $M$-fiber-density we associate a function on $G$ with values in $\CA^\KK$ which we will denote by $\int_M \nu$. 
This function is defined by the following procedure. 
Let $\mathcal{U}=\{\,U_a\mid a\in A\,\}$ be a cover of $M$ by coordinate charts with (local, even and odd) coordinates $(x^a, \xi^a)$, let $\rho_a$ be a partition of unity associated to this cover, and let $\mathcal V = \{V_b\mid b\in B\}$ be a cover of $G$ by coordinate charts with coordinates $(s,\sigma)$.  
For a point $g\in V_b$ with coordinates $g \cong (s^b,\sigma^b)$ we define the value $\bigl(\int_M \nu\bigr)(g) \equiv \bigl(\int_M \nu\bigr)(s,\sigma) \in \CA^\KK$ by
\begin{align*}
\Bigl(\int_M \nu\Bigr)(s^b,\sigma^b) 
&
= 
\sum_{a\in A} \int_{\{g\}\times U_a} \extder( x^{a}, \xi^{a})\ 
\rho_a(x^a,\xi^a) 
\\&
\kern5em
\cdot \nu\Bigl( \fracp{}{x^a}\bigrestricted_{(s^b,\sigma^b,x^a,\xi^a)}, \fracp{}{\xi^a}\bigrestricted_{(s^b,\sigma^b,x^a,\xi^a)} \Bigr) 
\mapob.
\end{align*}
Note that the $(s^b,\sigma^b)$-dependence enters via $\nu$ when we evaluate it on tangent vectors at $(s^b,\sigma^b,x^a,\xi^a)$. 
The following two results then are a direct application of \recalt{thefiberedchangeofvariablestheorem} (together with the definition).

\end{definition}

\begin{proclaim}{Lemma}
The smooth function $\int_M\nu : G\to \CA^\KK$ is well defined, independent of the choices of the local coordinate charts on $G$ and $M$ and independent of the chosen partition of unity. 

\end{proclaim}

\begin{proclaim}[resultsonintegratingMfiberdensities]{Lemma}
Let $G$ and $H$ be $\CA$-manifolds of dimension $p\vert q$, let $M$ and $N$ be $\CA$-manifolds of dimension $d\vert n$, let $\nu$ be a (smooth) fiber-density on $H\times N$, and let $\Phi:G\times M \to H\times N$ be a fiber-diffeomorphism,\ie, a diffeomorphism such that there exists a diffeomorphism $\Phi_1:G\to H$ making the following diagram commutative (see \recalf{fiberdiffeoascommutativediagram} 
\begin{moneq}
\begin{CD}
G\times M @>\Phi>> H\times N
\\ 
@V{\pi_G}VV @VV{\pi_H}V
\\
G @>>\Phi_1> H\rlap{\mapob.}
\end{CD}
\end{moneq}
\begin{enumerate}
\item
If $v=(v_i)_{i=1}^{d+n}$ is an $M$-frame at $(g,m)$, then $T_{(g,m)}\Phi(v)$ is an $N$-frame at $\Phi(g,m)\in H\times N$. 

\item
The map $\Phi^*\nu:\framebundle^M(G\times M) \to \CA^\KK$ defined by
$$
v\in \framebundle^M_{(g,m)}
\quad\Rightarrow\quad
(\Phi^*\nu)(v) = \nu\bigl(T_{(g,m)}\Phi(v)\bigr)
$$
is a fiber-density on $G\times M$. 

\item
\label{integratingpullbackofafiberdensity}
If the restriction of $\nu$ to any slice $\{h\}\times N$ has compact support, then the restriction of $\Phi^*\nu$ to any slice $\{g\}\times M$ has compact support and we have the equality (of smooth functions on $G$)
$$
\int_M \Phi^*\nu = \Bigl(\,\int_N \nu\Bigr) \scirc \Phi_1
\quad\text{or equivalently}\quad
\Bigl(\int_{M} \Phi^*\nu\Bigr)(g) = \Bigl(\int_{N} \nu\Bigr)\bigl(\Phi_1(g)\bigr)
\mapob.
$$

\end{enumerate}

\end{proclaim}

With the above preparations, we now can attack the question of invariant super metrics and the invariance of the associated super scalar products on $C^\infty_c(M;\CA^\KK)$. 
So let $M$ be an $\CA$-manifold of dimension $d\vert n$, $G$ a super Lie group and $\Phi:G\times M\to M$ a smooth left-action of $G$ on $M$. 
Let furthermore $\mfdmetric$ be a super metric on $M$ and let $\nu$ be a smooth density on $M$. 
There are several (more or less equivalent) ways to define the notion of invariance of $\nu$ and\slash or $\mfdmetric$ under the $G$ action. 
The most direct way is to use the generalized tangent map \cite[V.3.19]{Tu04}. 
For a fixed $g\in G$, the map $\Phi_g:M\to M$, $m\mapsto \Phi(g,m)$ is a bijection, but for $g\notin \body G$ it will (in general) not be smooth. 
Since we still want to define a tangent map $T\Phi_g$, we circumvent this problem by defining the (even linear) map $T_m\Phi_g : T_mM \to T_{\Phi(g,m)}M$ as follows:
we make the identification $T_gG\times T_mM \cong T_{(g,m)}(G\times M)$ and then we define 
$$
v\in T_mM
\quad\Rightarrow\quad
T_m\Phi_g(v) = T_{(g,m)}\Phi(0,v)
\mapob.
$$
With this definition of a generalized tangent map, it becomes easy to define when $\mfdmetric$ or $\nu$ is invariant. 

\noindent$\bullet$
We will say that \stresd{$\mfdmetric$ is invariant under the $G$-action} if for all possible choices we have (compare with \recalf{invariantmetricunderdiffeo})
\begin{moneq}[defofinvarianceunderdiffeos]
\mfdmetric_{\Phi(g,m)}\bigl( T_m\Phi_g(v), T_m\Phi_g(w)  \bigr)
=
\mfdmetric_m(v,w)
\mapob,
\end{moneq}
an equality that we could abbreviate to $(T_m\Phi_g)^*\mfdmetric_{\Phi(g,m)} = \mfdmetric_m$ or even shorter as $\Phi_g^*\mfdmetric = \mfdmetric$. 

\noindent$\bullet$
In the same way we will say that \stresd{$\nu$ is invariant under the $G$-action} if for all possible choices we have
$$
(v_i)_{i=1}^{d+n} \text{ a basis of } T_mM
\quad\Rightarrow\quad
\nu_{\Phi(g,m)}\Big( \bigl(T_m\Phi_g(v_i)\bigr)_{i=1}^{d+n}  \Bigr)
=
\nu_m\bigl(  (v_i)_{i=1}^{d+n} \bigr)
\mapob,
$$
an equality that we could abbreviate as $(T_m\Phi_g)^*\nu_{\Phi(g,m)} = \nu_m$ or even shorter as $\Phi_g^*\nu = \nu$ (see also \recalt{superintinvariantdiffeo}).

For readers that feel uncomfortable with the use of a generalized tangent map we can cast the definition of invariance of $\nu$ in terms of fiber-diffeomorphisms as follows (a similar definition is possible for the invariance of $\mfdmetric$). One starts by extending the density $\nu$ on $M$ in the obvious way to an $M$-fiber-density $\hat\nu$ on $G\times M$ (independent of $g\in G$) and then to require that we have the equality
$$
\Psi^*\hat\nu = \hat\nu
\mapob,
$$
where $\Psi:G\times M\to G\times M$ is the fiber-diffeomorphism defined as
$$
\Psi(g,m) = \bigl(g,\Phi(g,m)\bigr) \equiv \bigl(g,\Phi_g(m)\bigr)
\mapob.
$$
When one looks carefully at this definition, one will see that it boils down to exactly the same formula as the one given by the generalized tangent map.

\begin{proclaim}{Proposition}
Let $M$ be an \ood{} $\CA$-manifold, let $\mfdmetric$ be a super metric on $M$ and let $\Phi:G\times M\to M$ be a smooth left-action of a connected $\CA$-Lie group $G$ on $M$. 
If this $G$-action preserves $\mfdmetric$, then the density $\nu_\mfdmetric$ \recaltt{denstrivwithmetricandOOD}{canonicaltrivdensonOOD} is invariant under the $G$-action.

\end{proclaim}

\begin{preuve}
Let $(x^a, \xi^a)$ be an \ood{} coordinate system around $m\in M$ and let $(x^b, \xi^b)$ be a coordinate system around $\Phi_g(m)$. 
We thus can define the matrix of the even linear map $T_m\Phi_g$ in terms of the bases $(\partial_{x^a},\partial_{\xi^a})$ of $T_mM$ and $(\partial_{x^b},\partial_{\xi^b})$ of $T_{\Phi(g,m)}M$ (a matrix that can be seen as part of the matrix of $T_{(g,m)}\Phi$), which we can write more or less symbolically as
\begin{moneq}[matrixofatangentmapforgroupaction]
(T_m\Phi_g)(\partial_{x^a},\partial_{\xi^a}) = \mathrm{matrix}(T_m\Phi_g) \cdot (\partial_{x^b},\partial_{\xi^b})
\mapob.
\end{moneq}
We can also define the matrix of $\mfdmetric$ in these bases as in \recalf{defofmetricmatrixinaglobalbasis}, giving matrices $\mfdmetric_{ij}^a(m)$ and $\mfdmetric_{ij}^b\bigl(\Phi_g(m)\bigr)$. 
Using the invariance property \recalf{defofinvarianceunderdiffeos} of $\mfdmetric$ and applying the Berezinian, we then obtain the equality 
\begin{moneq}[relationBeronmetricwithdiffeo]
\Ber\bigl(\text{matrix}(T_m\Phi_g)\bigr)^2\cdot \Ber\bigl(\mfdmetric_{ij}^b(\Phi(g,m)\bigr)\,\bigr)
=
\Ber\bigl( \mfdmetric_{ij}^a(m)  \bigr)
\mapob.
\end{moneq}
Next we note that, by a continuity argument, connectedness of $G$ and the fact that $M$ is \ood, we have 
\begin{moneq}[positivityofBerpiundergroupaction]
\forall g\in G
\quad:\quad
\body \piBer\bigl(\text{matrix}(T_m\Phi_g)\bigr) >0
\mapob.
\end{moneq} 
And finally we recall that the definition of $\nu_\mfdmetric$ is given by
$$
\nu_\mfdmetric\Bigl( \fracp{}{x^a}\bigrestricted_{m}, \fracp{}{\xi^a}\bigrestricted_{m} \Bigr)
=
{\sqrt{\bigl\vert\Ber \bigl( \mfdmetric_{ij}^a(m) \bigr)\bigr\vert}} 
\mapob.
$$
We thus can compute:
\begin{align*}
\nu_\mfdmetric\restricted_m\Bigl( \fracp{}{x^a}\bigrestricted_{m}, \fracp{}{\xi^a}\bigrestricted_{m} \Bigr)
\kern0.7em
&
=
\kern0.7em
{\sqrt{\bigl\vert\Ber \bigl( \mfdmetric_{ij}^a(m) \bigr)\bigr\vert}} 
\\&
\oversetalign{\recalf{relationBeronmetricwithdiffeo}}\to=
\kern0.7em
\bigl\vert\,\Ber\bigl(\text{matrix}(T_m\Phi_g)\bigr)\,\bigr\vert
\cdot
{\sqrt{\bigl\vert\Ber \bigl( \mfdmetric_{ij}^b(\Phi(g,m)) \bigr)\bigr\vert}} 
\\&
\oversetalign{\recalf{positivityofBerpiundergroupaction}}\to=
\kern0.7em
\piBer\bigl(\text{matrix}(T_m\Phi_g)\bigr)
\cdot
{\sqrt{\bigl\vert\Ber \bigl( \mfdmetric_{ij}^b(\Phi(g,m)) \bigr)\bigr\vert}} 
\\&
=
\kern0.7em
\piBer\bigl(\text{matrix}(T_m\Phi_g)\bigr)
\cdot
\nu_\mfdmetric\restricted_m\Bigl( \fracp{}{x^b}\bigrestricted_{\Phi(g,m)}, \fracp{}{\xi^b}\bigrestricted_{\Phi(g,m)} \Bigr)
\\&
\oversetalign{\recalf{matrixofatangentmapforgroupaction}}\to=
\kern0.7em
\nu_\mfdmetric\restricted_{\Phi(g,m)}\Bigl( T_m\Phi_g\fracp{}{x^a}\bigrestricted_{m}, T_m\Phi_g\fracp{}{\xi^a}\bigrestricted_{m} \Bigr)
\mapob,
\end{align*}
which shows that $\nu_\mfdmetric$ is invariant under $\Phi_g$ as claimed. 
\end{preuve}

\begin{proclaim}{Lemma}
Let $G$ and $M$ be two $\CA$-manifolds and $F:G\times M\to \CA^\KK$ a smooth function. 
Then for any fixed $g\in G$ the map $F_g:M\to \CA^\KK$, $F_g(m) = F(g,m)$ belongs to $C^\infty(M; \CA^\KK) \otimes \CA^\KK$. 

\end{proclaim}

\begin{preuve}
Let $p\vert q$ be the dimension of $G$ and $d\vert n$ the dimension of $M$. 
Let $g\in G$ be fixed, choose a coordinate chart $V$ around $g$ with local coordinates $(s,\sigma)$ and choose a cover $\atlas=\{\,U_a\mid a\in A\,\}$ of $M$ by coordinate charts with local coordinates $(x^a,\xi^a)$ on $U_a$. 
Then, by definition of a smooth function, there exist ordinary smooth functions $F_{I,J}^a:\body V\times \body U_a\to \KK$ such that we have
$$
F(s,\sigma,x^a,\xi^a) = \sum_{I\subset\{1, \dots, n\} } \sum_{J\subset \{1, \dots, q\} } \sigma^J\cdot (\xi^a)^I\cdot (\Gextension F_{I,J}^a)(s,x^a)
\mapob.
$$
We now assume that $p=d=1$ in order to get less elaborate formul{\ae} (the general case is similar but simply far more elaborate to write). 

Writing $s=t+\alpha$ with $t=\body s\in \RR$ and $\alpha\in \CA_0$ even and nilpotent, and similarly $x^a = y^a + \beta^a$ with $y^a=\body x^a\in \RR$ and $\beta^a$ even and nilpotent, we have (by definition of the \Gextension{}):
\begin{align}
\notag
F(s,\sigma,x^a,\xi^a) 
&
= 
\sum_{I\subset\{1, \dots, n\} } \sum_{J\subset \{1, \dots, q\} } \sigma^J\cdot (\xi^a)^I\cdot 
\sum_{i=0}^\infty\sum_{j=0}^\infty 
\frac{\partial^{i+j}F^a_{I,J}}{(\partial s)^j\,(\partial x^a)^i}(t,y^a) \cdot \frac{\alpha^j\cdot \beta^i}{j! \, i!}
\\&
=
\sum_{J\subset \{1, \dots, q\} } 
\sigma^J\cdot 
\sum_{j=0}^\infty \frac{\alpha^j}{j!}  
\label{expansionfunctiononproduct}
\\&\kern4em
\cdot 
\sum_{I\subset\{1, \dots, n\} } 
(\xi^a)^I\cdot 
\sum_{i=0}^\infty 
\frac{\partial^i}{(\partial x^a)^i} \frac{\partial^{j}F^a_{I,J}}{(\partial s)^j}(t,y^a) \cdot \frac{\beta^i}{i!}
\mapob.
\notag
\end{align}
Looking at this formula suggests to introduce the ordinary smooth functions $f^a_{I,J,t,j}$ defined as
$$
f^a_{I,J,t,j}(y^a)
=
\frac{\partial^{j}F^a_{I,J}}{(\partial s)^j}(t,y^a)
$$
and the (super) smooth functions $f^a_{J,t,j}:U_a \to \CA^\KK$ as
$$
f^a_{J,t,j}(x^a,\xi^a) = 
\sum_{I\subset\{1, \dots, n\} } 
(\xi^a)^I\cdot 
\sum_{i=0}^\infty 
(\Gextension f^a_{I,J,t,j})(x^a) 
\mapob, 
$$
with which we can write
\begin{align*}
F(s,\sigma,x^a,\xi^a)
&
=
\sum_{J\subset \{1, \dots, q\} } 
\sigma^J\cdot 
\sum_{j=0}^\infty \frac{\alpha^j}{j!}  
\cdot 
\sum_{I\subset\{1, \dots, n\} } 
(\xi^a)^I\cdot 
\sum_{i=0}^\infty 
\frac{\partial^i f^a_{I,J,t,j}}{(\partial x^a)^i}(y^a) \cdot \frac{\beta^i}{i!}
\\&
=
\sum_{J\subset \{1, \dots, q\} } 
\sigma^J\cdot 
\sum_{j=0}^\infty \frac{\alpha^j}{j!}  
\cdot 
\sum_{I\subset\{1, \dots, n\} } 
(\xi^a)^I\cdot 
(\Gextension f^a_{I,J,t,j})(x^a) 
\\&
=
\sum_{J\subset \{1, \dots, q\} } 
\sum_{j=0}^\infty 
\frac{\sigma^J\cdot \alpha^j}{j!}  
\cdot 
f^a_{J,t,j}(x^a,\xi^a)
\mapob,
\end{align*}
which belongs, for fixed $g$ with coordinates $(s,\sigma)$, indeed to $C^\infty(U_a, \CA^\KK) \otimes \CA^\KK$ (remember, the infinite sum over $i$ is actually finite because $\alpha$ is nilpotent). 

In order to show that the functions $f^a_{J,t,j}(x^a,\xi^a)$ glue together to global smooth functions on $M$, we note that there is an intrinsic way to define these (local) functions. 
Looking at \recalf{expansionfunctiononproduct} it is immediate that when we derive with respect to the $\sigma$-coordinates and then evaluate at $\sigma=0$, we have
$$
\bigl((\partial_\sigma)^JF\bigr)(s,\sigma=0, x^a, \xi^a)
=
\sum_{I\subset\{1, \dots, n\} } (\xi^a)^I\cdot 
\sum_{i=0}^\infty\sum_{j=0}^\infty 
\frac{\partial^{i+j}F^a_{I,J}}{(\partial s)^j\,(\partial x^a)^i}(t,y^a) \cdot \frac{\alpha^j\cdot \beta^i}{j! \, i!}
\mapob.
$$
Similarly, deriving with respect to $s$ and then evaluating at $s=t$, \ie, at $\alpha=0$, we find
\begin{align*}
\shifttag{5em}
\bigl((\partial_s)^k(\partial_\sigma)^JF\bigr)(s=t,\sigma=0, x^a, \xi^a)
\\&
=
\sum_{I\subset\{1, \dots, n\} } (\xi^a)^I\cdot 
\sum_{i=0}^\infty\sum_{j=0}^\infty 
\frac{\partial^{i+j+k}F^a_{I,J}}{(\partial s)^{j+k}\,(\partial x^a)^i}(t,y^a) \cdot \frac{(\alpha=0)^j\cdot \beta^i}{j! \, i!}
\\&
=
\sum_{I\subset\{1, \dots, n\} } (\xi^a)^I\cdot 
\sum_{i=0}^\infty 
\frac{\partial^{i+k}F^a_{I,J}}{(\partial s)^k\,(\partial x^a)^i}(t,y^a) \cdot \frac{\beta^i}{i!}
\\&
=
\sum_{I\subset\{1, \dots, n\} } (\xi^a)^I\cdot 
\sum_{i=0}^\infty 
\frac{\partial^i}{(\partial x^a)^i}\frac{\partial^{k}F^a_{I,J}}{(\partial s)^k}(t,y^a) \cdot \frac{\beta^i}{i!}
=
f^a_{J,t,k}(x^a,\xi^a)
\mapob.
\end{align*}
Now $\partial_s$ and $\partial_\sigma$ are well defined vector fields on $V\times M\subset G\times M$. It follows that the functions $(\partial_s)^k(\partial_\sigma)^JF$ are well defined smooth functions on $V\times M$. 
We then invoke the property of super smooth functions that says that restriction of some of the variables to real values leaves a super smooth functions of the remaining variables. 
Hence the functions $f_{J,t,k} : M\to \CA^\KK$ defined by
$$
f_{J,t,k}(m) = \bigl( (\partial_s)^k(\partial_\sigma)^JF \bigr) (s=t, \sigma=0,m)
$$
are well defined smooth functions on $M$. 
It thus follows that for fixed $g\in G$ with coordinates $(s,\sigma)$ we have
$$
F_g(m) = 
F(s,\sigma, m) = 
\sum_{J\subset \{1, \dots, q\} } 
\sum_{j=0}^\infty 
\frac{\sigma^J\cdot \alpha^j}{j!}  
\cdot 
f_{J,t,j}(m)
\mapob,
$$
proving that the function $F_g$ belongs to $C^\infty(M; \CA^\KK) \otimes \CA^\KK$. 
\end{preuve}

\begin{proclaim}[therepresentationofagrouponCinftyM]{Corollary}
Let $G$ be a super Lie group acting on an $\CA$-manifold $M$. 
\begin{enumerate}
\item
For each fixed $g\in G$ the map $\rho(g):C_{(c)}^\infty(M; \CA^\KK)\to C_{(c)}^\infty(M; \CA^\KK) \otimes \CA^\KK$ defined by
$$
\bigl(\rho(g)\psi\bigr)(m) = \psi(g\mo m)
$$
extends to an even linear map $\rho(g) : C_{(c)}^\infty(M; \CA^\KK) \otimes \CA^\KK\to C_{(c)}^\infty(M; \CA^\KK) \otimes \CA^\KK$.

\item
The map $\rho:G\to \Aut\bigl( C_{(c)}^\infty(M; \CA^\KK) \otimes \CA^\KK\bigr)$, $g\mapsto \rho(g)$ is a group homomorphism. 

\end{enumerate}

\end{proclaim}

\begin{proclaim}[superscalarproductmfdmetricinvariantunderGaction]{Corollary}
Let $M$ be an \ood{} $\CA$-manifold, let $\mfdmetric$ be a super metric on $M$ and let $\Phi:G\times M\to M$ be the smooth action of an $\CA$-Lie group $G$ on $M$. 
If $\mfdmetric$ is invariant under the $G$-action, then the super scalar product $\supinprsym_\mfdmetric$ on $C_c^\infty(M;\CA^\KK)$ and extended to $C_c^\infty(M;\CA^\KK)\otimes \CA^\KK$ \recalt{superscalarbackandforthbodyextension} is invariant under the representation $\rho$ of $G$ on $C_c^\infty(M;\CA^\KK)\otimes \CA^\KK$ defined by $\bigl(\rho(g)\psi\bigr)(m) = \psi\bigl(\Phi(g\mo,m)\bigr)$. 
More precisely, for all $g\in G$ we have
$$
\homsuperinprod[2]{\rho(g)\chi}{\rho(g)\psi}{_\mfdmetric}
=
\homsuperinprod{\chi}{\psi}{_\mfdmetric}
\mapob.
$$

\end{proclaim}

\begin{preuve}
Instead of proving this property directly, let us work backwards. 
We first note that $\rho(g)\psi$ is the pull-back of $\psi$ by the bijection $\Phi_{g\mo}$. 
According to the definition we thus have:
$$
\homsuperinprod{\rho(g)\chi}{\rho(g)\psi}{_\mfdmetric}
=
\int_M \overline{\Phi_{g\mo}^*\chi} \cdot \Phi_{g\mo}^*\psi\ \nu_\mfdmetric
\mapob.
$$
But $\nu_\mfdmetric$ is invariant under the $G$-action, so we have $\nu_\mfdmetric = \Phi_{g\mo}^*\nu_\mfdmetric$. 
We \myquote{thus} can apply \recalt{superintinvariantdiffeo} to conclude that we have
\begin{align*}
\homsuperinprod{\rho(g)\chi}{\rho(g)\psi}{_\mfdmetric}
&
=
\int_M \overline{\Phi_{g\mo}^*\chi} \cdot \Phi_{g\mo}^*\psi\cdot \nu_\mfdmetric
\\&
=
\int_M \Phi_{g\mo}^*\bigl(\,\overline\chi \cdot \psi\cdot \nu\,\bigr)
\oversettext{\ref{superintinvariantdiffeo}}\to=
\int_M \overline\chi \cdot \psi\cdot \nu
=
\homsuperinprod{\chi}{\psi}{_\mfdmetric}
\mapob.
\end{align*}
However, \recalt{superintinvariantdiffeo} does not apply, as for generic $g$ the map $\Phi_{g\mo}$ is not a diffeomorphism, only for $g\in \body G$ we are sure it will be a diffeomorphism. 
The way out (as most of the time) is to consider all $g$ at the same time. 
We extend in the obvious way the density $\nu_\mfdmetric$ on $M$ to an $M$-fiber-density $\hat\nu_\mfdmetric$ on $G\times M$ and we extend the functions $\chi$ and $\psi$ in the obvious way to functions on $G\times M$. 
We furthermore define the fiber-diffeomorphism $\Psi:G\times M\to G\times M$ by $\Psi(g,m) = \bigl(g,\Phi(g\mo,m)\bigr)$. 
It follows immediately that we have
$$
(\Psi^*\chi)(g,m) = (\Phi_{g\mo}^*\chi)(m)
$$
and similarly for $\psi$. Moreover, the invariance of $\nu_\mfdmetric$ can be stated as the property $\Psi^*\hat\nu_\mfdmetric=\hat\nu_\mfdmetric$. 
With these preparations, the sought for equality is a direct consequence of \recaltt{integratingpullbackofafiberdensity}{resultsonintegratingMfiberdensities}. 
We start with the observation that the smooth function on $G$ defined by
$$
\int_{M} \overline\chi \cdot \psi\cdot \hat\nu_\mfdmetric
\mapob,
$$
(\ie, by integrating the $M$-fiber-density $\overline\chi \cdot \psi\cdot \hat\nu_\mfdmetric$, which is independent of $g\in G$) 
is constant equal to $\homsuperinprod{\chi}{\psi}{_\mfdmetric}$. 
On the other hand, as a function on $G$, the values $\homsuperinprod[2]{\rho(g)\chi}{\rho(g)\psi}{_\mfdmetric}$ can be obtained by integrating the $M$-fiber-density $\overline{\Psi^*\chi}\cdot \Psi^*\psi\cdot \hat\nu_\mfdmetric$. 
And then:
\begin{align*}
\homsuperinprod[2]{\rho(g)\chi}{\rho(g)\psi}{_\mfdmetric}
\quad
&
=
\quad
\Bigl(\int_{M} \overline{\Psi^*\chi}\cdot \Psi^*\psi\cdot \hat\nu_\mfdmetric\Bigr)(g)
\oversettext{inv.}\to=
\Bigl(\int_{M} \Psi^*\bigl(\overline{\chi}\cdot \psi\cdot \hat\nu_\mfdmetric\bigr)\Bigr)(g)
\\&
\oversetalign{\recaltt{integratingpullbackofafiberdensity}{resultsonintegratingMfiberdensities}}\to=
\quad
\Bigl(\int_{M} \overline{\chi}\cdot \psi\cdot \hat\nu_\mfdmetric\Bigr)(g)
=
\homsuperinprod{\chi}{\psi}{_\mfdmetric}
\mapob,
\end{align*}
which now is an equality between smooth functions on $G$. 
\end{preuve}

\masection{Berezin-Fourier, Hodge-star and super scalar products}
\label{BerezinFourierandHodgestarsection}

In this section we introduce the Berezin-Fourier transform (also called the super Fourier transform or Fermionic Fourier transform) and we show in \recalt{linksuperinprodHodgeandFFourier} that there is a close link between this operation, the Hodge-$*$ operation and the super scalar product $\supinprsym_{B,\mfdmetric}$. 
In \recals{backtomotivtingexamplesection} we will use the Berezin-Fourier transform in the decomposition of a super unitary representation into a family of super unitary representations depending upon odd parameters, providing a rigorous justification of our heuristic decomposition used in \recals{motivatingexamplesection}.

\firstofmysubsection
\mysubsection{BerezinFourierandHodgestarsection}{We start with the elementary Berezin-Fourier transform}

\begin{definition}[defclassicalBerezinFouriertransf]{Definition}
We define the \stresd{(elementary) Berezin-Fourier transform} (or \stresd{Fermi\-onic Fourier transform}, see \cite[\S7.1.5]{GuiSte1999}) as an application of Berezin integration.  
Let $(E, \inprodsym_E, \supinprsym_E)$ be a proto super Hilbert space over $\CC$ and $f:\CA_1^n\to E\otimes \CA^\CC$ a smooth map. 
There thus exist, for all $I\subset \{1, \dots, n\}$, elements $f_I\in E$ such that we have 
$$
f(\xi) = \sum_{I\subset \{1, \dots, n\}} \xi^I\,f_I
\mapob.
$$
We then define the \stresd{Berezin-Fourier transform} of the $f$ to be the smooth function $\Fourierodd f:\CA_1^n\to E\otimes \CA^\CC$ given by\footnote{Note the \myquote{subtle} difference with the notation for the ordinary Fourier transform: the Berezin-Fourier transform has a small \myquote{o} on top of the symbol $\Fourier$.}
$$
(\Fourierodd f)(\eta)
= 
\int_{\CA_1^n}\extder\xi^{(n)}\  \exp\Bigl({-i\,\sum_{p=1}^n\xi_p\eta_p}\Bigr)\,f(\xi) 
\mapob.
$$
We also define the odd maps $\Inv_j:C^\infty(\CA_1^n, E \otimes\CA^\CC) \to C^\infty(\CA_1^n, E \otimes\CA^\CC)$ by
$$
\Inv_j = \xi_j +i\, \partial_{\xi_j}
\mapob.
$$
And we finally equip $\CA_1^n$ with the standard super metric $\mfdmetric_o$ defined by
\begin{moneq}[standardmetriconCAsub1highn]
\mfdmetric_o(\partial_{\xi_j}, \partial_{\xi_k}) = i\cdot \delta_{jk}
\mapob.
\end{moneq}
Using this metric we can equip $C^\infty(\CA_1^n, E \otimes\CA^\CC) \equiv C^\infty_c(\CA_1^n, E \otimes\CA^\CC)$ with an ordinary metric $\inprodsym_{\mfdmetric_o}$ and a super scalar product $\supinprsym_{B,\mfdmetric_o}=\supinprsym_{\mfdmetric_o}$ according to \recalt{firstdefofordinarymetriconsuperfunctionsinE} and \recalt{defofsuperscalarprodforfunctionsonoodM}; they are given by 
$$
\inprod fg_{\mfdmetric_o} = \sum_{I\subset \{1, \dots, n\}} \inprod{f_I}{g_I}_E
\quad\text{and}\quad
\homsuperinprod fg{_{\mfdmetric_o}} = \sum_{I\subset \{1, \dots, n\}} (-1)^{\varepsilon(I,I^c)} \homsuperinprod{f_I}{g_{I^c}}{_E}
\mapob.
$$

\end{definition}

\begin{definition}[remarksondefofBerezinFourier]{Remarks}
$\bullet$ 
In defining the Berezin-Fourier transform as we did, we have made several arbitrary choices: we could have multiplied $\Fourierodd$ by an arbitrary factor (which might depend upon $n$!) and we could have used $\exp\bigl({+i\,\sum_{p=1}^n\xi_p\eta_p}\bigr)$ instead of $\exp\bigl({-i\,\sum_{p=1}^n\xi_p\eta_p}\bigr)$ (which is equivalent to a choice of order, as we have $\sum_{p=1}^n\xi_p\eta_p = -\sum_{p=1}^n\eta_p\xi_p$). 
However, if we want the Berezin-Fourier transform (or its generalization in the next subsection) to be an equivalence of proto super Hilbert spaces, then we have to use $\exp\bigl({-i\,\sum_{p=1}^n\xi_p\eta_p}\bigr)$ (see also \recalt{remarkonchoicesignforsupermetric} and \recalt{remarksonequivalenceofSHSandSUR}).\footnote{Leaving out the factor $i$ altogether, using $\exp\bigl({\sum_{p=1}^n\xi_p\eta_p}\bigr)$, is not a good idea, as then some of the results will depend on more then parity alone. 
Moreover, it is extremely useful, if not essential, in the Berezin-Fourier decomposition of a super unitary representation (see \recals{backtomotivtingexamplesection}).}

\medskip

$\bullet$ 
The attentive reader might have noticed that the definition of the Berezin-Fourier transform is not quite rigorous. 
In the first place it should be noted that it is an instance of integrating along fibers as explained in \recals{Integrationalongfiberssection}, as we integrate a function of the odd variables $\xi$ and $\eta$ only over the $\xi$ variables. 
But even then it is not quite rigorous, as in \recals{Integrationalongfiberssection} we only discussed integration of functions with values in $\CA^\KK$, not in an arbitrary graded vector space $E\otimes \CA^\KK$ (in which $E$ might be infinite dimensional). 
When one wants to integrate over even and odd variables, one needs to be (very) careful when the target space is infinite dimensional. 
On the other hand, when one integrates only over odd variables, the generalization to vector valued functions poses no problems. 
In order to show that our definition really is (or better, can be made) fully rigorous, we introduced the operators $\Inv_j$. 
Using these, we show in \recalt{BerezinFourierasproductofInvsNEW} that the Berezin-Fourier transform can be expressed in terms of these operators for \myquote{any} target space. 

\end{definition}

\begin{proclaim}[CliffordforoddInvoperators]{Lemma}
The odd maps $\Inv_j$ satisfy the relations
$$
\Inv_j\scirc \Inv_j = i\cdot\id
\qquad,\qquad
j\neq k
\Rightarrow
[\Inv_j,\Inv_k] \equiv \Inv_j\scirc \Inv_k + \Inv_k \scirc \Inv_j = 0
\mapob.
$$
Moreover, they are graded symmetric with respect to the super scalar product $\supinprsym_{\mfdmetric_o}$.

\end{proclaim}

\begin{preuve}
The first assertions are a direct computation. 
For graded skew-symmetry of $\Inv_j$, we take two smooth homogeneous functions $f,g:\CA_1^n\to E\otimes \CA^\CC$. 
Smoothness implies that there exists smooth functions $f_0,f_1, g_0, g_1$ of one odd variable less (the $\xi_j$!) such that we have
$$
f(\xi) = f_0(\xi') + \xi_j\,f_1(\xi')
\qquad\text{and}\qquad
g(\xi) = g_0(\xi') + \xi_j\,g_1(\xi')
\mapob,
$$
where $\xi'$ denotes the sequence of $n-1$ odd variables not containing $\xi_j$. 
Homogeneity means that $f_0$ has the same parity as $f$, whereas $f_1$ has the opposite parity (and similarly for $g_0$ and $g_1$). 
We then compute, using \recalt{leftlinearityofrightbilinear},
\begin{align*}
\homsuperinprod {\Inv_jf}g{_{\mfdmetric_o}} 
&
- 
(-1)^{\parity f}\,\homsuperinprod{f}{\Inv_jg}{_{\mfdmetric_o}}
\\&
=
\int_{\CA_1^n} \extder \xi^{(n)}\ 
\homsuperinprod[2]{(\Inv_jf)(\xi)}{g(\xi)}{_E} - (-1)^{\parity f}\,\homsuperinprod[2]{f(\xi)}{(\Inv_j g)(\xi)}{_E}
\\&
=
\int_{\CA_1^n} \extder \xi^{(n)}\ 
\homsuperinprod[2]{\xi_j\,f_0(\xi')+i\,f_1(\xi')}{g_0(\xi')+ \xi_j\,g_1(\xi')}{_E}
\\&
\kern8em
- (-1)^{\parity f}\,
\homsuperinprod[2]{f_0(\xi')+\xi_j\,f_1(\xi')}{\xi_j\,g_0(\xi') + i\, g_1(\xi')}{_E}
\\&
=
\int_{\CA_1^n} \extder \xi^{(n)}\ 
\homsuperinprod[2]{\xi_j\,f_0(\xi')}{g_0(\xi')}{_E} 
- (-1)^{\parity f}\,
\homsuperinprod[2]{f_0(\xi')}{\xi_j\,g_0(\xi')}{_E} 
\\&
\kern5em
+ \homsuperinprod[2]{i\,f_1(\xi')}{\xi_j\,g_1(\xi')}{_E}
- (-1)^{\parity f}\,\homsuperinprod[2]{\xi_j\,f_1(\xi')}{i\,g_1(\xi')}{_E}
\\&
=
0
\mapob,
\end{align*}
where for the third equality we used that the other terms do not contain $\xi_j$ and thus have zero Berezin integral. 
The final results follows directly from the anti-linearity property in the first variable \recaltt{defofantirightlinearity}{defofgradedhermitianforms} and the parities of $f_0$ and $f_1$.
\end{preuve}

\begin{definition}{Notation}
As in the sequel we will be confronted often with a specific combination of powers of $-1$ and of $i$, we will introduce an understandable shorthand for them. 
For any integer $k\in \ZZ$ we define $\oneori n$ as
$$
\oneori k = (-1)^{\frac12k(k-1)}\cdot i^k
=
\begin{cases}
1& \quad \text{if $k$ is even,}
\\
i & \quad\text{if $k$ is odd.}
\end{cases}
$$
Similarly we define $\oneormi n$ as
$$
\oneormi k = (-1)^k\,\oneori k
=
(-1)^{\frac12 k(k-1)}\cdot (-i)^k
=
\begin{cases}
1& \quad \text{if $k$ is even,}
\\
-i & \quad\text{if $k$ is odd.}
\end{cases}
$$
Note that we do have the conjugation property $\overline{\oneori k} = \oneormi k$, but that we do \stress{not} have any homomorphism like property such as $\oneori{ k+\ell} = \oneori k \cdot \oneori \ell$.

\end{definition}

\begin{proclaim}[BerezinFourierasproductofInvsNEW]{Lemma}
Unlike its ordinary counter part, the Berezin-Fourier transform has a simple operator description:
$$
\Fourierodd f = (-i)^n\,(\Inv_n \scirc \cdots \scirc \Inv_1)(f)
\equiv
\oneormi n \,\Inv^\nasset(f)
\mapob.
$$
More precisely, for any $I\subset \{1, \dots, n\}$ and any $\psi\in E$ we have
$$
\Fourierodd \xi^I\,\psi = \oneormi{\,\parity{I^c}\,}\, (-1)^{\varepsilon(I,I^c)}\,\xi^{I^c}\,\psi
\mapob.
$$

\end{proclaim}

\begin{preuve}
We first compare the two operations on monomials of the form $f(\xi) = \xi^I\,\psi$ for some $I\subset \{1, \dots, n\}$. 
We obtain:
\begin{align*}
(\Fourierodd \xi^I\,\psi)(\eta)
&
=
\int_{\CA_1^n} \extder \xi^{(n)} \ 
\exp\Bigl({-i\,\sum_{p=1}^n\xi_p\eta_p}\Bigr)\,\xi^I\,\psi
=
\int_{\CA_1^n} \extder \xi^{(n)} \ (-i)^{\parity{I^c}}\,
\xi^I\,\prod_{j\in I^c} \xi_j\eta_j\,\psi
\\&
=
(-1)^{\frac12 I^c(I^c-1)}\,(-i)^{\parity{I^c}}
\int_{\CA_1^n} \extder \xi^{(n)} \ \xi^I\,\xi^{I^c}\,\eta^{I^c}\,\psi
\\&
= (-1)^{\frac12 I^c(I^c-1)+\varepsilon(I,I^c)}\,(-i)^{\parity{I^c}}\,\eta^{I^c}\,\psi
\end{align*}
and
\begin{align*}
(-i)^n\,\Inv_n \scirc \cdots \scirc \Inv_1\, \xi^I\,\psi
&
=
(-1)^{\frac12 n(n-1)+\varepsilon(I^c,I)}\,(-i)^n\, \Inv^{I^c} \scirc \Inv^I\,\xi^I\,\psi
\\&
=
(-1)^{\frac12 n(n-1)+\frac12 I(I-1)+\varepsilon(I^c,I)}\,(-i)^n\, i^{\parity I}\, \xi^{I^c}\,\psi 
\\&
=
(-1)^{\frac12 n(n-1)+\frac12 I(I-1)+II^c+\varepsilon(I,I^c)}\, (-i)^{\parity{I^c}}\, \xi^{I^c}\,\psi 
\\&
=
(-1)^{\frac12 I^c(I^c-1)+\varepsilon(I,I^c)}\, (-i)^{\parity{I^c}}\, \xi^{I^c} \,\psi 
\mapob.
\end{align*}
As both operations are right-linear, we have equality for all $f\in C^\infty(\CA_1^n,E\otimes \CA^\CC)$.
\end{preuve}

\begin{proclaim}[squareofFermionicFourier]{Corollary \cite[Prop.~7.1.2]{GuiSte1999}}
$\Fourierodd\scirc\Fourierodd = \oneormi n\, \id$ and thus $\Fourierodd\mo = \oneori n\,\Fourierodd = \Inv^\nasset$.

\end{proclaim}

\mysubsection{BerezinFourierandHodgestarsection}{We then generalize to $\CA$-manifolds}

\begin{definition}[defofgeneralizedBFtransform]{Definition}
Let $(E, \inprodsym_E, \supinprsym_E)$ be a proto super Hilbert space over $\CC$, let $M$ be an \ood{} $\CA$-manifold, let $\mfdmetric$ be a super metric on $M$ and choose an \ood{} Batchelor atlas $\atlas = \{U_a\mid a\in A\}$. 
In any chart $U_a$ with coordinates $(x^a,\xi^a)$ we define the matrix $H^a(x,\xi) = -i\,(\Gextension\body \mfdmetric^{a,11})(x)=-i\,\mfdmetric^{a,11}(x,0)$ depending on the even coordinates $x^a$ only, \ie, 
$$
H^a_{pq}(x^a)
=
-i\cdot  
\mfdmetric\biggl(\, 
\fracp{}{\xi^a_p}\bigrestricted_{(x^a,\xi^a=0)}\ ,\ 
\fracp{}{\xi^a_q}\bigrestricted_{(x^a,\xi^a=0)}
\, \biggr)
\mapob.
$$
Comparing with \recalf{defofinducedmetriconBatchelorBundle} shows that $H^a$ is essentially the \Gextension{} of the inverse of the metric $\mfdmetric_{V\!M}$. 
Using this matrix we define for any smooth function $\psi:M\to E \otimes \CA^\CC$ its \stresd{generalized Berezin-Fourier transform} $\Fourierodd_{B,\mfdmetric}\psi$ in the chart $U_a$ by
$$
(\Fourierodd_{B,\mfdmetric} \psi)(x,\xi)
=
\int_{\CA_1^n} \frac{\extder\eta^{(n)}}{\sqrt{\det\bigl(H^a_{pq}(x)\bigr)}}\ 
\psi(x,\eta) \cdot \exp\Bigl(-i\,{\sum_{p,q=1}^n \eta_p \,H^a_{pq}(x) \,\xi_q}\Bigr)
\mapob.
$$

\end{definition}

\begin{proclaim}[BerezinFourierwithmetricwelldefined]{Lemma}
The Berezin-Fourier transform $\Fourierodd_{B,\mfdmetric}$ is a well defined operation $\Fourierodd_{B,\mfdmetric}:C^\infty(M;E \otimes\CA^\CC) \to C^\infty(M;E \otimes\CA^\CC)$ depending on the choice of an \ood{} Batchelor atlas.

\end{proclaim}

\begin{preuve}
The definition of the generalized Berezin-Fourier transform as given in \recalt{defofgeneralizedBFtransform} is rather incomplete, as it suppresses the use of the local charts, a fact that will be crucial for the proof. 
So let $\varphi_a:U_a \to O_a\subset F_0$ be the chart map with $O_a$ an open set in even part of a fixed graded vector space $F$ of dimension $d\vert n$. 
For any smooth function $\psi:M\to E\otimes\CA^\CC$ we thus have a family of smooth functions $\psi_a:O_a\to E\otimes \CA^\CC$ defined as $\psi_a = \psi\scirc \varphi_a\mo$. 
Obviously this family satisfies the condition
\begin{moneq}[compatibilitycondforfunctions]
\psi_b \scirc \varphi^{ba} = \psi_a
\mapob,
\end{moneq}
where the maps $\varphi^{ba}= \varphi_b \scirc \varphi_a\mo$ are the coordinate changes $(x^b, \xi^b) = \varphi^{ba}(x^a,\xi^a)$. 
The definition of the generalized Berezin-Fourier transform thus defines a family of maps
$$
(\Fourierodd_{B,\mfdmetric} \psi)_a(x^a,\xi^{a})
=
\int_{\CA_1^n} \frac{\extder\eta^{(n)}}{\sqrt{\det\bigl(H^a_{pq}(x^a)\bigr)}}\ 
\psi_a(x^a,\eta) \cdot
\exp\Bigl(-i\,{\sum_{p,q=1}^n \eta_p \,H^a_{pq}(x^a) \,\xi^a_q}\Bigr)
\mapob,
$$
for which we have to show that they satisfy the compatibility condition \recalf{compatibilitycondforfunctions}, \ie,
\begin{moneq}[compatibilityFourierforfunctions]
(\Fourierodd_{B,\mfdmetric} \psi)_b(x^b,\xi^b)
=
(\Fourierodd_{B,\mfdmetric} \psi)_a(x^a,\xi^a)
\mapob.
\end{moneq}
We now recall that our two charts are in a Batchelor atlas, meaning that the change of coordinates is of the form
$$
x^b=x^b(x^a)
\qquad\text{and}\qquad
\xi^b_p = \sum_{q=1}^n T^{ba}_{pq}(x^a)\,\xi^a_q
$$
for some (smooth) matrix valued function $T^{ba}(x^a)$.
Moreover, according to \recalf{changemetricbodycoordinates}, we have
$$
H^a_{pq}(x^a)
=
\sum_{j,k=1}^n
T^{ba}_{jp}(x^a) \cdot 
H^b_{jk}\bigl(x^b(x^a)\bigr) \cdot 
T^{ba}_{kq}(x^a)
$$
and thus, because  $\det\bigl(T_{ba}(x)\bigr)$ is positive by the assumption that the Batchelor atlas is \ood, we have
$$
\sqrt{\det\bigl(H^a_{pq}(x^a)\bigr)} = \det\bigl(T^{ba}(x^a)\bigr)\cdot \sqrt{\det\bigl(H^b_{jk}\bigl(x^b(x^a)\bigr)\,\bigr)}
\mapob.
$$
Using the explicit form of the coordinate change, we now simply compute the left hand side of \recalf{compatibilityFourierforfunctions}:
\begin{align*}
(\Fourierodd_{B,\mfdmetric} \psi)_b\bigl(x^b,\xi^b)
&
=
\int_{\CA_1^n} \frac{\extder\eta^{(n)}}{\sqrt{\det\bigl(H^b_{pq}(y)\bigr)}}\ \psi_b\bigl(x^b,\eta\bigr)
\\
\text{\footnotesize $\xi^b_p = \sum_{q=1}^n T^{ba}_{pq}(x^a)\,\xi^a_q$\quad}
&\kern5em
\cdot
\exp\Bigl({-i\sum_{j,k,q=1}^n \eta_j \,H^b_{jk}(x^b) \,T^{ba}_{kq}(x^a)\,\xi^a_q}\Bigr) 
\\
\text{\footnotesize(\recalt{changeofcoordinatesinsuperint} with $ \eta=T^{ba}(x)\,\zeta$)\quad}
&
=
\int_{\CA_1^n} \frac{\extder\zeta^{(n)}}{\det\bigl(T^{ba}(x^a)\bigr)\cdot\sqrt{\det\bigl(H^b_{pq}(x^b)\bigr)}}\ 
\psi_b\bigl(x^b,T^{ba}(x^a)\zeta\bigr)
\\&
\kern5em
\cdot 
\exp\Bigl({-i\sum_{j,k,p,q=1}^n T^{ba}_{jp}(x^a)\,\zeta_p \,H^b_{jk}(x^b) \,T^{ba}_{kq}(x^a)\,\xi^a_q}\Bigr)
\\&
=
\int_{\CA_1^n} \frac{\extder\zeta^{(n)}}{\sqrt{\det\bigl(H^a_{pq}(x^a)\bigr)}}\ 
\cdot \psi_b(x^b,T^{ba}(x^a)\zeta)
\\&\kern5em
\cdot 
\exp\Bigl({-i\sum_{p,q=1}^n \zeta_p \,H^a_{pq}(x^a) \,\xi^a_q}\Bigr)
\\&
=
(\Fourierodd_{B,\mfdmetric}\psi)_a(x^a,\xi^a)
\mapob,
\end{align*}
where for the last equality we used that the functions $\psi_a$ and $\psi_b$ satisfy the compatibility condition \recalf{compatibilitycondforfunctions}.
\end{preuve}

\begin{definition}{Remark}
When we apply the generalized Berezin-Fourier transform to the $\CA$-manifold $\CA_1^n$ equipped with the standard super metric $\mfdmetric_o$ \recalf{standardmetriconCAsub1highn}, we recover immediately the definition of the elementary Berezin-Fourier transform $\Fourierodd$ defined in the previous subsection. 
The generalized Berezin-Fourier transform $\Fourierodd_{B,\mfdmetric}$ thus is a genuine generalization. 

\end{definition}

\begin{proclaim}[generalizedBerezinFourierisequivalenceSHS]{Proposition}
The map $\Fourierodd_{B,\mfdmetric}:C^\infty(M;E \otimes\CA^\CC) \to C^\infty(M;E \otimes\CA^\CC)$ has the following properties. 
\begin{enumerate}
\item
$\Fourierodd$ is homogeneous of parity $\parity\Fourierodd = n$ and $\Fourierodd_{B,\mfdmetric} \scirc \Fourierodd_{B,\mfdmetric} = \oneormi n\, \id$.

\item
For $\chi,\psi\in C^\infty_c(M;E \otimes\CA^\CC)$ with $\chi$ homogeneous we have the equalities
\begin{align*}
\homsuperinprod{\Fourierodd_{B,\mfdmetric} \chi}{\Fourierodd_{B,\mfdmetric} \psi}{_{B,\mfdmetric}} 
&
=
\oneori{n}\,(-1)^{n \parity \chi}\,\homsuperinprod{\chi}{\psi}{_{B,\mfdmetric}} 
\\[2\jot]
\inprod{\Fourierodd_{B,\mfdmetric} \chi}{\Fourierodd_{B,\mfdmetric} \psi}_{B,\mfdmetric}
&
=
\inprod\chi\psi_{B,\mfdmetric}
\mapob,
\end{align*}
\ie, $\Fourierodd_{B,\mfdmetric} : C^\infty_c(M;E \otimes\CA^\CC) \to C^\infty_c(M;E \otimes\CA^\CC)$ is an equivalence of super Hilbert spaces \recalt{defofmainSHS}. 

\end{enumerate}

\end{proclaim}

\begin{preuve}
Let $\atlas=\{\,U_a\mid a\in A\,\}$ be an \ood{} Batchelor atlas adapted to the super metric $\mfdmetric$. 
Then according to \recalt{metriconsuperfunctionsinEinlocalcoordinates} we have
\begin{align*}
\inprod \chi\chi_\mfdmetric
&
=
\sum_{a\in A} 
\int_{\body U_a} \extder\Leb^{(d)}(x^a)\ 
\rho_a(x^a)\cdot
\sqrt{\bigl\vert\det\bigl(\body\mfdmetric^{a,00}(x^a)\bigr)\bigr\vert}
\\&
\kern7em
\cdot
\sum_{I\subset \{1, \dots, n\} } \inprod{\chi_{a,I}(x^a)}{\chi_{a,I}(x^a)}_E
\mapob,
\end{align*}
and according to \recalt{descriptionsuperscalarproductBginadaptedBatlas} we have 
\begin{align}
\homsuperinprod \chi\psi{_{B,\mfdmetric}}
&
=
\sum_{a\in A} 
\int_{\body U_a} \extder\Leb^{(d)}(x^a)\ 
\rho_a(x^a)\cdot
\sqrt{\bigl\vert\det\bigl(\body\mfdmetric^{a,00}_{ij}(x^a)\bigr)\bigr\vert}
\notag
\\&
\kern5em
\cdot
\sum_{I\subset \{1, \dots, n\} } (-1)^{\varepsilon(I,I^c)}\, \homsuperinprod[3]{\conjugate^{\parity I}\bigl(\chi_{a,I}(x^a)\bigr)}{\conjugate^n\bigl(\psi_{a,I^c}(x^a)\bigr)}{_E}
\mapob,
\label{expressionsupinprodforFourieroddisSHSequivalence}
\end{align}
At the same time, the Berezin-Fourier transform of $\chi$ (and similarly for $\psi$) is given in the same adapted atlas on the local chart $U_a$ by
$$
(\Fourierodd_{B,\mfdmetric} \chi)_a(x^a,\xi^{a})
=
\int_{\CA_1^n} \extder\eta^{(n)}\ 
\chi_a(x^a,\eta) \cdot
\exp\Bigl(-i\,{\sum_{p=1}^n \eta_p \,\xi^a_p}\Bigr)
\mapob,
$$
simply because in an adapted atlas the matrices $H^a$ are the identity. 
Using \recalt{BerezinFourierasproductofInvsNEW} we then find
\begin{moneq}[localexpressionBerezinFourierinadaptedchart]
(\Fourierodd_{B,\mfdmetric} \chi)_{a,I}(x^a)
=
\oneormi{\,\parity I\,}\,(-1)^{\varepsilon(I^c,I)}\,\chi_{a,I^c}(x^a)
\mapob.
\end{moneq}

\smallskip

In order to prove (i) we note that if $\chi$ has parity $\parity\chi$, then we have
\begin{moneq}[parityofchiaIintermsofparitychianda]
\parity{\chi_{a,I}} = \parity\chi+\parity I
\mapob.
\end{moneq}
This together with \recalf{localexpressionBerezinFourierinadaptedchart} shows that $\Fourierodd_{B,\mfdmetric}$ has parity $n$. 
Moreover, applying \recalf{localexpressionBerezinFourierinadaptedchart} twice gives us the equality
\begin{align*}
\bigl((\Fourierodd_{B,\mfdmetric} \scirc \Fourierodd_{B,\mfdmetric})(\chi)\bigr)_{a,I}(x^a)
&
=
\oneormi{\,\parity{I}\,}\,(-1)^{\varepsilon(I^c,I)}
\,(\Fourierodd_{B,\mfdmetric} \chi)_{a,I^c}(x^a)
\\&
=
\oneormi{\,\parity{I}\,}\,(-1)^{\varepsilon(I^c,I)}\,
\oneormi{\,\parity{I^c}\,}\,(-1)^{\varepsilon(I,I^c)}\,\chi_{a,I}(x^a)
\\&
=
\oneormi{n}\,\chi_{a,I}(x^a)
\mapob,
\end{align*}
from which the equality $\Fourierodd_{B,\mfdmetric} \scirc \Fourierodd_{B,\mfdmetric} = \oneormi n\, \id$ follows.

\smallskip

For (ii) we compute:
\begin{align*}
\inprod{\Fourierodd_{B,\mfdmetric} \chi}{\Fourierodd_{B,\mfdmetric} \psi}_\mfdmetric
&
=
\sum_{a\in A} 
\int_{\body U_a} \extder\Leb^{(d)}(x^a)\ 
\rho_a(x^a)\cdot
\sqrt{\bigl\vert\det\bigl(\body\mfdmetric^{a,00}(x^a)\bigr)\bigr\vert}
\\&
\kern7em
\cdot
\sum_{I\subset \{1, \dots, n\} } \inprod{\chi_{a,I^c}(x^a)}{\psi_{a,I^c}(x^a)}_E
\\&
=
\inprod\chi\psi_\mfdmetric
\mapob,
\end{align*}
and
\begin{align*}
\shifttag{3em}
\homsuperinprod{\Fourierodd_{B,\mfdmetric} \chi}{\Fourierodd_{B,\mfdmetric} \psi}{_{B,\mfdmetric}}
\\&
=
\sum_{a\in A} 
\int_{\body U_a} \extder\Leb^{(d)}(x^a)\ 
\rho_a(x^a)\cdot
\sqrt{\bigl\vert\det\bigl(\body\mfdmetric^{a,00}_{ij}(x^a)\bigr)\bigr\vert}
\\&
\kern2.5em
\cdot
\sum_{I\subset \{1, \dots, n\} } (-1)^{\varepsilon(I,I^c)}\, 
\homsuperinprod[3]{\conjugate^{\parity I}\bigl((\Fourierodd_{B,\mfdmetric} \chi)_{a,I}(x^a)\bigr)}{\conjugate^n\bigl((\Fourierodd_{B,\mfdmetric} \psi)_{a,I^c}(x^a)\bigr)}{_E}
\\&
=
\sum_{a\in A} 
\int_{\body U_a} \extder\Leb^{(d)}(x^a)\ 
\rho_a(x^a)\cdot
\sqrt{\bigl\vert\det\bigl(\body\mfdmetric^{a,00}_{ij}(x^a)\bigr)\bigr\vert}
\\&
\kern6em
\cdot
\sum_{I\subset \{1, \dots, n\} } (-1)^{\varepsilon(I,I^c)}\, 
\supinprsym_{E}\Bigl({\conjugate^{\parity I}\bigl(\oneormi{\,\parity{I}\,}\,(-1)^{\varepsilon(I^c,I)}\,\chi_{a,I^c}(x^a)\bigr)}\,,
\\&
\kern19em
{\conjugate^n\bigl(\oneormi{\,\parity{I^c}\,}\,(-1)^{\varepsilon(I,I^c)}\,\psi_{a,I}(x^a)\bigr)}\Bigr)
\\&
=
\sum_{a\in A} 
\int_{\body U_a} \extder\Leb^{(d)}(x^a)\ 
\rho_a(x^a)\cdot
\sqrt{\bigl\vert\det\bigl(\body\mfdmetric^{a,00}_{ij}(x^a)\bigr)\bigr\vert}
\\&
\kern3em
\cdot
\sum_{I\subset \{1, \dots, n\} } \oneormi{\,\parity{I^c}\,}\,\oneori{\,\parity{I}\,}\, (-1)^{\varepsilon(I^c,I)}\, 
\homsuperinprod[3]{\conjugate^{\parity I}\bigl(\chi_{a,I^c}(x^a)\bigr)}{\conjugate^n\bigl(\psi_{a,I}(x^a)\bigr)}{_E}
\\&
=
\sum_{a\in A} 
\int_{\body U_a} \extder\Leb^{(d)}(x^a)\ 
\rho_a(x^a)\cdot
\sqrt{\bigl\vert\det\bigl(\body\mfdmetric^{a,00}_{ij}(x^a)\bigr)\bigr\vert}
\\&
\kern3em
\cdot
\sum_{I\subset \{1, \dots, n\} } \oneori{n}\, (-1)^{\varepsilon(I,I^c)}\, 
\supinprsym_E\Bigl(
{(-1)^{n\parity I}\,\conjugate^{\parity{I^c}}\bigl(\chi_{a,I}(x^a)\bigr)}
\,,\,
{\conjugate^n\bigl(\psi_{a,I^c}(x^a)\bigr)}\Bigr)
\mapob.
\end{align*}
It now suffices to note that \recalf{parityofchiaIintermsofparitychianda} implies that, for homogenous $\chi$, we have
\begin{align*}
(-1)^{n\parity I}\,\conjugate^{\parity{I^c}}\bigl(\chi_{a,I}(x^a)\bigr)
&
=
(-1)^{n\parity I}\,(-1)^{\parity{I^c}\cdot (\parity\chi+\parity I)} \, \chi_{a,I}(x^a)
\\&
=
(-1)^{n\parity\chi}\,\conjugate^I{\bigl(\chi_{a,I}(x^a)\bigr)}
\mapob.
\end{align*}
Substituting this in our result for $\homsuperinprod{\Fourierodd_{B,\mfdmetric} \chi}{\Fourierodd_{B,\mfdmetric} \psi}{_{B,\mfdmetric}}$ and comparing with \recalf{expressionsupinprodforFourieroddisSHSequivalence} gives the announced result. 
\end{preuve}

\mysubsection{BerezinFourierandHodgestarsection}{An intermezzo on the Hodge-star operation}

We now recall that in \recals{BatchelorbundlewithHodgestarsection} we have used a Batchelor atlas to define a vector bundle $V\!M\to \body M$ and an isomorphism $\sigma$ between $C^\infty_{(c)}(M; \CA^\KK)$ and $\Gamma^\infty_{(c)}\bigl(\,\bigwedge V\!M\to \body M\,\bigr)$. 
Moreover, a super metric $\mfdmetric$ on $M$ induces an ordinary metric $\mfdmetric_{\body M}$ on $\body M$ and an ordinary metric $\mfdmetric_{V\!M}$ on the fibers of $V\!M$, which allowed us to define a metric $\inprodsym$ on $\Gamma^\infty_c\bigl(\,\bigwedge V\!M\to \body M\,\bigr)$. 
This metric was then transported to $C^\infty_c(M; \CA^\KK)$ and used to create the Hilbert space $\Hilbert$ as its completion. 
We now want to show that we can go one step further: we will prove that the Berezin-Fourier transform $\Fourierodd_{B,\mfdmetric}$ on $C^\infty_{(c)}(M; \CA^\KK)$ is closely related to the Hodge-star operation on $\Gamma^\infty\bigl(\,\bigwedge V\!M\to \body M\,\bigr)$. 
In order to do so, we start by recalling the definition of the Hodge-star operation on $\bigwedge V\!M$ (see also \cite[p79]{Wa83}). 

As we assume that $M$ is \ood{} and that the Batchelor atlas is \ood{}, it follows that the fibers of $V\!M$ are oriented. 
Then, if $v_1, \dots, v_n$ is an orthonormal oriented basis of $V_mM$ (the fiber of $V\!M$ above $m\in \body M$), then the Hodge-$*$ operation on $\bigwedge V_mM$ is defined by
\begin{moneq}[defofHodgestaronbasisNEW]
*(v^I)
=
(-1)^{\varepsilon(I,I^c)}\cdot v^{I^c}
\mapob,
\end{moneq}
for any $I\subset \{1, \dots, n\}$. 
As the monomials $v^I$ form an orthonormal basis of $\bigwedge V_mM$, the metric $\mfdmetric_{V\!M}(m)$ is given by
$$
\mfdmetric_{V\!M}(m)\Bigl(\ \sum_{I\subset \{1, \dots, n\}} a_I\,v^I\ ,\ \sum_{J\subset \{1, \dots, n\}} b_J\,v^J\ \Bigr)
=
\sum_{I\subset \{1, \dots, n\}} {a_I}\cdot b_J
\mapob.
$$
It then follows immediately that the $*$-operation is unitary with respect to this metric:
$$
\forall A,B\in \bigwedge V_mM
\quad:\quad
\mfdmetric_{V\!M}(m)(*A,*B) = \mfdmetric_{V\!M}(m)(A,B)
\mapob.
$$
The Hodge-$*$-operation extends directly to an operation on sections $s\in\Gamma^\infty\bigl(\,\bigwedge V\!M\to \body M\,\bigr)$ by
$$
(*s)(m) = *\bigl(s(m)\bigr)
\mapob.
$$
As $*$ is unitary with respect to $\mfdmetric_{V\!M}(m)$, the induced operation on $C^\infty_c(M; \CA)$ is unitary too with respect to the metric \recalf{firstdefofordinarymetriconsuperfunctionsinE}.

\bigskip

Now let $(E,\inprodsym_E, \supinprsym_E)$ be a proto super Hilbert space over $\KK$ and let (for $k=0, \dots, n$)  $J_k:E\to E$ be a continuous linear map (over $\KK$). 
With these ingredients we define the twisted $*$-operation $\Hstarh$ on $\bigwedge V\!M \otimes E$ by
$$
A\otimes e\in \bigwedge{}\!\!^k\,V\!M \otimes E
\quad\Rightarrow\quad
\Hstarh(A\otimes e) = *(A) \otimes J_k(e) \in \bigwedge{}\!\!^{n-k}\,V\!M \otimes E
\mapob.
$$
More precisely:
$$
\Hstarh(v^I\otimes e) = 
*(v^I) \otimes J_{\parity I}(e)
\equiv
(-1)^{\varepsilon(I,I^c)}\cdot v^{I^c} \otimes J_{\parity I}(e)
\mapob,
$$
where, as before, we assume that $v_1, \dots, v_n$ is an orthonormal oriented basis of $V_mM$.

\begin{proclaim}{Lemma}
The twisted $*$-operation $\Hstarh$ is well defined on the bundle $\bigwedge V\!M \otimes E\to M$. 
Moreover, if the maps $J_k$ are bijective, then so is $\Hstarh$. 

\end{proclaim}

\begin{preuve}
The transition functions for the tensor product bundle $\bigwedge V\!M \otimes E$ are the tensor products of the transition functions for the bundle $\bigwedge V\!M$ with the identity map. It then follows immediately from the linearity of the $J_k$ that the restriction of $\Hstarh$ is well defined on the subbundle $\bigwedge^kV\!M\otimes E$, hence on the whole bundle.  
As the restriction of $\Hstarh$ to the subspace $\bigwedge^kV\!M\otimes E$ is the (tensor product) map $*\otimes J_k$, it follows immediately that $\Hstarh$ is bijective when the $J_k$ are. 
\end{preuve}

As for the ordinary Hodge-$*$ operation, the twisted $*$ operation induces an operation on sections $\psi\in \Gamma(\bigwedge V\!M\otimes E\to M)$. 
And because we assume that the maps $J_k$ are continuous, it follows that $\Hstarh$ preserves the subspace of smooth sections: 
$$
\Hstarh:\Gamma^\infty_{(c)}(\bigwedge V\!M\otimes E\to \body M) \to \Gamma^\infty_{(c)}(\bigwedge V\!M\otimes E\to \body M)
\mapob.
$$

\mysubsection{BerezinFourierandHodgestarsection}{Berezin-Fourier versus (twisted) Hodge-$*$}

\begin{proclaim}[linksuperinprodHodgeandFFourier]{Proposition}
Let $(E,\inprodsym_E, \supinprsym_E)$ be a proto super Hilbert space over $\KK$, let $M$ be an \ood{} $\CA$-manifold of graded dimension $d\vert n$ and let $\mfdmetric$ be a super metric on $M$. 
Suppose furthermore that $J:E\to E$ is a continuous linear map such that we have
$$
\forall e,f\in E
\quad:\quad
\homsuperinprod ef{_E} = \inprod[2]{J(e)}f_E
\mapob.
$$
Then for $\chi,\psi\in C^\infty_c(M;E\otimes\CA^\KK)$ we have
$$
\homsuperinprod\chi\psi{_{B,\mfdmetric}}
=
\inprod[2]{\Hstarh\bigl(\sigma(\chi)\bigr)}{\sigma(\psi)}_\mfdmetric
\mapob,
$$
where $\supinprsym_{B,\mfdmetric}$ is the super scalar product introduced in \recalt{defofsuperscalarprodforfunctionsonoodM}, where $\inprodsym_\mfdmetric$ is the metric on $\Gamma^\infty_c(M;E\otimes\CA^\KK)$ introduced in \recalf{firstdefofordinarymetriconsuperfunctionsinE} (see also \recalt{defofordinarmetriconCinftyMbymfdmetric}), and where $\Hstarh$ is the twisted Hodge-$*$ operation associated to the sequence of linear maps
$
J_k = \conjugate^n \scirc J\scirc \conjugate^{k}
$. 
If in addition $J$ preserves $\inprodsym_E$, \ie, $\inprod{J(e)}{J(f)}_E = \inprod ef_E$ for all $e,f\in E$, then $\Hstarh$ preserves $\inprodsym_\mfdmetric$ on $\Gamma^\infty_c(M;E\otimes\CA^\KK)$. 

In the special case $E=\KK$ we have $\Hstarh=*$ and in the (more) special case $E=\CC$ we have for homogeneous $\chi\in \Gamma^\infty_c(M;\CA^\CC)$ the equality
$$
\sigma\bigl( \Fourierodd_{B,\mfdmetric}(\chi)\bigr) = \oneormi{n-\parity\chi\,}\cdot *\bigl( \sigma(\chi)\bigr) 
\mapob.
$$

\end{proclaim}

\begin{preuve}
Let $U_a$ be a chart with local coordinates $(x^a, \xi^a)$ in an \ood{} Batchelor atlas adapted to the super metric $\mfdmetric$ \recalt{Batcheloratlasadaptedtometric}. 
For any $\chi\in C^\infty(M; \CA^\KK)$ with local expression on the chart $U_a$
$$
\chi_a(x^a,\xi^a) = \sum_{I\subset \{1, \dots, n\} } (\xi^a)^I\, \chi_{a,I}(x^a)
$$
the local expression of $(\sigma\mo\scirc\Hstarh\scirc\sigma)(\chi)$ in the chart $U_a$ is then given by
\begin{align*}
\bigl((\sigma\mo\scirc\Hstarh \scirc\sigma)(\chi)\bigr)_a(x^a,\xi^{a})
&
=
\sum_{I\subset \{1, \dots, n\} } (-1)^{\varepsilon(I,I^c)}\,(\xi^a)^{I^c} \, J_{\parity I}\bigl(\chi_{a,I}(x^a)\bigr)
\mapob.
\end{align*}
Using \recalt{metriconsuperfunctionsinEinlocalcoordinates} it then follows that we have
\begin{align*}
\inprod[2]{\Hstarh\bigl(\sigma(\chi)\bigr)}{\sigma(\psi)}_\mfdmetric
&
=
\sum_{a\in I} 
\int_{\body U_a} \extder\Leb^{(d)}(x^a)\ 
\rho_a(x^a)\cdot
\sqrt{\bigl\vert\det\bigl(\body\mfdmetric^{a,00}(x^a)\bigr)\bigr\vert}
\\&
\kern7em
\cdot
\sum_{I\subset \{1, \dots, n\} } (-1)^{\varepsilon(I,I^c)}\, 
\inprod{J_{\parity I}\bigl(\chi_{a,I}(x^a)\bigr)}{\psi_{a,I^c}(x^a)}_E
\end{align*}
We now recall that we have $\inprod{E_0}{E_1}_E=0$ and thus, as in the proof of \recalt{superscalarBgiscontinuous}, the conjugation map $\conjugate$ preserves $\inprodsym_E$. 
We thus have:
\begin{align*}
\inprod[2]{J_{\parity I}\bigl(\chi_{a,I}(x^a)\bigr)}{\psi_{a,I^c}(x^a)}_E
&
=
\inprod[2]{(\conjugate^{n}\scirc J\scirc \conjugate^{\parity I})\bigl(\chi_{a,I}(x^a)\bigr)}{\psi_{a,I^c}(x^a)}_E
\\&
=
\inprod[2]{(J\scirc \conjugate^{\parity I})\bigl(\chi_{a,I}(x^a)\bigr)}{\conjugate^{n}\bigl(\psi_{a,I^c}(x^a)\bigr)}_E
\\&
=
\homsuperinprod[3]{\conjugate^{\parity I}\bigl(\chi_{a,I}(x^a)\bigr)}{\conjugate^{n}\bigl(\psi_{a,I^c}(x^a)\bigr)}{_E}
\mapob.
\end{align*}
Comparison with \recalt{descriptionsuperscalarproductBginadaptedBatlas} then gives the desired result. 

Let us now assume that $J$ preserves $\inprodsym_E$. 
Then we find, still using our \ood{} Batchelor atlas adapted to $\mfdmetric$:
\begin{align*}
\inprod{\Hstarh\chi}{\Hstarh\psi}_\mfdmetric
&
=
\sum_{a\in I} 
\int_{\body U_a} \extder\Leb^{(d)}(x^a)\ 
\rho_a(x^a)\cdot
\sqrt{\bigl\vert\det\bigl(\body\mfdmetric^{a,00}(x^a)\bigr)\bigr\vert}
\\&
\kern5em
\cdot
\sum_{I\subset \{1, \dots, n\} } 
\inprod{J_{\parity I}\bigl(\chi_{a,I}(x^a)\bigr)}{J_{\parity I}\bigl(\psi_{a,I}(x^a)\bigr)}_E
\\&
=
\sum_{a\in I} 
\int_{\body U_a} \extder\Leb^{(d)}(x^a)\ 
\rho_a(x^a)\cdot
\sqrt{\bigl\vert\det\bigl(\body\mfdmetric^{a,00}(x^a)\bigr)\bigr\vert}
\\&
\kern5em
\cdot
\sum_{I\subset \{1, \dots, n\} } 
\inprod{(\conjugate^{n}\scirc J\scirc \conjugate^{\parity I})\bigl(\chi_{a,I}(x^a)\bigr)}{(\conjugate^{n}\scirc J\scirc \conjugate^{\parity I})\bigl(\psi_{a,I}(x^a)\bigr)}_E
\\&
=
\sum_{a\in I} 
\int_{\body U_a} \extder\Leb^{(d)}(x^a)\ 
\rho_a(x^a)\cdot
\sqrt{\bigl\vert\det\bigl(\body\mfdmetric^{a,00}(x^a)\bigr)\bigr\vert}
\\&
\kern5em
\cdot
\sum_{I\subset \{1, \dots, n\} } 
\inprod{\chi_{a,I}(x^a)}{\psi_{a,I}(x^a)}_E
=
\inprod{\chi}{\psi}_\mfdmetric
\end{align*}
as claimed. 

\medskip

In the special case $E=\KK$ we have $J=\mathrm{id}$ and, because then $E$ has only even elements, $\conjugate=\mathrm{id}$, so all $J_k=\mathrm{id}$ and thus $\Hstarh = *$ as claimed. 
Remains the link between $\Fourierodd_{B,\mfdmetric}$ and $*$ in the special case $E=\CC$. 
Still using our \ood{} Batchelor atlas adapted to $\mfdmetric$, we apply \recalf{localexpressionBerezinFourierinadaptedchart} to obtain
\begin{moneq}[localexpressionFourierBgpsi]
(\Fourierodd_{B,\mfdmetric}\chi)_a(x,\xi) = 
\sum_{I\subset \{1, \dots, n\} } 
\oneormi{\,\parity{I^c}\,}\, (-1)^{\varepsilon(I,I^c)}\,\xi^{I^c}\,\chi_{a,I}(x)
\mapob.
\end{moneq}
On the other hand, according to the definition of $\sigma$ \recalf{tempdefofsigmasubflinkfuncwithBatchNEW}, the section $\sigma(\psi)$ is given in the local trivialization associated to the local chart $U_a$ by
\begin{moneq}
\bigl(\sigma(\psi)\bigr)_a(x)
=
\Bigl(\ x
\ ,\ 
\sum_{I\subset \{1,\dots, n\}} \chi_{a,I}(x)\,e^I
\ \Bigr)
\mapob,
\end{moneq}
where $e_1, \dots, e_n$ denotes the canonical basis of the typical fiber $\RR^n$. 
Moreover, in this same trivialization the metric $\mfdmetric_{V\!M}$ is given by the inverse of the matrix $H^a$, which is the identity matrix. 
In other words, the local trivializing sections given by the canonical basis $e_1, \dots, e_n$ is orthonormal. 
It then follows that the Hodge-$*$ operation is given on this local trivialization by
$$
\Bigl(*\bigl(\sigma(\chi)\bigr)\Bigr)_a(x) 
=
\Bigl(\ x
\ ,\ 
\sum_{I\subset \{1,\dots, n\}} (-1)^{\varepsilon(I,I^c)}\,\chi_{a,I}(x)\,e^{I^c}
\ \Bigr)
\mapob.
$$
And when we compute the local expression for $\sigma\bigl(\Fourierodd_{B,\mfdmetric}(\chi)\bigr)$ using \recalf{localexpressionFourierBgpsi} we find
$$
\Bigl(\sigma\bigl(\Fourierodd_{B,\mfdmetric}(\chi)\bigr)\Bigr)_a(x)
=
\Bigl(\ x
\ ,\ 
\sum_{I\subset \{1,\dots, n\}} \oneormi{\,\parity{I^c}\,}\,(-1)^{\varepsilon(I,I^c)}\,\chi_{a,I}(x)\,e^{I^c}
\ \Bigr)
\mapob.
$$
Comparing these two expressions, we see that the only difference is a factor $\oneormi{\,\parity{I^c}\,}$ under the summation sign. 
When we assume that $\psi$ has parity $\alpha$, this means that in the summation over $I\subset \{1, \dots, n\}$ only subsets $I$ appear with $\parity I = \alpha$ (more precisely, the coefficients of the $\xi^I$ will be zero for those $I$ with $\parity I\neq\alpha$), and thus $\parity{I^c} = n-\alpha$. 
It follows that in that case the factor $\oneori{\,\parity{I^c}\,}$ is independent of $I$, showing that for $\psi$ with fixed parity $\alpha$ we have the equality
$$
\Bigl(\sigma\bigl(\Fourierodd_{B,\mfdmetric}(\psi)\bigr)\Bigr)_b(x)
=
\oneormi{n-\alpha} \cdot \Bigl(*\bigl(\sigma(\psi)\bigr)\Bigr)_b(x)
\mapob.
$$
As the local charts $U_b$ cover $M$, it follows that we have the global equality $\sigma\scirc \Fourierodd_{B,\mfdmetric} = \oneormi{n-\alpha}\cdot *\scirc \sigma$ on $C^\infty(M;\CA^\CC)_\alpha$ as announced.
\end{preuve}

\begin{definition}[linkingsupscprtometricviaKrein]{Discussion}
We can summarize a part of these results as follows. 
The spaces $\RR$ and $\CC$ are naturally super Hilbert spaces with $\inprodsym=\supinprsym$. 
Moreover the metric and the super scalar product are related by a bijective linear map preserving $\inprodsym$, viz. the identity map. 
And then we have shown that if $(E,\inprodsym_E, \supinprsym_E)$ is a proto super Hilbert space in which $\inprodsym_E$ and $\supinprsym_E$ are related by a bijective linear map preserving $\inprodsym_E$, then so is the space $\Gamma^\infty_c(M; E\otimes \CA^\KK)$ of compactly supported smooth function on $M$ with values in $E\otimes \CA^\KK$ (where $M$ is supposed to be equipped with a super metric $\mfdmetric$). 
As the bijective linear map preserves the metric, it is bi-continuous. 
It then follows that the extension of the super scalar product $\supinprsym_{B,\mfdmetric}$ on $\Gamma^\infty_c(M; E\otimes \CA^\KK)$ to its Hilbert space completion (with respect to the metric $\inprodsym_\mfdmetric$) automatically remains non-degenerate. 
In other words, we have created a new super Hilbert space out of the space of compactly supported smooth functions on $M$ with values in a super Hilbert space in which again (as for $\KK$) the super scalar product and the metric are related by bijective  linear map preserving the metric. 
The subcategory of this kind of super Hilbert spaces is very close to the definition of a super Hilbert space given by A.~de Goursac in \cite[\S2.4]{BDGT:2012} (a definition inspired by the notion of a Krein space) and further elaborated in \cite{deGoursacMichel:2015}. 
Unfortunately, already for the left-regular representation of a super Lie group the natural invariant super scalar product is not $\supinprsym_{B,\mfdmetric}$ but $\supinprsym_\mfdmetric$ and there is no obvious natural way to create such a link between the metric $\inprodsym_\mfdmetric$ and this super scalar product $\supinprsym_\mfdmetric$. 
For $\CA$-Lie groups $G$ there actually is a (homogeneous) bi-continuous map $\superspmatrix :C^\infty_c(G; \CA^\KK)\to C^\infty_c(G; \CA^\KK)$ (see \recalt{linkingsupscprtometriconsuperLiegroupBUTnotunit}) that relates the metric $\inprodsym_\mfdmetric$ to the super scalar product 
$\supinprsym_\mfdmetric$ (for a left-invariant metric $\mfdmetric$ on $G$), but in general this map does not preserve the metric, as can be seen in the example of the group $\OSp(1,2)$ \recals{OSp12section}. 

\end{definition}

\masection{Left-regular representations}
\label{leftregularrepresentationssection}

For an ordinary (non super) Lie group $G$, the left-regular representation can be defined as follows. 
One starts with a metric at $T_eG\cong \Liealg g$, extends it to a left-invariant metric $\mfdmetric$ on $G$ and determines the associated left-invariant (Haar) measure $\nu_o$ on $G$. 
This determines the pre-Hilbert space $C_c^\infty(G;\KK)$ of compactly supported smooth functions on $G$ equipped with the scalar product
\begin{moneq}[usualmetriconL2G]
\inprodd{\chi}{\psi} = \int_G \overline\chi\cdot \psi\ \nu_o
\mapob.
\end{moneq}
The associated Hilbert space is its completion and gives $\Hilbert = L^2(G;\KK;\nu_o)$. 
If no confusion is possible, we will abbreviate this as $\Hilbert = L^2(G)$. 
The representation $\rho$ on $\Hilbert$ is defined as usual by the formula 
$$
\bigl(\rho(g)\psi\bigr)(h) = \psi(g\mo h)
\mapob,
$$
and the fact that $\mfdmetric$ and especially $\nu_o$ are left-invariant then guarantees that this representation is unitary. 

For a super Lie group $G$ of dimension $d\vert n$ we will follow the same scheme. 
We start with a super metric on $T_eG\cong \Liealg g$, we extend it by left-translation to a left-invariant super metric $\mfdmetric$ on $G$ and we determine the associated trivializing density $\nu_\mfdmetric$ on $G$, as well as the ordinary metrics $\mfdmetric_{\body G}$ on $\body G$ and $\mfdmetric_{VG}$ on the fibers of the vector bundle $VG\to \body G$ (see \recals{BatchelorbundlewithHodgestarsection}). 
The first allows us to define the super scalar product $\supinprsym_\mfdmetric$ on $C_c^\infty(G; \CA^\KK)$ and the last two allow us to define the metric $\inprodsym_\mfdmetric$ on $C_c^\infty(G; \CA^\KK)$, turning the latter into a pre-Hilbert space. 
We also have the representation $\rho$ of $G$ on $C_{(c)}^\infty(G;\CA^\KK) \otimes \CA^\KK$ defined by (see \recalt{therepresentationofagrouponCinftyM})
$$
\bigl(\rho(g)\psi\bigr)(h) = \psi(g\mo h)
$$
and according to \recalt{superscalarproductmfdmetricinvariantunderGaction} the super scalar product $\supinprsym_\mfdmetric$ (extended to $C_{c}^\infty(G;\CA^\KK) \otimes \CA^\KK$) is invariant under this action. 

In order to prove that this can be turned into a super unitary representation of $G$, we need some more information on the structure of $C_{(c)}^\infty(G;\CA^\KK)$. 
For that we will rely heavily upon the identification\slash diffeomorphism $\Phi:\wod G \times \oddp{\Liealg g}_0 \to G$, $\Phi(g,X)=g\cdot \exp(X)$ given in \recalf{diffeoGwithGwodtimesoddLiealgg1}. 
And in order to \myquote{simplify} notation, we will most of the time not write the symbol $\Phi$. 
Th first consequence is that $G$ is \ood{} and thus that we indeed are allowed to apply \recalt{canonicaltrivdensonOOD} to obtain $\nu_\mfdmetric$. 

We now fix (once and for all) a basis $e_1, \dots, e_d$ of $\body\Liealg g_0$ and a basis $f_1, \dots, f_n$ of $\body \Liealg g_1$. 
This will give us global odd coordinates $\xi_1, \dots, \xi_n$ on $\oddp{\Liealg g}_0$ via $X=\sum_{i=1}^n \xi_i\,f_i$. 
As we have a direct product, it follows immediately that any smooth function $f:G\cong \wod G \times \oddp{\Liealg g}_0\to \CA^\CC$ determines $2^n$ unique smooth functions $f_I:\body G\to \CC$, $I\subset \{1, \dots, n\}$ such that we gave
$$
f(g,\xi) = \sum_{I\subset \{1, \dots, n\} } \xi^I\cdot (\Gextension f_I)(g)
\mapob.
$$
We thus have an identification
\begin{moneq}[identificationCinftyGwithproductsofCinftybodyG]
C^\infty_{(c)}(G;\CA^\CC) \cong \bigl(\,C^\infty_{(c)}(\body G)\,\bigr)^{2^n}
\qquad,\qquad
f \cong (f_I)_{I\subset \{1, \dots, n\} }
\mapob.
\end{moneq}
In order to determine the (ordinary) scalar product associated to a left-invariant metric $\mfdmetric$ on $G$, we need a more precise control over the left and right-invariant vector fields on $G$ in terms of our product structure $G\cong \wod G \times \oddp{\Liealg g}_0$. 
For that we will need a refinement of \recalt{separatingevenandoddinproductofexponentials} in the special case when $Y$ contains a single odd parameter.

\begin{definition}{Definition}
We define the functions $f,h,b, b_-, b_+:\RR\to\RR$ by
\begin{gather*}
f(t) = \frac{\eexp^t - 1}t
\qquad,\qquad
h(t) = \frac{\eexp^t - 1}{\eexp^t+1}
\qquad,\qquad
b(t)
=
\frac{t(\eexp^t +1)}{2(\eexp^t - 1)}
\\
b_-(t) = b(t) -\tfrac12\,t\,h(t) = \frac t{\sinh(t)}
\qquad,\qquad
b_+(t) = b(t) +\tfrac12\,t\,h(t) = \frac{t\,\cosh(t)}{\sinh(t)} 
\mapob.
\end{gather*}
All five have convergent power series in a neighborhood of $0\in \RR$. Moreover $h$ is an odd function, whereas $b$, $b_\pm$ are even functions.
Their leading terms are given as
\begin{gather*}
h(t) = \tfrac12\,t - \tfrac1{24}\,t^3 + \cdots
\qquad,\qquad
b(t) = 1+ \tfrac1{12}\,t^2 - \tfrac1{720}\,t^4+\cdots
\\ 
b_-(t) = 1 - \tfrac16\,t^2+ \tfrac7{360}\,t^4 + \cdots
\qquad,\qquad
b_+(t) = 1+\tfrac13\,t^2 - \tfrac1{45}\,t^4 +\cdots
\mapob.
\end{gather*}
For future use we note that we have the following relations among these functions:
\begin{moneq}[relationsamongbfandhoft]
b(t) = \frac1{f(t)} +\tfrac12 t 
\qquad,\qquad
b(t)\cdot h(t) = \tfrac12\,t
\qquad,\qquad
b_+(t) = b(2t)
\mapob.
\end{moneq}

\end{definition}

\begin{proclaim}[separatinggwodandoddpginexponentials]{Lemma}
Let $X\in \Liealg g_0$, $Y\in \Liealg g_1$ and $\tau\in \CA_1$. Then we have the following equalities:
\begin{align*}
\exp(X+\tau Y) 
&
=
\exp(X)\cdot \exp\Bigl(\tau\, f\bigl(-\ad(X)\bigr)Y \Bigr) 
\\&
= 
\exp\Bigl(\tau\, f\bigl(\ad(X)\bigr)Y \Bigr) \cdot \exp(X)
\mapob,
\\
\noalign{\vskip3\jot}
\exp(X) \cdot \exp(\tau Y) 
&= 
\exp\Bigl( X +  \tfrac12\,\tau\,[X,Y] + \tau\, b\bigl(\ad(X)\bigr)Y \Bigr)
\\&
=
\exp\Bigl( \tau\, h\bigl(\ad(X)\bigr)Y \Bigr)   \cdot
\exp\Bigl( X + \tau\, b_+\bigl(\ad(X)\bigr)Y  \Bigr )
\\
\noalign{\vskip3\jot}
\exp(\tau Y) \cdot \exp(X) 
&= 
\exp\Bigl( X - \tfrac12 \,\tau\,[X,Y] + \tau\, b\bigl(\ad(X)\bigr)Y \Bigr)
\\&
=
\exp\Bigl( -\tau\, h\bigl(\ad(X)\bigr)Y\Bigr)\  \cdot
\exp\Bigl(X + \tau\, b_-\bigl(\ad(X)\bigr)Y \Bigr )
\mapob.
\end{align*}

\end{proclaim}

\begin{preuve}
These results can be proven by using the Baker-Campbell-Hausdorff formula, but it is faster to use the derivative of the exponential map (see \cite[VI.3.15]{Tu04} or \cite[p.24]{DK00} for proof in the non-super setting), simply because any Taylor expansion with an odd parameter stops at order two ($\tau^2=0$). 
We thus can compute:
\begin{align*}
\exp(X+\tau Z) 
&= 
\exp(X) + (T_X\exp)(\tau Z)
\\&
=
\exp(X) + (T_e L_{\exp(X)})\Bigl( \frac{1-\eexp^{-\ad(X)}}{\ad(X)} \Bigr)(\tau Z)
\\&
=
(L_{\exp(X)}\scirc \exp)(0) + (T_e L_{\exp(X)}\scirc T_{0}\exp)\Bigl( \bigl(f(-\ad(X)\bigr)(\tau Z) \Bigr)
\\&
=
(L_{\exp(X)}\scirc \exp)\Bigl( 0+\tau\bigl(f(-\ad(X)\bigr) Z \Bigr)
\\&
=
\exp(X)\cdot \exp\Bigl( \tau f\bigl(-\ad(X)\bigr) Z \Bigr)
\mapob,
\end{align*}
which is the first equality to be proven.
Taking inverses on both sides (and replacing $X$ and $Z$ by their opposites) one finds the second equality.

We now recall that if $A$ is an endomorphism of a finite dimensional vector space and if $f$ is a polynomial or a converging power series in a single variable, \ie, $f(t)=\sum_k c_k\,t^k$, then we can define the endomorphism $f(A)$ simply by $f(A)v= \sum_k c_k\,A^k(v)$. Moreover, if $f$ and $g$ are two such functions, we have the equality
$$
f(A)\scirc g(A) = (f\cdot g)(A)
\mapob.
$$
We will use this property repeatedly in the upcoming computations.

Writing $f\bigl(-\ad(X)\bigr)Z = Y$ we have (using \recalf{relationsamongbfandhoft})
$$
Z = \Bigl(\,\frac1{f\bigl(-\ad(X)\bigr)}\,\Bigr)Y = \Bigl(\tfrac12 \,\ad(X)+ b\bigl(\ad(X)\bigr) \Bigr)Y
$$
and thus
\begin{align}
\exp(X) \cdot \exp(\tau Y)
&=
\exp\Bigl(X + \tfrac12\, \tau\, \ad(X)(Y) + \tau\, b\bigl(\ad(X)\bigr)(Y)\Bigr)
\notag
\\&
=
\exp\Bigl(X + \tfrac12 \,\tau\, [X,Y] + \tau\, b\bigl(\ad(X)\bigr)(Y)\Bigr)
\mapob,
\label{firstBCHwithonetau}
\end{align}
which is the third equality to be proven. 
Taking inverses and using that $b$ is an even function we also find
\begin{moneq}[secondBCHwithonetau]
\exp(\tau Y) \cdot \exp(X)
=
\exp\Bigl(X - \tfrac12\, \tau\, \ad(X)(Y) + \tau\, b\bigl(\ad(X)\bigr)(Y)\Bigr)
\mapob,
\end{moneq}
which is the fifth equality.

Writing, with $\varepsilon=\pm1$,
\begin{align*}
X_\varepsilon 
&
= X + \tau \Bigl(\, b\bigl(\ad(X)\bigr)(Y) - \tfrac12\,\varepsilon \, h\bigl(\ad(X)\bigr) (\,[X,Y]\,)\,\Bigr)
\\
Y_\varepsilon 
&
= \varepsilon\, h\bigl(\ad(X)\bigr)(Y)
\end{align*} 
we have, according to \recalf{firstBCHwithonetau},
$$
\exp(X_\varepsilon)
\cdot \exp(\tau Y_\varepsilon)
=
\exp(Z_\varepsilon)
$$
with
\begin{align*}
Z_\varepsilon
&
=
X_\varepsilon + \tfrac12 \,\tau\, \ad(X_\varepsilon)Y_\varepsilon + \tau \,b(\ad(X_\varepsilon))Y_\varepsilon
\\[2\jot]
\text{\tiny ($\tau^2=0$)}\qquad
&
=
X_\varepsilon + \tfrac12\, \tau\, \ad(X)Y_\varepsilon + \tau\, b(\ad(X))Y_\varepsilon
\\[2\jot]
&
=
\Bigl(X + \tau\, b\bigl(\ad(X)\bigr)Y -  \tfrac12\, \varepsilon\, \tau\, h\bigl(\ad(X)\bigr) [X,Y]\Bigr)
\\
&
\qquad\qquad
+ 
\tfrac12\,\varepsilon \, \tau\, \ad(X) h\bigl(\ad(X)\bigr)Y
+
\varepsilon\, \tau\, b\bigl(\ad(X)\bigr) h\bigl(\ad(X)\bigr)Y
\\[2\jot]&
=
X + \tau\Bigl( b\bigl(\ad(X)\bigr) 
+
\varepsilon\,  b\bigl(\ad(X)\bigr)\scirc h\bigl(\ad(X)\bigr)\Bigr)Y
\\[2\jot]&
\oversetalign{\recalf{relationsamongbfandhoft}}\to=
\ 
X + \tau\Bigl( b\bigl(\ad(X)\bigr)Y + \tfrac12\, \varepsilon\, \ad(X)Y \Bigr) 
\mapob.
\end{align*}
Comparing this with \recalf{firstBCHwithonetau} and \recalf{secondBCHwithonetau} we thus find:
\begin{align*}
\exp(X_{+1})\cdot\exp(\tau Y_{+1})
&
=
\exp(X)\cdot \exp(\tau Y)
\\
\exp(X_{-1})\cdot\exp(\tau Y_{-1})
&
=
\exp(\tau Y)\cdot \exp(X)
\mapob.
\end{align*}
The second of these equalities is the sixth equality to be proven, whereas taking inverses in the first gives us the fourth equality to be proven.
\end{preuve}

\begin{proclaim}[leftinvariantvfonGwodtimesgsub1]{Lemma}
In the decomposition $G = \wod G \times \oddp{\Liealg g}_0$ given by the map $(g,v)\mapsto g\cdot \exp(v)$, we have the natural identification of the tangent space as 
$$
T_{(g,v)}\bigl(\wod G \times \oddp{\Liealg g}_0\bigr)\cong T_{g}\wod G \times \oddp{\Liealg g}
\mapob.
$$
In this identification the left- and right-invariant vector fields take the following form.
\begin{enumerate}
\item\label{leftandrightinvvfintrivializationofGintoGwodalgopddp}
For $X\in \body T_e\wod G\equiv \body \Liealg g_0$ (and thus even) we have
\begin{align*}
TL_{(g,v)} X\caprestricted_{(e,0)} 
&
= 
\bigl(\,TL_{g}X\caprestricted_{e} \,,\, [v,X]\caprestricted_v \,\bigr)
\cong
TL_{g}X\caprestricted_{e} + [v,X]\caprestricted_v 
\\
TR_{(g,v)} X\caprestricted_{(e,0)} 
&
= 
\bigl(\,TR_{g}X\caprestricted_{e} \,,\, 0\caprestricted_v \,\bigr)
\cong
TR_{g}X\caprestricted_{e} + 0\caprestricted_v 
\mapob.
\end{align*}

\item
For $Y\in \body \oddp{\Liealg g} \equiv \body \Liealg g_1$ (and thus odd) we have
\begin{align*}
TL_{(g,v)} Y\caprestricted_{(e,0)} 
&
= 
\Bigl(\  TL_{g}\bigl(h\bigl(\ad(v)\bigr)Y\caprestricted_e\,\bigr) \ ,\  b_+\bigl(\ad(v)\bigr)Y  \ \Bigr)
\\&
\cong
TL_{g}\Bigl(h\bigl(\ad(v)\bigr)Y\caprestricted_e\Bigr) + \Bigl(b_+\bigl(\ad(v)\bigr)Y \Bigr)\bmidrestricted_v
\\
TR_{(g,v)} Y\caprestricted_{(e,0)} 
&
= 
\Bigl(\  -TR_{g}\bigl(h\bigl(\ad(\Ad(g)v)\bigr)Y\caprestricted_e\,\bigr) \ ,\  \bigl(b_-\bigl(\ad(v)\bigr)\Ad(g\mo)Y \,\bigr) \ \Bigr)
\\&
\cong
-{TR_{g}}\Bigl(h\bigl(\ad(\Ad(g)v)\bigr)Y\caprestricted_e\Bigr)+ \Bigl(b_-\bigl(\ad(v)\bigr)\Ad(g\mo)Y\Bigr)\bmidrestricted_v
\mapob.
\end{align*}

\end{enumerate}
Note that these formul{\ae} make sense: for $X\in \body\Liealg g_0$ and $v\in \oddp{\Liealg g}_0$ we have $[v,X]\in \oddp{\Liealg g}_0$. Moreover, as $h$ is odd and $b_\pm$ even, we have indeed, for any $Y\in \body\Liealg g_1\subset \oddp{\Liealg g}$, that $h\bigl(\ad(v)\bigr)Y$ belongs to $\wod{\Liealg g}$ and $b_\pm\bigl(\ad(v)\bigr)Y$ to $\oddp{\Liealg g}$.

\end{proclaim}

\begin{preuve}
\def\ddt{\frac{\extder}{\extder t}\bigrestricted_{t=0}\,}
For (i) we compute:
\begin{align*}
{TL_{(g,v)}} X\caprestricted_{(e,0)}
&
=
\ddt\Bigl( (g,v)\cdot\bigl(\exp(tX),0\bigr)\Bigr)
\\[2\jot]&
=
\ddt \Bigl(g \cdot \exp(v) \cdot \exp(tX)\Bigr)
\\[2\jot]&
=
\ddt\Bigl( g\cdot \exp(tX) \cdot \exp\bigl(\Ad\bigl(\exp(-tX)\bigr)v\bigr)\Bigr)
\\[2\jot]&
=
\ddt \Bigl(g\cdot \exp(tX)\ ,\ \exp\bigl(-\ad(tX)\bigr)v \Bigr)
\\[2\jot]&
=
{TL_{g}}X\caprestricted_e - \ad(X)v\restricted_v
=
{TL_{g}}X\caprestricted_e + [v,X]\caprestricted_v
\end{align*}
and
\begin{align*}
{TR_{(g,v)}} X\caprestricted_{(e,0)}
&
=
\ddt\Bigl( \bigl(\exp(tX),0\bigr)\cdot(g,v)\Bigr)
\\[2\jot]&
=
\ddt \Bigl(\bigl(\exp(tX)\cdot g,v\bigr)\Bigr)
\\[2\jot]&
=
\bigl({TR_{g}}X\caprestricted_e , 0\bigr)
\cong
{TR_{g}}X\caprestricted_e + 0\caprestricted_v
\end{align*}

For (ii) we compute:
\def\ddtau{{\frac{\extder}{\extder \tau}\,}}
\begin{align*}
{TL_{(g,v)}}Y\caprestricted_{(e,0)}
&
=
\ddtau \bigl((g,v)\cdot (e,{\tau Y})\bigr)
=
\ddtau \bigl(g\cdot \exp(v)\cdot \exp(\tau Y)\bigr)
\\[2\jot]&
\oversetalign{\recalf{separatinggwodandoddpginexponentials}}\to=
\ddtau \Bigl(g\cdot \exp\bigl(\tau h\bigl(\ad(v)\bigr)Y\bigr)  
\cdot \exp\bigl(v + \tau\,b_+\bigl(\ad(v)\bigr)Y \bigr)  \Bigr)
\\[2\jot]&
=
\ddtau\Bigl( g\cdot \exp(\tau h(\ad(v))Y) \ ,\ 
v + \tau\,b_+\bigl(\ad(v)\bigr)Y  \Bigr)
\\&
=
{TL_{g}}\Bigl(h\bigl(\ad(v)\bigr)Y\Bigr)\bmidrestricted_e 
+ \Bigl( b_+\bigl(\ad(v)\bigr)Y\Bigr)\bmidrestricted_v
\mapob.
\end{align*}
and 
\begin{align*}
\shifttag{3.5em}
{TR_{(g,v)}}Y\caprestricted_{(e,0)}
=
\ddtau \bigl( (e,{\tau Y})\cdot (g,v)\bigr)
=
\ddtau \bigl( \exp(\tau Y)\cdot g\cdot \exp(v)\bigr)
\\[2\jot]&
=
\ddtau \bigl( g\cdot \exp(\tau \Ad(g\mo)Y)\cdot \exp(v)\bigr)
\\[2\jot]&
\oversetalign{\recalf{separatinggwodandoddpginexponentials}}\to=
\ddtau \Bigl( \exp\bigl(-\tau h\bigl(\ad(\Ad(g)v)\bigr)Y\bigr)  
\cdot g\cdot \exp\bigl(v+\tau\,b_-\bigl(\ad(v)\bigr)\Ad(g\mo)Y \bigr)  \Bigr)
\\[2\jot]&
=
\ddtau\Bigl( \exp(-\tau h(\ad(\Ad(g)v))Y) \cdot g \ ,\ 
v+\tau \,b_-\bigl(\ad(v)\bigr)\Ad(g\mo)Y  \Bigr)
\\&
=
-{TR_{g}}\Bigl(h\bigl(\ad(\Ad(g)v)\bigr)Y\Bigr)\bmidrestricted_e 
+ \Bigl(b_-\bigl(\ad(v)\bigr)\Ad(g\mo)Y\Bigr)\bmidrestricted_v
\mapob.
\QEDici
\end{align*}
\smallskip
\end{preuve}

\begin{definition}{Notation}
In general the Lie algebra of a Lie group $G$ \myquote{is} the space of left-invariant vector fields on $G$, which can be identified with the tangent space at the identity. 
However, in the sequel we will want to distinguish this dual nature of an element of the Lie algebra: a (left-invariant) vector field on $G$ or an (abstract) element of an (abstract) vector space $\Liealg g \cong T_eG$. 
In order to distinguish these two aspects, we will denote by $\vec X$ the left-invariant vector field on $G$ whose value at the identity $e\in G$ is $X$: $\vec X\caprestricted_e = X$. 
And we will denote by $X^R$ the \emph{right-invariant} vector field on $G$ whose value at $e$ is $X$. 

\end{definition}

Now let $x_1, \dots, x_d$ be local even coordinates on $\wod G$, then, together with our global odd coordinates $\xi_1, \dots, \xi_n$ on $\oddp{\Liealg g}_0$, \ie, $(\xi_1, \dots, \xi_n)$, the couple $(x,\xi)$ becomes a local system of coordinates on $G\cong \wod G \times \oddp{\Liealg g}_0$. 
In terms of these local coordinates we now can express the left-invariant vector fields as follows. 
We start with the observation that there exists an invertible matrix $C(x)$ on the local chart for $\wod G$ such that for the basis $e_i$ of $\Liealg g_0$ we have (with $x$ being the local coordinates of $g$)
$$
\vec e_i\restricted_g
\equiv
TL_g(\vec e_i\restricted_e) = \sum_j \partial_{x_j}\,C_{ji}(x)
\mapob.
$$
There also exist matrices $A(\xi)$, $B(\xi)$ and $H(\xi)$ (only $B$ being necessarily square) depending only upon the odd coordinates such that
\begin{gather*}
[v,e_i] = \sum_j f_j\,A_{ji}(\xi)
\qquad,\qquad
b_+\bigl(\ad(v)\bigr)f_j = \sum_k f_k\,B_{kj}(\xi)
\\
\qquad\text{and}\qquad
h\bigl(\ad(v)\bigr)f_j = \sum_i e_i\,H_{ij}(\xi)
\mapob.
\end{gather*}
We now recall that in the identification $T_{(g,v)}\bigl(\wod G \times \oddp{\Liealg g}_0\bigr)\cong T_{g}\wod G \times \oddp{\Liealg g}$ the basis vector $f_i$ gets identified with the tangent vector $\partial_{\xi_i}\caprestricted_{(g,v)}$. 
Applying \recalt{leftinvariantvfonGwodtimesgsub1} then immediately gives
\begin{align*}
TL_{(g,v)}\vec e_i\restricted_{(e,0)} 
&
= \sum_j \partial_{x_j}\,C_{ji}(x) + \sum_k \partial_{\xi_k}\,A_{ki}(\xi)
\\
TL_{(g,v)}\vec f_j\restricted_{(e,0)} 
&
= 
\sum_j \partial_{x_k} \, C_{ki}(x)\,H_{ij}(\xi) + \sum_k \partial_{\xi_k} \, B_{kj}(\xi)
\mapob.
\end{align*}
We can write this in matrix form as
\begin{align}
TL_{(g,v)}
\begin{pmatrix} \vec{e}_\bullet\restricted_{(e,0)} & \partial_{\xi}\restricted_{(e,0)}
\end{pmatrix}
&
=
\begin{pmatrix} \partial_{x}\caprestricted_{(g,v)} & \partial_{\xi}\caprestricted_{(g,v)} \end{pmatrix}
\cdot
\begin{pmatrix} C(x) & C(x)\cdot H(\xi) \\ A(\xi) & B(\xi) \end{pmatrix}
\notag
\\&
=
\begin{pmatrix} \partial_{x}\caprestricted_{(g,v)} & \partial_{\xi}\caprestricted_{(g,v)} \end{pmatrix}
\cdot
\begin{pmatrix} C(x) & \mathbf0 \\ \mathbf0 & \oneasmatrix \end{pmatrix}
\cdot
\begin{pmatrix} \oneasmatrix & H(\xi) \\ A(\xi) & B(\xi) \end{pmatrix}
\label{expressionleftinvvfonGasproduct}
\\&
=
\begin{pmatrix} \vec{e}_\bullet\restricted_g & \partial_{\xi}\caprestricted_{(g,v)} \end{pmatrix}
\cdot
\begin{pmatrix} \oneasmatrix & H(\xi) \\ A(\xi) & B(\xi) \end{pmatrix}
\notag
\mapob.
\end{align}
We now use the \myquote{standard} super metric on $\Liealg g \cong T_{(e,0)}G$, \ie, $\mfdmetric_{(e,0)}(e_j,e_k) = \delta_{jk}$ and $\mfdmetric_{(e,0)}(f_j,f_k) = i\delta_{jk}$, which we extend by left-translation to a global left-invariant super metric $\mfdmetric$ on $G$. 
It then follows from the definition of the matrix $\mfdmetric_{ij}\restricted_{(x,\xi)} = \mfdmetric(\partial_i\restricted_{(x,\xi)},\partial_j\restricted_{(x,\xi)})$ (where $\partial_i$ denotes \myquote{all} partial derivatives\slash basis tangent vectors at $(x,\xi)\in G$) and \recalf{expressionleftinvvfonGasproduct} that this matrix satisfies the equation (compare with \recalf{changeofbasisforametricformula}, but beware: the definition of the matrix defining the change of basis is not the same, there is a change from left to right; see also \cite[II.4.10]{Tu04})
\begin{moneq}[expressionleftinvariantmetriconG]
\begin{gathered}
\begin{pmatrix} \oneasmatrix & \mathbf0 \\ \mathbf0 & i\,\oneasmatrix \end{pmatrix}
=\kern4em
\\[2\jot]
\begin{pmatrix} {}\trans C(x) & -\ {}\trans\! A(\xi) \\ {}\trans\! H(\xi)\cdot {}\trans C(x) & {}\trans\! B(\xi) \end{pmatrix}
\cdot
\begin{pmatrix}
\mfdmetric^{00}_{(x,\xi)} & \mfdmetric^{01}_{(x,\xi)} 
\\[2\jot]
\mfdmetric^{10}_{(x,\xi)} & \mfdmetric^{11}_{(x,\xi)}
\end{pmatrix}
\cdot
\begin{pmatrix} C(x) & C(x)\cdot H(\xi) \\ A(\xi) & B(\xi) \end{pmatrix}
\mapob.
\end{gathered}
\end{moneq}
As we have $\body B(\xi) = \oneasmatrix$, it follows easily that the metric $\body \mfdmetric^{00}$ on $\body G$ and the metric $-i\,\body \mfdmetric^{11}$ on the Batchelor bundle $B=\body G \times \body\Liealg g_1 \to \body G$ (see \recalf{changemetricbodycoordinates}) are as follows: $\body \mfdmetric^{00}$ is a left-invariant metric on $\body G$ that is the standard metric on $\Liealg g_0$ in terms of the basis $e_\bullet$ and $-i\,\body \mfdmetric^{11}$ is the standard metric on $\Liealg g_1$ in terms of the basis $f_\bullet$.

\begin{proclaim}[variousdensitiesonGorbodyGassotoinvariantmfdmetric]{Corollary}
Let $\nu_\mfdmetric$ and $\nu_{B,\mfdmetric}$ be the two densities on $G$ associated to the left-invariant metric $\mfdmetric$ as defined in \recalt{canonicaltrivdensonOOD}. 
Then in the local coordinate system $x,\xi$ we have
\begin{align*}
\nu_\mfdmetric\Bigl( \fracp{}{x}\bigrestricted_{(x,\xi)}, \fracp{}{\xi}\bigrestricted_{(x,\xi)} \Bigr)
&
=
\bigl\vert \Det\bigl(C(x)\bigr)\bigr\vert\mo  \cdot \frac{\Det\bigl( B(\xi) \bigr)}{\Det\bigl( \oneasmatrix - H(\xi)\cdot B(\xi)\mo\cdot A(\xi) \bigr)}
\\
\nu_{B,\mfdmetric}\Bigl( \fracp{}{x}\bigrestricted_{(x,\xi)}, \fracp{}{\xi}\bigrestricted_{(x,\xi)} \Bigr)
&
=
\bigl\vert \Det\bigl(C(x)\bigr)\bigr\vert\mo  
\mapob.
\end{align*}
Moreover, the density $\nu_o$ on $\body G$ defined in the local coordinates $x$ (now with $x\in \RR^d$) by
$$
\nu_{o}\Bigl( \fracp{}{x}\bigrestricted_{(x)} \Bigr)
=
\bigl\vert \Det\bigl(C(x)\bigr)\bigr\vert\mo  
$$
is the left-invariant density (volume form or Haar Measure) on $\body G$ associated to the left-invariant metric $\body \mfdmetric^{00}$ on $\body G$.

\end{proclaim}

\begin{preuve}
According to \recalt{canonicaltrivdensonOOD} these densities are determined by the Berezinian of the matrix associated to the super metric $\mfdmetric$. 
But according to \recalf{expressionleftinvariantmetriconG} we have the relation
$$
\Ber(\begin{pmatrix} \oneasmatrix & \mathbf0 \\ \mathbf0 & i\,\oneasmatrix \end{pmatrix})
=
\Ber(\begin{pmatrix}
\mfdmetric^{00}_{(x,\xi)} & \mfdmetric^{01}_{(x,\xi)}
\\[2\jot]
\mfdmetric^{10}_{(x,\xi)} & \mfdmetric^{11}_{(x,\xi)}
\end{pmatrix})
\cdot
\left(\Ber(\begin{pmatrix} C(x) & C(x)\cdot H(\xi) \\ A(\xi) & B(\xi) \end{pmatrix})\right)^2
\mapob.
$$
It follows immediately that we have
\begin{align*}
\Ber(\begin{pmatrix}
\mfdmetric^{00}_{(x,\xi)} & \mfdmetric^{01}_{(x,\xi)}
\\[2\jot]
\mfdmetric^{10}_{(x,\xi)} & \mfdmetric^{11}_{(x,\xi)}
\end{pmatrix})
&
=
\left(\frac{\det\bigl(B(\xi)\bigr)}{\det\bigl(C(x)-C(x)\cdot H(\xi)\cdot B(\xi)\mo\cdot A(\xi)\bigr)}\right)^2
\\
&
=
\left(\frac{\det\bigl(B(\xi)\bigr)}{\det\bigl(C(x)\bigr)\cdot \det\bigl(\oneasmatrix- H(\xi)\cdot B(\xi)\mo\cdot A(\xi)\bigr)}\right)^2
\mapob.
\end{align*}
Applying \recalt{canonicaltrivdensonOOD} immediately gives the desired results (using that $B(\xi)$ starts with the identity, and thus we don't need absolute values). 
\end{preuve}

\begin{proclaim}{Lemma}
Changing the initial \myquote{standard} super metric on $T_eG\cong \Liealg g$ in \recalf{expressionleftinvariantmetriconG} changes the densities $\nu_\mfdmetric$, $\nu_{B,\mfdmetric}$ and $\nu_o$ in \recalt{variousdensitiesonGorbodyGassotoinvariantmfdmetric} with a constant real multiple. 
These densities are thus unique up to multiplication by real numbers.

\end{proclaim}

\begin{definition}{Definitions}
For any $\CA$-Lie group $G$ equipped with a super metric $\mfdmetric$, we define $L^2(\body G, \KK, \nu_o)$ as the completion of $C^\infty_c(\body G; \KK)$ with respect to the metric $\inprodd\cdot\cdot$ induced by the invariant volume form $\nu_o \equiv \mathrm{Vol}_{\mfdmetric_{\body G}}$ determined by the metric $\mfdmetric_{\body G}$ \recalt{firstdefofordinarymetriconsuperfunctionsinE}, \recalt{variousdensitiesonGorbodyGassotoinvariantmfdmetric} (a volume form which is an invariant Haar measure), \ie, we have for $\chi,\psi\in L^2(\body G, \KK, \nu_o)$ (see also \recalf{usualmetriconL2G}):
$$
\inprodd\chi\psi = \int_{\body G} \overline{\chi(g)}\cdot \psi(g) \ \Vol_{\mfdmetric_{\body G}}
\mapob.
$$
As the metric $\inprodd\cdot\cdot$ is unique up to a real multiple, and as the field $\KK$ is most of the time understood implicitly, we will shorten the name of this Hilbert space $L^2(\body G, \KK, \nu_o)$ to $L^2(\body G)$. 
We finally denote by $\rho_{\body G}$ the left-regular representation of $\body G$ on $L^2(\body G)$ defined (as usual) by $\bigl(\rho_{\body G}(g)\psi\bigr)(h) = \psi(g\mo h)$. 

\end{definition}

\begin{proclaim}[metricandcompletionforleftregularrep]{Corollary}
The metric $\inprodsym_\mfdmetric$ and the super scalar product $\supinprsym_{B,\mfdmetric}$ on $C^\infty_{c}(G;\CA^\KK)$ induced by the left-invariant super metric $\mfdmetric$ on $G$ are given in the identification $C^\infty_{c}(G;\CA^\KK) \cong \bigl(\,C^\infty_{c}(\body G; \KK)\,\bigr){}^{2^n}$ 
\recalf{identificationCinftyGwithproductsofCinftybodyG} by
\begin{align*}
\inprod \chi\psi_\mfdmetric 
&
= \sum_{I\subset \{1, \dots, n\} } \int_{\body G} \overline{\chi_I(m)}\cdot \psi_I(m)\ \nu_o(m)
\equiv
\sum_{I\subset \{1, \dots, n\} } \inprodd{\chi_I}{\psi_I}
\end{align*}
and
\begin{align}
\homsuperinprod\chi\psi{_{B,\mfdmetric}}
&
\sum_{I\subset \{1, \dots, n\} } (-1)^{\varepsilon(I,I^c)}\,\inprodd{\chi_I}{\psi_{I^c}}
\mapob.
\label{explicitformofsuperinprodBmfdmetriconleftregrep}
\end{align}
The completion $\Hilbert$ of $C_c^\infty(G;\CA^\KK)$ thus is (isomorphic to)
$$
\Hilbert \cong \bigl(\,L^2(\body G, \KK, \nu_o)\,\bigr)^{2^n}
\equiv \bigl(\,L^2(\body G)\,\bigr)^{2^n}
\mapob.
$$

\end{proclaim}

\begin{preuve}
The expression for $\inprodsym_\mfdmetric$ is a direct consequence of \recalt{metriconsuperfunctionsinEinlocalcoordinates}. 
For $\supinprsym_{B,\mfdmetric}$ we invoke \recalt{descriptionsuperscalarproductBginadaptedBatlas}, which gives us:
\begin{align*}
\homsuperinprod\chi\psi{_{B,\mfdmetric}}
&
=
\int_M \nu_{B,\mfdmetric}(m) \cdot \overline{\chi(m)}\cdot \psi(m)
\notag
\\&
=
\sum_{a\in A} \int_{ U_a} \extder\Leb^{(d)}(x^a)\, \extder\xi^{(n)}\ \rho_a(x^a) \cdot 
\bigl\vert \Det\bigl(C(x)\bigr)\bigr\vert\mo 
\cdot \overline{\chi_a(x,\xi)}\cdot \psi_a(x,\xi) 
\notag
\\&
=
\sum_{I,J\subset \{1, \dots, n\} }\sum_{a\in A} 
\int_{\body U_a} \extder\Leb^{(d)}(x^a)\,\extder\xi^{(n)} \  \bigl\vert \Det\bigl(C(x)\bigr)\bigr\vert\mo  \cdot \rho_a(x^a)
\notag
\\&
\kern10em
\cdot \overline{\chi_{a,I}(x^a)}\cdot \psi_{a,J}(x^a)\cdot \xi^I\cdot \xi^J
\notag
\\&
=
\sum_{I\subset \{1, \dots, n\} } (-1)^{\varepsilon(I,I^c)}\cdot \int_{\body G} \overline{\chi_I(g)}\cdot \psi_{I^c}(g)\ \nu_o(g)
\notag
\\&
\equiv
\sum_{I\subset \{1, \dots, n\} } (-1)^{\varepsilon(I,I^c)}\,\inprodd{\chi_I}{\psi_{I^c}}
\end{align*}
as claimed. 
\end{preuve}

\begin{definition}{Definition}
We will denote by $L^2(G; \CA^\KK; \mfdmetric)$ the completion of the pre-Hilbert space $C^\infty_c(G, \CA^\KK)$ with respect to the metric $\inprodsym_\mfdmetric$. 
And as for the non-super case (and with the same arguments), we will shorten this to $L^2(G)$ if no confusion is possible. 

\end{definition}

The explicit formula for $\inprodsym_\mfdmetric$ in the identification $C^\infty_c(G; \CA^\KK) \cong \bigl(\,C^\infty_{c}(\body G; \KK)\,\bigr){}^{2^n}$ tells us immediately that we have a natural identification 
$$
L^2(G) \cong \bigl(\,L^2(\body G)\,\bigr)^{2^n}
\mapob.
$$
Looking at \recalf{explicitformofsuperinprodBmfdmetriconleftregrep} confirms that $\supinprsym_{B,\mfdmetric}$ is continuous with respect to the metric on $C^\infty_c(G; \CA^\KK) \subset L^2(G)$ \recalt{superscalarBgiscontinuous} and thus extends to $L^2(G)$. 
Moreover, it is also immediate from \recalf{explicitformofsuperinprodBmfdmetriconleftregrep} that the extension to $\Hilbert$ is (remains) non-degenerate, so it defines a super scalar product on $\Hilbert$.

In order to get a good understanding of the super scalar product $\supinprsym_{\mfdmetric}$ \recalt{defofsuperscalarprodforfunctionsonoodM} we define the function $\Delta:\oddp{\Liealg g}_0\to \CA_0$ by
$$
\Delta(\xi) = 
\frac{\Det\bigl( B(\xi) \bigr)}{\Det\bigl( \oneasmatrix - H(\xi)\cdot B(\xi)\mo\cdot A(\xi) \bigr)}
\mapob.
$$
As it depends only upon the odd coordinates $\xi$, there exist real constants $\Delta_I\in \RR$, $I\subset\{1, \dots, n\}$ such that
$$
\Delta(\xi) = \sum_{I\subset \{1, \dots, n\} } \xi^I\,\Delta_I
\mapob.
$$
Moreover, because $\Delta$ takes values in $\CA_0$, $\Delta_I=0$ whenever ${I}$ contains an odd number of elements. 
The definition of $\supinprsym_\mfdmetric$ then tells us that we have
\begin{align*}
\homsuperinprod\chi\psi{_\mfdmetric}
&
=
\int_M \nu_\mfdmetric(m) \cdot \overline{\chi(m)}\cdot \psi(m)
\\&
=
\sum_{a\in A} \int_{ U_a} \extder\Leb^{(d)}(x^a)\, \extder\xi^{(n)}\ \rho_a(x^a) \cdot 
\frac{\bigl\vert \Det\bigl(C(x)\bigr)\bigr\vert\mo  \cdot \Det\bigl( B(\xi) \bigr)}{\Det\bigl( \oneasmatrix - H(\xi)\cdot B(\xi)\mo\cdot A(\xi) \bigr)}
\\&
\kern7em
\cdot \overline{\chi_a(x,\xi)}\cdot \psi_a(x,\xi) 
\\&
=
\sum_{I,J,K\subset \{1, \dots, n\} }\ \sum_{a\in A} 
\int_{ U_a} \extder\Leb^{(d)}(x^a)\,\extder\xi^{(n)} \  \bigl\vert \Det\bigl(C(x)\bigr)\bigr\vert\mo  \cdot \rho_a(x^a)
\\&
\kern10em
\cdot \overline{\chi_{a,I}(x^a)}\cdot \psi_{a,J}(x^a)\cdot
\Delta_K\cdot \xi^I\cdot \xi^J \cdot \xi^K
\\&
=
\sum_{I,J,K\subset \{1, \dots, n\} }\ \biggl(\sum_{a\in A} 
\int_{\body U_a} \extder\Leb^{(d)}(x^a) \  \bigl\vert \Det\bigl(C(x)\bigr)\bigr\vert\mo  \cdot \rho_a(x^a)
\\&
\kern7em
\cdot \overline{\chi_{a,I}(x^a)}\cdot \psi_{a,J}(x^a)\biggr)\cdot
\int_{\CA_1^n}\extder\xi^{(n)} \ 
\Delta_K\cdot \xi^I\cdot \xi^J \cdot \xi^K 
\mapob.
\end{align*}
This suggests that we define the numbers $\superspmatrix_{IJ}\in \RR$ by
$$
\superspmatrix_{IJ} = \int \extder\xi^{(n)}\ \Delta_{(I\cup J)^c}\cdot \xi^I\cdot \xi^J \cdot \xi^{(I\cup J)^c} 
\mapob,
$$
which allows us to obtain finally
\begin{moneq}[explicitformulasuperinproductonGwithsuperspmatrix]
\homsuperinprod\chi\psi{_\mfdmetric}
=
\sum_{I,J\subset \{1, \dots, n\} } \superspmatrix_{IJ}\cdot \inprodd{\chi_I}{\psi_J}
\mapob.
\end{moneq}

\begin{definition}[linkingsupscprtometriconsuperLiegroupBUTnotunit]{Remark}
If we define the map $\superspmatrix:L^2(G) \to L^2(G)$ by
$$
(\superspmatrix \psi)_I = \sum_J \superspmatrix_{JI}\,\psi_J
\mapob,
$$
then we automatically have the equality
$$
\homsuperinprod\chi\psi{_\mfdmetric} = \inprod{\superspmatrix \chi}\psi
\mapob,
$$
linking the super scalar product $\supinprsym_\mfdmetric$ to the (ordinary) metric $\inprodsym_\mfdmetric$. 
However, the map $\superspmatrix$ will in general not preserve the metric $\inprodsym_\mfdmetric$ (see \recals{axiplusbetagroupsection} for an explicit example, see also \recalt{linkingsupscprtometricviaKrein}).

\end{definition}

\begin{proclaim}[superspmatrixisinvertibleandmore]{Corollary}
The matrix $(\superspmatrix_{IJ})_{I,J\subset \{1, \dots, n\}}$ (of size $2^n \times 2^n$) is invertible and we have the implication
$$
\bigl[\ I\cap J\neq \emptyset \text{ or }
\parity I + \parity J =  n + 1 \ (\mathrm{ mod }\,2)
\ \bigr]
\qquad\Longrightarrow\qquad
\superspmatrix_{IJ}=0
\mapob. 
$$
Moreover,the super scalar product $\supinprsym_\mfdmetric$ is continuous with respect to the topology induced by $\inprodsym_\mfdmetric$ on $C^\infty_c(G; \CA^\KK)\subset L^2(G)$ and extends to a continuous super scalar product of parity $ n$ on $L^2(G)$. 

\end{proclaim}

\begin{preuve}
Looking at the definition of $\superspmatrix_{IJ}$ shows that it is necessarily zero whenever $I\cap J\neq \emptyset$ or whenever $\Delta_{(I\cup J)^c}=0$. 
But $\Delta$ is an even function, so $\Delta_K=0$ whenever $\parity K =1$. 
And thus $\superspmatrix_{IJ}=0$ when $I\cap J\neq \emptyset$ and when $I\cap J=\emptyset$ it will also be zero when $\parity{(I\cup J)^c}=1$, \ie, whenever $\parity{I} + \parity{ J}=n+1$ as claimed. 

To prove invertibility, we suppose that we have a vector of numbers $c_I\in \RR$, not all zero, such that $\sum_{J\subset \{1, \dots, n\}} \superspmatrix_{IJ}\,c_J=0$ for all $I$. 
We then choose $I\subset \{1, \dots, n\}$ such that $c_I\neq0$ and $c_J=0$ for all $J\subset \{1, \dots, n\}$ with a smaller cardinal than $I$ (such $I$ need not be unique). 
By hypothesis we have in particular
$$
0 =
\sum_{J\subset \{1, \dots, n\}} \superspmatrix_{I^cJ}\,c_J
=
\sum_{J\subset I} \superspmatrix_{I^cJ}\,c_J
\mapob,
$$
where the second equality is a consequence of the fact that $\superspmatrix_{I^cJ}=0$ whenever $I^c\cap J\neq \emptyset$. 
But by assumption $c_J=0$ when the cardinal of $J$ is smaller than the cardinal of $I$. 
It thus follows that we must have $\superspmatrix_{I^cI}\,c_I=0$, \ie, $\superspmatrix_{I^cI}=0$. 
But $\Delta(0)=\Delta_\emptyset = 1$ and thus $\superspmatrix_{I^cI}=(-1)^{\varepsilon(I^c,I)}\neq0$. 
This contradiction shows that the kernel of the matrix $\superspmatrix$ is zero and thus this matrix is invertible. 

The only part concerning the super scalar product $\supinprsym_\mfdmetric$ that is not directly obvious is non-degeneracy of the extension of $\supinprsym$ to $L^2(G)$. 
To do so, we take $0 \neq \chi\in L^2(G)$ and we define $\psi$ by
$$
\psi_I = \sum_{J\subset \{1, \dots, n\}} (\superspmatrix\mo)_{IJ} \chi_J
\mapob.
$$
We then have 
\begin{align*}
\homsuperinprod{\chi}{\psi}{_\mfdmetric}
&
=
\sum_{I,K\subset \{1, \dots, n\}} \superspmatrix_{IK} \,\inprodd{\chi_I}{\psi_K}
=
\sum_{I,J,K\subset \{1, \dots, n\}} \superspmatrix_{IK}\,(\superspmatrix\mo)_{KJ}\, \inprodd{\chi_I}{\chi_J}
\\&
=
\sum_{I\subset \{1, \dots, n\}} \inprodd{\chi_I}{\chi_I} 
=
\inprod{\chi}{\chi}_\mfdmetric
\neq
0
\mapob,
\end{align*}
proving non-degeneracy of $\supinprsym_\mfdmetric$ on the whole of $L^2(G)$. 
\end{preuve}

\begin{proclaim}{Corollary}
The two triples $\bigl(L^2(G), \inprodsym_\mfdmetric, \supinprsym_{B,\mfdmetric}\bigr)$ and $\bigl(L^2(G), \inprodsym_\mfdmetric, \supinprsym_{\mfdmetric}\bigr)$ are super Hilbert spaces.

\end{proclaim}

We now turn our attention to the representation $\rho$ of $G$ on $C_{(c)}^\infty(G;\CA^\KK) \otimes \CA^\KK$ defined by (see \recalt{therepresentationofagrouponCinftyM})
$$
\bigl(\rho(g)\psi\bigr)(h) = \psi(g\mo h)
\mapob.
$$
In terms of the identification $G\cong \wod G \times \oddp{\Liealg g}_0$ the left-translation by an element $g\in \wod G$ is given by
$$
L_g(h,X) = (gh,X)
\mapob.
$$
It follows immediately that the maps $L_g$ with $g\in \body G\subset \wod G$ are diffeomorphisms that are linear in the odd coordinates (in the Batchelor atlas defined by the above identification). 
We thus can apply \recalt{linearinoddcoordinateshencepullbackunitary} (also because $\mfdmetric$ is left-invariant) and conclude that the maps $\rho(g)$ with $g\in \body G$ determine unitary maps on $C_{c}^\infty(G;\CA^\KK) \otimes \CA^\KK$. 
Extending them to the completion $L^2(G)$, we thus obtain a unitary representation $\rho_o$ of $\body G$ on $L^2(G)$ \recalt{DenseinsideCinftyrhoo}. 
In terms of the identification $\psi\in L^2(G) \cong (\psi_I)_{I\subset \{1, \dots, n\}} \in \bigl(L^2(\body G)\bigr){}^{2^n}$ the representation $\rho_o$ is given by
$$
\bigl(\rho_o(g)\psi\bigr)_I(h) = \psi_I(g\mo h)
\mapob,
$$
\ie, $\rho_o$ consists of $2^n$ copies of the left-regular representation $\rho_{\body G}$ of $\body G$.

Before turning our attention to the (full) representation $\rho$ on $G$, we first recall the definition of a fundamental vector field associated to a smooth left-action $\Phi:G\times M\to M$ of a Lie group $G$ on a manifold $M$ (a definition that is also valid in the context of $\CA$-manifolds). 
For any $X\in \Liealg g$, the associated \stresd{fundamental vector field $X^M$ on $M$} is defined by
$$
X^M\caprestricted_m = 
\frac{d}{dt}\bigrestricted_{t=0} \Phi(\eexp^{-tX},m)
\mapob.
$$
Because of the minus sign, this map is a Lie algebra morphism: for any $X,Y\in \Liealg g$ we have $[X,Y]^M = \bigl[ X^M, Y^M\bigr]$. 
When we apply this to the multiplication map on $G$: $m:G\times G\to G$, the associated fundamental vector field $X^G$ is the right-invariant vector field on $G$ whose value at $e$ is $-X$ (which confirms that we need the minus sign to get a Lie algebra morphism). 
When we apply it to the map\slash left-action $\Phi:G\times G\to G$ defined by $\Phi(g,h) = h\cdot g\mo$ (which is essentially the right-action of $G$ on itself), the associated fundamental vector field $X^G$ is exactly the left-invariant vector field $\vec X$. 

Looking at the (full) representation $\rho$ of $G$ on $C^\infty(G;\CA^\KK) \otimes \CA^\KK$, we know that its infinitesimal form is given by the fundamental vector fields of the left-action of $G$ on itself. 
In particular for $X\in \body \Liealg g$ homogeneous we have (with $\parity t=\parity X$)
$$
\frac{\extder}{\extder t}\bigrestricted_{t=0}\,\bigl(\rho\bigl(\exp(Xt)\bigr)\psi\bigr)(h)
=
\frac{\extder}{\extder t}\bigrestricted_{t=0}\,\psi\bigl(\exp(-tX)h\bigr)
=
-(X^R\psi)(h)
\mapob.
$$
We thus define the map $\tau:\body \Liealg g\to \End\bigl(C^\infty(G;\CA^\KK) \bigr)$ by
$$
\tau(X)\psi = -X^R\psi
\mapob,
$$
which is an even graded Lie algebra morphism preserving the subspace $C_c^\infty(G;\CA^\KK)$. 

We now recall a result by Poulsen \cite[p114]{Poulsen:1972} that says that $C^\infty(\rho_{\body G})$ is given by\footnote{There is a slight abuse of notation in identifying $\body \Liealg g_0$ with the Lie algebra of $\body G$ and by identifying the right-invariant vector fields $X^R$, $X\in \body \Liealg g_0$ with the right-invariant vector fields on $\body G$.}
$$
C^\infty(\rho_{\body G}) = 
\{\,\psi\in C^\infty(\body G;\KK) \mid
\forall k\in \NN\ \forall X_i\in \body \Liealg g_0: X_1^R\cdots X_k^R\psi \in L^2(\body G)\,\}
\mapob.
$$
Inspired by this result, we define the space $\Dense_\rho$ by
\begin{align*}
\Dense_{\rho}
&
=
\bigl\{\,\psi\in C^\infty(G; \CA^\KK) \mid
\forall k\in \NN\ \forall X_i\in \body \Liealg g: X_1^R\cdots X_k^R\psi \in L^2( G)\,\bigr\}
\mapob.
\end{align*}
Using the identifications
$$
C_{(c)}^\infty(G;\CA^\KK) \cong \bigl(C_{(c)}^\infty(\body G; \KK)\bigr)^{2^n}
\qquad\text{and}\qquad
L^2(G) \cong \bigl(L^2(\body G)\bigr)^{2^n}
\mapob,
$$
and knowing that $\rho_o$ is $2^n$ copies of the left-regular representation $\rho_{\body G}$, it follows easily that we have
$$
C_c^\infty(G; \CA^\KK)
\subset
C^\infty(\rho_o)
\cong
\bigl(C^\infty(\rho_{\body G})\bigr)^{2^n}
\subset
L^2(G)
\mapob.
$$
It then follows easily that we have the alternative description
\begin{moneq}[alternatesimpifieddescriptionDensesubrho]
\Dense_\rho = 
\{\,\psi\in C^\infty(\rho_o) \mid \forall k\in \NN\ \forall X_i\in \body \Liealg g_1: X_1^R\cdots X_k^R\psi \in L^2( G)\,\}\mapob,
\end{moneq}
and that we have the inclusions
$$
C^\infty_c(G; \CA^\KK) \subset \Dense_\rho \subset C^\infty(\rho_o)
\mapob,
$$
proving that $\Dense_\rho$ is dense in $\Hilbert$. 
Before stating our main result, we note that obviously $\Dense_\rho$ is the maximal domain that is invariant under all maps $\tau(X)$, $X\in \body \Liealg g$.

\begin{proclaim}[maintheoremLeftRegularRep]{Theorem}
The couple $\bigl(C^\infty_c(G; \CA^\KK), \rho\bigr)$ is a \psur{} of $G$ on the super Hilbert space $\bigl(L^2(G), \inprodsym_\mfdmetric, \supinprsym_\mfdmetric\bigr)$ whose infinitesimal form is the triple $\bigl(\rho_o, C^\infty_c(G; \CA^\KK), \tau\bigr)$. 
Moreover, the couple $(\Dense_\rho, \rho)$ is the unique super unitary representation of $G$ on $L^2(G)$ that extends $\bigl(C^\infty_c(G; \CA^\KK), \rho\bigr)$; its infinitesimal form is the triple $(\rho_o, \Dense_\rho, \tau)$. 

\end{proclaim}

The proof of \recalt{maintheoremLeftRegularRep} will be broken down into a series of smaller results in such a way that at the end it suffices to put these smaller results together to obtain a full proof.

\begin{proclaim}[CinftysubcisPSUR]{Lemma}
The triple $\bigl(\rho_o, C^\infty_c(G; \CA^\KK), \tau\bigr)$ is a \psur{} in its infinitesimal form of $G$ on the super Hilbert space $\bigl(L^2(G), \inprodsym_\mfdmetric, \supinprsym_\mfdmetric\bigr)$. 

\end{proclaim}

\begin{preuve}
Looking at the \myquote{definition} of a \psur{} in infinitesimal form \recalt{equivalentDefSuperUnitaryRepNEW}, it is immediate that $\rho_o$ and the maps $\tau(X)$ preserve the space $C^\infty_c(G; \CA^\KK)$. 
We thus have to prove the properties (\ref{alternateSUR1}--\ref{alternateSUR3}). 

For the first property we recall that $\rho_o$ consists of $2^n$ copies of $\rho_{\body G}$ and that the right-invariant vector fields $X^R$ (for $X\in \body\Liealg g_0$) \myquote{are} $2^n$ copies of the corresponding right-invariant vector fields on $\body G$ \recaltt{leftandrightinvvfintrivializationofGintoGwodalgopddp}{leftinvariantvfonGwodtimesgsub1}. 
As these right-invariant vector fields on $\body G$ are the generators of the left-regular representation $\rho_{\body G}$, it follows immediately that the maps $\tau(X)$, $X\in \body\Liealg g_0$ are the infinitesimal generators of $\rho_o$, proving \recaltt{alternateSUR1}{equivalentDefSuperUnitaryRepNEW}. 

By the same argument it follows that the maps $\tau(X)$ with $X\in \body \Liealg g_0$ are skew-symmetric with respect to $\inprodsym_\mfdmetric$.
It then suffices to look at the form \recalf{explicitformulasuperinproductonGwithsuperspmatrix} of $\supinprsym_\mfdmetric$ (and to use the fact that $X^R$ consists of $2^n$ copies of the corresponding right-invariant vector field on $\body G$) to see that $\tau(X)$ is (graded) skew-symmetric with respect to $\supinprsym_\mfdmetric$. 
On the other hand, for $X\in \body \Liealg g_1$ we can take $\xi\in \CA_1$ and consider the map $F:\CA_1\times G \to \CA^\KK$ defined by 
$$
F(\xi,g)=\psi\bigl(\exp(-\xi X)\,g\bigr)
$$
for any $\psi\in C^\infty_c(G; \CA^\KK)$. As this is a smooth map and as $\xi^2=0$, we have the equality
$$
F(\xi,g) = F(0,g) + \xi\cdot \fracp F\xi(0,g)
\mapob.
$$
But by definition $(\partial_\xi F)(0,g)$ is given by $(X^R\psi)(g)$, which also belongs to $C^\infty_c(G; \CA^\KK)$. 
We thus have, for any $\psi\in C^\infty_c(G; \CA^\KK)$ the equality
$$
\rho\bigl(\exp(\xi X)\bigr)\psi = \psi+\xi\,\tau(X)\psi
\mapob.
$$
Using that $\rho$ preserves $\supinprsym_\mfdmetric$ on $C^\infty_c(G; \CA^\KK)$ \recalt{superscalarproductmfdmetricinvariantunderGaction}, we deduce from this, just as in the proof of (the direct part of) \recalt{equivalentDefSuperUnitaryRepNEW} (where we used \refmetnaam{\labelMainSUR}{MainSURlabelrhopsiissmooth} to obtain the analogous equality \recalf{deftauforoddbyexpansion}) that $\tau(X)$ is graded skew-symmetric with respect to $\supinprsym_\mfdmetric$. 
And thus all maps $\tau(X)$, $X\in \body\Liealg g$ are graded skew-symmetric with respect to $\supinprsym_\mfdmetric$, proving \recaltt{alternateSUR2}{equivalentDefSuperUnitaryRepNEW}. 

As an immediate consequence of the properties of a fundamental vector field associated to a smooth action (here the left-action of $G$ on itself) we have
$$
\forall g\in \body g
\quad:\quad
\tau\bigl(\Ad(g)X\bigr) = \rho(g)\scirc \tau(X) \scirc \rho(g\mo)
\mapob,
$$
which is condition \recaltt{alternateSUR3}{equivalentDefSuperUnitaryRepNEW}. 
\end{preuve}

\begin{proclaim}[tauonDenseintegratestorho]{Lemma}
If a \psur{} $(\rho_o, \Dense, \tau)$ is an extension of the (infinitesimal) \psur{} $\bigl(\rho_o, C^\infty_c(G; \CA^\KK), \tau\bigr)$ then its integrated form is $(\Dense, \rho)$. 

\end{proclaim}

\begin{preuve}
If $(\rho_o, \Dense, \tau)$ is an infinitesimal \psur, then according to \recalt{equivalentDefSuperUnitaryRepNEW} there exists a unique \psur{} $(\Dense, \rhoh)$ whose infinitesimal form is $(\rho_o, \Dense, \tau)$. 
We thus have to show that we necessarily have $\rhoh(g)=\rho(g)\restricted_\Dense$. 
In order to do so, we follow more or less the proof of (the converse part of) \recalt{equivalentDefSuperUnitaryRepNEW} in which the representation $\rhoh$ is constructed out of $\tau$. 

By definition of $\rhoh$ we have, for $g\in \body G$, $\rhoh(g) = \rho_o(g) = \rho(g)$ and thus $\rhoh=\rho$ when restricted to $\body G$. 
The next step in the construction of $\rhoh$ is to extend $\rho_o$ by  \Gextension{} to $\wod G$. 
Here we have to be careful, because for fixed $\psi\in \Dense$, the \Gextension{} is done with respect to the smooth map $\FgroupHilbert_\psi:\body G\to L^2(G)$, $g\mapsto \rho_o(g)\psi$, not with respect to the smooth map $(g,h)\mapsto \bigl(\rho_o(g)\psi\bigr)(h) = \psi(g\mo h)$, which takes values in $\CA^\KK$. 
We thus first recall the above cited result of Poulsen that says that $C^\infty(\rho_{\body G})$ is given by
$$
C^\infty(\rho_{\body G}) = 
\{\,\psi\in C^\infty(\body G;\KK) \mid
\forall k\in \NN\ \forall X_i\in \body \Liealg g_0: X_1^R\cdots X_k^R\psi \in L^2(\body G)\,\}
\mapob,
$$
and in particular that the elements of $C^\infty(\rho_{\body G})$ are smooth functions. 
It then is easy to deduce that (still in the context of ordinary non-super Lie groups and their left-regular representation) we have, for $X\in \body\Liealg g_0$ and $\psi\in C^\infty(\rho_{\body G})$:
\begin{moneq}[derivativeofFgroupHilbertforbodyLiealgg0]
\vec X \FgroupHilbert_\psi = \FgroupHilbert_{\tau(X)\psi}
\mapob,
\end{moneq}
simply because $\tau(X) = -X^R$ is the fundamental vector field associated to the left-action of $G$ on itself. 
We next recall that $\rho_o$, the representation of $\body G$ on $L^2(G) \cong \bigl(L^2(\body G)\bigr){}^{2^n}$ is just $2^n$ copies of $\rho_{\body G}$ and that the right-invariant vector field $X^R$ boils down to $2^n$ copies of the restriction of $X^R$ to $L^2(\body G)$ (see \recaltt{leftandrightinvvfintrivializationofGintoGwodalgopddp}{leftinvariantvfonGwodtimesgsub1}). 
It follows that \recalf{derivativeofFgroupHilbertforbodyLiealgg0} remains valid for $\psi\in C^\infty(\rho_o)$, \ie, we have for all $\psi\in C^\infty(\rho_o)$ and all $g,h\in G$ the equality
$$
\Bigl(\frac{\extder}{\extder t} \bigrestricted_{t=0} \rhoh\bigl(g\,\exp(tX)\bigr)\psi\Bigr)(h)
=
\Bigl(\rhoh(g)\bigl(\tau(X)\psi\bigr)\Bigr)(h)
\mapob.
$$

Now there is a general statement concerning the flow of a smooth vector field $X$ on a manifold $M$ that says that if $\phi_t$ denotes the flow of $X$ and if, for a smooth function $f:M\to V$ (with $V$ some normed vector space) we define the function $F:\RR\to V$ by $F(t) = f\bigl(\phi_t(m)\bigr)$ for some $m\in M$, then we have for the $k$-th order derivative:
$$
F^{(k)}(t) = (X^kf)\bigl(\phi_t(m)\bigr)
\mapob.
$$
Moreover, this statement remains valid for $\CA$-manifolds and smooth even homogeneous vector fields. 

In order to apply this in our situation, we consider a left-invariant vector field $\vec X$ on $\wod G$ with $X\in \body\Liealg g_0$, whose flow is given by $g\mapsto g\,\exp(tX)$. 
We also consider the smooth function $\FgroupHilbert_\psi:\wod G\to \Dense \otimes \CA^\KK$ obtained by \Gextension{} from the smooth function with the same name $\FgroupHilbert_\psi:\body G\to \Dense$. 
We thus define, for any $g\in \wod G$, the function $F:\RR\to \Dense \otimes \CA^\KK$ by $F(t) = (\FgroupHilbert_\psi)\bigl(g\,\exp(tX)\bigr)$, and we find
$$
F^{(k)}(t) = \bigl((\vec X)^k\FgroupHilbert_\psi\bigr)\bigl(g\,\exp(tX)\bigr)
=
\FgroupHilbert_{\tau(X)^k\psi}\bigl(g\,\exp(tX)\bigr)
\mapob.
$$
By Taylor expansion and the definition\slash construction of $\rhoh$ we thus can make the computation for even nilpotent $n\in \CA_0$ (which yields a finite Taylor series):
\begin{align*}
\rhoh\bigl(g\,\exp(nX)\bigr)\psi
&
=
\FgroupHilbert_\psi\bigl(g\,\exp(nX)\bigr)
=
F(n) 
=
\sum_{k=0}^\infty \frac{n^k}{k!}\,F^{(k)}(0)
\\&
=
\sum_{k=0}^\infty \frac{n^k}{k!}\,\FgroupHilbert_{\tau(X)^k\psi}(g)
=
\sum_{k=0}^\infty \frac{n^k}{k!}\,\rhoh(g)\bigl(\tau(X)^k\psi\bigr)
\mapob.
\end{align*}
We now suppose that we know that, for this $g\in \wod G$, we have the equality $\rhoh(g) = \rho(g)$. 
Then we may conclude that we have, for all $h\in G$:
\begin{align*}
\Bigl(\rhoh\bigl(g\,\exp(nX)\bigr)\psi\Bigr)(h)
&
=
\sum_{k=0}^\infty \frac{n^k}{k!}\,\Bigl(\rhoh(g)\bigl(\tau(X)^k\psi\bigr)\Bigr)(h)
=
\sum_{k=0}^\infty \frac{n^k}{k!}\,\bigl(\tau(X)^k\psi\bigr)(g\mo h)
\\&
=
\psi\bigl(\exp(-nX)\,g\mo h\bigr)
=
\Bigl(\rho\bigl(g\,\exp(nX)\bigr)\psi\Bigr)(h)
\mapob,
\end{align*}
simply because the flow of the vector field $\tau(X) = -X^R$ on $G$ is given by $h\mapsto \exp(-tX)h$. 
The conclusion is that, if we know $\rhoh(g)=\rho(g)$, then we also know $\rhoh\bigl(g\,\exp(nX)\bigr) = \rho\bigl(g\,\exp(nX)\bigr)$. 
Starting with $g\in \body G \subset \wod G$ for which we know it is true, we conclude (by taking successive basis elements $X=e_i$ for $\body \Liealg g_0$) that it is true for all $g\in \wod G$. 

To finish, we take $X\in \body\Liealg g_1$ and $\xi\in \CA_1$. 
Then by definition of $\rhoh$ we have 
$$
\rhoh\bigl(\exp(\xi X)\bigr)\psi = \psi+\xi\tau(X)\psi
$$
and thus for any $g\in G$:
\begin{align*}
\Bigl(\rhoh\bigl(\exp(\xi X)\bigr)\psi\Bigr)(g) 
&
= \bigl(\psi+\xi\tau(X)\psi\bigr)(g)
=
\psi(g) -\xi (X^R\psi)(g)
\\&
=
\psi\bigl(\exp(-\xi X)\,g\bigr)
=
\Bigl(\rho\bigl(\exp(\xi X)\bigr)\psi\Bigr)(g)
\mapob,
\end{align*}
where the third equality is a consequence of the property of the flow of the even vector field $\xi X^R$ whose flow is given by $g\mapsto \exp(t\xi X)\,g$. 
We thus may conclude that for any $X\in \body\Liealg g_1$ and any $\xi\in \CA_1$ we have $\rhoh\bigl(\exp(\xi X)\bigr) = \rho\bigl(\exp(\xi X)\bigr)$. 
But $\wod G$ and elements of the form $\exp(\xi X)$ generate the whole group $G$ and thus $\rhoh=\rho$ as claimed. 
\end{preuve}

\begin{proclaim}[propertiesofsupinprsymmfdmetricwithsupportrestrictions]{Lemma}
Let $\psi\in L^2(G) \cap C^\infty(G; \CA^\KK)$.
\begin{enumerate}
\item
If $\chi\in L^2(G)$ is such that $\supp(\chi) \cap \supp(\psi) = 0$, then $\homsuperinprod\chi\psi{_\mfdmetric} = \inprod\chi\psi_\mfdmetric=0$.

\item
If for some open subset $U\subset G$ we have the property
$$
\forall \chi\in C^\infty_c(G;\CA^\KK)
\quad:\quad
\supp(\chi)\subset U
\quad\Rightarrow\quad
\homsuperinprod\chi\psi{_\mfdmetric} = 0
\mapob,
$$
then $\psi\restricted_U$, the restriction of $\psi$ to $U$, is zero.

\end{enumerate}

\end{proclaim}

\begin{preuve}
If we have $\supp(\chi) \cap \supp(\psi) = \emptyset$, then in particular for any $I,J\subset \{1, \dots, n\}$ we also have $\supp(\chi_I) \cap \supp(\psi_J) = \emptyset$. 
It follows directly that we have
$$
\inprodd{\chi_I}{\psi_J}
=
\int_{\body G} \overline{\chi_I(g)}\cdot \psi_J(g) \ \Vol_{\mfdmetric_{\body G}}
=
0
\mapob.
$$
Applying the explicit expressions for $\inprodsym_\mfdmetric$ \recalt{metricandcompletionforleftregularrep} and $\supinprsym_\mfdmetric$ \recalf{explicitformulasuperinproductonGwithsuperspmatrix} then shows the first point. 

\medskip

For the second point we fix $I_o\subset \{1, \dots, n\}$ and we take an arbitrary $\chi_o\in C^\infty_c(\body G; \KK)$ with $\supp(\chi_o)\subset U$. 
We then define $\chi\in C^\infty_c(G; \CA^\KK) \cong \bigl( C^\infty_c(\body G; \KK)\bigr){}^{2^n}$ by $\chi_{I_o} = \chi_o$ and $\chi_I=0$ for $I\neq I_o$. 
We then have $\supp(\chi)\subset U$ and we have
$$
0 = \homsuperinprod{\chi}{\psi}{_\mfdmetric} = \sum_{I,J\subset \{1, \dots, n\} } \superspmatrix_{IJ}\cdot \inprodd{\chi_I}{\psi_J}
=
\sum_{J\subset \{1, \dots, n\}} \superspmatrix_{I_oJ}\cdot \inprodd{\chi_o}{\psi_J}
\mapob.
$$
As this is true for an arbitrary $I_o$ and as the matrix $\superspmatrix_{IJ}$ is invertible \recalt{superspmatrixisinvertibleandmore}, we may conclude that we have
$$
\forall I_o,J\subset \{1, \dots, n\}
\quad:\quad
0 = \inprodd{\chi_o}{\psi_J} = \int_{\body G} \overline{\chi_o(g)}\cdot \psi_J(g) \ \Vol_{\mfdmetric_{\body G}}
\mapob.
$$
The results then follows from standard measure theory and the fact that the functions are smooth. 
\end{preuve}

\begin{proclaim}[tauisuniqueforleftregularrep]{Lemma}
Let $(\Dense, \tauh)$ be such that $\Dense\subset L^2(G)$ is a graded subspace satisfying $C_c^\infty(G; \CA^\KK) \subset \Dense \subset C^\infty(G; \CA^\KK)$ and $\tauh:\body\Liealg g\to \End(\Dense)$ an even graded Lie algebra morphism satisfying the conditions
\begin{enumerate}
\item
for all $X\in \body\Liealg g$ the restriction of $\tauh(X)$ to $C_c^\infty(G; \CA^\KK)$ equals $\tau(X)$ and

\item
for all $X\in \body \Liealg g$ the map $\tauh(X)$ is graded skew-symmetric with respect to $\supinprsym_\mfdmetric$.

\end{enumerate}
Then $\tauh=\tau$ and in particular $\Dense$ is invariant under all maps $\tau(X)$, $X\in \body \Liealg g$.

\end{proclaim}

\begin{preuve}
Let us start with the remark that by \recalt{CinftysubcisPSUR} the maps $\tau(X)$ are graded skew-symmetric with respect to $\supinprsym_\mfdmetric$ on $C^\infty_c(G; \CA^\KK)$. 
It thus makes sense to require that the restriction of $\tauh$ to $C_c^\infty(G; \CA^\KK)$ equals $\tau$. 

We now fix $X\in \body \Liealg g$ homogeneous, $\psi\in \Dense$ and $g_o\in G$ and we choose $\phi\in C_c^\infty(G; \CA)$ and an open neighborhood $U\subset G$ of $g_o$ such that $\phi(g)=1$ for all $g\in U$. Such $\phi$ (with compact support!) and $U$ exist because $G$ is locally compact. 
As $\phi$ has compact support and $\psi$ is smooth, the product $\phi\cdot \psi$ belongs to $C_c^\infty(G; \CA^\KK)$. 
Because $\tauh$ and $\tau$ are graded skew-symmetric with respect to $\supinprsym_\mfdmetric$ on $C^\infty_c(G; \CA^\KK)$, we have for any $\chi\in C^\infty_c(G; \CA^\KK)$ with $\supp(\chi)\subset U$ the equality 
$$
\homsuperinprod[2]{\chi}{\tau(X)(\phi\cdot\psi) - \tauh(X)\psi}{_\mfdmetric}
=
-\homsuperinprod[2]{\tau(X)\conjugate^{\parity X}(\chi)}{\phi\cdot\psi - \psi}{_\mfdmetric}
\mapob.
$$
But for all $h\in \supp(\chi)$ we have $\phi(h)=1$ and thus, because $\tau(X)=-X^R$ is a differential operator, we have 
$$
\supp\bigl( \tau(X)\conjugate^{\parity X}(\chi) \bigr) \cap \supp(\phi\cdot\psi - \psi)
\subset 
\supp(\chi) \cap \supp(\phi\cdot\psi - \psi)=\emptyset
\mapob.
$$
It then follows from \recalt{propertiesofsupinprsymmfdmetricwithsupportrestrictions} that we have
$$
\homsuperinprod[2]{\tau(X)\conjugate^{\parity X}(\chi)}{\phi\cdot\psi - \psi}{_\mfdmetric}
=
0
\mapob.
$$
We thus have shown that 
$$
\forall \chi\in C_c^\infty(G; \CA^\KK)
\quad:\quad
\supp(\chi)\subset U
\ \Rightarrow\ 
\homsuperinprod[2]{\chi}{\tau(X)(\phi\cdot\psi) - \tauh(X)\psi}{_\mfdmetric}
= 0
\mapob.
$$
Invoking \recalt{propertiesofsupinprsymmfdmetricwithsupportrestrictions} again we obtain that $\tau(X)(\phi\cdot\psi) - \tauh(X)\psi$ is identically zero on $U$. 
For all $g\in U$ we thus have
\begin{align*}
\bigl(\tauh(X)\psi\bigr)(g)
&
=
\bigl(\tau(X)(\phi\cdot\psi)\bigr)(g)
=
-\bigl(X^R(\phi\cdot\psi)\bigr)(g)
\\&
=
-(X^R\phi)(g)\cdot \psi(g) + \phi(g)\cdot \bigl(\tau(X)\psi\bigr)(g)
\mapob.
\end{align*}
But $\phi$ is identically $1$ on $U$ so $X^R\phi$ is identically $0$ on $U$. 
And thus in particular $\bigl(\tauh(X)\psi\bigr)(g_o) = \bigl(\tau(X)\psi\bigr)(g_o)$. 
As $g_o$ was arbitrary, we have shown $\tauh(X)\psi = \tau(X)\psi$. 
And as $\psi\in \Dense$ was arbitrary, we have shown $\tauh(X)=\tau(X)$ for all homogeneous $X\in \body\Liealg g$. 
\end{preuve}

\begin{proclaim}{Lemma}
The maps $\tau(X)$, $X\in \body\Liealg g$ are graded skew-symmetric on $\Dense_\rho$. 

\end{proclaim}

\begin{preuve}
It is obvious from the definition of $\Dense_\rho$ that this graded subspace is invariant under all maps $\tau(X)$, so it makes sense to investigate whether they are graded skew-symmetric. 
By linearity it suffices to prove this only for our fixed basis $e_1, \dots, e_d$, $f_1, \dots, f_n$ of $\body \Liealg g$. 
As we already have seen, the maps $\tau(e_i)$ are the generators of $\rho_o$ and are thus (graded) skew symmetric on $C^\infty(\rho_o) \supset \Dense_\rho$ (see the proof of \recalt{CinftysubcisPSUR}). 

In order to prove that a $\tau(f_j)$ is graded skew-symmetric, we use a change of coordinates that will allow us to \myquote{rectify} the right-invariant vector field $f_j^R$. 
Besides the identification $\Phi : \wod G \times \oddp{\Liealg g}_0 \to G$ introduced in \recals{superunitarydefandequivalencesection}, the map $\Phih_j : \CA_1^n \times \wod G \to G$ defined by
\begin{align*}
\shifttag{4em}
\Phih\bigl((\xi_1, \dots, \xi_n), g\bigr) =
\\&
\exp(\xi_j \, f_j)\cdot \exp(\xi_1\,f_1) \cdots \exp(\xi_{j-1}\, f_{j-1}) \cdot \exp(\xi_{j+1}\, f_{j+1}) \cdot \exp(\xi_{n}\, f_{n}) \cdot g
\end{align*}
also is a diffeomorphism. 
And in these coordinates (the odd $\xi_i$ and even coordinates on $\wod G$) the right-invariant vector field $f_j^R$ is given by $f_j^R = \partial_{\xi_j}$. 
We then choose a sequence of smooth maps with compact support $\zeta_n:\body G \to [0,1]$ such that for all $g\in \body G$ we have $\limn \zeta_n(g) = 1$. 
Such a sequence exists because any ordinary Lie group is $\sigma$-compact. 
Denoting by the same symbol $\zeta_n$ the extension of these maps first to $\wod G$, then to $\CA_1^n\times \wod G$ and finally to $G$ via the map $\Phih_j$, we obtain the property, 
$$
\limn \zeta_n(g)\cdot \psi(g) = \psi(g)
\mapob,
$$
where $\psi:G\to \CA^\KK$ is any function. 
The important points for us are the following: 
\begin{enumerate}
\item
if $\psi$ is smooth, then $\zeta_n\cdot \psi$ belongs to $C^\infty_c(G; \CA^\KK)$, 

\item
if $\psi$ belongs to $L^2(G)$, then $\zeta_n\cdot \psi$ converges to $\psi$ in $L^2(G)$,

\item
if $\psi$ is smooth, then $f_j^R(\zeta_n\cdot \psi) = \zeta_n\cdot (f_j^R\psi)$. 

\end{enumerate}
We then choose $\chi,\psi\in \Dense_\rho \subset L^2(G)\cap C^\infty(G; \CA^\KK)$ and we compute:
\begin{align*}
\homsuperinprod{f_j^R\chi}{\psi}{_\mfdmetric} 
&
=
\homsuperinprod{\lim_{m\to\infty} \zeta_m\,f_j^R\chi}{\limn \zeta_n \,\psi}{_\mfdmetric}
=
\lim_{m\to\infty}\limn 
\homsuperinprod{\zeta_m\,f_j^R\chi}{ \zeta_n \,\psi}{_\mfdmetric}
\\&
=
\lim_{m\to\infty}\limn 
\homsuperinprod{f_j^R(\zeta_m\,\chi)}{ \zeta_n \,\psi}{_\mfdmetric}
=
\lim_{m\to\infty}\limn 
-\homsuperinprod{\zeta_m\,\conjugate\chi}{f_j^R( \zeta_n \,\psi)}{_\mfdmetric}
\\&
=
\lim_{m\to\infty}\limn 
-\homsuperinprod{\zeta_m\,\conjugate\chi}{\zeta_n \,f_j^R\psi)}{_\mfdmetric}
= 
-\homsuperinprod{\lim_{m\to\infty}\zeta_m\,\conjugate\chi}{\limn\zeta_n \,f_j^R\psi)}{_\mfdmetric}
\\&
=
-\homsuperinprod{\conjugate\chi}{f_j^R\psi)}{_\mfdmetric}
\mapob,
\end{align*}
where for the first and last equality we used (ii) above, for the second and sixth we used the continuity of $\supinprsym_\mfdmetric$, for the third and fifth we used (iii) above, and for the fourth equality we used (i) above and the fact that $f_j^R$ is graded skew-symmetric with respect to $\supinprsym_\mfdmetric$ on $C^\infty_c(G; \CA^\KK)$ \recalt{CinftysubcisPSUR} (and that $f_j^R$ is odd). 
The final result then says that $f_j^R\equiv -\tau(f_j)$ is graded skew-symmetric with respect to $\supinprsym_\mfdmetric$ on $\Dense_\rho$ as claimed. 
\end{preuve}

\begin{preuve}[Proof of \recalt{maintheoremLeftRegularRep}]
By \recalt{CinftysubcisPSUR} we know that $\bigl(\rho_o, C^\infty_c(G; \CA^\KK), \tau\bigr)$ is a \psur{} in infinitesimal form and then by \recalt{tauonDenseintegratestorho} we know that its integrated form is $\bigl( C^\infty_c(G; \CA^\KK), \rho\bigr)$, proving the first part. 

For the second part, we start by showing that $\bigl(\rho_o, \Dense_\rho, \tau\bigr)$ is a \psur{} in infinitesimal form \recalt{equivalentDefSuperUnitaryRepNEW}. 
For that we note that we already know that $\tau:\body \Liealg g \to \End(\Dense_\rho)$ is an even graded Lie algebra morphism that preserves $\supinprsym_\mfdmetric$ and thus in particular condition \recaltt{alternateSUR2}{equivalentDefSuperUnitaryRepNEW} is satisfied. 
And as the $\tau(X)$ are fundamental vector fields, we also have condition \recaltt{alternateSUR3}{equivalentDefSuperUnitaryRepNEW}. 
As $\Dense_\rho\subset C^\infty(\rho_o)$, it also follows immediately that  \recaltt{alternateSUR1}{equivalentDefSuperUnitaryRepNEW} is satisfied. 
It thus remains to prove that $\Dense_\rho$ is invariant under $\rho_o$. 
For that we recall the description \recalf{alternatesimpifieddescriptionDensesubrho} of $\Dense_\rho$, the fact that $C^\infty(\rho_o)$ is invariant under $\rho_o$ and the property \recaltt{alternateSUR3}{equivalentDefSuperUnitaryRepNEW} just proven to conclude that indeed $\Dense_\rho$ is invariant under $\rho_o$. 

Once we know that $\bigl(\rho_o, \Dense_\rho, \tau\bigr)$ is a \psur{} in infinitesimal form, we apply \recalt{tauonDenseintegratestorho} again to conclude that $(\Dense_\rho, \rho)$ is a \psur. 
Now suppose $(\rho_o,\Dense, \tauh)$ is any maximal extension of $\bigl(\rho_o, C^\infty_c(G; \CA^\KK), \tau\bigr)$. 
Then by \recalt{tauisuniqueforleftregularrep} we must have $\tauh=\tau$ and $\Dense$ is invariant under all maps $\tau(X)$, $X\in \body\Liealg g$. 
But $\Dense_\rho$ is the maximal graded subspace invariant under all these maps, so we must have $\Dense\subset \Dense_\rho$. 
By maximality of $(\Dense, \tauh)$ we conclude that we have $\Dense=\Dense_\rho$, proving that $\bigl(\rho_o, \Dense_\rho, \tau\bigr)$ is the unique maximal extension of $\bigl(\rho_o, C^\infty_c(G; \CA^\KK), \tau\bigr)$. 
And thus $(\Dense_\rho, \rho)$ is the unique super unitary representation extending $\bigl( C^\infty_c(G; \CA^\KK), \rho\bigr)$. 
\end{preuve}

\masection{Super direct integrals}
\label{directintegralsofrepssection}

\begin{definition}[defofdirectintegraloverRkofsuperunireps]{Definition}
Let $(\Hilbert, \inprodsym, \supinprsym)$ be a super Hilbert space. 
We then define the space $L^2(\RR^d; \Hilbert, \Leb^{(d)})$ by
\begin{align*}
\shifttag{17em}
L^2(\RR^d; \Hilbert, \Leb^{(d)}) = \Bigl\{\,\psi:\RR^d\to \Hilbert \ \Bigm\vert\  \psi \text{ is measurable and } 
\\&
\int_{\RR^d} \inprod{\psi(k)}{\psi(k)}\ \extder\Leb^{(d)}(k)<\infty\,\Bigr\}
\mapob,
\end{align*}
which we equip with the metric $\inprodsym_{\mathrm{DI}}$ defined as
\begin{moneq}[metricondirectintegralofsuperHilbert]
\inprod{\chi}{\psi}_{\mathrm{DI}}
=
\int_{\RR^d} \inprod{\chi(k)}{\psi(k)}\ \extder\Leb^{(d)}(k)
\mapob.
\end{moneq}
So far this is a particular case of the direct integral of a family of Hilbert spaces (in this case all the same $\Hilbert$). 
But on $L^2(\RR^d; \Hilbert, \Leb^{(d)})$ we also define the tentative super scalar product $\supinprsym_{\mathrm{DI}}$ by
\begin{moneq}[superscalarprodondirectintegralofsuperHilbert]
\homsuperinprod\chi\psi{_{\mathrm{DI}}}
=
\int_{\RR^d} \superinprod[2]{\chi(k)}{\psi(k)} \ \extder\Leb^{(d)}(k)
\mapob.
\end{moneq}
As $\supinprsym$ is continuous with respect to $\inprodsym$, it is easy to establish first that $\supinprsym_{\mathrm{DI}}$ is well defined on $L^2(\RR^d; \Hilbert, \Leb^{(d)})$ (the defining integral is convergent) and then that it is continuous with respect to $\inprodsym_{\mathrm{DI}}$. 
We then make the assumption that it also is non-degenerate, which will turn it into a bona-fide super scalar product on $L^2(\RR^d; \Hilbert, \Leb^{(d)})$. 

An element $\psi\in L^2(\RR^d; \Hilbert, \Leb^{(d)})$ will be homogeneous of degree $\alpha$ if (and only if) for all $k\in \RR^d$ we have $\psi(k)\in \Hilbert_\alpha$. 
It then follows easily that with this definition the triple $\bigl(L^2(\RR^d; \Hilbert, \Leb^{(d)}), \inprodsym_{\mathrm{DI}}, \supinprsym_{\mathrm{DI}}\bigr)$ becomes a super Hilbert space. 

\medskip

Following \cite[\S8.4]{Kirillov:1976} we will say that a super unitary representation (in its infinitesimal form, see \recalt{equivalentDefSuperUnitaryRepNEW}) $(\rho_o, \Dense, \tau)$ of a super Lie group $G$ on the super Hilbert space $(L^2(\RR^d; \Hilbert, \Leb^{(d)}), \inprodsym_{\mathrm{DI}}, \supinprsym_{\mathrm{DI}})$ \stresd{decomposes as the direct integral of a family of super unitary representations $(\rho_{k,o}, \Dense_k, \tau_k)$, $k\in \RR^d$ on $(\Hilbert, \inprodsym, \supinprsym)$} if the following conditions are satisfied. 
\begin{enumerate}
\item
For each $k\in \RR^d$ the triple $(\rho_{k,o}, \Dense_k, \tau_k)$ is a super unitary representation of $G$ on $(\Hilbert, \inprodsym, \supinprsym)$ in its infinitesimal form.

\item
For each $k\in \RR^d$, each $g\in \body G$ and each $\psi\in L^2(\RR^d; \Hilbert, \Leb^{(d)})$ we have
$$
\bigl(\rho_o(g)\psi\bigr)(k) = \rho_{k,o}(g)\psi(k) \equiv \bigl(\rho_{k,o}(g)\bigr)\bigl(\psi(k)\bigr)
\mapob.
$$

\item
For each $k\in \RR^d$, each $X\in \body \Liealg g$ and each $\psi\in \Dense \subset L^2(\RR^d; \Hilbert, \Leb^{(d)})$ we have
$$
\bigl(\tau(X)\psi\bigr)(k)
=
\tau_k(X)\psi(k)
\mapob,
$$
which presupposes that, for fixed $\psi\in \Dense$ the elements $\psi(k)$ belong to $\Dense_k$ for almost all $k\in \RR^d$ (with respect to the Lebesgue measure). 

\end{enumerate}

\end{definition}

In order to generalize this procedure to direct integrals over odd parameters, we note that the construction of $L^2(\RR^d; \Hilbert, \Leb^{(d)})$ can be seen as a special case of the construction of a super Hilbert space out of $C^\infty_c(M;\CA^\KK)$ described in \recals{BatchelorbundlewithHodgestarsection}.3. 
Taking $M=\CA_0^d$, we have $\body M = \RR^d$ and $M=\wod M = \Gextension\body M \equiv \Gextension \RR^d$ and (by definition of smooth functions of only even coordinates)
$$
C^\infty_{(c)}(M; \Hilbert\otimes \CA^\KK) = C^\infty_{(c)}(\RR^d; \Hilbert)
\mapob.
$$
Taking the standard (\myquote{euclidean}) metric $\mfdmetric$ on $M$, it follows easily that the metric volume form $\mathrm{Vol}_{\mfdmetric_{\body M}}$ is the Lebesgue volume form\slash measure and that, according to \recalt{metriconsuperfunctionsinEinlocalcoordinates}, the metric \recalf{firstdefofordinarymetriconsuperfunctionsinE} induced on $C^\infty_c(\RR^d; \Hilbert)$ is given by \recalf{metricondirectintegralofsuperHilbert}. 
As $C^\infty_c(\RR^d; \KK)$ is dense in $L^2(\RR^d; \Hilbert, \Leb^{(d)})$ it follows immediately that the Hilbert space associated to the construction described in \recals{BatchelorbundlewithHodgestarsection}.3 is exactly $L^2(\RR^d; \Hilbert, \Leb^{(d)})$ as defined above. 

Moreover, by nearly the same argument, the trivializing density $\nu_\mfdmetric$ is the Lebesgue volume form\slash measure and the super scalar product on $C^\infty_c(M; \Hilbert\otimes \CA^\KK)$ defined in \recalt{superscalarprforHilbertvaluedfunctions}\slash \recalt{descriptionsuperscalarproductBginadaptedBatlas} is exactly \recalf{superscalarprodondirectintegralofsuperHilbert}. 
And thus we can see the direct integral of the super Hilbert space $(\Hilbert, \inprodsym, \supinprsym)$ over $\RR^d$ as a special case of the procedure described earlier (including the assumption that the super scalar product $\supinprsym_\mfdmetric$ remains non-degenerate when extended to $L^2(\RR^d; \Hilbert, \Leb^{(d)})$).

\begin{definition}[defofdirectintegraloverCA1nofsuperunireps]{Definitions}
Let $(\Hilbert, \inprodsym, \supinprsym)$ be a super Hilbert space over $\KK$. 
By analogy with the \myquote{ordinary} case we define \stresd{the direct integral over $\CA_1^n$ of $(\Hilbert, \inprodsym, \supinprsym)$} as the Hilbert space $C^\infty(\CA_1^n; \Hilbert\otimes \CA^\KK)$ equipped with the super Hilbert space structure defined by the standard metric $\mfdmetric_o$ \recalf{standardmetriconCAsub1highn} on $\CA_1^n$ as described in \recalt{metriconsuperfunctionsinEinlocalcoordinates} and \recalt{descriptionsuperscalarproductBginadaptedBatlas} (note that in this case, just as in the case above, we have $\nu_{\mfdmetric_o}=\nu_{B,\mfdmetric_o}$). 
We thus have the identification
$$
C^\infty(\CA_1^n, \Hilbert \otimes \CA^\KK)
\cong
\Hilbert^{2^n}
\mapob,
$$
with, for $\psi\in C^\infty(\CA_1^n, \Hilbert \otimes \CA^\KK)$:
$$
\psi\cong (\psi_I)_{I\subset \{1, \dots, n\}}\in \Hilbert^{2^n}
\qquad\Longleftrightarrow\qquad
\psi(\xi) = \sum_{I\in \{1, \dots, n\} } \xi^I\,\psi_I 
\mapob.
$$
Moreover, the metric and super scalar product are given by
$$
\inprod\chi\psi_{\mfdmetric_o} = \sum_{I\subset \{1, \dots, n\} } \inprod{\chi_I}{\psi_I}
$$
and 
$$
\homsuperinprod\chi\psi{_{\mfdmetric_o}}
=
\int_{\CA_1^n}\extder \xi^{(n)}\  \superinprod[2]{\chi(\xi)}{\psi(\xi)} 
=
\sum_{I\subset \{1, \dots, n\} } (-1)^{\varepsilon(I,I^c)}\, \superinprod[2]{\conjugate^I(\chi_{I})}{\conjugate^n(\psi_{I^c})}
\mapob.
$$

In order to mimic the definition of a decomposition of a \psur{} as a direct integral over $\CA_1^n$ of a family of \psur{}s, we have to define this notion. 
What we want to say is that it is (in infinitesimal form) a family $(\rhob_\kappa, \Dense_\kappa, \tau_\kappa)$, $\kappa\in \CA_1^n$ of \psur{}s on $\Hilbert$, but once one thinks about that, several obvious problems come to mind with respect to the fact that we have an odd parameter. 
And simply saying that we should replace $\Hilbert$ by $\Hilbert \otimes \CA^\KK$ does not solve all our problems. 
We thus start essentially from the other side, \ie, from the representation on $C^\infty(\CA_1^n; \Hilbert\otimes \CA^\KK)$ and we adapt the definition in such a way that it is compatible with a \psur{} on $C^\infty(\CA_1^n; \Hilbert\otimes \CA^\KK)$ in the obvious way. 
This leads us to the following definition.\footnote{Even though this definition might look complicated, contrived and tailor made for our purpose, heuristic arguments in \cite{Tuynman:2009} strongly suggest that it is more or less the only way to do it.} 
A \stresd{smooth family of \psur{}s (in infinitesimal form) of $G$ on $\Hilbert$} is a triple $(\rhob, (\Dense_I)_{I\subset \{1, \dots, n\}}, (\tau_I)_{I\subset \{1, \dots, n\}})$ with the following $9$ (nine) properties.
\begin{enumerate}
\item
$\rhob$ is an even unitary representation of $\body G$ on $\Hilbert$ (we have no $\kappa$-dependence here).

\item
$\Dense_I\subset C^\infty(\rhob)$ is a dense graded subspace of $\Hilbert$ invariant under the action of $\rhob$.

\item\label{DIoddpropertyofdomainofthetausubI}
$\forall X\in \body \Liealg g$ we have a map $\tau_I(X): \sum_{J\cap I=\emptyset} \Dense_J \to \Hilbert$.

\item
For homogeneous $X\in \body \Liealg g$ the map $\tau_I$ is homogeneous of parity $\parity{\tau_I} = \parity X + \parity I$. 

\end{enumerate}
With these ingredients we can define the sets $\Dense_\kappa\subset \Hilbert \otimes \CA^\KK$, which are graded submodules (but not necessarily graded subspaces of the form $V\otimes \CA^\KK$ for some graded subspace $V\subset \Hilbert$) by 
$$
\Dense_\kappa = \sum_{I\subset \{1, \dots, n\}} \kappa^I \cdot \Dense_I \otimes \CA^\KK
\mapob.
$$
We then can use property (\ref{DIoddpropertyofdomainofthetausubI}) to show that we can define, for all $X\in \body \Liealg g$, a map $\tau_\kappa(X)$ on $\Dense_\kappa$ by
$$
\tau_\kappa(X) = \sum_{I\subset \{1, \dots, n\}} \kappa^I \cdot \tau_I(X)
\mapob,
$$
simply because for $\psi\in \Dense_\kappa$ we have $\psi=\sum_J \kappa^J \psi_J$ and thus
$$
\tau_\kappa(X)\psi
=
\sum_{I,J} \kappa^I\,\tau_I(X)\kappa^J\,\psi_J
=
\sum_{I,J, I\cap J=\emptyset} (-1)^{\varepsilon(I,J)}\,\kappa^{I\cup J} \, \tau_I(X)\psi_J
\mapob.
$$
And then we require in addition the following properties.
\begin{enumerate}
\setcounter{enumi}{4}
\item
$\forall X\in \body \Liealg g$ the image $\bigl(\tau_\kappa(X)\bigr)(\Dense_\kappa)$ is contained in $\Dense_\kappa$.

\item
The maps $\tau_\kappa:\body \Liealg g \to \End(\Dense_\kappa)$ are even graded Lie algebra morphisms.

\item
For each $X\in \body \Liealg g_0$ the map $\tau_\kappa(X)$ is the restriction of the infinitesimal generator of $\rhob\bigl(\exp(tX)\bigr)$ to $\Dense_\kappa$.

\item
For all $X\in \body \Liealg g$ the map $\tau_\kappa(X)$ is graded skew-symmetric with respect to $\supinprsym_{\mfdmetric_o}$.

\item
For all $g\in \body G$ and all $X\in \body \Liealg g_1$ we have
$$
\tau_\kappa\bigl(\Ad(g)X\bigr) = \rhob(g)\scirc \tau_\kappa(X) \scirc \rhob(g\mo)
\mapob.
$$

\end{enumerate}
We will say that a smooth family of \psur{}s $(\rhob, (\Dense_I)_{I\subset \{1, \dots, n\}},\allowbreak (\tau_I)_{I\subset \{1, \dots, n\}})$ is a \stresd{smooth family of super unitary representations} when it is maximal within the given constraints (in the sense that there is no non-trivial extension).

When we compare the conditions (i) and (v)--(ix)  with the properties of an infinitesimal \psur{} \recalt{equivalentDefSuperUnitaryRepNEW}, we see that these conditions essentially say that the triple $(\rhob,\Dense_\kappa, \tau_\kappa)$ is an infinitesimal \psur. 
Except of course that $\Dense_\kappa$ is not a dense subspace of $C^\infty(\rhob)$ but of $C^\infty(\rhob) \otimes \CA^\KK$. 
And that it is not a graded subspace, but simply a graded submodule. 
And of course, $\Dense_\kappa$ is not just any graded submodule, but of a particular form, as are the maps $\tau_\kappa$. 
In the sequel we will indicate such a smooth family of \psur{}s as a triple $(\rhob,\Dense_\kappa, \tau_\kappa)$, $\kappa\in \CA_1^n$, but it should be understood that these objects are defined as above in terms of the families $(\Dense_I)_{I\subset \{1, \dots, n\}}$ and $(\tau_I)_{I\subset \{1, \dots, n\}}$ satisfying the conditions (ii)--(iv).

As for the case of a direct integral of a super Hilbert space over $\RR^d$, we now say that a super unitary representation (in its infinitesimal form) $(\rho_o, \Dense, \tau)$ of a super Lie group $G$ on $\bigl(C^\infty(\CA_1^n, \Hilbert \otimes \CA^\KK), \inprodsym_{\mfdmetric_o}, \supinprsym_{\mfdmetric_o}\bigr)$ \stresd{decomposes as the direct integral of a smooth family of super unitary representations $(\rhob, \Dense_\kappa, \tau_\kappa)$, $\kappa\in \CA_1^n$ on $(\Hilbert, \inprodsym, \supinprsym)$} if the following conditions are satisfied. 
\begin{enumerate}
\item
For all $\kappa\in \CA_1^n$ and all $\psi\in \Dense$ we have $\psi(\kappa)\in \Dense_\kappa$. 

\item
For all $\kappa\in \CA_1^n$, all $g\in \body G$ and all $\psi\in C^\infty(\CA_1^n, \Hilbert \otimes \CA^\KK)$ we have
$$
\bigl(\rho_o(g)\psi\bigr)(\kappa) = \rhob(g)\psi(\kappa) \equiv \bigl(\rhob(g)\bigr)\bigl(\psi(\kappa)\bigr)
\mapob.
$$

\item
For all $\kappa\in \CA_1^n$, all $X\in \body \Liealg g$ and all $\psi\in \Dense$ we have
$$
\bigl(\tau(X)\psi\bigr)(\kappa)
=
\tau_\kappa(X)\psi(\kappa)
\mapob.
$$

\end{enumerate}
The only difference with a decomposition as a direct integral of \psur{}s over $\RR^d$ is that the unitary representation of $\body G$ is not allowed to depend upon the odd parameter $\kappa\in \CA_1^n$, which is allowed with the real parameter $k\in \RR^d$. 

\end{definition}

\begin{definition}{Remarks}
$\bullet$
A tedious but straightforward verification of the proof of \recalt{equivalentDefSuperUnitaryRepNEW} shows that we can integrate the infinitesimal version of a smooth family $(\rhob, \Dense_\kappa, \tau_\kappa)$ of \psur{}s depending upon an odd parameter $\kappa\in \CA_1^n$ to a family $\rho_\kappa$ of group homomorphisms 
$$
\rho_\kappa:G\to \Aut(\Dense_\kappa)
$$
with the following properties:
\begin{enumerate}
\item
For all $g\in  G$ and all $\kappa\in \CA_1^n$, $\rho_\kappa(g)$ preserves $\supinprsym$ and

\item
for all $g\in \body G$ and all $\kappa\in \CA_1^n$ the restriction $\rhob(g)\caprestricted_{C^\infty(\rho_o)}$ equals $\rho_\kappa(g)$.

\end{enumerate}

$\bullet$
The fact that for $X\in \body\Liealg g_0$ the map $\tau_\kappa(X)$ must be the generator of $\rho_o\bigl( \exp(tX)\bigr)$ implies that the restriction of $\tau_\kappa$ to $\Liealg g_0$ is independent of $\kappa$. 
The dependence of $\kappa$ thus can only concern the odd part of $\Liealg g$.

\end{definition}

\masection{Back to our motivating example}
\label{backtomotivtingexamplesection}

Recall that our super Lie group is given as $G=E_0$ with $E$ a graded vector space of dimension $1\vert n$ and multiplication defined by
$$
(x,\xi)\cdot (x',\xi') = (\,x+x'+\tfrac12\,\inprod\xi{\xi'}, \xi+\xi' \,)
\mapob,
$$
where $(x,\xi_1, \dots, \xi_n)$ denote global coordinates on $E_0$ and where $\inprod\xi{\xi'} = \sum_{j=1}^n\xi_j\,\xi'_j$. 
For this $G$, the natural basis $e, f_1, \dots, f_n$ of $\Liealg g \cong T_eG \cong E$ gives us the following left-invariant vector fields and associated right-invariant vector fields:
$$
\vec e = \partial_x
\quad,\quad
\vec f_j = \partial_{\xi_j} - \tfrac12\xi_j\,\partial_x
\qquad\text{and}\qquad
e^R = \partial_x
\quad,\quad
f_j^R
=
\partial_{\xi_j} + \tfrac12\xi_j\,\partial_x
\mapob.
$$
The associated basis of the left-invariant $1$-forms is given by
$$
\extder \xi_j
\quad,\quad
\extder x - \tfrac12\sum_j \xi_j\,\extder\xi_j
\mapob.
$$
It is not hard to show that the exponential map on $\oddp{\Liealg g}_0$ is given by
$$
\exp\Bigl(\ \sum_{j=1}^n \xi_j\,f_j\,\Bigr) = (0,\xi)
\mapob,
$$
which shows that the coordinates $(x,\xi)$ are already the coordinates adapted to the \myquote{decomposition} $G \cong \wod G \times \oddp{\Liealg g}_0$.

Using the \myquote{standard} super metric on $\Liealg g \cong T_eG$ (see \recals{leftregularrepresentationssection}), it it not hard to show that the extension to a left-invariant super metric on $G$ is given by
\begin{gather*}
\mfdmetric\Bigl( \fracp{}x\bigrestricted_{(x,\xi)}\ ,\ \fracp{}x\bigrestricted_{(x,\xi)}\,\Bigr)
=
1
\qquad,\qquad
\mfdmetric\Bigl( \fracp{}{\xi_j}\bigrestricted_{(x,\xi)}\ ,\ \fracp{}{\xi_k}\bigrestricted_{(x,\xi)}\,\Bigr)
=
i\,\delta_{jk} + \tfrac14\xi_j\xi_k
\\
\mfdmetric\Bigl( \fracp{}x\bigrestricted_{(x,\xi)}\ ,\ \fracp{}{\xi_j}\bigrestricted_{(x,\xi)}\,\Bigr)
=
\tfrac12\,\xi_j
=
\mfdmetric\Bigl( \fracp{}{\xi_j}\bigrestricted_{(x,\xi)}\ ,\ \fracp{}{x}\bigrestricted_{(x,\xi)}\,\Bigr)
\mapob.
\end{gather*}
Its matrix thus is
$$
\mathrm{matrix}(\mfdmetric) = \begin{pmatrix}
1 & \tfrac12\xi_k \\ \tfrac12\xi_j & i\delta_{jk}+\tfrac14\xi_j\xi_k
\end{pmatrix}
\mapob,
$$
whose graded determinant is given by $\Ber(\mfdmetric_{jk}) = i^n$, which is independent of the odd coordinates. 
It follows that the trivializing densities $\nu_\mfdmetric$ and $\nu_{B,\mfdmetric}$ are equal and left-invariant. 
Moreover, the induced volume form $\nu_o$ on $\body G=\RR$ \recalt{variousdensitiesonGorbodyGassotoinvariantmfdmetric} is exactly the Lebesgue measure $\Leb$. 
According to \recalt{metricandcompletionforleftregularrep} we thus have the (natural) identification
$$
L^2(G) \cong \bigl(L^2(\RR)\bigr)^{2^n}
$$
with the representation $\rho_o$ being $2^n$ copies of the standard left-regular representation $\rho_\RR$ of $\body G=\RR$ (with addition as group law); it is the Hilbert space we heuristically \myquote{imposed} in \recals{motivatingexamplesection}. 

Writing $\chi(x,\xi) = \sum_I \xi^I\,\chi_I(x)$ (and similarly for $\psi$), it follows directly from \recalt{metricandcompletionforleftregularrep}, \recalf{explicitformofsuperinprodBmfdmetriconleftregrep} and the fact that $\nu_{B,\mfdmetric} = \nu_\mfdmetric$ that the metric and super scalar product on $\Hilbert$ are given by
\begin{align*}
\inprod\chi\psi_\mfdmetric
&
=
\sum_{I\subset \{1, \dots, n\} }
\int_\RR \extder x\ \overline{\chi_I(x)}\,\psi_I(x)
\equiv
\sum_{I\subset \{1, \dots, n\} }
\inprodd{\chi_I}{\psi_I}
\\
\homsuperinprod\chi\psi{_\mfdmetric}
&
=
\int_\RR\extder x\ \int_{\CA_1^n} \extder\xi^{(n)}\ 
\overline{\chi(x,\xi)}\cdot \psi(x,\xi)
\\&
=
\sum_I (-1)^{\varepsilon(I,I^c)}\cdot
\int_\RR\extder x\ \overline{\chi_I(x)}\cdot \psi_{I^c}(x)
\equiv
\sum_I (-1)^{\varepsilon(I,I^c)}\cdot
\inprodd{\chi_I}{\psi_{I^c}}
\mapob,
\end{align*}
the last of which was heuristically introduced in \recals{motivatingexamplesection} (and simply denoted by $\supinprsym$). 

In order to determine the graded subspace $\Dense_\rho$ of the super unitary (left-regular) representation $(\Dense_\rho,\rho)$ of $G$, we first note that, according to the result by Poulsen, we have
$$
C^\infty(\rho_\RR) = \{\,\psi\in C^\infty(\RR) \mid \forall k\in \NN : \psi^{(k)} \in L^2(\RR)\,\}
\mapob.
$$
Looking at the description \recalf{alternatesimpifieddescriptionDensesubrho} of $\Dense_\rho$, we have to determine the differential operators
$$
\tau(f_j)
=
-f_j^R
\equiv
-\partial_{\xi_j} - \tfrac12\xi_j\,\partial_x
\mapob.
$$
In the decomposition $\psi\cong (\psi_I)_{I\subset\{1, \dots, n\}}$, it is easy to show that we have
$$
\bigl(\tau(f_j)\psi\bigr)_I
=
\begin{cases}
-\tfrac12\, (-1)^{\varepsilon(\{j\},I\setminus\{j\})}\,\psi'_{I\setminus\{j\}} & \  j\in I
\\
-(-1)^{\varepsilon(\{j\},I)}\,\psi_{I\cup\{j\}} & \  j\notin I
\mapob.
\end{cases}
$$
In other words, $\tau(f_j)$ permutes the elements of $\bigl(C^\infty(\rho_\RR)\bigr){}^{2^n} \cong C^\infty(\rho_o)$, sometimes \myquote{adding} a derivative and\slash or a coefficient. 
But $C^\infty(\rho_\RR)$ is invariant under taking derivatives, and thus we find
$$
\Dense_\rho = C^\infty(\rho_o) \equiv \bigl(C^\infty(\rho_\RR)\bigr){}^{2^n}
\mapob.
$$

\bigskip

\noindent\textbf{Fourier decomposition via the central coordinate}
\medskip

\noindent\textbf{Nota Bene.}
Till now we could have chosen arbitrarily $\KK=\RR$ or $\KK=\CC$. 
However, once we apply a Fourier transform, the choice $\KK=\CC$ is imposed. 

\bigskip

We denote by $\Fourier:L^2(\RR)\to L^2(\RR)$ the standard Fourier transform defined by
$$
\bigl(\Fourier(\psi_I)\bigr)(k) = 
\frac1{\sqrt{2\pi}} \int_\RR \eexp^{-ikx}\,\psi_I(x)\ \extder x
$$
and we extend it to a map $\Fourier:L^2(G)\cong \bigl(L^2(\RR)^{2^n}\bigr) \to L^2(G)$ (slight abuse of notation) by applying $\Fourier$ to each component separately:
$$
\bigl(\Fourier(\psi)\bigr)_I = \Fourier(\psi_I)
\mapob.
$$
As, on $L^2(\RR)$, the Fourier transform is unitary, this extension to $L^2(G)$ preserves both the metric $\inprodsym_\mfdmetric$ and the super scalar product $\supinprsym_\mfdmetric$ on $L^2(G)$, \ie, it is an even equivalence of super Hilbert spaces. 
It follows immediately that the triple $(\rhoh_o,\hat\Dense, \tauh)$ defined by
$$
\rhoh_o(g) = \Fourier\scirc \rho_o(g)\scirc \Fourier\mo
\qquad\text{and}\qquad
\tauh(X) = \Fourier\scirc \tau(X)\scirc \Fourier\mo
$$
for $g\in \body G$ and $\body X\in \Liealg g$ and
$$
\hat\Dense = \Fourier(\Dense) = C^\infty(\rhoh_o)
$$
is a super unitary representation on $(L^2(G), \inprodsym_\mfdmetric, \supinprsym_\mfdmetric)$. 
When we now remember that (forgetting to denote the \Gextension{}) a family $(\psi_I)_{I\subset \{1, \dots, n\}}\in \Hilbert$ represents the function $\psi:G\to \CA^\CC$ defined by
$
\psi(x,\xi) = \sum_{I\subset \{1, \dots, n\} } \xi^I\,\psi_I(x)
$,
it is immediate that the image $\psih=\Fourier\psi\in \Hilbert$ represents the function $\psih:G\to \CA^\CC$ defined by
\begin{align*}
\psih(k,\xi) 
&
= 
\sum_{I\subset \{1, \dots, n\} }\xi^i\,\psih_I(k)
=
\sum_{I\subset \{1, \dots, n\} }\xi^i\,(\Fourier\psi_I)(k)
\\&
=
\frac1{\sqrt{2\pi}} \int_\RR \eexp^{-ikx}\sum_{I\subset \{1, \dots, n\} }\xi^i\,\psi_I(x)\ \extder x
=
\frac1{\sqrt{2\pi}} \int_\RR \eexp^{-ikx}\psi(x)\ \extder x
\mapob.
\end{align*}
In other words, apart from the fact that we forgot about taking the \Gextension{}, the map $\Fourier:\Hilbert\to \Hilbert$ is given by the standard formula for the Fourier transform. 

As $\rho_o$ is ($2^n$ copies of) the left-regular representation of $\RR$, it follows immediately that, for $y\in \RR$, \ie $(y,0)\in \body G$ and $\psih \cong(\psih_I)_{I\subset \{1, \dots, n\}}\in \Hilbert$ we have:
$$
\bigl(\rhoh_o(y,0)\psih\bigr)_I(k) = \eexp^{-iky}\cdot \psih_I(k)
\mapob.
$$
Moreover, for $\psih\in C^\infty(\rhoh_o) = \Fourier\bigl(C^\infty(\rho_o)\bigr)$ we have
$$
\bigl(\tauh(e)\psih\bigr)_I(k) = -ik\,\psih_I(k)
$$
and, for $j=1, \dots, n$:
$$
\bigl(\tauh(f_j)\psih\bigr)_I(k) = 
\begin{cases}
-\tfrac12\,ik\, (-1)^{\varepsilon(\{j\},I\setminus\{j\})}\,\psih_{I\setminus\{j\}}(k) & \  j\in I
\\[2\jot]
-(-1)^{\varepsilon(\{j\},I)}\,\psih_{I\cup\{j\}}(k) & \  j\notin I
\mapob.
\end{cases}
$$
Combining the components $\psih_I$ into a single function $\psih:G\to \CA^\CC$, these formul{\ae} can be written as
$$
\bigl(\rhoh_o(y,0)\psih\bigr)(k,\xi)
=
\eexp^{-iky}\,\psih(k,\xi)
$$
and
$$
\bigl(\tauh(e)\psih\bigr)(k,\xi) 
=
-ik\,\psih(k,\xi)
\quad,\quad
\bigl(\tauh(f_j)\psih\bigr)(k,\xi) = 
-(\partial_{\xi_j}\psih)(k,\xi) - \tfrac12\, ik\,\xi_j\,\psih(k,\xi)
\mapob.
$$
According to \recalt{equivalentDefSuperUnitaryRepNEW} this infinitesimal form of a super unitary representation should integrate to a super unitary representation $\bigl(C^\infty(\rhoh_o),\rhoh\bigr)$. 
Using the results heuristically found in \recals{motivatingexamplesection}, it is not hard to show that $\rhoh$ is (indeed) given by
$$
\bigl(\rhoh(y,\eta)\psih\bigr)(k,\xi) = \eexp^{-ik(y+\frac12\inprod\eta\xi)} \, \psih(k,\xi-\eta)
\mapob.
$$

\medskip

We now note that the super Hilbert space $(L^2(G), \inprodsym_\mfdmetric, \supinprsym_\mfdmetric)$ is the direct integral over $\RR$ of the super Hilbert space $\Hilbert = C^\infty(\CA_1^n; \CA^\CC)$ or more precisely  the super Hilbert space associated to the $\CA$-manifold $\CA_1^n$ equipped with the standard constant metric $\mfdmetric_o$ \recalf{standardmetriconCAsub1highn} according to the procedure described in \recalt{metriconsuperfunctionsinEinlocalcoordinates} and \recalt{descriptionsuperscalarproductBginadaptedBatlas}; its metric and super scalar product are given by
$$
\inprod{f}{g}_{\mfdmetric_o} = \sum_{I\subset \{1, \dots, n\} } \overline{f_I}\cdot g_I
\quad\text{and}\quad
\homsuperinprod{f}{g}{_{\mfdmetric_o}} = \sum_{I\subset \{1, \dots, n\} } (-1)^{\varepsilon(I,I^c)}\,\overline{f_I}\cdot g_{I^c}
\mapob,
$$
for $f(\xi) = \sum_{I\subset \{1, \dots, n\}} \xi^I\,f_I$, $f_I\in \CC$ (and similarly for $g$). 
More important is that the family of representations $(\Hilbert,\rhoh_k)$ of $G$ on this (finite-dimensional) super Hilbert space $\Hilbert=C^\infty(\CA_1^n; \CA^\CC)$ defined by
$$
\bigl(\rhoh_{k}(y,\eta)f\bigr)(\xi)
=
\eexp^{-ik(y+\frac12\inprod\eta\xi)} \, f(\xi-\eta)
$$
is a family of super unitary representations (invariance of $\supinprsym_{\mfdmetric_o}$ can be checked using \recalt{superscalarproductmfdmetricinvariantunderGaction}). 
It then is immediate that the super unitary representation $(\hat\Dense,\rhoh)$ decomposes as the direct integral over $\RR$ of the family $(\Hilbert, \rhoh_k)$, $k\in \RR$ of super unitary representations (see \recalt{defofdirectintegraloverRkofsuperunireps}). 
And thus we have fully justified our heuristic (Fourier) decomposition of the (left-regular) super unitary representation $(\Dense,\rho)$ as a direct integral over $\RR$ of the family of super unitary representations $(\Hilbert, \rhoh_k)$, $k\in \RR$.

\bigskip

\noindent\textbf{Decomposing the 0-Fourier mode}
\medskip

The super unitary representation $\rhoh_{k=0}$ on $\Hilbert \equiv C^\infty(\CA_1^n; \CA^\CC)$ is given by
$$
\bigl(\rhoh_{k=0}(y,\eta)f\bigr)(\xi) = f(\xi-\eta)
\mapob,
$$
or in its infinitesimal form by
$$
\tauh_{k=0}(e) = 0
\qquad\text{and}\qquad
\tauh_{k=0}(f_j) = -\partial_{\xi_j}
\mapob.
$$
We now want to show that the Berezin-Fourier decomposition of this representation is the direct integral over $\CA_1^n$ of 1-dimensional (irreducible) super unitary representations $\rhob_\kappa$, $\kappa\in \CA_1^n$ on the super Hilbert space $(\RR, \inprodsym, \supinprsym=\inprodsym)$ (see \recalt{defofdirectintegraloverCA1nofsuperunireps}). 
To do so, we follow exactly the same procedure as for the ordinary Fourier decomposition. 

We first recall the elementary Berezin-Fourier transform $\Fourierodd : \Hilbert \to \Hilbert$ introduced in \recals{BerezinFourierandHodgestarsection} and which is an equivalence of super Hilbert spaces \recalt{generalizedBerezinFourierisequivalenceSHS} of parity $\parity n$:
$$
(\Fourierodd f)(\kappa)
= 
\int_{\CA_1^n}\extder\xi\  \exp\Bigl({-i\,\sum_{p=1}^n\xi_p\kappa_p}\Bigr)\,f(\xi) 
\mapob.
$$
We then define the map $\rhob$ by
$$
\rhob(X) = \Fourierodd\scirc \rhoh_{k=0}(X)\scirc \Fourierodd\mo
\mapob.
$$
Using \recalt{BerezinFourierasproductofInvsNEW} and \recalt{generalizedBerezinFourierisequivalenceSHS} it is elementary to show that $\bigl( \Hilbert ,\rhob\bigr)$ is a super unitary representation of $G$ on $C^\infty(\CA_1^n; \CA^\CC)$. 
Using the explicit formula for $\Fourierodd$, it is elementary to show that we have
$$
\bigl(\rhob(y,\eta)\bar f\bigr)(\kappa) = \eexp^{i\,\inprod\kappa\eta}\,\bar f(\kappa)
\mapob,
$$
or in infinitesimal form
$$
\taub(e) = 0
\qquad\text{and}\qquad
\taub(f_j) = -i\kappa_j
\mapob,
$$
from which one also can deduce immediately that it is a super unitary representation. 

Combining the various definitions, it is then not hard to see that $\Hilbert \equiv C^\infty(\CA_1^n; \CA^\CC)$ is the direct integral over $\CA_1^n$ of the super Hilbert space $(\CC; \inprodsym, \supinprsym=\inprodsym)$, that the family $(\id, \Dense_\kappa, \taub_\kappa)$ defined by $\Dense_\kappa = \CA^\CC$ (or equivalently, $\Dense_I = \CC$ for all $I\subset \{1, \dots, n\}$) as well as
$$
\taub_\kappa(e) = 0
\qquad\text{and}\qquad
\taub_\kappa(f_j) = -i\kappa_j
\mapob,
$$
is a smooth family of ($1$-dimensional, irreducible) super unitary representations of $G$ in infinitesimal form on $(\CC; \inprodsym, \supinprsym=\inprodsym)$ and that these infinitesimal representations integrate to the family $(\Dense_\kappa,\rhob_\kappa)$, $\kappa\in \CA_1^n$ defined by
$$
\rhob_\kappa(y,\eta) v = \eexp^{i\,\inprod\kappa\eta}\,v
\mapob.
$$
It then is immediate that the super unitary representation $(\Hilbert,\rhob)$ decomposes as the direct integral of the smooth family of super unitary representations $(\id, \Dense_\kappa, \taub_\kappa)$.

\bigskip

\noindent\textbf{The case $k\neq0$ for $n=2$}
\medskip

The computations made in \recals{motivatingexamplesection} for this case are done in finite dimensions and thus need no additional analytical justifications. 
However, only algebraic computations were made there, no mention has been made of the super Hilbert space structures. 
In terms of the isomorphism $C^\infty(\CA_1^2; \CA^\CC) \cong \CC^4$ given by
$$
\psi(\xi_1,\xi_2) = \psi_0+\xi_1\,\psi_1 + \xi_2\,\psi_2+\xi_1\xi_2\,\psi_{12}
\qquad\mapsto\qquad
(\psi_0,\psi_1,\psi_2,\psi_{12})\in \CC^4
$$
the metric $\inprodsym_{\mfdmetric_o}$ and the super scalar product $\supinprsym_{\mfdmetric_o}$ are given by
\begin{align*}
\inprod{\chi}{\psi}_{\mfdmetric_o} 
&
= \overline{\chi_0}\cdot \psi_0 + \overline{\chi_1}\cdot \psi_1 + \overline{\chi_2}\cdot \psi_2 + \overline{\chi_{12}}\cdot \psi_{12}
\\[2\jot]
\homsuperinprod{\chi}{\psi}{_{\mfdmetric_o}} 
&
=
\psi_0\omega_{12} + \psi_{12}\omega_0 + \psi_1\omega_2 - \psi_2\omega_1
\mapob.
\end{align*}
It is elementary to show that the restrictions of  $\inprodsym_{\mfdmetric_o}$ and $\supinprsym_{\mfdmetric_o}$ to each subspace $\Hilberth_\varepsilon$ turn these two subspaces into super Hilbert spaces (non-degeneracy of $\supinprsym_{\mfdmetric_o}$ on $\Hilberth_\varepsilon$ is a consequence of the condition $k\neq0$). 
Moreover, we have $\homsuperinprod{\Hilberth_1}{\Hilberth_{-1}}{_{\mfdmetric_o}}=0$, \ie, these two subspaces are {orthogonal} with respect to $\supinprsym_{\mfdmetric_o}$. 
On the other hand, except for the case $1-\tfrac14k^2=0$, they are not orthogonal with respect to $\inprodsym_{\mfdmetric_o}$.

\begin{definition}{Remarks}
$\bullet$
For fixed $k$ the representations obtained by restricting $\rhoh_k$ to $\Hilberth_\varepsilon$ are non-isomorphic irreducible. 
On the other hand, complex conjugation intertwines $\rhoh_k$ on $\Hilberth_\varepsilon$ with $\rhoh_{-k}$ on $\Hilbert_{-\varepsilon}$. 
To be more precise, we have a family of representations $\rhoh_k$ on $\Hilberth^k$ with $\Hilberth^k = C^\infty(\CA_1^2; \CA^\CC)$ for all $k\in \RR$. 
We thus can consider the map (an anti-linear isomorphism) $\mathcal{C}:\Hilberth^k\to \Hilberth^{-k}$ defined as $\bigl(\mathcal{C}(\psi)\bigr)(\xi) = \overline{\psi(\xi)}$. 
And a direct computation shows that we have the commutative diagram
$$
\begin{CD}
\Hilberth^k_\varepsilon @>\rhoh_k>> \Hilberth^k_\varepsilon
\\
@V\mathcal{C}VV @VV\mathcal{C}V
\\
\Hilberth^{-k}_{-\varepsilon} @>>\rhoh_{-k}> \Hilberth^{-k}_{-\varepsilon}\rlap{\mapob.}
\end{CD}
$$

\bigskip

$\bullet$
It is not only tempting, but even recommended to interpret the subspaces $\Hilberth_\varepsilon$ as spaces of \myquote{holomorphic} functions of the single odd-complex variable $\zeta={\xi_1+i\varepsilon\,\xi_2}$ completed with the \myquote{Gaussian weight} $\eexp^{-\frac12\varepsilon k\xi_1\xi_2}$. 
They can be obtained in a straightforward manner by geometric quantization of a coadjoint orbit of $G$ with the {(anti-)} holomorphic polarization as explained in \cite{Tuynman:2009}. They are the direct analogues of the Bargmann representation of the space $L^2(\RR^n)$ as the space of (anti-) holomorphic functions on $\CC^n$ with the Gaussian weight $\eexp^{-z\zb}$, a space that is obtained by geometric quantization with the (anti-) holomorphic polarization. 
The only difference is that in the \myquote{real} case the choice of the (anti-) holomorphic polarization is fixed by the \myquote{sign} of the symplectic form, whereas in the super case, in which the spaces are finite dimensional, both choices are possible. 

\end{definition}

\masection{The \texorpdfstring{$a\xi+\beta$}{axi+beta} group}
\label{axiplusbetagroupsection}

The group $G$ of affine transformations of the odd line $\CA_1$ consists of the transformations $\zeta\mapsto a\zeta+\beta$ with $(a,\beta)\in \CA_0\times \CA_1$ with $\body a>0$. 
For computational reasons we change coordinates and we realize this group as $G=E_0$ with $E$ a graded vector space of dimension $1\vert 1$ with group law
$$
(x,\xi) \cdot (y,\eta) = (x+y, \eta + \eexp^{-y}\,\xi)
\mapob.
$$
This corresponds to the matrix representation (acting on a vector $(\zeta,1)\,$)
$$
\begin{pmatrix}
\eexp^x & \eexp^x \xi
\\
0 & 1 
\end{pmatrix}
\cdot
\begin{pmatrix}
\eexp^y & \eexp^y \eta
\\
0 & 1 
\end{pmatrix}
=
\begin{pmatrix}
\eexp^{x+y} & \eexp^{x+y}\,(\eta + \eexp^{-y}\, \xi)
\\
0 & 1 
\end{pmatrix}
\mapob.
$$
This group is interesting because it is an example in which the dense subspace $\Dense_\rho$ in the left-regular representation is smaller than $C^\infty(\rho_o)$. 

In order to establish the left-regular representation in detail, we start with a basis of the left-invariant vector fields $e,f\in \body \Liealg g$ given by
$$
\vec e\,\restricted_{(x,\xi)} = \partial_x\caprestricted_{(x,\xi)} - \xi\,\partial_\xi\caprestricted_{(x,\xi)}
\qquad\text{and}\qquad
\vec f\,\caprestricted_{(x,\xi)} = \partial_\xi\caprestricted_{(x,\xi)}
\mapob.
$$
The corresponding right-invariant vector fields are given by
$$
e^R\caprestricted_{(x,\xi)} = \partial_x\caprestricted_{(x,\xi)}
\qquad\text{and}\qquad
f^R\caprestricted_{(x,\xi)} = \eexp^{-x}\,\partial_\xi\caprestricted_{(x,\xi)}
\mapob.
$$
A direct computation shows that, in the coordinates $(x,\xi)$, the matrix of the invariant metric $\mfdmetric$ is given by
$$
\mathrm{matrix}(\mfdmetric_{(x,\xi)}) = 
\begin{pmatrix}
1 & i\xi \\ -i\xi & i 
\end{pmatrix}
\mapob.
$$
It follows immediately that the invariant measure on $\body G = \RR$ is the Lebesgue measure and that the we have $\Delta(\xi) \equiv 1$. 
We thus have $L^2(G)\cong L^2(\RR)^2$ with identification (for smooth elements) 
$$
\psi\in L^2(G)
\qquad\leftrightarrow\qquad
(\psi_0,\psi_1)\in L^2(\RR)^2
\qquad,\qquad
\psi(x,\xi) = \psi_0(x) + \xi\,\psi_1(x)
\mapob.
$$
In terms of this identification the metric $\inprodsym_\mfdmetric$ and the super scalar product $\supinprsym_\mfdmetric$ are given by
$$
\inprod\chi\psi_\mfdmetric = \int_\RR \overline{\chi_0(x)}\,\psi_0(x) + \overline{\chi_1(x)}\,\psi_1(x) \ \extder\Leb^{(1)}(x)
\equiv
\inprodd{\chi_0}{\psi_0} + \inprodd{\chi_1}{\psi_1}
$$
and
$$
\homsuperinprod\chi\psi{_\mfdmetric} 
=
\inprodd{\chi_0}{\psi_1} + \inprodd{\chi_1}{\psi_0}
\mapob.
$$

The left-regular representation is given by
\begin{align*}
\bigl(\rho(y,\eta)\psi\bigr)(x,\xi)
&
=
\psi\bigl((y,\eta)\mo\cdot (x,\xi)\bigr)
=
\psi(x-y,\xi-\eexp^{y-x}\,\eta)
\\&
=
\psi_0(x-y) + (\xi- \eexp^{y-x}\,\eta) \psi_1(x-y)
\\&
=
\bigl( \psi_0(x-y) - \eexp^{y-x}\,\eta\, \psi_1(x-y)\bigr) + \xi\,\psi_1(x-y)
\mapob.
\end{align*}
Infinitesimally we get (taking derivatives with respect to $y$ and $\eta$ at $(y,\eta)=(0,0)\,$):
\begin{align*}
\bigl(\tau(e)\psi\bigr)(x,\xi) 
&
= 
-\psi_0'(x) - \xi\,\psi_1'(x)
=-(\partial_x\psi)(x,\xi) = (-e^R\psi)(x,\xi)
\\
\bigl(\tau(f)\psi\bigr)(x,\xi) 
&
= 
-\eexp^{-x}\,\psi_1(x)
=
-\eexp^{-x}\,(\partial_\xi\psi)(x,\xi)
=
(-f^R\psi)(x,\xi)
\mapob,
\end{align*}
which confirms that these operations are the action of minus the right-invariant vector fields. 
In terms of the components we can write the action of $\tau(e)$ and $\tau(f)$ in matrix form (acting on vectors $(\psi_0, \psi_1)\,$) as
$$
\tau(e) = 
\begin{pmatrix}
-\partial_x & 0 \\ 0 & -\partial_x
\end{pmatrix}
\qquad\text{and}\qquad
\tau(f) =
\begin{pmatrix}
0 & -\eexp^{-x} \\ 0 & 0
\end{pmatrix}
\mapob.
$$
Now the standard left-regular representation of $\RR$ has as set of smooth vectors all those smooth functions on $\RR$ for which all order derivatives belong to $L^2(\RR)$:
$$
C^\infty(\rho_{\body G}) \equiv C^\infty(\rho_\RR) = 
\{\,f\in C^\infty(\RR) \mid \forall k\in \NN: f^{(k)}\in L^2(\RR)\,\}
\mapob.
$$
This set is not invariant under multiplication by the function $\phi(x)=\eexp^{-x}$, which leads us to the definition of the set $C^\infty_\phi(\rho_\RR)$ as
$$
C^\infty_\phi(\rho_\RR)
=
\{\,f\in C^\infty(\rho_\RR) \mid \phi\cdot f\in C^\infty(\rho_\RR)\,\}
\mapob.
$$
It then is not hard to show that the maximal domain $\Dense_\rho$ for the super unitary (left-regular) representation is given by
$$
\Dense_\rho = C^\infty(\rho_\RR) \times C^\infty_\phi(\rho_\RR) \subset L^2(\RR) \times L^2(\RR) \cong \Hilbert
\mapob.
$$

\medskip

But, just as in our motivating example, we can go one step further: we can use the Berezin-Fourier transform to decompose this representation as an odd direct integral of super unitary representations. 
In order to apply the Berezin-Fourier transform, we have two points of view: we can either apply the generalized Berezin-Fourier transform $\Fourierodd_{B,\mfdmetric} : C^\infty_c(G; \CA^\CC)\to C^\infty_c(G; \CA^\CC)$ and then extend it by continuity to its completion $L^2(G;\mfdmetric)$, or we can view $L^2(G)$ as the space
$$
L^2(G) \cong C^\infty\bigl(\CA_1; L^2(\RR)\bigr)
\mapob,
$$
\ie, as the space of (smooth) functions of a single odd coordinate with values in the super Hilbert space $L^2(\RR)$. 
In both cases the map $\Fourierodd : L^2(G) \to L^2(G)$ is given (in the first case formally) by
$$
(\Fourierodd\psi)(x,\kappa) = \int \extder\xi\ \eexp^{-i\xi\kappa} \, \psi(x,\xi)
\mapob,
$$
or in terms of the representation $\psi\cong (\psi_0,\psi_1)\in L^2(\RR)^2$:
$$
\Fourierodd(\psi_0,\psi_1) = (\psi_1, -i\psi_0)
\mapob.
$$
According to \recalt{generalizedBerezinFourierisequivalenceSHS} the map $\Fourierodd$ is an odd equivalence of super Hilbert spaces, a fact that one can easily confirm by an explicit computation. 
As in \recals{backtomotivtingexamplesection} we now define the representation $\rhob$ on $L^2(G)$ by
$$
\forall g\in G : \rhob(g) = \Fourierodd \scirc \rho(g) \scirc \Fourierodd\mo
\mapob.
$$
which gives us explicitly
$$
\bigl(\rhob(y,\eta)\psib\bigr)(x,\kappa)
=
\bigl(1-i\eta\kappa\,\phi(x-y)\bigr)\psib(x-y,\kappa) 
\equiv
\eexp^{-i\eta\kappa\,\phi(x-y)}\psib(x-y,\kappa)
\mapob,
$$
or in infinitesimal form
$$
\bigl(\taub(e)\psib\bigr)(x, \kappa)
=
-(\partial_x\psib)(x,\kappa)
\qquad\text{and}\qquad
\bigl(\taub(f)\psib\bigr)(x, \kappa)
=
-i\kappa\,{\phi(x)}\,\psib(x,\kappa)
\mapob.
$$
Its domain of definition $\bar\Dense = \Fourier(\Dense)$ is given by
$$
\bar\Dense = C^\infty_\phi(\rho_\RR) \times C^\infty(\rho_\RR)
\mapob.
$$
It then is easy to see that this representation decomposes as the direct integral over $\CA_1$ of the smooth family of super unitary representations $(\rho_\RR, \Dense_\kappa, \taub_\kappa)$ on the super Hilbert space $\bigl(L^2(\RR), \inprodsym=\inprodd\cdot\cdot = \supinprsym\bigr)$, where $\rho_\RR$ is the left-regular representation of $\RR$ on $L^2(\RR)$ (and thus defined by $\bigl(\rho_\RR(y)\psi\bigr)(x) = \psi(x-y)$), where 
$$
\Dense_\kappa = \Dense_\emptyset \otimes \CA^\CC + \kappa\Dense_{\{1\}}\otimes \CA^\CC
\mapob,
$$
with $\Dense_\emptyset = C^\infty_\phi(\rho_\RR)$ and $\Dense_{\{1\}} = C^\infty(\rho_\RR)$, and where $\taub_\kappa$ is defined by
$$
\bigl(\taub_\kappa(e)\psi\bigr)(x)
=
-(\partial_x\psi)(x)
\qquad\text{and}\qquad
\bigl(\taub_\kappa(f)\psi\bigr)(x)
=
-i\kappa\,\phi(x)\,\psi(x)
\mapob,
$$
or, in the form of differential or multiplication operators:
$$
\taub_\kappa(e) = -\partial_x
\qquad\text{and}\qquad
\taub_\kappa(f) = -i\kappa\,\phi
\mapob.
$$
This smooth family of super unitary representations in infinitesimal form integrates to the family $\rhob_\kappa$ defined by
\begin{align*}
\bigl(\rhob_\kappa(y,\eta)\psib\bigr)(x) 
&
= 
\bigl(1-i\eta\kappa\,\phi(x-y)\bigr)\psib(x-y) 
\equiv
\eexp^{-i\eta\kappa\,\phi(x-y)}\psib(x-y)
\\
&
=
\bigl(\rho_\RR(y) (\eexp^{-i\eta\kappa\phi}\psib)\bigr)(x)
\equiv
\bigl(\rho_\RR(y) (\psib - i\eta\kappa\,\phi\,\psib)\bigr)(x)
\mapob,
\end{align*}
or, as translation and multiplication operators:
$$
\rhob_\kappa(y,\eta) = \rho_\RR(y) \scirc \eexp^{-i\eta\kappa\,\phi}
\mapob.
$$

\masection{The group \texorpdfstring{$\OSp(1,2)$}{OSp1,2}}
\label{OSp12section}

The group $\OSp(1,2)$ is an interesting example because the function $\Delta(\xi)$ is not constant, implying that the two super scalar products $\supinprsym_{\mfdmetric}$ and $\supinprsym_{B,\mfdmetric}$ are not the same. 

\medskip

$\OSp(1,2)$ is the (connected) subgroup of the group $\Gl(1,2)=\Aut(V)$ (with $V$ a graded vector space over $\RR$ of dimension $\dim(V)=1\vert 2$) of automorphisms that preserve the graded symmetric bilinear form 
$$
\begin{pmatrix} 1 & 0 & 0 \\ 0 & 0 & 1 \\ 0 & -1 & 0 \end{pmatrix}
\mapob.
$$ 
It can be described as
$$
\OSp(1,2) = \left\{\ 
\left.
\begin{pmatrix}
1+\alpha\beta & \alpha & \beta
\\
a\beta - b\alpha & a & b  
\\
c\beta - d\alpha  & c & d  
\end{pmatrix}
\hskip0.5em
\right\vert
\hskip0.5em
\biggl\{\ 
\begin{matrix}
a,b,c,d\in \CA_0\ \&\ \alpha,\beta\in \CA_1
\\[1\jot]
ad-bc = 1-\alpha\beta
\end{matrix}
\ 
\right\}
\mapob.
$$
Its dimension thus is $3\vert 2$ with $\body \OSp(1,2)=\Sl(2,\RR)$. 
A basis $e_1,e_2,e_3, f_1,f_2$ for its Lie algebra $\Liealg g$ is given by the matrices (see \cite{Mus12})
$$
e_1 = 
\begin{pmatrix} 0&0&0 \\ 0 & 1&0 \\ 0&0&-1 \end{pmatrix}
\qquad,\qquad
e_2 = 
\begin{pmatrix} 0&0&0 \\ 0 & 0&1 \\ 0&0&0\end{pmatrix}
\qquad,\qquad
e_3 = 
\begin{pmatrix} 0&0&0 \\ 0 & 0&0 \\ 0&1&0\end{pmatrix}
$$
and
$$
f_1 = 
\begin{pmatrix} 0&0&1 \\ -1&0&0 \\ 0 & 0 & 0 \end{pmatrix}
\qquad\text{and}\qquad
f_2 = 
\begin{pmatrix} 0 & 1&0 \\ 0 & 0&0 \\ 1 & 0 & 0 \end{pmatrix}
\mapob.
$$
An elementary computations shows that the exponential of $\xi f_1+\eta f_2$ is given by
$$
\exp(\xi f_1+\eta f_2)
=
\exp(
\begin{pmatrix} 0 & \eta&\xi \\ \xi & 0&0 \\ -\eta & 0 & 0 \end{pmatrix}
)
=
\begin{pmatrix} 1-\xi\eta & \eta&\xi \\ \xi & 1+\tfrac12\xi\eta & 0 \\ -\eta & 0 & 1+\tfrac12\xi\eta \end{pmatrix}
\mapob.
$$
It follows directly that the identification $\OSp(1,2) \cong \wod{\OSp(1,2)} \times \oddp{\Liealg g}_1 \equiv \Gextension\Sl(2,\RR) \times \CA_1^2$ is given by (but remember, we have the constraint $ad-bc=1$)
\begin{align*}
\Bigl(\ 
\begin{pmatrix} a & b \\ c & d \end{pmatrix} 
\ ,\ 
(\xi,\eta)
\ \Bigr)
&
\cong
\begin{pmatrix}
1 & 0 & 0
\\
0 & a & b  
\\
0 & c & d  
\end{pmatrix}
\cdot
\begin{pmatrix} 1-\xi\eta & \eta&\xi \\ \xi & 1+\tfrac12\xi\eta & 0 \\ -\eta & 0 & 1+\tfrac12\xi\eta \end{pmatrix}
\\[3\jot]&
=
\begin{pmatrix}
1-\xi\eta & \eta & \xi
\\
a\xi - b\eta & (1+\tfrac12\xi\eta)\,a & (1+\tfrac12\xi\eta)\,b  
\\
c\xi - d\eta & (1+\tfrac12\xi\eta)\,c & (1+\tfrac12\xi\eta)\,d  
\end{pmatrix}
\mapob.
\end{align*}
According to the general theory given in \recals{leftregularrepresentationssection}, we have to compute 
$\ad(v)$ for $v=\xi f_1+\eta f_2 \in \oddp{\Liealg g}_0$ and $\Ad(g)v$ for $g\in \body G \equiv \Sl(2,\RR)$, as well as the series $h\bigl(\ad(v)\bigr)$ and $b_\pm\bigl(\ad(v)\bigr)$. 
We thus start by computing the commutators between our basis elements:
$$
\begin{aligned}{}
[e_1,e_2] &=  2e_2
&\qquad
[e_1,e_3] &=  -2e_3
&\qquad
[e_2,e_3] &=  e_1
\\
[e_1,f_1] &= f_1
&
[e_2,f_1] &= 0
&
[e_3,f_1] &= -f_2
\\
[e_1,f_2] &= -f_2
&
[e_2,f_2] &= -f_1
&
[e_3,f_2] &= 0
\\
[f_1,f_1] &= -2e_2
&
[f_1,f_2] &= -e_1
&
[f_2,f_2] &= 2e_3
\mapob.
\end{aligned}
$$
It follows easily that the matrix of $\ad(\xi f_1+\eta f_2)$ with respect to this basis is given by
$$
\mathrm{matrix}\bigl(\ad(\xi f_1+\eta f_2)\bigr)
=
\begin{pmatrix} 
0&0&0&-\eta & -\xi 
\\ 
0&0&0& -2\xi & 0 
\\ 
0 &0&0& 0 & 2\eta 
\\ 
\xi & -\eta & 0 & 0 & 0 
\\ 
-\eta & 0 & -\xi &0&0 \end{pmatrix}
\mapob.
$$
Its square is given by
\begin{align*}
\ad(\xi f_1+\eta f_2)^2
&
=
\begin{pmatrix} 2\xi\eta & 0 & 0 & 0 & 0 
\\ 
0 & 2\xi\eta & 0 & 0 & 0 
\\ 
0 & 0 & 2\xi\eta & 0 & 0 
\\ 
0 & 0 & 0 & -3\xi\eta &0 
\\ 
0 & 0 & 0 & 0 & -3\xi\eta 
\end{pmatrix}
\end{align*}
and all higher order powers are zero, simplifying the computations of the series $h\bigl(\ad(v)\bigr)$ and $b_\pm\bigl(\ad(v)\bigr)$. 
We now recall that the matrices $A$, $B$ and $H$ are defined as follows: $A$ is the lower-left block of $\ad(\xi f_1+\eta f_2)$, $B$ is the lower right block of $b_+\bigl(\ad(v)\bigr) = \oneasmatrix+\tfrac13 \ad(\xi f_1+\eta f_2)^2$ and $H$ is the upper-right block of $h\bigl(\ad(v)\bigr) = \tfrac12 \ad(\xi f_1+\eta f_2)$. 
We thus get
$$
\begin{pmatrix}
\oneasmatrix & H \\ A & B \end{pmatrix}
\cong
\begin{pmatrix}
1 & 0 & 0 & -\tfrac12\eta & -\tfrac12 \xi
\\
0 & 1 & 0 & -\xi & 0
\\
0 & 0 & 1 & 0 & \eta
\\
\xi & -\eta & 0 & 1-\xi\eta & 0
\\
-\eta & 0 & -\xi & 0 & 1-\xi\eta
\end{pmatrix}
\mapob.
$$
This means that the function $\Delta(\xi)$ is given as
$$
\Delta(\xi,\eta) = 
\frac{\Det\bigl( B(\xi,\eta) \bigr)}{\Det\bigl( \oneasmatrix - H(\xi,\eta)\cdot B(\xi)\mo\cdot A(\xi,\eta) \bigr)}
=
1+\xi\eta
\mapob.
$$
In the identification $\Hilbert\cong \bigl(L^2(\Sl(2,\RR)\bigr)^4$, $\psi(x,\xi,\eta) = \sum_{I\subset \{1, 2\} }\xi^I\,\psi_I(x) \cong\break (\psi_\emptyset, \psi_\xi,\psi_\eta,\psi_{\xi\eta})$ we thus find (see \recalf{explicitformulasuperinproductonGwithsuperspmatrix})
\begin{align*}
\homsuperinprod{\chi}{\psi}{_\mfdmetric}
=
\inprodd{\chi_\emptyset}{\psi_\emptyset}
+
\inprodd{\chi_\emptyset}{\psi_{\xi\eta}}
+
\inprodd{\chi_\xi}{\psi_\eta}
-
\inprodd{\chi_\eta}{\psi_\xi}
+
\inprodd{\chi_{\xi\eta}}{\psi_\emptyset}
\mapob,
\end{align*}
where $\inprodd\cdot\cdot$ denotes the usual metric on $L^2\bigl(\Sl(2,\RR)\bigr)$ given by an invariant Haar measure. 

Our next task is to establish the right-invariant vector fields in terms of the decomposition $\OSp(1,2)\cong \Gextension \,\Sl(2,\RR) \otimes \CA_1^2$ as given by \recalt{leftinvariantvfonGwodtimesgsub1}:
$$
TR_{(g,v)} Y\caprestricted_{(e,0)} 
=
-{TR_{g}}\Bigl(h\bigl(\ad(\Ad(g)v)\bigr)Y\caprestricted_e\Bigr)+ \Bigl(b_-\bigl(\ad(v)\bigr)\Ad(g\mo)Y\Bigr)\bmidrestricted_v
$$
For $g\in \wod{\OSp(1,2)} \cong \Gextension \,\Sl(2,\RR)$ with 
$$
g = \begin{pmatrix}
1 & 0 & 0
\\
0 & a & b  
\\
0 & c & d  
\end{pmatrix}
$$
one easily computed
$$
\Ad(g)
=
\begin{pmatrix}
ad+bc & -ac & bd & 0 & 0
\\
-2ab & a^2 & -b^2 & 0 & 0
\\
2cd & -c^2 & d^2 & 0 & 0
\\
0 & 0 & 0 & a & -b
\\
0 & 0 & 0 & -c & d
\end{pmatrix}
\mapob.
$$
With $Y=y_1f_1+y_2f_2 \cong \bigl( \begin{smallmatrix} y_1 \\ y_2 \end{smallmatrix}\bigr)$, $y_i\in\RR$ and $b_-\bigl(\ad(v)\bigr) = \oneasmatrix - \tfrac16\ad(v)^2$ we thus find (in terms of the basis $e_\bullet$, $f_\bullet$)
\begin{align*}
b_-\bigl(\ad(v)\bigr)\Ad(g\mo)Y
&
=
\begin{pmatrix}
1+\tfrac12\xi\eta & 0 \\ 0 & 1+\tfrac12\xi\eta
\end{pmatrix}
\cdot
\begin{pmatrix} d & b \\ c & a \end{pmatrix}\cdot \begin{pmatrix} y_1 \\ y_2 \end{pmatrix}
\\&
=
(1+\tfrac12\xi\eta)\cdot\begin{pmatrix} dy_1 + by_2 \\ cy_1 + ay_2  \end{pmatrix}
\end{align*}
and (only noting the components with respect to the basis vectors $e_\bullet$ because the components with respect to $f_\bullet$ are necessarily zero)
$$
h\bigl(\ad(\Ad(g)v)\bigr)Y
=
\begin{pmatrix}
-\tfrac12(d\eta - c\xi) & -\tfrac12 (a\xi - b\eta)
\\
-(a\xi - b\eta) & 0
\\
0 & d\eta - c\xi
\end{pmatrix}\cdot \begin{pmatrix} y_1 \\ y_2 \end{pmatrix}
\mapob,
$$
which means that we have
\begin{align*}
-h\bigl(\ad(\Ad(g)v)\bigr)Y
&
=
\tfrac12 \bigl( (d\eta - c\xi)y_1 + (a\xi - b\eta)y_2 \bigr)e_1
\\&
\kern3em
+
(a\xi - b\eta)y_1 e_2
-
(d\eta - c\xi)y_2 e_3
\mapob.
\end{align*}
Combining these results we find for the right-invariant vector fields associated to the basis $f_1,f_2$ of $\body\Liealg g_1$ at the point $g$ with coordinates $\bigl((a,b,c,d), (\xi,\eta)\bigr)\in \Gextension\Sl(2,\RR) \times \body \Liealg g_1$:
\begin{align*}
f_1^R\caprestricted_g
&
=
\tfrac12(d\eta - c\xi)\, e_1^R\caprestricted_g
+ (a\xi - b\eta)\, e_2^R\caprestricted_g 
+
(1+\tfrac12\xi\eta)d\,\partial_{\xi} + (1+\tfrac12\xi\eta)c\,\partial_{\eta}
\\
f_2^R\caprestricted_g
&
=
\tfrac12(a\xi - b\eta)\, e_1^R\caprestricted_g
+ (c\xi - d\eta)\, e_3^R\caprestricted_g
+
(1+\tfrac12\xi\eta)b\,\partial_{\xi} + (1+\tfrac12\xi\eta)a\,\partial_{\eta}
\mapob.
\end{align*}
In the decomposition $\Hilbert \cong \bigl(L^2\bigl(\Sl(2,\RR)\bigr)\,\bigr)^4$, $\psi \cong (\psi_I)_{I\subset \{1, 2\}} \equiv (\psi_\emptyset, \psi_{\xi\eta}, \psi_\xi, \psi_\eta)$ where the first two components are even and the last two odd, we can write the action of $f_i^R$ on such a function in matrix form as
\begin{align*}
f_1^R 
&
\cong
\begin{pmatrix}
0 & 0 & d & c
\\
0 & 0 
& \tfrac12d - \tfrac12d\,e_1^R + b\,e_2^R  
& \tfrac12c - \tfrac12c\,e_1^R + a\,e_2^R  
\\
-\tfrac12c\,e_1^R + a\,e_2^R & -c & 0 & 0 
\\
\tfrac12d\,e_1^R - b\,e_2^R & d & 0 & 0
\end{pmatrix}
\\[3\jot]
f_2^R
&
\cong
\begin{pmatrix}
0 & 0 & b & a
\\
0 & 0 
& \tfrac12 b + \tfrac12 b\,e_1^R + d\,e_3^R  
& \tfrac12 a + \tfrac12a\,e_1^R + c\,e_3^R 
\\
\tfrac12a\,e_1^R - c\,e_3^R & -a & 0 & 0 
\\
-\tfrac12 b\,e_1^R + d\,e_3^R & b & 0 & 0
\end{pmatrix}
\mapob.
\end{align*}
A direct computation (using that the operators\slash right-invariant vector fields $e_i^R$ are skew-symmetric on $L^2\bigl( \Sl(2,\RR)\bigr)\,$) shows that the operators $f_i^R$ indeed preserve $\supinprsym_\mfdmetric$ (on the space of compactly supported functions $C^\infty_c(G; \CA^\KK)\,$) and that the term $\inprodd{\chi_\emptyset}{\psi_\emptyset}$ (due to the function $\Delta(\xi)$) is essential for this result.

\masection{Appendix: an overview of \texorpdfstring{$\CA$}{A}-manifold theory}
\label{appendixonAmanifoldssection}

This appendix\footnote{It is a truncated version of the appendix I wrote in \cite{Tuynman:2015}} is a very short overview of $\CA$-manifold theory, intended for readers with some familiarity with supermanifold theory, but not with this approach. It will not always be complete, and sometimes it might be slightly besides the truth, but it intends to give the gist of the theory, not the most precise formulation (which can be found in \cite{Tu04}).

\firstofmysubsection
\mysubsection{appendixonAmanifoldssection}
{The basic graded commutative ring}

The starting point of $\CA$-manifold theory is a graded-commutative ring $\CA$. We will fix it as being the exterior algebra of an infinite dimensional (real) vector space $V$: 
$$
\CA=\bigwedge V = \bigoplus_{k=0}^\infty \bigwedge{}{\kern-3pt\raise1.5ex\hbox{$\scriptstyle k$}\kern2pt} V =  \Bigl(\, \bigoplus_{k=0}^\infty \bigwedge {\kern-3pt\raise1.5ex\hbox{$\scriptstyle 2k$}\kern2pt} V \,\Bigr) \oplus \Bigl(\, \bigoplus_{k=0}^\infty \bigwedge{\kern-3pt\raise1.5ex\hbox{$\scriptstyle 2k+1$}\kern2pt}V \,\Bigr)
=
\CA_0\oplus \CA_1
\mapob.
$$
We denote by $\body:\CA\to \CA$ (sic) the canonical projection onto the direct summand $\RR\equiv\bigwedge^0V\subset \bigoplus_{k=0}^\infty\bigwedge^kV$; we will call $\body$ the \stresd{body map} and the image $\body(a)\in \RR$ the \stresd{body of $a\in \CA$}.
The set of nilpotent elements $\nilpotent=\bigoplus_{k=1}^\infty \bigwedge^kV$ (an ideal) is a supplement (in $\CA$) to $\RR=\body(\CA)=\bigwedge^0V$. 
We thus can identify the quotient $\CA/\nilpotent$ (which is a field) with $\body(\CA)\subset \CA$ and the canonical projection $\CA\mapsto \CA/\nilpotent$ with the body map $\body$.
Another feature of (this) $\CA$ (due to the fact that $V$ is infinite dimensional) is that for any $n\in \NN^*$ there exist elements $\xi_1, \dots, \xi_n\in \CA_1=\bigoplus_{k=0}^\infty \bigwedge^{2k+1}V$ such that the product $\xi_1\cdots\xi_n \equiv \xi_1\wedge\cdots\wedge \xi_n$ is non-zero.

\begin{definition}{Remark}
Actually, the particular form of $\CA$ is not that important, as long as it has the above mentioned features: a supplement to the nilpotent elements which is isomorphic to $\RR$ and for any $n\in \NN^*$ elements $\xi_1,\dots, \xi_n\in \CA_1$ whose product is non-zero. And even the last condition is slightly stronger than strictly needed: if an $\CA$-manifold has odd dimension $q$, then there should exist $q+1$ odd elements whose product is non-zero (and thus $V$ should have at least dimension $q+1$). However, as one does not wish to be restricted in the choice of the odd dimension, it is preferable to have this condition for any $n\in \NN^*$, not only for $q+1$. Hence the choice of an infinite dimensional vector space $V$.
\end{definition}

In the sequel we will never have any use for the particular form of $\CA$, so the symbol $V$ will no longer be reserved for the vector space whose exterior algebra is $\CA$.

\mysubsection{appendixonAmanifoldssection}
{\texorpdfstring{$\CA$}{A}-vector spaces}

As $\CA$ is not commutative, there is a difference between left- and right-modules over $\CA$. We will first concentrate on graded bi-modules, meaning an abelian group $M$ which is at the same time a left- and a right-module over $\CA$, and which splits as a direct sum $M=M_0\oplus M_1$ satisfying the conditions
$$
a\in \CA_\alpha \text{ and } m\in M_\beta
\quad\Rightarrow\quad
a\cdot m =(-1)^{\alpha\beta}m\cdot a \in M_{\alpha+\beta}
\mapob,
$$
where the gradings $\alpha,\beta$ should be seen as belonging to $\ZZ/2\ZZ$ and thus $1+1=0$. 
(A slight abuse of naming, as by a \myquote{graded bi-module} one usually means a graded left- and right-module over $\CA$, without the condition that the left- and right-actions are related by $a\cdot m =(-1)^{\alpha\beta}m\cdot a$.)

Again because $\CA$ is not commutative, there is a difference between left-linear homomorphisms and right-linear ones, even for graded bi-modules. 
A left-linear morphism that is even (meaning that it maps homogeneous elements of the source module to homogeneous elements of the target module of the same parity) is automatically also right-linear (and vice-versa), but a similar statement for non-even morphisms is not true. 
One thus has two morphism modules between two graded bi-modules, the left-linear ones and the right-linear ones. 
These morphism modules are in a natural way graded bi-modules over $\CA$ and they are in a natural way isomorphic.

\begin{definition}{Convention}
As we wish to adhere to the convention that interchanging two homogeneous elements (of whatever nature) induces a minus sign whenever both are odd, we are led to denote evaluation of left-linear morphisms on the left instead of on the right as is usual for maps. More precisely, if $f:M\to N$ is a left-linear morphism between the graded bi-modules $M$ and $N$, we will denote the image of the element $m\in M$ by the map $f$ as $\contrf{m}f$. In that way, left-linearity gets the form
$$
\contrf{a\cdot m}f = a\cdot(\contrf{m}f)
\mapob,
$$
instead of
$$
f(am) = af(m)
\mapob,
$$
which would violate the convention in case $f$ is not even.
This notational convention is already used in ordinary differential geometry when one evaluates\slash contracts a $k$-form with a tangent vector (to yield a $k-1$-form).

\end{definition}

The (for us) important constructions on vector spaces can be carried out also for graded bi-modules over $\CA$, such as direct sums, direct products, morphism modules (left or right), tensor products (over the graded-commutative ring $\CA$), exterior algebras and free bi-modules on a set of homogeneous generators. 

\medskip

We next turn our attention to a special case of graded bi-modules over $\CA$: those that are of the form $V\otimes_\RR \CA$ with $V=V_0\oplus V_1$ a graded vector space over $\RR$, which we call \stresd{$\CA$-vector spaces}. 
Inside such an $\CA$-vector space we have the (real) subspace $V\otimes_\RR \body(\CA)\cong V$ of elements that can be written in the form $v\otimes 1$ for some $v\in V$ and the projection (body map) $\body_V:V\otimes_\RR \CA\to V\otimes_\RR \body(\CA)\cong V$ defined by $\body_V(v\otimes a) = v\otimes \body(a)\cong \body(a)\cdot v$.

If $M$ is any (left) module over $\CA$, we can define the subset $\nilpotent_M\subset M$ of those elements $m$ for which there exists $a\in \CA$, $a\neq0$ such that $a\cdot m=0$.
As the reals are included in $\CA$, any module $M$ over $\CA$ is in particular a real vector space and one can show that $\nilpotent_M$ is a vector subspace of $M$ (as a real vector space).
For an $\CA$-vector space $M=V\otimes_\RR \CA$ we have
$$
\nilpotent_M = \{\,v\otimes n\mid v\in V\,,\, n\in\nilpotent\subset \CA\,\}
\mapob.
$$
Another way to define an $\CA$-vector space thus is as those graded bi-modules $M$ for which there exists a supplement $V\subset M$ for $\nilpotent_M$ in the category of graded vector spaces over $\RR$.
The body map then becomes the projection onto this summand $M=V\oplus \nilpotent_M$. 
In this way the definition of an $\CA$-vector space obtains the same flavor as the conditions imposed on our graded commutative ring $\CA$.

\begin{definition}{Remark}
The map $V\mapsto V\otimes_\RR \CA$ from the category of graded vector spaces over $\RR$ to the category of $\CA$-vector spaces is a functor (nearly an isomorphism of categories). In particular the constructions one can perform on these categories (direct sums, tensor products, exterior algebras) are preserved by this map. However, one has to make a choice what to do with morphisms, as in the category of graded vector spaces over $\RR$ there is no difference between left- and right-linear, whereas in the category of $\CA$-vector spaces there is.

\end{definition}

\mysubsection{appendixonAmanifoldssection}
{Smooth functions}
\label{subsecsmoothfunc}

In order to give a definition of smooth super functions from first principles, we first give an alternate description of smooth functions in the ordinary (non-super) case. We then just copy this alternate description to the super case to obtain our definition of super smooth functions.
The starting point is the formula
\begin{equation}\label{fminfequalsintoverder}
f(x) - f(y) = \sum_{i=1}^m (x_i-y_i)\cdot \left( \int_0^1 (\partial_if)\bigl(sx + (1-s)y\bigr)\,\extder s\right)
\mapob,
\end{equation}
valid for any function $f:O\subset \RR^m\to F$ of class $C^1$ defined on a convex set $O$ with values in a normed vector space $F$. If we define the functions $g_i:O\times O\to F$ by
$$
g_i(x,y) = \int_0^1 (\partial_if)\bigl(sx + (1-s)y\bigr)\,\extder s
\mapob,
$$
then we have:
\begin{equation}\label{smoothnesstrick}
\forall x,y\in O : 
f(x)-f(y) = \sum_{i=1}^m (x_i-y_i)\cdot g_i(x,y)
\mapob.
\end{equation}
It then it is easy to show that, for any $k\in \NN$, $f$ is of class $C^{k+1}$ if and only if there exist functions $g_i$ of class $C^k$ such that \recalf{smoothnesstrick} is valid.
Moreover, the partial derivatives of $f$ are given by $(\partial_if)(x) = g_i(x,x)$. As is well known, these partial derivatives are unique, while the functions $g_i$ themselves are not (apart from the case $n=1$).

\begin{definition}{\myquote{Classical} definitions}
A \stresd{smooth system} is an assignment $\smooths$ that associates, for any dimension $m\in\NN$, to any open set $O\subset \RR^m$ and any target (normed vector) space $F$ a collection of \stress{continuous} functions $\smooths(O,F)\subset C^0(O,F)$ verifying the property
\begin{multline*}
\forall f\in \smooths(O,F)
\quad
\exists g_i\in \smooths(O^2,F)
\quad
\forall x,y\in O\subset \RR^m
\quad:\quad
\\
f(x) - f(y) = \sum_{i=1}^m (x_i-y_i)\cdot g_i(x,y)
\mapob.
\end{multline*}
A smooth system $\smooths_1$ is said to be \stresd{smaller than} a smooth system $\smooths_2$, denoted as $\smooths_1\le \smooths_2$, if we have
$$
\forall O: \smooths_1(O,F)\subset \smooths_2(O,F)
\mapob.
$$

\end{definition}

\begin{proclaim}{Theorem}
$C^\infty$ is the (not a) maximal (with respect to the order $\le$) smooth system.
\end{proclaim}

\begin{definition}{Remark}
The above definitions and statements should not be taken literally, as they are wrong as stated. Formula \recalf{fminfequalsintoverder} is valid for convex sets, but not in general for arbitrary (open) sets. \stress{If} one disposes of smooth partitions of unity, one can re-establish such a result, but that would create, in the case of $\RR$, a circular definition, and in the case of $\CC$ it would be impossible. The solution is to use covers at every stage, but that complicates the notation (not the idea), so we left it out in this summary.

\end{definition}

In order to mimic these definitions in the super case, we have to start with a topology.

\begin{definition}[defoftheDeWittTopologyinAppendix]{Definition}
Let $V$ be any normed vector space over $\RR$ and let $E=V\otimes_\RR \CA$ be the associated $\CA$-vector space, \ie, $V=\body E$. 
The \stresd{DeWitt topology} on $E$ is the coarsest topology on $E$ for which the projection $\body:E\to \body E$ is continuous (when $\body E$ is equipped with the topology induced by the given norm; for finite dimensional $V$ this is the usual euclidean topology).
In particular $\CA$ itself is a (finite dimensional) $\CA$-vector space with $\body \CA=\RR\oplus\{0\}$ (i.e., $(\body\CA)_0=\RR$ and $(\body\CA)_1=\{0\}$), and thus is equipped with the DeWitt topology.
It should be noted that the DeWitt topology is not separated\slash Hausdorff and that for any open set $U\subset E$ we have $\body U\subset U$ and $U=\body\mo(\body U)$.

\end{definition}

Now let $E$ be a finite dimensional $\CA$-vector space and let $e_1, \dots, e_p,f_1, \dots, f_q\in \body E$ be a homogeneous basis, i.e., $e_1, \dots, e_p$ is a basis of the vector space (over $\RR$) $(\body E)_0\subset E_0$ and $f_1, \dots, f_q$ a basis of $(\body E)_1\subset E_1$. It follows that any element $e\in E$ is described uniquely by $p+q$ elements $x_1, \dots, x_p, \xi_1, \dots, \xi_q\in\CA$ according to
$$
e = \sum_{i=1}^p x_ie_i + \sum_{j=1}^q \xi_jf_j
\mapob.
$$
Now if $e$ belongs to the \stress{even} part of $E$, then the $x_i$ belong to $\CA_0$ and the $\xi_j$ to $\CA_1$. In other words, an element of $E_0$ is described by $p$ even \myquote{coordinates} and $q$ odd \myquote{coordinates.}
Moreover, a subset $U\subset E_0$ of the even part of $E$ is open (in the DeWitt topology) if and only if there exists an open set $O\subset \RR^p$ such that
\begin{equation}\label{defofopeninE0}
e = \sum_{i=1}^p x_ie_i + \sum_{j=1}^q \xi_jf_j \in U
\qquad\Longleftrightarrow\qquad
(\body x_1, \dots, \body x_p) \in O
\mapob.
\end{equation}
Another way to state this equivalence is the equality $U = (\body\caprestricted_{E_0})\mo(\body U)$ with $\body U = O\oplus \{0\}\subset (\body E)_0 \oplus (\body E)_1=\RR^p \oplus \RR^q$.

\begin{definition}[superdefinitions]{\myquote{Super} definitions}
$\bullet$
A \stresd{super smooth system} is an assignment $\smooths$ that associates to any open set $U\subset E_0$ in the even part of any finite dimensional $\CA$-vector space $E$ and to any $\CA$-vector space $F$ (the target space, equipped with the DeWitt topology induced by a norm on $\body F$) a collection of \stress{continuous} functions $\smooths(U,F)\subset C^0(U,F)$ verifying the two properties
$$
f(\body U) \subset \body F
$$
and
\begin{multline*}
\forall f\in \smooths(U,F)
\quad
\exists g_i\in \smooths(U^2,F)
\quad
\forall x,y\in U\subset E_0
\quad:\quad
\\
f(x) - f(y) = \sum_{i=1}^{p+q} (x_i-y_i)\cdot g_i(x,y)
\mapob.
\end{multline*}

$\bullet$
A super smooth system $\smooths_1$ is said to be \stresd{smaller than} a super smooth system $\smooths_2$, denoted as $\smooths_1\le \smooths_2$, if we have
$$
\forall U,F: \smooths_1(U,F)\subset \smooths_2(U,F)
\mapob.
$$

$\bullet$
$C^\infty$ is \stresd{the} maximal super smooth system with respect to the order $\le$.

\end{definition}

We now have a definition of super smooth functions that is a look-alike to a possible definition of smooth functions in the non-super case, but two questions remain: (1) why the additional condition $f(\body U)\subset \body F$ and (2) can we say anything interesting about these smooth functions?
One possible answer to (1) is that it allows us to give a positive answer to (2). 
Another possible answer to (1) is that we wish to stay as close as possible to non-super smooth functions and that this conditions is trivially satisfied in the non-super case.

\begin{proclaim}[Gextension]{Lemma}
Let $U\subset E_0$ be an open set in the even part of the $p\vert q$-dimensional $\CA$-vector space $E$, let $O=\body U\subset \RR^p$ be the corresponding open subset in $\RR^p$ (see \recalf{defofopeninE0}) and let $f:O\to \body F$ be an ordinary smooth function with values in the real vector space $\body F$ (where we thus ignore the grading).

With these data, we define the function $\mathbf{G}f:U\to F$ by the formula
$$
(\mathbf{G}f)(x_1, \dots, x_p, \xi_1, \dots, \xi_q)
=
\sum_{k=0}^\infty \frac1{k!} \bigl((D^kf)(r_1, \dots, r_p)\bigr)(\overbrace{n,\vrule width0pt height2.3ex\dots, n}^{k \text{ terms}})
\mapob,
$$
where $r_i=\body x_i\in \RR$, $n_i=x_i-r_i\in \nilpotent_0$, $n=(n_1, \dots, n_p)$, $(D^kf)(r_1, \dots, r_p)$ the $k$th order derivative of $f$ at $(r_1, \dots, r_p)\in O$ as $k$-linear symmetric map with values in $\body F$ and $\bigl((D^kf)(r_1, \dots, r_p)\bigr)({n,\dots, n})$ the formal evaluation of $(D^kf)(r_1, \dots, r_p)$ in the even nilpotent coefficients $n_i$.

This $\mathbf{G}f$ is a super smooth function.

\end{proclaim}

\begin{proclaim}[structureofsmoothfunctions]{Theorem}
Let $U\subset E_0$ be an open set in the even part of the $p\vert q$-dimensional $\CA$-vector space $E$, let $O=\body U\subset \RR^p$ be the corresponding open subset in $\RR^p$ and let $F$ be an  
$\CA$-vector space.
Then a function $f:U\to F$ belongs to $C^\infty(U,F)$ (i.e., is super smooth) if and only if there exist (ordinary) smooth functions $f_{I}:O\to \body F$ such that we have
$$
f(x_1, \dots, x_p, \xi_1, \dots, \xi_q)
=
\sum_{I\subset \{1, \dots, q\}} (\mathbf{G}f_I)(x_1, \dots, x_p)\cdot \xi^I
\mapob,
$$
where $\xi^I$ denotes the product
$$
I = \{j_1, \dots, j_k\} \text{ with } 1\le j_1<\cdots <j_k\le q
\qquad\Longrightarrow\qquad
\xi^I = \xi_{j_1} \cdots \xi_{j_k}
\mapob,
$$
with $\xi^\emptyset = 1$.

\end{proclaim}

\begin{proclaim}{Corollary}
Let $U\subset E_0$ be an open set in the even part of the $p\vert q$-dimensional $\CA$-vector space $E$, let $O=\body U\subset \RR^p$ be the corresponding open subset in $\RR^p$ and let $F$ be an  
$\CA$-vector space. Then we have the equality
$$
C^\infty(U,F)
=
C^\infty(O,\body F) \otimes \bigwedge \RR^q
\mapob.
$$

\end{proclaim}

\begin{definition}{Remark}
Let $E$ and $F$ be two $\CA$-vector spaces (equipped with their DeWitt topologies induced by norms on $\body E$ and $\body F$ respectively). 
For any (left- or right-) linear map $f:E\to F$, we have to ask whether $f$ is continuous, and if so, whether it is smooth. 
Now if $E$ is finite dimensional (and we will always assume that the source space is!), $f$ will automatically be continuous, but it need not be smooth. 
It will be smooth if and only if we have the inclusion $f(\body E)\subset \body F$. 
This implies in particular that the space of smooth (left- or right-) linear maps from $E$ to $F$ is isomorphic to the space of continuous $\RR$-linear maps from $\body E$ to $\body F$. 

\end{definition}

Once we know the structure of super smooth functions, we can address the question of partial derivatives. For that we imitate the non-super case and want to define the partial derivative $\partial_if$ as
$$
(\partial_if)(x) = g_i\bigl(x,x)
\mapob,
$$
where the super smooth functions $g_i$ are given by the definition of super smoothness as
$$
f(x) - f(y) = \sum_{i=1}^{p+q} (x_i-y_i)\cdot g_i(x,y)
\mapob.
$$
But for this to be a coherent definition, one must show that the diagonal $g_i(x,x)$ is uniquely determined by $f$, even when the functions $g_i(x,y)$ are not unique. 
And it is here that one needs the fact that there exist $q+1$ odd elements in $\CA_1$ whose product is non-zero: if that condition is satisfied, the $\partial_if$ are well defined; if not, the diagonal functions $g_i(x,x)$ will not be unique and thus our super smooth functions will not have derivatives. 
(Of course one could define the partial derivatives by hand in terms of the functions $f_I$ of \recalt{structureofsmoothfunctions}, but then the whole purpose of giving an intrinsic definition of super smoothness would be superfluous.)

\medskip

With our choice of $\CA$, the partial derivatives are well defined and behave exactly as expected. If the index $i$ is associated to an odd coordinate, $\partial_i$ is an odd derivation of $C^\infty(U,F)$, whereas it is an even derivation when associated to an even coordinate. 
One should note that the $\partial_i$ are \stress{right}-derivations (as opposed to left-derivations), meaning that they are right-linear, which is a consequence of the choice in \recalt{superdefinitions} to write the \myquote{coefficients} $(x_i-y_i)$ to the left of $g_i(x,y)$ instead of to the right.

Once we know what smooth functions and their derivatives are, nearly all of the usual results of real analysis, and in particular the inverse function theorem, remain true (in most cases even the proofs are close copies).

\mysubsection{appendixonAmanifoldssection}
{\texorpdfstring{$\CA$}{A}-manifolds}
\label{subsecAmfds}

Once we know what open subset are and what (local) diffeomorphisms are, we can copy (nearly) all standard differential geometric constructions, such as manifolds (using an atlas and local diffeomorphisms for chart-changing), fiber bundles, vector bundles, the tangent bundle of a manifold, its cotangent bundle etc{\ae}tera. 
There is however one important difference that concerns vector bundles. 
An $\CA$-manifold is modelled by charts that are open subsets in the \stress{even} part of a finite dimensional $\CA$-vector space $E$; as such it has an even and an odd dimension (the dimensions of $(\body E)_0$ and $(\body E)_1$). 
On the other hand, the typical fiber of a vector bundle is \stress{all} of a finite dimensional $\CA$-vector space $F$. 
There are multiple reasons for doing so, one of them being that the standard constructions such as tensor products and exterior powers can be performed on $\CA$-vector spaces, but not on their even parts (at least not in a way that gives satisfactory results). 
Now a full $\CA$-vector space $F$ can be seen in a natural way as the even part of another $\CA$-vector space: the direct sum of $F$ with its parity dual $\prod F$ (the even part of $\prod F$ is the odd part of $F$ and the odd part of $\prod F$ is the even part of $F$), so vector bundles are still $\CA$-manifolds.

Concerning these constructions, one can show the following results. 
\begin{itemize}
\item
The body map can be extended to $\CA$-manifolds and their smooth maps, the result being an ordinary (non-super) manifold and an ordinary smooth map (taking the body map is essentially mapping all nilpotent elements in $\CA$ to zero).

\item
If $M$ is an $\CA$-manifold of graded dimension $p\vert q$, then $\body M$ is an ordinary manifold of dimension $p$ and for any (sufficiently small) open set $O\subset \body M$ with $U=\body\mo O\subset M$ we have the equality
$$
C^\infty(U, \CA) = C^\infty(O)\otimes \bigwedge \RR^q
\mapob,
$$
providing the link with the sheaf theoretic/ringed spaces approach (in this way we create a ringed space on the ordinary manifold $\body M$, and conversely, every appropriate ringed space appears this way).

\item
The set $C^\infty(M)$ of all smooth maps $f:M\to \CA$ on an $\CA$-manifold $M$ is a graded $\RR$-algebra and the set $\Gamma(B)$ of all smooth sections of a vector bundle $B\to M$ is a graded bi-module over $C^\infty(M)$.

\item
The tangent map $Tf:TM\to TN$ associated to a smooth map $f:M\to N$ between $\CA$-manifolds is naturally \stress{left}-linear and even. In particular we have, for local coordinate systems $x$ in $M$ and $y$ in $N$ with $y=f(x)$, the formula
$$
\iota\Bigl({\sum_{i=1}^{p+q}X_i \cdot \fracp{}{x_i}\bigrestricted_x}\Bigr){Tf} = \sum_{i=1}^{p+q}\sum_{j=1}^{p'+q'} X_i \cdot (\partial_if_j)(x)\cdot \fracp{}{y_j}\bigrestricted_y
\mapob,
$$
where $p\vert q$ is the graded dimension of $M$ and $p'\vert q'$ that of $N$.

\item
Any smooth and \stress{even} vector field $X$ on an $\CA$-manifold $M$ can be integrated to produce a flow, i.e., a $1$-parameter group of local diffeomorphisms $\Phi_t$ with $t\in \CA_0$, which has the usual properties of a flow.
According to \recalt{structureofsmoothfunctions} it has a local expression $\Phi_t(x,\xi)\equiv \Phi(t,x,\xi) = \sum_I (\mathbf{G}\Phi_I)(t,x)\cdot \xi^I$ in terms of ordinary smooth functions of real variables and products of odd coordinates. 
These functions can be computed by induction on the number of elements in $I$ (details can be found in \cite[Ch.V\S4]{Tu04}); the equation for $\Phi_\emptyset$ is a \myquote{standard} first order differential equation for the flow of a non-super vector field (it is the flow equation for the ordinary smooth vector field $\body X$ on the ordinary manifold $\body M$) and the other terms are determined by first order inhomogeneous linear differential equations (and as such they do not restrict the domain of definition of the flow). 

\item
A smooth \stress{odd} vector field $X$ can be integrated to a flow $\Phi_\tau$ with an odd time parameter $\tau\in \CA_1$ (but otherwise with the same properties of a flow) if and only if the auto-commutator $[X,X]$ is zero.
In local coordinates the expression for the flow essentially boils down to the formula
$$
\Phi_\tau(x) = x+\tau\cdot X(x)
\mapob.
$$

\item
Lie's third theorem and its converse are true: to any $\CA$-Lie group $G$ (i.e., an $\CA$-manifold equipped with a group structure such that multiplication and inversion are smooth maps) is associated an $\CA$-Lie algebra $\Liealg g$ (isomorphic to the tangent space at the identity and isomorphic to the space of all left-invariant vector fields). 
Conversely, for any finite dimensional $\CA$-Lie algebra $\Liealg g$ there exists a (unique up to isomorphisms) simply connected $\CA$-Lie group $G$ whose $\CA$-Lie algebra is (isomorphic to) $\Liealg g$.

\item
The notion of an $\CA$-Lie group is equivalent to the notion of a \myquote{Lie supergroup pair} or \myquote{super Harish-Chandra pair,} \ie, an \myquote{ordinary} Lie group $H$, an $\CA$-Lie algebra $\Liealg g$ and an action of $H$ on $\Liealg g$ that satisfy the two conditions that the Lie algebra $\Liealg h$ of $H$ is isomorphic to $\body\Liealg g_0$ and such that the restriction to $\Liealg h\cong\body\Liealg g_0$ of the action of $H$ on $\Liealg g$ is the adjoint action.
Starting from an $\CA$-Lie group $G$, one obtains the Lie supergroup pair directly by taking $H=\body G$ and as the action of $H$ on $\Liealg g$ the restriction to $H\subset G$ of the adjoint action of $G$ on $\Liealg g$.
Conversely, if we have a Lie supergroup pair $(H,\Liealg g)$, we first construct the simply connected $\CA$-Lie group $\widetilde G$ whose $\CA$-Lie algebra is $\Liealg g$. It then follows that $\body \widetilde G$ is a covering group of $H$.
If we denote by $Z\subset \body \widetilde G\subset \widetilde G$ the kernel of this covering map, then the condition that the action of $H$ on $\Liealg g$ extends the adjoint action of $H$ on $\Liealg h\cong \body\Liealg g_0$ implies that $Z$ is central (and discrete) in $\widetilde G$. 
One then shows easily that $G=\widetilde G/Z$ is an $\CA$-Lie group whose associated Lie supergroup pair (via the previous construction) is the given Lie supergroup pair.

\end{itemize}

\begin{definition}{Remark}
If we replace (everywhere) the field of real numbers $\RR$ by the field of complex numbers $\CC$ (which means in particular that we replace $\CA$ by $\CA^\CC = \CA\otimes_\RR\CC = \CA \oplus i\CA$), we get complex smooth super functions. As ordinary smooth complex differentiable functions (from open sets in $\CC^p$ to complex vector spaces) are holomorphic, this means that for $U\subset E_o$ open in an $\CA^\CC$-vector space of dimension $p\vert q$ with $O= \body U\subset \CC^p$ and $F$ any $\CA^\CC$-vector space, we have the equality
$$
C^\infty(U,F) = \operatorname{Hol}(O,\body F) \otimes \bigwedge \CC^q
\mapob,
$$
where $\operatorname{Hol}(O,\body F)$ denotes the set of complex analytic\slash holomorphic functions on $O$ with values in $F$.
As for real analytic superfunctions (and despite the fact that I don't have an intrinsic definition in the same spirit as for smoothness), it is easy to define them as a subset of super smooth real functions, just by requiring the functions $f_I$ in \recalt{structureofsmoothfunctions} to be real analytic.

Applying this to manifold theory, we obtain the categories of real and complex analytic $\CA^{(\CC)}$-manifolds (split and non-split). 
As the results on integration of vector fields depends (only) upon the corresponding statement in ordinary non-super analysis, one obtains directly that an analytic vector field on a real or complex analytic $\CA^{(\CC)}$-manifold admits an analytic flow.
It follows that an $\CA$-Lie group admits a unique structure of a real analytic $\CA$-manifold, simply because in the non-super setting one can prove it by means of the flow of an analytic vector field (see \cite{DK00}).

\end{definition}

\section*{Acknowledgements}

I thank my colleagues S.~De Bièvre, H.~Queffélec, F.~Ziegler and in particular V.~Thilliez for their patience, support and crucial help during the preparation of this paper. 
I also thank A.~De Goursac and J.-Ph.~Michel for enlightening discussions.

This work was supported in part by the Labex CEMPI  (ANR-11-LABX-0007-01).

\providecommand{\bysame}{\leavevmode\hbox to3em{\hrulefill}\thinspace}
\providecommand{\MR}{\relax\ifhmode\unskip\space\fi MR }
\providecommand{\MRhref}[2]{%
  \href{http://www.ams.org/mathscinet-getitem?mr=#1}{#2}
}
\providecommand{\href}[2]{#2}

\end{document}